\documentstyle[amssymb, amsfonts]{amsart}

\newenvironment{proof}{\noindent {\bf Proof} }{\endprf\par}
\def \endprf{\hfill  {\vrule height6pt width6pt depth0pt}\medskip}

\theoremstyle{plain}
  \newtheorem{theorem}[subsection]{Theorem}
  
  \newtheorem{proposition}[subsection]{Proposition}
  \newtheorem{lemma}[subsection]{Lemma}
  \newtheorem{corollary}[subsection]{Corollary}

\theoremstyle{remark}

\theoremstyle{definition}
  \newtheorem{definition}[subsection]{Definition}

\begin{document}

  \title{Stability of Spherically Symmetric Wave Maps. }
  \author{J.
  Krieger}\thanks{ Author partially supported by NSF grant
  DMS-0401177.}
\date{}
\maketitle
\begin{abstract}
We study Wave Maps from ${\mathbf{R}}^{2+1}$ to the hyperbolic
plane ${\mathbf{H}}^{2}$ with smooth compactly supported initial
data which are close to smooth spherically symmetric initial data
with respect to some $H^{1+\mu}$, $\mu>0$. We show that such Wave
Maps don't develop singularities in finite time and stay close to
the Wave Map extending the spherically symmetric data(whose
existence is ensured by a theorem of
Christodoulou-Tahvildar-Zadeh) with respect to all $H^{1+\delta},
\delta<\mu_{0}$ for suitable $\mu_{0}(\mu)>0$. We obtain a similar
result for Wave Maps whose initial data are close to geodesic
ones. This strengthens a theorem of Sideris for this context.
\end{abstract}

 \begin{section}{Introduction}\markboth{Joachim
 Krieger}{Introduction}

 Let $\{({\bf{x}},\,{\bf{y}})|{\bf{y}}>0\}$ be the hyperbolic
 plane, equipped with metric
 $dg=\frac{d{\bf{x}}^{2}+d{\bf{y}}^{2}}{{\bf{y}}^{2}}$, and let
 ${\mathbf{R}}^{n+1},\,n\geq 1$ denote the standard Minkowski space equipped
 with the metric $dh=-dx_{0}^{2}+\sum_{i=1}^{n}d_{x_{i}}^{2}$. We shall
 also use the identifications $x_{0}=t$ (time), $\partial_{x_{\nu}}=\partial_{\nu}$, $\nu=0,1,\ldots,n$.
 A Wave Map from
 Minkowski space to ${\mathbf{H}}^{2}$ is a map $u:{\mathbf{R}}^{n+1}\rightarrow {\mathbf{H}}^{2}$
 which is critical with respect to the functional
 \begin{equation}\nonumber
 u\rightarrow
 \int_{{\mathbf{R}}^{n+1}}<\partial_{\alpha}u,\partial^{\alpha}u>_{g}dx_{0}d_{x_{1}}\ldots d_{x_{n}},
 \end{equation}
where $\partial_{\alpha}u = u_{*}(\partial_{\alpha})$,
$\partial^{\alpha}=h^{\alpha\beta}\partial_{\beta}$ and Einstein's
summation convention is in force. The Euler Lagrange equations of
this problem read as follows:
\begin{equation}\label{Euler1}
\Box\ln{\bf{y}}=-\frac{\partial_{\nu}{\bf{x}}}{{\bf{y}}}\frac{\partial^{\nu}{\bf{x}}}{\bf{y}}
\end{equation}
\begin{equation}\label{Euler2}
\Box(\frac{{\bf{x}}}{\bf{y}})=\frac{{\bf{x}}}{\bf{y}}\frac{\partial_{\nu}{\bf{y}}\partial^{\nu}{\bf{y}}+\partial_{\nu}{\bf{x}}\partial^{\nu}{\bf{x}}}{{\bf{y}}^{2}}
\end{equation}
If $n=2$, the fundamental Conjecture associated with this problem
is the following, which flows from the intuition that the negative
curvature of the target should prevent a focusing of energy in
small spatial regions:
\\

{\bf{Conjecture}}(e. g. Klainerman \cite{Kl}) Let $n=2$. Given
smooth initial data $({\bf{x}},{\bf{y}})$,
$(\partial_{t}{\bf{x}},\partial_{t}{\bf{y}}):{\mathbf{R}}^{2}\times\{0\}\rightarrow
{\mathbf{H}}^{2}\times T{\mathbf{H}}^{2}$, there exists a
global-in-time smooth Wave Map extending them.
\\

This is expected to be generalizable to arbitrary targets of
negative curvature and satisfying some geometric niceness
conditions. We stick in this paper to the ${\mathbf{H}}^{2}$ model
on account of its simplicity.\\
The difficulty in establishing the above conjecture stems from the
fact that the problem is energy critical, i. e. the natural scale
invariant Sobolev space is exactly the energy space $\dot{H}^{1}$
(the energy
$\sum_{\alpha=0}^{n}||\frac{\partial_{\alpha}{\bf{x}}}{{\bf{y}}}||_{L_{x}^{2}}^{2}+||\frac{\partial_{\alpha}{\bf{y}}}{{\bf{y}}}||_{L_{x}^{2}}^{2}$
is preserved under the Wave Map flow). Establishing global
regularity for such problems consists customarily of showing that
smooth small data imply global regularity, as well as
non-concentration of energy in physical space. The latter needs to
depend subtly on the geometry of the target, since a priori
analytic reasons cannot rule out a rapid shift of the energy from
low to high frequency modes, resulting in sudden focusing.\\
In the case $n=3$, one expects breakdown of solutions for large
data for analytic reasons (the scale invariant Sobolev space
$\dot{H}^{\frac{3}{2}}$ which in some sense controls the local
well-posedness behavior is not controlled by the energy). We can
formulate
\\

{\bf{Conjecture:}} {\it{Let $n=3$. There exist (large) smooth
initial data
$({\bf{x}},{\bf{y}}),\,(\partial_{t}{\bf{x}},\partial_{t}{\bf{y}}):{\mathbf{R}}^{3}\times\{0\}\rightarrow
{\mathbf{H}}^{2}\times T{\mathbf{H}}^{2}$, which lead to breakdown
in finite time.}}
\\

Breakdown solutions are known in $3+1$ dimensions, but only for
special targets \cite{Caz} not including the hyperbolic plane. \\
The best result known at this point pertaining to the first
Conjecture is the following theorem of the author \cite{Kr-3},
\cite{Kr-4}:
\begin{theorem}\label{small data} Let $n=2,3,\ldots$. Then there exists $\epsilon>0$
such that for smooth initial data
$({\bf{x}},{\bf{y}}),\,(\partial_{t}{\bf{x}},\partial_{t}{\bf{y}}):{\mathbf{R}}^{n}\times\{0\}\rightarrow
{\mathbf{H}}^{2}\times T{\mathbf{H}}^{2}$ satisfying
\begin{equation}\nonumber
\int_{{\mathbf{R}}^{2}}\sum_{\nu=0}^{n}||\frac{\partial_{\nu}{\bf{x}}}{\bf{y}}||_{\dot{H}^{\frac{n}{2}}}
+||\frac{\partial_{\nu}{\bf{y}}}{\bf{y}}||_{\dot{H}^{\frac{n}{2}}}<\epsilon,
\end{equation}
there exists a smooth global-in-time Wave Map extending them.
\end{theorem}
This is similar to earlier results of Tao \cite{Tao 2} when the
target is a sphere, as well as of Tataru \cite{Tatu1} when the
target is uniformly isometrically embeddable into a Euclidean
space. Similar results in dimensions $n\geq 4$ for quite general
targets were achieved by Klainerman-Rodninanski\cite{Kl-R},
Shatah-Struwe\cite{Str-Sh} as well as
Nahmod-Stefanov-Uhlenbeck\cite{N-St-Uh} after Tao's initial breakthrough \cite{Tao 1}, \cite{Tao 2}. \\
Thus the preceding theorem does not yet exhibit behavior
reflecting the geometric nature\footnote{Paradoxically, the proof
of this result involved extra complications over the case of
target a sphere, on account of the fact that one needs to work
with the derivative of the Wave Map, losing one degree of
smoothness.} of ${\mathbf{H}}^{2}$.
\\

What we set out to do in this paper is to try to exploit features
which appear to hinge on geometric properties of this target and
set it apart from positively curved targets such as the sphere
$S^{2}$. We need the following definition:
\\
\begin{definition}: We call a Wave Map
$u:{\mathbf{R}}^{n+1}\rightarrow M$ 'spherically symmetric'
provided $u(t,\rho x)=u(t,x)$ $\forall x\in{\mathbf{R}}^{2}$ and
$\rho:S^{1}\rightarrow SO(2)$ the standard representation of
$S^{1}$ on ${\mathbf{R}}^{2}$.
\end{definition}

We shall use the deep results of
Christodoulou-Tahvildar-Zadeh\cite{Chr-Ta1} on the asymptotic
behavior of spherically symmetric Wave Maps, valid for certain
targets which amongst other things have no conjugate points, to
conclude the following:
\begin{theorem}\label{large data}(Stability) Let $n=2$ and $u[0]=(u(0),\partial_{t}u(0)):{\mathbf{R}}^{2}\times\{0\}
\rightarrow {\mathbf{H}}^{2}\times T{\mathbf{H}}^{2}$ be a smooth
spherically symmetric Wave Map with compactly supported initial
data. Then for any $\sigma>0$ there exists a number $\epsilon>0$
such that for all initial data $\tilde{u}[0]$ which are
$\epsilon$-close to $u[0]$ in $H^{1+\sigma}$, there exists a
smooth global in time Wave Map $\tilde{u}$ extending
$\tilde{u}[0]$. Moreover, $\tilde{u}$ will stay close to $u$ in
the energy topology\footnote{Indeed, even in a certain range of
subcritical spaces $H^{1+\lambda}$.} (in a suitable sense)globally
in time.
\end{theorem}

This is a type of 'large data result', although of course there is
still a smallness assumption present. As far as the case of target
$S^{2}$ is concerned, a surprising result of M. Struwe \cite{Str2}
asserts that smooth radial data lead to global Wave Maps. This
suggests the important question of whether these solutions are
unstable:
\\

{\bf{Question}}{\it{(Instability?). Let
$u[0]:{\mathbf{R}}^{2}\times\{0\}\rightarrow S^{2}\times TS^{2}$
be large generic spherically symmetric initial data. Is it true
that for any $\sigma>0$, $\epsilon>0$, there exist smooth initial
data $\tilde{u}[0]:{\mathbf{R}}^{2}\times\{0\}\rightarrow
S^{2}\times TS^{2}$ with the property that
$||u[0]-\tilde{u}[0]||_{H^{1+\sigma}}<\epsilon$ while\footnote{To
define this norm, use $S^{2}\hookrightarrow {\mathbf{R}}^{3}$ and
use standard coordinates} the smooth Wave Map $\tilde{u}$
extending $\tilde{u}[0]$ locally in time breaks down after finite
time? More precisely for any $\delta>0$
\begin{equation}\nonumber
\exists T<\infty \rightarrow \forall\delta>0\,\lim_{t\rightarrow
T_{-}}||\tilde{u}[t]||_{H^{1+\delta}\times H^{\delta}}=\infty
\end{equation}}}

Unfortunately, our techniques appear to have no bearing on this
question. For example, we don't know what the asymptotic behavior
of Struwe's solutions is.
\\

The key ingredient to prove Theorem~\ref{large data} is the
boundedness of a range of subcritical Sobolev norms for large
spherically symmetric Wave Maps:
\begin{theorem}: There exists $\delta_{0}>0$ such that $\forall
0\leq\delta<\delta_{0}$ and spherically symmetric smooth Wave Maps
$u:{\mathbf{R}}^{2+1}\rightarrow {\mathbf{H}}^{2}$ we have
\begin{equation}\nonumber
\sup_{t}||u[t]||_{H^{1+\delta}\times H^{\delta}}<\infty
\end{equation}
\end{theorem}
The proof of this will follow from the asymptotic results of
\cite{Chr-Ta1}, which in turn rely on a careful analysis of
conservation laws
associated with \eqref{Euler1}, \eqref{Euler2}. \\
We shall then rely on the setup of \cite{Kr-4}, using the
intrinsic derivative formulation (by differentiating
\eqref{Euler1}, \eqref{Euler2}) and passing to the Coulomb Gauge.
The new difficulties by comparison with \cite{Kr-4} concern
nonlinear terms which are linear in the difference
$\nabla[\tilde{u}-u]$. Working in the Coulomb Gauge, this
corresponds to perturbing the flat d'Alembertian $\Box$ with a
potential term $V$ which is in some sense quadratic in the
derivatives of the spherically symmetric Wave Map. We shall show
that the good decay behavior of the spherically symmetric Wave Map
allows us to treat these terms as source terms, instead of having
to modify the linear operator. However, the fact that we cannot
just work with mixed Lebesgue type norms but complicated
null-frame spaces will force simultaneous localizations in
physical and frequency space on us, which make the argument quite
intricate. These types of estimates might be useful when working
on the general large data problem. Our analysis shall have as
simple corollary a generalization of a result of Sideris\cite{Si}
to $n=2$: we define a geodesic Wave Map $u(t,x)$ to be of the form
$u=\gamma(v)$ where $\gamma(.):{\mathbf{R}}\rightarrow
{\mathbf{H}}^{2}$ is a geodesic and $\Box v=0$. Then we have the
following:
\begin{theorem} Let $u(t,x):{\mathbf{R}}^{2+1}\rightarrow {\mathbf{H}}^{2}$ be a smooth geodesic Wave
Map. Then there exists $\epsilon>0$ such that for all initial data
$\tilde{u}[0]$ $\epsilon$-close to $u[0]$ in $H^{1+\sigma}$, there
exists a global Wave Map extending $\tilde{u}[0]$. Also ,
$\tilde{u}$ will stay close to $u$ in the energy topology in a
suitable sense.
\end{theorem}
We also point out that due to a result of Shatah-Tahvildar-Zadeh
\cite{Sh-Z} on the asymptotic behavior of equivariant Wave Maps,
one expects a similar result for perturbations of large
equivariant Wave Maps to hyperbolic targets.
\end{section}
\begin{section}{A priori estimates for spherically symmetric Wave Maps.}

For a Wave Map
$u=({\bf{x}},\,{\bf{y}}):{\mathbf{R}}^{2+1}\rightarrow
{\mathbf{H}}^{2}$, we define the norm $||u(t)||_{H^{s}}$ as
\begin{equation}\nonumber
||({\bf{x}},\,{\bf{y}})||_{H^{s}}:=\sum_{\nu=0}^{2}||\frac{\partial_{\nu}{\bf{x}}}{\bf{y}}(t)||_{H^{s-1}}+||\frac{\partial_{\nu}{\bf{y}}}{\bf{y}}||_{H^{s-1}}
\end{equation}

We also introduce the following notation:
$r=\sqrt{x_{1}^{2}+x_{2}^{2}}$. Now let
$u(t,x)=({\bf{x}},{\bf{y}})$ be a spherically symmetric Wave Map
with compactly supported smooth initial data. Then we have
\begin{lemma}\label{bounded}The image of the Wave Map belongs to a bounded
subset of ${\mathbf{H}}^{2}$. More precisely, we have
\begin{equation}\nonumber
||\ln{\bf{y}}||_{L_{t}^{\infty}L_{x}^{\infty}}<\infty,\,||\frac{{\bf{x}}}{{\bf{y}}}||_{L_{t}^{\infty}L_{x}^{\infty}}<\infty
\end{equation}
The bounds depend (at most) on the size of the support as well as
some norm $||u[0]||_{H^{1+\delta}}$, $\delta>0$.
\end{lemma}
\begin{proof}: We shall rely on the following Proposition in \cite{Chr-Ta1}:
\begin{proposition}\label{Chr1}(Chr.-Tah) Under the previous assumptions, the following
inequalities hold:
\begin{equation}\nonumber
\text{'Good
derivative'}:\,|\frac{\partial_{v}{\bf{x}}}{{\bf{y}}}|+|\frac{\partial_{v}{\bf{y}}}{\bf{y}}|\lesssim
(t+r)^{-\frac{3}{2}},\,\partial_{v}=\partial_{t}+\partial_{r}
\end{equation}
\begin{equation}\nonumber
\text{'Bad
derivative'}:\,|\frac{\partial_{u}{\bf{x}}}{{\bf{y}}}|+\frac{\partial_{u}{\bf{y}}}{\bf{y}}|\lesssim
(|t-r|+1)^{-1}(t+r)^{-\frac{1}{2}},\,\partial_{u}=\partial_{t}-\partial_{r}
\end{equation}
\end{proposition}

By local well-posedness of \eqref{Euler1}, $\eqref{Euler2}$ in
$H^{s}$, $s>1$, there exists a time interval $[-T,T]$ with
$T=T(||u[0]||_{H^{s}})$ on which we uniformly control
$||\ln{\bf{y}}(t)||_{H^{s}}$, $t\in [-T,T]$. Given an arbitrary
point $(t,x)\in{\mathbf{R}}^{2+1}$ at distance
$<\frac{T}{\sqrt{2}}$ from the forward light cone (say), connect
it to a point in the strip $[-T,T]\times{\mathbf{R}}^{2}$ by means
of a null-geodesic $\gamma$ given by $t-r=u=\text{const}$. We have
\begin{equation}\nonumber
|\frac{\partial_{v}{\bf{y}}}{\bf{y}}|\lesssim (t+r)^{-\frac{3}{2}}
\end{equation}
on $\gamma$, whence
\begin{equation}\nonumber
\int_{\gamma\cap \{t\geq
T\}}|\frac{\partial_{v}{\bf{y}}}{\bf{y}}|dt<\infty
\end{equation}
Combining this with the embedding $H^{s}\subset L_{x}^{\infty}$
yields the claim for $\ln{\bf{y}}$ for such points $(t,\,x)$.
Given a point $(t,x)$ at distance $\geq \frac{T}{\sqrt{2}}$ from
the forward light cone, connect it via a geodesic $\gamma$:
$t+r=\text{const}$ to a point in the strip of thickness
$\frac{T}{\sqrt{2}}$ around the light cone. Using
\begin{equation}\nonumber
|\frac{\partial_{u}{\bf{y}}}{\bf{y}}|\lesssim
u^{-1}v^{-\frac{1}{2}}
\end{equation}
yields
\begin{equation}\nonumber
\int_{\gamma}|\frac{\partial_{u}{\bf{y}}}{\bf{y}}|dt<\infty
\end{equation}
Thus the claim follows in general for $\ln{\bf{y}}$. With this,
one proceeds similarly for ${\bf{x}}$.
\end{proof}
The following is the main result of this section:
\begin{proposition}\label{subcritical} Let the assumptions be as in the preceding
lemma. There exists $\epsilon_{0}>0$ such that $\forall
0<\epsilon<\epsilon_{0}$, there exists a global bound
\begin{equation}\nonumber
\sum_{\nu=0}^{2}||\partial_{\nu}\ln{\bf{y}}||_{L_{t}^{\infty}H^{\epsilon}}<C_{\epsilon},\,\sum_{\nu=0}^{2}||\frac{\partial_{\nu}\bf{x}}{\bf{y}}||_{L_{t}^{\infty}H^{\epsilon}}<C_{\epsilon},
\end{equation}
where $C_{\epsilon}$ depends on the size of the support as well as
the $L^{2}$-mass of finitely many derivatives of
$u[0]$\footnote{We are being imprecise here. All that matters to
us is the global bound as stated.}
\end{proposition}
\begin{proof}: We shall need the following lemma, which is also
due to Christodoulou-Tahvildar-Zadeh:
\begin{lemma}\label{Chr2}(Chr.-Tah) Introduce the vector fields
$S=t\partial_{t}+r\partial_{r}$,
$\Omega=r\partial_{t}+t\partial_{r}$. For a smooth spherically
symmetric Wave Map $u(t,x):{\mathbf{R}}^{2+1}\rightarrow
{\mathbf{H}}^{2}$ and $\epsilon>0$, we have
\begin{equation}\nonumber
r^{\frac{1}{2}+\epsilon}v^{1-\epsilon}[|\partial_{v}S{\bf{x}}|+|\partial_{v}S{\bf{y}}|]<\infty
\end{equation}
\begin{equation}\nonumber
r^{\frac{1}{2}+\epsilon}v^{1-\epsilon}[|\partial_{v}\Omega{\bf{x}}+|\partial_{v}\Omega{\bf{y}}|]<\infty
\end{equation}
As a consequence, we conclude that
\begin{equation}\nonumber
|\nabla_{x,t}\partial_{v}{\bf{x}}|+|\nabla_{x,t}\partial_{v}{\bf{y}}|\lesssim
v^{-(1-\epsilon)}u^{-1}r^{-\frac{1}{2}-\epsilon}
\end{equation}
\end{lemma}
We shall in fact prove that the quantity
\begin{equation}\nonumber\begin{split}
&A(t):=\sum_{k\in{\mathbf{Z}}}2^{\delta|k|}[||P_{k}(\ln{\bf{y}}(t,.))||_{\dot{H}_{x}^{1}}+||P_{k}(\frac{\nabla_{x}{\bf{x}}}{{\bf{y}}}(t,.))||_{L_{x}^{2}}\\&\hspace{4cm}+||P_{k}\partial_{t}(\ln{\bf{y}})(t,.)||_{L_{x}^{2}}
+||P_{k}(\frac{\partial_{t}{\bf{x}}}{{\bf{y}}})(t,.)||_{L_{x}^{2}}]\\
\end{split}\end{equation}
is bounded globally in time, provided $\delta>0$
is chosen sufficiently small. We have introduced here the
Littlewood-Paley multipliers $P_{k},\,k\in{\mathbf{Z}}$, which
localize the (spatial) Fourier support to dyadic size $|\xi|\sim
2^{k}$. More precisely, let $\chi_{0}(.)\in
C_{0}^{\infty}({\mathbf{R}}_{>0})$ have support contained in
$(\frac{1}{2},2)$ and satisfy
\begin{equation}\nonumber
\sum_{j\in{\mathbf{Z}}}\chi_{0}(\frac{x}{2^{j}})=1\,\forall
x\in{\mathbf{R}}_{>0}.
\end{equation}
Then we define (see also \cite{St}) $P_{k}f$ for any $f\in
L^{1}({\mathbf{R}})$ via
\begin{equation}\nonumber
\widehat{P_{k}f}(\xi)=\chi_{0}(\frac{|\xi|}{2^{k}})\hat{f}(\xi)
\end{equation}
We need to show that $A(T)\lesssim 1+\int_{0}^{T}A(t)\phi(t)dt$
where $\phi(t)$ is integrable, in order to be able to apply
Gronwall's inequality.

We first establish this for the contribution  from $\ln{\bf{y}}$.
Frequency-localize the equation \eqref{Euler1}, resulting in
\begin{equation}\label{Euler'}
\Box
P_{k}\ln{\bf{y}}=-P_{k}[\frac{\partial_{\nu}{\bf{x}}}{{\bf{y}}}\frac{\partial^{\nu}{\bf{x}}}{\bf{y}}]
\end{equation}
We observe that
$\frac{\partial_{\nu}{\bf{x}}}{{\bf{y}}}\frac{\partial^{\nu}{\bf{y}}}{{\bf{y}}}=\frac{(\partial_{r}-\partial_{t}){\bf{x}}}{{\bf{y}}}
\frac{(\partial_{r}+\partial_{t}){\bf{x}}}{{\bf{y}}}$, and apply a
Littlewood-Paley trichotmomy:
\begin{equation}\label{LP1}\begin{split}
&\Box
P_{k}\ln{\bf{y}}=-P_{k}[P_{<k-10}\frac{(\partial_{r}-\partial_{t}){\bf{x}}}{{\bf{y}}}
P_{[k-5,k+5]}\frac{(\partial_{r}+\partial_{t}){\bf{x}}}{{\bf{y}}}]\\
&\hspace{2cm}-P_{k}[P_{[k-10,k+10]}\frac{(\partial_{r}-\partial_{t}){\bf{x}}}{{\bf{y}}}
P_{<k+15}\frac{(\partial_{r}+\partial_{t}){\bf{x}}}{{\bf{y}}}]\\&\hspace{2cm}-P_{k}[P_{>k+10}\frac{(\partial_{r}-\partial_{t}){\bf{x}}}{{\bf{y}}}
P_{>k+5}\frac{(\partial_{r}+\partial_{t}){\bf{x}}}{{\bf{y}}}]\\
\end{split}\end{equation}
We use the following terminology: $P_{<a}=\sum_{k<a}P_{k}$,
$a,\,k\in{\mathbf{Z}}$, $P_{[a,b]}=\sum_{k\in [a,b]}P_{k}$ etc. We
restrict ourselves to time interval $[0,T]$, and let $T\rightarrow
\infty$. The general case will follow from time reversal symmetry.
Applying Duhamel's formula, we see that we need to control the
norm
$\sum_{k\in{\mathbf{Z}}}2^{\delta|k|}||P_{k}(.)||_{L_{t}^{1}L_{x}^{2}}$
of the right hand side by an expression
$1+\int_{0}^{T}A(t)\phi(t)dt$. We may restrict ourselves to a time
interval $[c,\infty]$ for some $c>0$ (depending on the initial
data), on account of local-in-time well-posedness and finite
propagation speed. We estimate each of the terms on the right-hand
side of \eqref{LP1}: {\it{first assume}} $k\geq 0$.
\\

(i) {\it{The first term}}: We would like to place the 2nd input
$P_{[k-5,k+5]}\frac{(\partial_{r}+\partial_{t}){\bf{x}}}{{\bf{y}}}$
into $L_{t}^{\infty}L_{x}^{2}$ and the first input
$P_{<k-10}\frac{(\partial_{r}-\partial_{t}){\bf{x}}}{{\bf{y}}}$
into $L_{t}^{1}L_{x}^{\infty}$. This doesn't integrate up,
however. Placing the 2nd input into $L_{t}^{\infty}L_{x}^{\infty}$
will work provided $t$ is much larger than $2^{|k|}$, but not in
the opposite case: Thus we subdivide
\begin{equation}\nonumber\begin{split}
&P_{k}[P_{<k-10}\frac{(\partial_{r}-\partial_{t}){\bf{x}}}{{\bf{y}}}
P_{[k-5,k+5]}\frac{(\partial_{r}+\partial_{t}){\bf{x}}}{{\bf{y}}}]
\\&=\phi_{\geq
2^{\frac{k}{C}}}(t))P_{k}[P_{<k-10}\frac{(\partial_{r}-\partial_{t}){\bf{x}}}{{\bf{y}}}
P_{[k-5,k+5]}\frac{(\partial_{r}+\partial_{t}){\bf{x}}}{{\bf{y}}}]\\&+\phi_{<2^{\frac{k}{C}}}(t)P_{k}[P_{<k-10}\frac{(\partial_{r}-\partial_{t}){\bf{x}}}{{\bf{y}}}
P_{[k-5,k+5]}\frac{(\partial_{r}+\partial_{t}){\bf{x}}}{{\bf{y}}}],\\
\end{split}\end{equation}
where $\phi_{\geq a}(t)$, $\phi_{<a}(t)$ are smooth cutoffs  to
dilates of the regions $t\geq a$, $t<a$, adding up to $1$. Also,
$C$ is a large number to be chosen. We can immediately estimate
\begin{equation}\nonumber\begin{split}
&2^{\delta k}||\phi_{\geq
2^{\frac{k}{C}}}(t))P_{k}[P_{<k-10}\frac{(\partial_{r}-\partial_{t}){\bf{x}}}{{\bf{y}}}
P_{[k-5,k+5]}\frac{(\partial_{r}+\partial_{t}){\bf{x}}}{{\bf{y}}}]||_{L_{t}^{1}L_{x}^{2}}
\\&\lesssim 2^{\delta
k}||P_{<k-10}\frac{(\partial_{r}-\partial_{t}){\bf{x}}}{{\bf{y}}}||_{L_{t}^{\infty}L_{x}^{2}}
||\phi_{\geq
2^{\frac{k}{C}}}(t))P_{[k-5,k+5]}\frac{(\partial_{r}+\partial_{t}){\bf{x}}}{{\bf{y}}}||_{L_{t}^{1}L_{x}^{\infty}}\\
&\lesssim 2^{\delta k}2^{-\frac{k}{2C}}\\
\end{split}\end{equation}
This can be summed over $k\geq 0$ provided $\delta<\frac{1}{2C}$.
\\
Now we proceed to the case in which time is dominated by
frequency, $t\lesssim 2^{\frac{k}{C}}$. We shall distinguish
between the region separated from the light cone, where we use
lemma~\ref{Chr2}, as well as the region very close to the light
cone, where we use Proposition~\ref{Chr1} as well as Hoelder's
inequality: we decompose
\begin{equation}\nonumber\begin{split}
&P_{[k-5,k+5]}\frac{(\partial_{r}+\partial_{t}){\bf{x}}}{{\bf{y}}}=
P_{[k-5,k+5]}(\psi_{>\frac{t}{2}}(|t|-|x|)\frac{(\partial_{r}+\partial_{t}){\bf{x}}}{{\bf{y}}})\\
&\hspace{3.5cm}+P_{[k-5,k+5]}(\psi_{\frac{t}{2}>.\geq 2^{-\mu
k}}(|t|-|x|)\frac{(\partial_{r}+\partial_{t}){\bf{x}}}{{\bf{y}}})\\
&\hspace{3.5cm}+P_{[k-5,k+5]}(\psi_{< 2^{-\mu
k}}(|t|-|x|)\frac{(\partial_{r}+\partial_{t}){\bf{x}}}{{\bf{y}}})\\
\end{split}\end{equation}
where the smooth cutoffs $\psi_{>\frac{t}{2}}(.)$,
$\psi_{\frac{t}{2}>.\geq 2^{-\mu k}}(.)$, $\psi_{< 2^{-\mu k}}(.)$
add up to $1$ and localize, respectively, to dilates of the
regions indicated in their subscripts. We let $\mu$ be a small
positive number to be chosen. We have
\begin{equation}\nonumber\begin{split}
&2^{\delta
k}||\phi_{<2^{\frac{k}{C}}}(t)P_{k}[P_{<k-10}\frac{(\partial_{r}-\partial_{t}){\bf{x}}}{{\bf{y}}}
\\&\hspace{4.5cm}P_{[k-5,k+5]}(\psi_{>\frac{t}{2}}(|t|-|x|)\frac{(\partial_{r}+\partial_{t}){\bf{x}}}{{\bf{y}}})]||_{L_{t}^{1}L_{x}^{2}([c,T]\times{\mathbf{R}}^{2})}\\
&\lesssim 2^{(\delta
-1)k}\int_{c}^{T}||P_{<k-10}\frac{(\partial_{r}-\partial_{t}){\bf{x}}}{{\bf{y}}}(t)||_{L_{x}^{\infty}}\\&\hspace{5.3cm}||P_{[k-5,k+5]}\nabla_{x}(\psi_{>\frac{t}{2}}(|t|-|x|)\frac{(\partial_{r}+\partial_{t}){\bf{x}}}{{\bf{y}}})||_{L_{x}^{2}}dt\\
\end{split}\end{equation}
Now using lemma~\ref{Chr2} as well as Proposition~\ref{Chr1}
\begin{equation}\nonumber
||P_{[k-5,k+5]}\nabla_{x}(\psi_{>\frac{t}{2}}(|t|-|x|)\frac{(\partial_{r}+\partial_{t}){\bf{x}}}{{\bf{y}}})||_{L_{x}^{2}}
\lesssim
t^{-(2-\epsilon)}\sqrt{\int_{0}^{\frac{t}{2}}r^{-1-2\epsilon}rdr}\lesssim
t^{-\frac{3}{2}}.
\end{equation}
This implies, reiterating application of Proposition~\ref{Chr1}
\begin{equation}\nonumber\begin{split}
&2^{\delta
k}||\phi_{<2^{\frac{k}{C}}}(t)P_{k}[P_{<k-10}\frac{(\partial_{r}-\partial_{t}){\bf{x}}}{{\bf{y}}}
\\&\hspace{4.5cm}P_{[k-5,k+5]}(\psi_{>\frac{t}{2}}(|t|-|x|)\frac{(\partial_{r}+\partial_{t}){\bf{x}}}{{\bf{y}}})]||_{L_{t}^{1}L_{x}^{2}([c,T]\times{\mathbf{R}}^{2})}\\
&\lesssim 2^{(\delta-1)k}\int_{c}^{T}t^{-\frac{3}{2}}dt,\\
\end{split}\end{equation}
which can be summed over $k\geq 0$, provided $\delta<1$. Next, we
estimate
\begin{equation}\nonumber\begin{split}
&2^{\delta
k}||\phi_{<2^{\frac{k}{C}}}(t)P_{k}[P_{<k-10}\frac{(\partial_{r}-\partial_{t}){\bf{x}}}{{\bf{y}}}
\\&\hspace{2cm}P_{[k-5,k+5]}(\psi_{\frac{t}{2}>.\geq 2^{-\mu
k}}(|t|-|x|)\frac{(\partial_{r}+\partial_{t}){\bf{x}}}{{\bf{y}}})]||_{L_{t}^{1}L_{x}^{2}([c,T]\times{\mathbf{R}}^{2})}\\
&\lesssim
2^{(\delta-1)k}||P_{<k-10}\frac{(\partial_{r}-\partial_{t}){\bf{x}}}{{\bf{y}}}||_{L_{t}^{\infty}L_{x}^{2}}
\\&\hspace{3cm}||P_{[k-5,k+5]}\nabla_{x}(\psi_{\frac{t}{2}>.\geq 2^{-\mu
k}}(|t|-|x|)\frac{(\partial_{r}+\partial_{t}){\bf{x}}}{{\bf{y}}})||_{L_{t}^{1}L_{x}^{\infty}([c,T]\times{\mathbf{R}}^{2})}\\
&\lesssim
2^{(\delta-1+\mu)k}||P_{<k-10}\frac{(\partial_{r}-\partial_{t}){\bf{x}}}{{\bf{y}}}||_{L_{t}^{\infty}L_{x}^{2}}
\int_{c}^{T}t^{-\frac{3}{2}}dt.\\
\end{split}\end{equation}
We can sum here over $k$ provided $\mu+\delta<1$. Finally, we
calculate using Hoelder's inequality as well as
Proposition~\ref{Chr1}
\begin{equation}\nonumber\begin{split}
&2^{\delta
k}||\phi_{<2^{\frac{k}{C}}}(t)P_{k}[P_{<k-10}\frac{(\partial_{r}-\partial_{t}){\bf{x}}}{{\bf{y}}}
\\&\hspace{3cm}P_{[k-5,k+5]}(\psi_{<2^{-\mu
k}}(|t|-|x|)\frac{(\partial_{r}+\partial_{t}){\bf{x}}}{{\bf{y}}})||_{L_{t}^{1}L_{x}^{2}([c,T]\times{\mathbf{R}}^{2})}\\
&\lesssim \min\{T,2^{\frac{k}{C}}\}2^{\delta
k}||P_{<k-10}\frac{(\partial_{r}-\partial_{t}){\bf{x}}}{{\bf{y}}}||_{L_{x}^{\infty}}||P_{[k-5,k+5]}(\psi_{<2^{-\mu
k}}(|t|-|x|)\frac{(\partial_{r}+\partial_{t}){\bf{x}}}{{\bf{y}}})||_{L_{x}^{2}}\\
&\lesssim 2^{(\delta+\frac{1}{C}-\frac{\mu}{2})k}\\
\end{split}\end{equation}
This can be summed over $k\geq 0$ provided we have
$\delta+\frac{1}{C}<\frac{\mu}{2}$. Combining with the conditions
obtained earlier, namely $\delta<\frac{1}{2C}$ as well as
$\delta+\mu<1$, we get $\delta<\frac{1}{7}$.
\\

(ii) {\it{The 2nd term of \eqref{LP1}}} . This term appears
immediate on account of Proposition~\ref{Chr1}. Formally
\begin{equation}\nonumber\begin{split}
&2^{\delta
k}||P_{k}[P_{[k-10,k+10]}\frac{(\partial_{r}-\partial_{t}){\bf{x}}}{{\bf{y}}}
P_{<k+15}\frac{(\partial_{r}+\partial_{t}){\bf{x}}}{{\bf{y}}}]||_{L_{t}^{1}L_{x}^{2}([c,T]\times{\mathbf{R}}^{2})}
\\&\hspace{7cm}\lesssim
\int_{c}^{T}A(t)(1+t)^{-\frac{3}{2}}dt\\
\end{split}\end{equation}
We have to argue more carefully here since $\partial_{u}$ involves
$\partial_{r}=\frac{x_{1}}{r}\partial_{x_{1}}+\frac{x_{2}}{r}\partial_{x_{2}}$.
Decompose
\begin{equation}\label{LP2}\begin{split}
&P_{k}(\frac{\partial_{r}{\bf{x}}}{\bf{y}})=P_{k}[P_{<k-10}(\frac{x_{1}}{r})P_{[k-10,k+10]}(\frac{\partial_{1}{\bf{x}}}{\bf{y}})\\
&+P_{[k-10,k+10]}(\frac{x_{1}}{r})P_{<k+15}(\frac{\partial_{1}{\bf{x}}}{\bf{y}})+P_{>k+10}(\frac{x_{1}}{r})P_{>k+5}(\frac{\partial_{1}{\bf{x}}}{\bf{y}})]\\
&+\text{similar terms}\\
\end{split}\end{equation}
The first term in this Littlewood-Paley trichotomy is estimated
exactly as before, so we treat the 2nd and third term. Let
$\chi_{0}(.)\in C_{0}^{\infty}({\mathbf{R}}_{>0})$ be the cutoff
used for the Littlewood-Paley localizers $P_{k}$. We note that
\begin{equation}\nonumber
P_{k}(\frac{x_{1}}{r})(x) =
2^{2k}\int_{{\mathbf{R}}^{2}}\widehat{\chi_{0}}(2^{k}(x-y))\frac{y_{1}}{|y|}dy=
2^{2k}\int_{{\mathbf{R}}^{2}}\widehat{\chi_{0}}(2^{k}y)[\frac{y_{1}-x_{1}}{|y-x|}-\frac{x_{1}}{|x|}]dy
\end{equation}
On account of the inequality
\begin{equation}\nonumber
|\frac{y_{1}-x_{1}}{|y-x|}-\frac{x_{1}}{|x|}|\lesssim
\min\{\frac{|y|}{|x|},\,1\},
\end{equation}
we get, using the rapid decay of $y\rightarrow
\widehat{\chi_{0}}(2^{k}y)$ outside of a disc of radius $\sim
2^{-k}$:
\begin{equation}\nonumber
|P_{k}(\frac{x_{1}}{r})(x)|\lesssim \min\{\frac{2^{-k}}{|x|},1\}.
\end{equation}
We introduce another cutoff $\chi_{2^{-\frac{k}{2}}}(x)$ which
smoothly localizes to a dilate of the disc
$B_{2^{-\frac{k}{2}}}({\bf{0}})$ centered at ${\bf{0}}=(0,0)$. We
then decompose
\begin{equation}\nonumber\begin{split}
&P_{k}[P_{[k-10,k+10]}(\frac{x_{1}}{r})P_{<k+15}(\frac{\partial_{1}{\bf{x}}}{\bf{y}})P_{<k+15}\frac{\partial_{v}{\bf{x}}}{\bf{y}}]=\\
&P_{k}[\chi_{2^{-\frac{k}{2}}}(x)P_{[k-10,k+10]}(\frac{x_{1}}{r})P_{<k+15}(\frac{\partial_{1}{\bf{x}}}{\bf{y}})P_{<k+15}\frac{\partial_{v}{\bf{x}}}{\bf{y}}]\\
&+P_{k}[(1-\chi_{2^{-\frac{k}{2}}}(x))P_{[k-10,k+10]}(\frac{x_{1}}{r})P_{<k+15}(\frac{\partial_{1}{\bf{x}}}{\bf{y}})P_{<k+15}\frac{\partial_{v}{\bf{x}}}{\bf{y}}]\\
\end{split}\end{equation}
Using Proposition~\ref{Chr1} and Hoelder's inequality, we get
\begin{equation}\nonumber\begin{split}
&2^{\delta k}||P_{k}[\chi_{2^{-\frac{k}{2}}}(x)P_{[k-10,k+10]}(\frac{x_{1}}{r})P_{<k+15}(\frac{\partial_{1}{\bf{x}}}{\bf{y}})P_{<k+15}\frac{\partial_{v}{\bf{x}}}{\bf{y}}]||_{L_{1}^{1}L_{x}^{2}}\\
&\lesssim 2^{\delta
k}||\chi_{2^{-\frac{k}{2}}}(x)P_{[k-10,k+10]}(\frac{x_{1}}{r})||_{L_{t}^{\infty}L_{x}^{2}}||P_{<k+15}(\frac{\partial_{1}{\bf{x}}}{\bf{y}})P_{<k+15}\frac{\partial_{v}{\bf{x}}}{\bf{y}}||_{L_{t}^{1}L_{x}^{\infty}}\\
&\lesssim 2^{(\delta-\frac{1}{2})k}.\\
\end{split}\end{equation}
On the other hand, using the preceding calculations as well as
Proposition~\ref{Chr1} we get
\begin{equation}\nonumber\begin{split}
&2^{\delta k}||P_{k}[(1-\chi_{2^{-\frac{k}{2}}}(x))P_{[k-10,k+10]}(\frac{x_{1}}{r})P_{<k+15}(\frac{\partial_{1}{\bf{x}}}{\bf{y}})P_{<k+15}\frac{\partial_{v}{\bf{x}}}{\bf{y}}]||_{L_{t}^{1}L_{x}^{2}}\\
&\lesssim 2^{\delta
k}||(1-\chi_{2^{-\frac{k}{2}}}(x))P_{[k-10,k+10]}(\frac{x_{1}}{r})||_{L_{t}^{\infty}L_{x}^{\infty}}\\&\hspace{5cm}||P_{<k+15}(\frac{\partial_{1}{\bf{x}}}{\bf{y}})||_{L_{t}^{\infty}L_{x}^{2}}
||P_{<k+15}\frac{\partial_{v}{\bf{x}}}{\bf{y}}||_{L_{t}^{1}L_{x}^{\infty}}\\
&\lesssim 2^{(\delta-\frac{1}{2})k}.\\
\end{split}\end{equation}
Since we have to choose $\delta<\frac{1}{7}$, both can be summed
over $k\geq 0$. The case corresponding to the third term in
\eqref{LP2} as well as the remaining terms are handled
analogously.
\\

(iii) {\it{The third term of \eqref{LP1}}}: We can write
\begin{equation}\nonumber\begin{split}
&P_{k}[P_{>k+10}(\frac{\partial_{u}{\bf{x}}}{\bf{y}})P_{>k+5}(\frac{\partial_{v}{\bf{x}}}{\bf{y}})]
=\sum_{l_{1}>k+10,\,|l_{1}-l_{2}|<5}P_{k}[P_{l_{1}}(\frac{\partial_{u}{\bf{x}}}{\bf{y}})P_{l_{2}}(\frac{\partial_{v}{\bf{x}}}{\bf{y}})]\\
\end{split}\end{equation}
Next, we estimate, using Proposition~\ref{Chr1}
\begin{equation}\nonumber\begin{split}
&\sum_{l_{1}>k+10,\,|l_{1}-l_{2}|<5}2^{\delta k}||P_{k}[P_{l_{1}}(\frac{\partial_{u}{\bf{x}}}{\bf{y}})P_{l_{2}}(\frac{\partial_{v}{\bf{x}}}{\bf{y}})]||_{L_{t}^{1}L_{x}^{2}([c,T]\times{\mathbf{R}}^{2})}\\
&\lesssim\sum_{l_{1}>k+10,\,|l_{1}-l_{2}|<5}2^{\delta(k-l_{1})}\int_{c}^{T}[2^{\delta
l_{1}}||P_{l_{1}}(\frac{\partial_{u}{\bf{x}}}{\bf{y}})(t,.)||_{L_{x}^{2}}]t^{-\frac{3}{2}}dt,
\end{split}\end{equation}
and one can sum here over both $l_{1},\,k$ to obtain  the upper
bound $\lesssim 1+\int_{c}^{T}A(t)t^{-\frac{3}{2}}dt$.  We have to
argue for $\frac{\partial_{u}{\bf{x}}}{\bf{y}}$ as in (ii). This
completes the estimates for case $k\geq 0$. For the case $k<0$, we
have for $M>6$
\begin{equation}\nonumber\begin{split}
&||P_{k}[\frac{\partial_{u}{\bf{x}}}{{\bf{y}}}\frac{\partial_{v}{\bf{x}}}{{\bf{y}}}]||_{L_{t}^{1}L_{x}^{2}}
\lesssim
2^{\frac{2k}{M}}||\frac{\partial_{u}{\bf{x}}}{{\bf{y}}}||_{L_{t}^{\infty}L_{x}^{2}}||\frac{\partial_{v}{\bf{x}}}{{\bf{y}}}||_{L_{t}^{1}L_{x}^{M}}
\lesssim 2^{\frac{2k}{M}},\\
\end{split}\end{equation}
We  have used {\it{Bernstein's inequality}} which states that for
any rectangle $R\subset{\mathbf{R}}^{2}$  and smooth cutoff
$\chi_{R}$ supported in $R$ we have\footnote{We denote the spatial
Fourier transform of $f(x)$ either by $\hat{f}$ or
${\mathcal{F}}f$.}
\begin{equation}\nonumber
||{\mathcal{F}}^{-1}(\chi_{R}{\mathcal{F}}f)||_{L_{x}^{q}}\lesssim
|R|^{\frac{1}{p}-\frac{1}{q}}||f||_{L_{x}^{p}},\,p\leq q.
\end{equation}
Also, the the estimate for
$||\frac{\partial_{v}{\bf{x}}}{{\bf{y}}}||_{L_{t}^{1}L_{x}^{M}}$
follows from  interpolating between the decay estimate for
$||\frac{\partial_{v}{\bf{x}}}{{\bf{y}}}(t)||_{L_{x}^{\infty}}$
and energy conservation.\\  The estimates for
$\frac{\partial_{\nu}{\bf{x}}}{\bf{y}}$ are similar: we have by
the same reasoning as before
\begin{equation}\nonumber
\sum_{k\in{\mathbf{Z}}}2^{\delta
k}||P_{k}[\frac{\partial_{\nu}{\bf{y}}\partial^{\nu}{\bf{y}}+\partial_{\nu}{\bf{x}}\partial^{\nu}{\bf{x}}}{{\bf{y}}^{2}}]||_{L_{t}^{1}L_{x}^{2}([c,T]\times{\mathbf{R}}^{2})}\lesssim
1+\int_{c}^{T}A(t)t^{-\frac{3}{2}}dt
\end{equation}
Now as for the nonlinearity on the right hand side of
\eqref{Euler2}, the small frequency case $k<0$ follows exactly as
above from the boundedness of $\frac{{\bf{x}}}{{\bf{y}}}$, see
lemma~\ref{bounded}. As for the large frequency case, we have the
usual frequency trichotomy
\begin{equation}\nonumber\begin{split}
&P_{k}[\frac{{\bf{x}}}{\bf{y}}\frac{\partial_{\nu}{\bf{y}}\partial^{\nu}{\bf{y}}+\partial_{\nu}{\bf{x}}\partial^{\nu}{\bf{x}}}{{\bf{y}}^{2}}]
=P_{k}[P_{[k-5,k+5]}(\frac{{\bf{x}}}{\bf{y}})P_{<k-10}[\frac{\partial_{\nu}{\bf{y}}\partial^{\nu}{\bf{y}}+\partial_{\nu}{\bf{x}}\partial^{\nu}{\bf{x}}}{{\bf{y}}^{2}}]]\\
&\hspace{4cm}+P_{k}[P_{<k+15}(\frac{{\bf{x}}}{\bf{y}})P_{[k-10,k+10]}[\frac{\partial_{\nu}{\bf{y}}\partial^{\nu}{\bf{y}}+\partial_{\nu}{\bf{x}}\partial^{\nu}{\bf{x}}}{{\bf{y}}^{2}}]]\\
&\hspace{4cm}+P_{k}[P_{>k+5}(\frac{{\bf{x}}}{\bf{y}})P_{>k+10}[\frac{\partial_{\nu}{\bf{y}}\partial^{\nu}{\bf{y}}+\partial_{\nu}{\bf{x}}\partial^{\nu}{\bf{x}}}{{\bf{y}}^{2}}]]\\
\end{split}\end{equation}
We need
\begin{lemma} The following inequality holds:
\begin{equation}\nonumber
\sum_{k\geq 0}2^{\delta
k}||P_{k}\nabla_{x}(\frac{{\bf{x}}}{\bf{y}})(t)||_{L_{x}^{2}}\lesssim
A(t)+1
\end{equation}
\end{lemma}
\begin{proof}: Call the left hand side $B(t)$. Note that
\begin{equation}\nonumber
B(t)\leq A(t)+\sum_{k\geq 0}2^{\delta
k}||P_{k}(\frac{{\bf{x}}}{\bf{y}}\frac{\nabla{\bf{y}}}{{\bf{y}}})(t)||_{L_{x}^{2}}
\end{equation}
We have the frequency trichotomy
\begin{equation}\nonumber\begin{split}
&P_{k}[\frac{\nabla{\bf{y}}}{\bf{y}}\frac{{\bf{x}}}{\bf{y}}]=P_{k}[P_{<k-10}(\frac{\nabla{\bf{y}}}{\bf{y}})P_{[k-5,k+5]}(\frac{{\bf{x}}}{\bf{y}})]\\
&+P_{k}[P_{[k-10,k+10]}(\frac{\nabla{\bf{y}}}{\bf{y}})P_{<k+15}(\frac{{\bf{x}}}{\bf{y}})]+P_{k}[P_{>k+10}(\frac{\nabla{\bf{y}}}{\bf{y}})P_{>k+5}(\frac{{\bf{x}}}{\bf{y}})]\\
\end{split}\end{equation}
The estimate is immediate for the 2nd term on the right hand side.
As to the first, we have
\begin{equation}\nonumber\begin{split}
&2^{\delta
k}||P_{k}[P_{<k-10}(\frac{\nabla{\bf{y}}}{\bf{y}})P_{[k-5,k+5]}(\frac{{\bf{x}}}{\bf{y}})](t)||_{L_{x}^{2}}
\\&\lesssim
2^{(\delta-1)k}||P_{<k-10}(\frac{\nabla{\bf{y}}}{\bf{y}})||_{L_{x}^{\infty}}||\nabla_{x}P_{[k-5,k+5]}(\frac{{\bf{x}}}{\bf{y}})(t)||_{L_{x}^{2}}\\
&\lesssim 2^{(\delta-1)k}B(t)\\
\end{split}\end{equation}
One can also estimate this term by $\lesssim 2^{\delta k}$ from
energy conservation and lemma~\ref{bounded}. The estimate for the
third term in the preceding trichotomy is similar. We conclude
that
\begin{equation}\nonumber
B(t)\lesssim A(t)+\sum_{0\leq k\leq C}2^{\delta
k}+\sum_{k>C}2^{(\delta-1)k}B(t)
\end{equation}
Choosing $C$ large enough, one obtains the claim of the lemma.
\end{proof}
Armed with this, we now have (we may assume $k\geq 10$)
\begin{equation}\nonumber\begin{split}
&2^{\delta k}||P_{k}[P_{[k-5,k+5]}(\frac{{\bf{x}}}{\bf{y}})P_{<k-10}[\frac{\partial_{\nu}{\bf{y}}\partial^{\nu}{\bf{y}}+\partial_{\nu}{\bf{x}}\partial^{\nu}{\bf{x}}}{{\bf{y}}^{2}}]]||_{L_{t}^{1}L_{x}^{2}([c,T]\times{\mathbf{R}}^{2})}\\
&\lesssim
2^{(\delta-1)k}\int_{c}^{T}||\nabla_{x}P_{[k-5,k+5]}(\frac{{\bf{x}}}{\bf{y}})(t)||_{L_{x}^{2}}||P_{<k-10}[\frac{\partial_{\nu}{\bf{y}}\partial^{\nu}{\bf{y}}+\partial_{\nu}{\bf{x}}\partial^{\nu}{\bf{x}}}{{\bf{y}}^{2}}(t)||_{L_{x}^{\infty}}dt\\
\end{split}\end{equation}
Using Proposition~\ref{Chr1} as well as the preceding lemma and
summing over $k\geq 10$, we bound this by $\lesssim
1+\int_{c}^{T}A(t)t^{-\frac{3}{2}}dt$. The estimate for the third
term in the frequency trichotomy preceding the last lemma is more
of the same. Thus we get
\begin{equation}\nonumber
\sum_{k\in{\mathbf{Z}}}2^{\delta|k|}||P_{k}[\frac{{\bf{x}}}{\bf{y}}\frac{\partial_{\nu}{\bf{y}}\partial^{\nu}{\bf{y}}+\partial_{\nu}{\bf{x}}\partial^{\nu}{\bf{x}}}{{\bf{y}}^{2}}]||_{L_{t}^{1}L_{x}^{2}([c,T]\times{\mathbf{R}}^{2})}
\lesssim 1+\int_{c}^{T}A(t)t^{-\frac{3}{2}}dt
\end{equation}
Using Duhamel's formula, we get
\begin{equation}\nonumber
\sum_{k\in{\mathbf{Z}}}2^{\delta|k|}[||P_{k}\nabla_{x}(\frac{{\bf{x}}}{{\bf{y}}})(T)||_{L_{x}^{2}}+||P_{k}\partial_{t}(\frac{{\bf{x}}}{{\bf{y}}})(T)||_{L_{x}^{2}}]\lesssim
1+\int_{c}^{T}A(t)t^{-\frac{3}{2}}dt
\end{equation}
We need to estimate $\frac{\nabla{\bf{x}}}{{\bf{y}}}(T)$, which
differs from the preceding by
$\frac{\bf{x}}{{\bf{y}}}\frac{\nabla{\bf{y}}}{\bf{y}}(T)$. For
frequencies $\geq 0$, this is estimated as in the preceding lemma,
observing that we already improved the estimate for
$||\frac{\nabla{\bf{y}}}{\bf{y}}(T)||_{L_{x}^{2}}$ from the
preceding estimates (i)-(iii). The only case not yet covered
concerns small frequencies. However, we have for $k<0$
\begin{equation}\nonumber\begin{split}
&P_{k}[\frac{\bf{x}}{{\bf{y}}}\frac{\nabla{\bf{y}}}{\bf{y}}]=P_{k}[P_{[k-5,k+5]}(\frac{\bf{x}}{{\bf{y}}})P_{<k-10}(\frac{\nabla{\bf{y}}}{\bf{y}})]\\
&+P_{k}[P_{<k+15}(\frac{\bf{x}}{{\bf{y}}})P_{[k-10,k+10]}(\frac{\nabla{\bf{y}}}{\bf{y}})]+P_{k}[P_{>k+5}(\frac{\bf{x}}{{\bf{y}}})P_{>k+10}(\frac{\nabla{\bf{y}}}{\bf{y}})]\\
\end{split}\end{equation}
Then
\begin{equation}\nonumber\begin{split}
&\sum_{k<0}2^{\delta
|k|}||P_{k}[P_{[k-5,k+5]}(\frac{\bf{x}}{{\bf{y}}})P_{<k-10}(\frac{\nabla{\bf{y}}}{\bf{y}})](T)||_{L_{x}^{2}}
\\&\hspace{3.5cm}\lesssim
\sum_{k<0}2^{\delta|k|}||P_{[k-5,k+5]}(\frac{\bf{x}}{{\bf{y}}})||_{L_{x}^{\infty}}||P_{<k-10}(\frac{\nabla{\bf{y}}}{\bf{y}})(T)||_{L_{x}^{2}},\\
\end{split}\end{equation}
which in turn is bounded by $\lesssim
1+\int_{c}^{T}A(t)t^{-\frac{3}{2}}dt$, as is easily\footnote{Use
Bernstein's inequality.} verified. The estimate for the 2nd term
is immediate and the estimate for the third term as follows:
\begin{equation}\nonumber\begin{split}
&\sum_{k\in{\mathbf{Z}}_{<0}}2^{\delta|k|}||P_{k}[P_{>k+5}(\frac{\bf{x}}{{\bf{y}}})P_{>k+10}(\frac{\nabla{\bf{y}}}{\bf{y}})](T)||_{L_{x}^{2}}\\
&\lesssim
\sum_{k\in{\mathbf{Z}}_{<0}}\sum_{l_{1}>k+10,\,|l_{1}-l_{2}|<5}||P_{k}[P_{l_{2}}(\frac{\bf{x}}{{\bf{y}}})P_{l_{1}}(\frac{\nabla{\bf{y}}}{\bf{y}})](T)||_{L_{x}^{2}}\\
&\lesssim
\sum_{k\in{\mathbf{Z}}_{<0}}\sum_{l_{1}>k+10,\,|l_{1}-l_{2}|<5}2^{\delta
|k|}2^{k-l_{1}}||P_{l_{2}}\nabla_{x}(\frac{\bf{x}}{{\bf{y}}})(T)||_{L_{x}^{2}}||P_{l_{1}}(\frac{\nabla{\bf{y}}}{\bf{y}})||_{L_{x}^{2}}\\
\end{split}\end{equation}
One verifies easily from the preceding estimates that this is
$\lesssim 1+\int_{c}^{\infty}A(t)t^{-\frac{3}{2}}dt$. Putting all
of these ingredients together, we obtain
\begin{equation}\nonumber
A(t)\lesssim 1+\int_{c}^{T}A(t)t^{-\frac{3}{2}}dt.
\end{equation}
The desired upper bound now follows from Gronwall's inequality.
\end{proof}

\begin{corollary}\label{crux} Let
$N(\nabla{\bf{x}},\,\nabla{\bf{y}},\,{\bf{x}},\,{\bf{y}})$ denote
any one of the nonlinearities occuring on the right hand side of
\eqref{Euler1}, \eqref{Euler2}. Then for $\delta$ as in
Proposition~\ref{subcritical}, we have for $\delta<\frac{1}{7}$
\begin{equation}\nonumber
\sum_{k\in{\mathbf{Z}}}2^{\delta
|k|}||P_{k}N(\nabla{\bf{x}},\,\nabla{\bf{y}},\,{\bf{x}},\,{\bf{y}})||_{L_{t}^{1}L_{x}^{2}([-T,T]\times{\mathbf{R}}^{2})}<\infty
\end{equation}
\end{corollary}
This follows from the preceding proof and time reversal symmetry.
In the same vein, we have the following lemma:
\begin{lemma}\label{technical25} Choosing $\delta>0$ small enough, we have the
inequality
\begin{equation}\nonumber
\sum_{k\in{\mathbf{Z}}}2^{\delta
|k|}||P_{k}N(\nabla{\bf{x}},\,\nabla{\bf{y}},\,{\bf{x}},\,{\bf{y}})||_{L_{t}^{2}\dot{H}^{-\frac{1}{2}}}<\infty
\end{equation}
\end{lemma}
\begin{proof}: We work with
$N(...)=\frac{\partial_{\nu}{\bf{x}}}{{\bf{y}}}\frac{\partial^{\nu}{\bf{x}}}{{\bf{y}}}$,
the other cases being similar. Divide into the cases $k\geq 0$ and
$k<0$. In the first case, estimate
\begin{equation}\nonumber\begin{split}
&2^{\delta
k}||P_{k}N(\nabla{\bf{x}},\,\nabla{\bf{y}},\,{\bf{x}},\,{\bf{y}})||_{L_{t}^{2}\dot{H}^{-\frac{1}{2}}}\\&\lesssim
2^{(\delta-\frac{1}{2})k}||\frac{(\partial_{t}-\partial_{r}){\bf{x}}}{{\bf{y}}}||_{L_{t}^{\infty}L_{x}^{2}}
||\frac{(\partial_{t}+\partial_{r}){\bf{x}}}{{\bf{y}}}||_{L_{t}^{2}L_{x}^{\infty}}\lesssim
2^{(\delta-\frac{1}{2})k}\\
\end{split}\end{equation}
In the 2nd case, estimate
\begin{equation}\nonumber\begin{split}
&2^{-\delta
k}||P_{k}N(\nabla{\bf{x}},\,\nabla{\bf{y}},\,{\bf{x}},\,{\bf{y}})||_{L_{t}^{2}\dot{H}^{-\frac{1}{2}}}\\&\lesssim
2^{(-\delta
+\frac{2}{4-}-\frac{1}{2})k}||\frac{(\partial_{t}-\partial_{r}){\bf{x}}}{{\bf{y}}}||_{L_{t}^{\infty}L_{x}^{2}}||\frac{(\partial_{t}+\partial_{r}){\bf{x}}}{{\bf{y}}}||_{L_{t}^{2}L_{x}^{4-}}\lesssim
2^{(-\delta +\frac{2}{4-}-\frac{1}{2})k}\\
\end{split}\end{equation}
We have used here that
\begin{equation}\nonumber
||\frac{(\partial_{t}+\partial_{r}){\bf{x}}}{{\bf{y}}}(t)||_{L_{x}^{4-}}\lesssim
t^{-\frac{3}{4+}}
\end{equation}
which follows from interpolating between Proposition~\ref{Chr1}
and energy conservation. Choosing $\delta
<\frac{2}{4-}-\frac{1}{2}$ results in the claim of the lemma.
\end{proof}
\end{section}

\begin{section}{The perturbation argument}\markboth{Joachim
 Krieger}{The perturbation argument}

\subsection{Precise statement of theorem. Outline of the procedure}

The formulation \eqref{Euler1}, \eqref{Euler2}, while good enough
for the purposes of the last section, will not suffice for us
here\footnote{It appears that the fact that we impose stronger
control over $\tilde{u}$ than just the energy (indeed stronger
than a Besov norm) should allow us to work with the original
coordinate formulation, see e. g. \cite{Tatu}. However, it appears
that the bilinear null-structure in \eqref{Euler1}, \eqref{Euler2}
is not good enough to obtain the gains in time we shall need, see
Proposition~\ref{refined bootstrap}. Indeed, proving an equivalent
of this Proposition for the bilinear expressions appears to
require time decay (in the sense that the norm evaluated on the
function truncated to large times decays) for norms such as
$||u||_{\dot{X}_{k}^{1,\frac{1}{2},\infty}}$, which already fails
for free waves. Moreover, our proof will actually reveal that one
gets an honest $H^{1}$-stability result provided one restricts
oneself to large enough times.}. Instead, following the procedure
in \cite{Kr-4}, we shall pass to the derivative formulation of the
problem, and translate everything into the Coulomb Gauge. More
precisely, introduce the variables
$\phi^{1}_{\nu}=-\frac{\partial_{\nu}{\bf{x}}}{{\bf{y}}}$,
$\phi^{2}_{\nu}=-\frac{\partial_{\nu}{\bf{y}}}{{\bf{y}}}$, pass to
complex notation $\phi_{\nu}=\phi^{1}_{\nu}+i\phi^{2}_{\nu}$, and
revert to the Coulomb Gauge by introducing the variables
$\psi_{\nu}=\phi_{\nu}e^{-i\triangle^{-1}\sum_{i=1,2}\partial_{i}\phi^{1}_{i}}$.
One gets the following remarkable self-contained divergence curl
system:
\begin{equation}\label{Div-Curl1}
\partial_{\alpha}\psi_{\beta}-\partial_{\beta}\psi_{\alpha}
=i\psi_{\beta}\triangle^{-1}\sum_{j=1,2}\partial_{j}(\psi^{1}_{\alpha}\psi^{2}_{j}-\psi^{2}_{\alpha}\psi^{1}_{j})
-i\psi_{\alpha}\triangle^{-1}\sum_{j=1,2}\partial_{j}(\psi^{1}_{\beta}\psi^{2}_{j}-\psi^{2}_{\beta}\psi^{1}_{j})
\end{equation}
\\
\begin{equation}\label{Div-Curl2}
\partial_{\nu}\psi^{\nu}=i\psi^{\nu}\triangle^{-1}\sum_{j=1,2}\partial_{j}(\psi^{1}_{\nu}\psi^{2}_{j}-\psi^{2}_{\nu}\psi^{1}_{j}).
\end{equation}
From these one easily deduces the following system of wave
equations:
\begin{equation}\label{Wave Equation}\begin{split}
\Box\psi_{\alpha} =
&i\partial^{\beta}[\psi_{\alpha}\triangle^{-1}\sum_{j=1}^{2}\partial_{j}[\psi^{1}_{\beta}\psi^{2}_{j}-\psi^{2}_{\beta}\psi^{1}_{j}]]\\
&-i\partial^{\beta}[\psi_{\beta}\triangle^{-1}\sum_{j=1}^{2}\partial_{j}[\psi^{1}_{\alpha}\psi^{2}_{j}-\psi^{2}_{\alpha}\psi^{1}_{j}]]\\
&+i\partial_{\alpha}[\psi_{\nu}\triangle^{-1}\sum_{j=1}^{2}\partial_{j}[\psi^{1\nu}\psi^{2}_{j}-\psi^{2\nu}\psi^{1}_{j}]].\\
\end{split}\end{equation}
As in \cite{Kr-4}, these in conjunction with the underlying
first-order system \eqref{Div-Curl1}, \eqref{Div-Curl2} shall form
the basis for our estimates. We can now give the precise version
of Theorem~\ref{large data}:
\begin{theorem}\label{Main} Let
$u:{\mathbf{R}}^{2+1}\rightarrow {\mathbf{H}}^{2}$ be a smooth
spherically symmetric Wave Map with compactly supported (large)
initial data. Let $\{\psi_{\nu}\}_{\nu=0}^{2}$ be the derivative
components in the Coulomb Gauge. Then for any $\mu>0$ there exists
$\epsilon =\epsilon(u,\mu)>0$ such that for all smooth initial
data $\tilde{u}[0]=(\tilde{u}(0),\,\partial_{t}\tilde{u}(0))$ with
$||(u-\tilde{u})[0]||_{H^{1+\mu}\times H^{\mu}}<\epsilon$, there
exists a smooth Wave Map $\tilde{u}$ extending $\tilde{u}[0]$.
Also, $\tilde{u}$ will stay close to $u$ in the sense that
\begin{equation}\nonumber
\sup_{t}||(\psi_{\nu}-\tilde{\psi}_{\nu})(t)||_{L_{x}^{2}}\lesssim
\epsilon
\end{equation}
\end{theorem}

The proof of this shall consist in analyzing the wave equation
satisfied by the difference
$\delta\psi_{\nu}:=\tilde{\psi}_{\nu}-\psi_{\nu}$. Subtracting the
wave equations for $\tilde{\psi}_{\nu},\,\psi_{\nu}$, and
eliminating the $\tilde{\psi}_{\nu}$ results in terms linear,
quadratic and cubic in $\delta\psi_{\nu}$. As these expressions
have no apparent null-structure in them, we shall revert to the
device of a Hodge-type decomposition used already in \cite{Kr-3},
\cite{Kr-4}: we shall write $\psi_{\nu}=R_{\nu}\psi+\chi_{\nu}$
and similarly for $\tilde{\psi}_{\nu}$, where we impose the
condition $\sum_{i=1,2}\partial_{i}\chi_{i}=0$.  Note that this
results in a similar decomposition for $\delta\psi_{\nu}$. One
easily deduces an elliptic div-curl system for $\chi_{\nu}$,
$\tilde{\chi}_{\nu}$, from which one deduces the schematic
identities $\chi_{\nu}=\nabla^{-1}(\psi\nabla^{-1}(\psi^{2}))$
etc., where the operators $\nabla^{-1}$ stand for linear
combinations of operators of the form
$\triangle^{-1}\partial_{j}$. Plugging these ingredients back into
the wave equations satisfied by the $\delta\psi_{\nu}$ and
eliminating all $\tilde{\psi}_{\nu}$ results in trilinear
null-form terms as well terms of higher degree of linearity,
either linear or of higher degree in the $\delta\psi_{\nu}$. All
of this is just like in \cite{Kr-4}. Terms which are at least
quadratic in the $\delta\psi_{\nu}$ can be treated just as there,
using the fact that Corollary~\ref{crux} shall allow us to
retrieve all the necessary estimates about $\psi_{\nu}$. The only
added difficulty comes from the terms linear in
$\delta\psi_{\nu}$. One way to think of these is as an extra
driving term added to the flat operator $\Box$. However, the very
good decay estimates satisfied by the $\psi_{\nu}$ shall allow us
to treat these terms as source terms instead. The added difficulty
over \cite{Kr-4} we encounter here has to do with the fact that we
need to gain explicitly in time. This will force us to localize
simultaneously in physical and frequency space. In fact, we shall
use a kind of wave packet decomposition to get the necessary
estimates. The next two subsections provide the technical setup.
In the same vein as the preceding theorem, we have
\begin{theorem}\label{Main'} Let $u:{\mathbf{R}}^{2+1}\rightarrow
{\mathbf{H}}^{2}$ be a smooth geodesic Wave Map with compactly
supported initial data. Then for any $\mu>0$, there exists
$\epsilon=\epsilon(u,\mu)>0$ such that for all smooth initial data
$\tilde{u}[0]=(\tilde{u}(0),\partial_{t}\tilde{u}(0))$, with
$||(u-\tilde{u})(0)||_{H^{1+\mu}\times H^{\mu}}<\epsilon$, there
exists a smooth Wave Map $\tilde{u}:{\mathbf{R}}^{2+1}\rightarrow
{\mathbf{H}}^{2}$ extending $\tilde{u}[0]$. $\tilde{u}$ will stay
close to $u$ in the sense that
$\sup_{t}||(\psi_{\nu}-\tilde{\psi}_{\nu})(t)||_{L_{x}^{2}}\lesssim
\epsilon$.
\end{theorem}

\subsection{Sobolev type spaces}

We commence by introducing the functional analytic framework of
\cite{Tatu}, \cite{Tao 2}, \cite{Kr-4} which we have to rely on to
run the perturbation argument. We recall the Littlewood-Paley
multipliers $P_{k}$ introduced in the previous section:
\begin{equation}\nonumber
\widehat{P_{k}f}(\xi)=\chi_{0}(\frac{|\xi|}{2^{k}})\hat{f}(\xi),
\end{equation}
for a suitable cutoff $\chi_{0}(.)$. These are not flexible
enough, and we also introduce the multipliers $Q_{j}$ which
localize the space-time Fourier support to dyadic distance $\sim
2^{j}$ from the light cone: letting
\begin{equation}\nonumber
\tilde{\phi}(\xi,\tau)=\int_{{\mathbf{R}}^{2+1}}e^{-i(\tau
t+\xi\cdot x)}\phi(t,\,x)dt dx
\end{equation}
denote the space-time Fourier transform, we let
\begin{equation}\nonumber
\widetilde{Q_{j}\phi}(\tau,\xi):=\chi_{0}(\frac{||\tau|-|\xi||}{2^{j}})\tilde{\phi}(\tau,\xi)
\end{equation}
where $\chi_{0}(.)$ is as for the $P_{k}$'s. We note that these
definitions entail the identities
\begin{equation}\nonumber
\sum_{k\in{\mathbf{Z}}}P_{k}\phi=\phi,\,\sum_{j\in{\mathbf{Z}}}Q_{j}\phi=\phi,\,\phi\in{\mathcal{S}}({\mathbf{R}}^{2+1}).
\end{equation}
We have the basic inhomogeneous Sobolev spaces $H^{s}$, and their
homogeneous counterparts $\dot{H}^{s}$:
\begin{equation}\nonumber
||\phi||_{H^{s}}=\sqrt{\int_{{\mathbf{R}}^{2}}(1+|\xi|^{2})^{s}|\hat{\phi}(\xi)|^{2}d\xi},\,
||\phi||_{\dot{H}^{s}}=\sqrt{\int_{{\mathbf{R}}^{2}}|\xi|^{2s}|\hat{\phi}(\xi)|^{2}d\xi}
\end{equation}

Note that the space $H^{s}$ is defined as completion of
${\mathcal{S}}({\mathbf{R}}^{2})$ with respect to the first norm.
Trying to do the same for $\dot{H}^{s}$ leads to difficulties (one
gets not necessarily locally integrable functions). We shall only
work with smooth functions anyways, so we only care about
$||.||_{\dot{H}^{s}}$. These norms are not flexible enough, and we
also need the $X^{s,\theta}$ spaces of Klainerman-Machedon as well
as their ('frequency localized') homogeneous Besov analogs(again
only as norms):
\begin{equation}\nonumber
||\phi||_{X^{s,\theta}}:=\sqrt{\int_{{\mathbf{R}}^{2+1}}(1+|\xi|^{2})^{s}(1+||\tau|-|\xi||)^{2\theta}|\tilde{\phi}(\tau,\xi)|^{2}d\tau
d\xi},\,\theta>\frac{1}{2}
\end{equation}
\begin{equation}\nonumber\begin{split}
&||\phi||_{\dot{X}_{k}^{a,b,c}}:=2^{ka}(\sum_{j\in{\mathbf{Z}}}[2^{b
j}||Q_{j}\phi||_{L_{t}^{2}L_{x}^{2}}]^{c})^{\frac{1}{c}},\,c<\infty\\&\hspace{6cm}||\phi||_{\dot{X}_{k}^{a,b,\infty}}:=
2^{ak}\sup_{j\in{\mathbf{Z}}}[2^{b
j}||Q_{j}\phi||_{L_{t}^{2}L_{x}^{2}}]\\
\end{split}\end{equation}
We shall always have $b=\frac{1}{2}$. The latter norms can be
assembled to 'global versions', most naturally via
\begin{equation}\nonumber
||\phi||_{\dot{X}^{a,b,c}}:=\sqrt{\sum_{k\in{\mathbf{Z}}}||P_{k}\phi||_{\dot{X}_{k}^{a,b,c}}^{2}}
\end{equation}
The most intuitive way to think about the $X^{s,\theta}$ etc is to
view them as superpositions of 'twisted free waves', gotten by
foliating space-time by cones $||\tau|-|\xi||=\lambda$. One has
the representation (see \cite{Kl-Se})
\begin{equation}\nonumber
\phi=\int_{\lambda\in{\mathbf{R}}}\phi_{\lambda}e^{it\lambda}d\lambda,
\end{equation}
where $\Box\phi_{\lambda}=0$ and
\begin{equation}\nonumber
\int_{\lambda}||\phi_{\lambda}||_{H^{s}}d\lambda\lesssim
||\phi||_{X^{s,\theta}}.
\end{equation}
At the homogeneous level, we have the embedding\footnote{The way
to think about these is in the sense of inequalities between the
associated norms: $A\subset B\rightarrow ||u||_{B}\lesssim
||u||_{A}$.}
\begin{equation}\nonumber
\dot{X}_{k}^{0,\frac{1}{2},1}\subset L_{t}^{\infty}L_{x}^{2}.
\end{equation}
More generally, the Strichartz estimates (see e. g. \cite{Ke-T})
imply that the following embeddings hold:
\begin{equation}\nonumber
\dot{X}_{k}^{0,\frac{1}{2},1}\subset
2^{k(1-\frac{1}{p}-\frac{2}{q})}L_{t}^{p}L_{x}^{q},
\end{equation}
where $\frac{1}{p}+\frac{1}{2q}\leq\frac{1}{4}$. Similar
embeddings hold for the 'subcritical spaces' $X^{s,\theta}$. We
shall need slightly shrunk versions of the spaces $X^{s,\theta}$
etc. which give stronger control for the 'elliptic regions' far
away from the light cone. For example, we have the norms (see
\cite{Kl-Se}) $||.||_{{\mathcal{X}}^{s,\theta}}$, which are
defined via
\begin{equation}\nonumber
||\phi||_{{\mathcal{X}}^{s,\theta}}:=||\phi||_{X^{s,\theta}}+||\partial_{t}\phi||_{X^{s-1,\theta}}
\end{equation}
Similarly, we introduce the space ${\mathcal{H}}^{s}$ defined as
the completion of ${\mathcal{S}}({\mathbf{R}}^{2})$ under the norm
\begin{equation}\nonumber
||\psi||_{{\mathcal{H}}^{s}}:=||\psi||_{H^{s}}+||\partial_{t}\psi||_{H^{s-1}}
\end{equation}

\subsection{Tataru's null-frame spaces}

This subsection also summarizes material expounded in greater
detail elsewhere (e. g. \cite{Tao 2}, \cite{Kr-3}, \cite{Kr-4}).
The spaces $X^{s,\theta}$ and their homogeneous Besov counterparts
are unfortunately only part of the story. This has to do with the
fact that even the strongest homogeneous versions of these norms
(the norms $||.||_{X^{a,\frac{1}{2},1}}$) do not yield good
algebra type estimates, due to logarithmic divergences in low
frequencies. A solution to this problem is given by 'spaces'
incorporating Tataru's null-frame spaces. We present here a first
version of spaces that overcome this difficulty.  We shall
construct norms $||.||_{S}$ assembled from a family of 'frequency
localized' norms $||.||_{S[k]}$:
\begin{equation}\nonumber
||\phi||_{S}:=\sqrt{\sum_{k\in{\mathbf{Z}}}||P_{k}\phi||_{S[k]}^{2}}
\end{equation}
The norms $||.||_{S[k]}$ in turn are gotten as in \cite{Kr-4}:
they are constructed to satisfy
$||.||_{\dot{X}_{k}^{0,\frac{1}{2},\infty}}\leq
||.||_{S[k]}\lesssim ||.||_{\dot{X}_{k}^{0,\frac{1}{2},1}}$. We
arrange that the norms are invariant under the natural scaling
operation associated with derivatives of Wave Maps in $2+1$
dimensions, since we shall be working at the level of the
derivative. The precise definition of $S[k]$ is complicated: we
first construct norms $||.||_{S[k,\kappa]}$ associated with every
integer $k$ and cap $\kappa\subset S^{1}$. To do so, we introduce
{\it{null-frame coordinates}} $(t_{\omega},\,x_{\omega})$,
$\omega\in S^{1}$, on space-time, whose definition is as follows:
\begin{equation}\nonumber
t_{\omega}=\frac{1}{\sqrt{2}}(1,\omega)\cdot(t,x)
\end{equation}
\begin{equation}\nonumber
x_{\omega}=(t,x)-\frac{t_{\omega}}{\sqrt{2}}(1,\omega)
\end{equation}
Thus these are Cartesian coordinates with respect to a tilted
reference frame, whose 'time axis' with direction
$\frac{1}{\sqrt{2}}(1,\omega)$ lies along the light cone. Now we
introduce the space $PW[\kappa]$ defined as the atomic Banach
space whose atoms are Schwartz functions $\psi\in
{\mathcal{S}}({\mathbf{R}}^{2+1})$ satisfying
\begin{equation}\nonumber
\inf_{\omega\in
\tilde{\kappa}}||\psi||_{L_{t_{\omega}}^{2}L_{x_{\omega}}^{\infty}}\leq
1,
\end{equation}
where $\tilde{\kappa}$ is a slightly grown version of $\kappa$
(say by a factor $\frac{11}{10}$) concentric with it. Thus for
$\psi\in {\mathcal{S}}({\mathbf{R}}^{2+1})$,we have
\begin{equation}\nonumber
||\psi||_{PW[\kappa]}:=\inf_{\int_{\tilde{\kappa}}\psi_{\omega}d\omega=\psi}\int_{\tilde{\kappa}}||\psi_{\omega}||_{L_{t_{\omega}}^{2}L_{x_{\omega}}^{\infty}}d\omega
\end{equation}
Moreover, we put for $\psi$ as above
\begin{equation}\nonumber
||\psi||_{NFA[\kappa]^{*}}=\sup_{\omega\notin
2\kappa}\text{dist}(\omega,\kappa)||\psi||_{L_{t_{\omega}}^{\infty}L_{x_{\omega}}^{2}}
\end{equation}
Now we put
\begin{equation}\nonumber
||\psi||_{S[k,\kappa]}=||\psi||_{L_{t}^{\infty}L_{x}^{2}}+2^{-\frac{k}{2}}|\kappa|^{-\frac{1}{2}}||\psi||_{PW[\kappa]}
+|\psi||_{NFA[\kappa]^{*}}
\end{equation}
This definition immediately entails the following fundamental
first bilinear inequality
\begin{equation}\label{bilinear1}
||\phi\psi||_{L_{t}^{2}L_{x}^{2}}\lesssim
\frac{|\kappa'|^{\frac{1}{2}}2^{\frac{k'}{2}}}{\text{dist}(\kappa,\kappa')}||\phi||_{S[k,\kappa]}||\psi||_{S[k',\kappa']},
\end{equation}
provided $2\kappa\cap 2\kappa'=\emptyset$. We now construct the
norms $||\psi||_{S[k]}$ by evaluating suitably microlocalized
pieces of $\psi$ with respect to the $||.||_{S[k,\kappa]}$, taking
a suitable mean and combining this with
$||.||_{\dot{X}_{k}^{a,b,c}}$ type norms. The null-frame norms may
be thought of as controlling the 'free wave-like' character of
$\psi$, while the remaining norms may be thought of
as controlling the 'elliptic character' of $\psi$. \\
For every integer $l<-10$, subdivide $S^{1}$ into a uniformly
finitely overlapping collection $K_{l}$ of caps $\kappa$ of
diameter $2^{l}$. Also, for every integer $\lambda$ with
$-10\geq\lambda\geq l$, we subdivide the angular sector
$\{\xi\in{\mathbf{R}}^{2}|\frac{\xi}{|\xi|}\in\kappa,\,|\xi|\sim
2^{k}\}$ into a uniformly finitely overlapping collection
$C_{k,\kappa,\lambda}$ of slabs $R$ of width $2^{k+\lambda}$. We
introduce various localization operators associated with these
regions: for each $\kappa\in K_{l}$, choose a smooth cutoff
$a_{\kappa}:S^{1}\rightarrow {\mathbf{R}}_{\geq 0}$ supported on a
dilate of $\kappa$. These are to be chosen such that
$\sum_{\kappa\in K_{l}}a_{\kappa}=1$. We also introduce cutoffs
$m_{R}(.):{\mathbf{R}}_{>0}\rightarrow {\mathbf{R}}_{\geq 0}$ such
that the cutoff $m_{R}(|\xi|)a_{\kappa}(\frac{\xi}{|\xi|})$
localizes to a dilate of the slab $R$. Also, we require that
$\sum_{R\in C_{k,\kappa,\lambda}}m_{R}(|\xi|)=
\chi_{0}(\frac{|\xi|}{2^{k}})$. We have the associated pseudo
differential operator $\tilde{P}_{R}\psi$:
\begin{equation}\nonumber
\widehat{\tilde{P}_{R}\psi}(t,\,\xi)=m_{R}(|\xi|)a_{\kappa}(\frac{\xi}{|\xi|})\hat{\psi}(\xi)
\end{equation}
We also have the $\Psi$DO's $P_{k,\kappa}$ associated with
multiplier
$a_{\kappa}(\frac{\xi}{|\xi|})\chi_{0}(\frac{|\xi|}{2^{k}})$.
Then, almost\footnote{The original definition also contained a
norm $||P_{k}Q_{\geq
k}\partial_{t}\psi||_{L_{t}^{M}\dot{H}^{-1+\frac{1}{M}}}$.
Inspection of the proof there yields that this is superfluous,
though, since the norm $||.||_{\dot{X}_{k}^{-\frac{1}{2},-1,2}}$
suffices for the elliptic estimates.} as in \cite{Kr-4} we define
\begin{equation}\label{original definition}\begin{split}
||\psi||_{S[k]}: =
&||\psi||_{L_{t}^{\infty}L^{2}}+||\psi||_{\dot{X}_{k}^{0,\frac{1}{2},\infty}}+||\psi||_{\dot{X}_{k}^{-\frac{1}{2},1,2}}
\\
&+\sup_{\pm}\sup_{l<-10}\sup_{-10\geq\lambda\geq
 l}|\lambda|^{-1}(\sum_{\kappa\in
K_{l}}\sum_{R\in C_{k,\kappa,\lambda}}||\tilde{P}_{R}Q^{\pm}_{<k+2l}\psi||^{2}_{S[k,\pm\kappa]})^{\frac{1}{2}}.\\
\end{split}\end{equation}
This norm looks very complicated, but it isn't too hard to get
control over its ingredients. A fundamental inequality \cite{Kr-4}
for example states that
\begin{equation}\label{tame}
||P_{k}Q_{<k+O(1)}\psi||_{S[k]}\lesssim
||P_{k}\psi||_{\dot{X}_{k}^{0,\frac{1}{2},1}}
\end{equation}
We also need a norm $||.||_{N}$ with respect to which we evaluate
the nonlinearities of our wave equations . We shall put
\begin{equation}\nonumber
||\psi||_{N}:=\sqrt{\sum_{k\in{\mathbf{Z}}}||P_{k}\psi||_{N[k]}^{2}},
\end{equation}
where the $N[k]$ will be constructed as atomic Banach spaces. More
precisely, we let $N[k]$ be the atomic Banach space whose atoms
are Schwartz functions $F\in{\mathcal{S}}({\mathbf{R}}^{2+1})$
with spatial Fourier support contained in the region $|\xi|\sim
2^{k}$ and
\begin{enumerate}
\item{$||F||_{L_{t}^{1}\dot{H}^{-1}}\leq 1$ and $F$ has modulation
$<2^{k+100}$.}
\\
\item{$F$ is at modulation $\sim 2^{j}$ and satisfies
$||F||_{L_{t}^{2}L_{x}^{2}}\leq 2^{\frac{j}{2}}2^{k}$.}
\\
\item{$F$ satisfies $||F||_{\dot{X}_{k}^{-\frac{1}{2},-1,2}}\leq
1$.}
\\
\item{There exists an integer $l<-10$, and Schwartz functions
$F_{\kappa}$ with Fourier support in the region
\\
\begin{equation}\nonumber
\{(\tau,\xi)|\pm\tau>0,\,\mid|\tau|-|\xi|\mid\leq
2^{k-2l-100},\,\frac{\xi}{|\xi|}\in \pm\kappa\}
\end{equation}
\\
with the properties
\\
\begin{equation}\nonumber
F=\sum_{\kappa\in K_{l}}F_{\kappa},\,(\sum_{\kappa\in
K_{l}}||F_{\kappa}||_{NFA[\kappa]}^{2})^{\frac{1}{2}}\leq 2^{k}
\end{equation}
\\
In the last inequality, $NFA[\kappa]$ denotes the dual of
$NFA[\kappa]^{*}$ (the completion of
${\mathcal{S}}({\mathbf{R}}^{2+1})$ with respect to
$||.||_{NFA[\kappa]^{*}}$) used in the definition of
$S[k,\kappa]$: Thus $NFA[\kappa]$ is the atomic Banach space whose
atoms $F$ satisfy
\begin{equation}\nonumber
\frac{1}{\text{dist}(\omega,\,\kappa)}||F||_{L_{t_{\omega}}^{1}L_{x_{\omega}}^{2}}\leq
1
\end{equation}
for some $\omega\notin 2\kappa$.}
\end{enumerate}
Observe that $\dot{X}_{k}^{-1,-\frac{1}{2},1}\subset N[k]$,
$P_{k}(L_{t}^{1}\dot{H}^{-1})\subset N[k]$.

This definition immediately entails the fundamental 2nd bilinear
inequality (again $2\kappa\cap 2\kappa'=\emptyset$)
\begin{equation}\label{bilinear2}
||\phi\psi||_{NFA[\kappa]}\lesssim
\frac{2^{\frac{k'}{2}}|\kappa'|^{\frac{1}{2}}}{\text{dist}(\kappa,\kappa')}||\phi||_{L_{t}^{2}L_{x}^{2}}||\psi||_{S[k',\kappa']}
\end{equation}
We quickly summarize here the main properties of these spaces we
shall need, all proved in \cite{Kr-4}: all functions $\phi,\,\psi$ etc. below shall be in ${\mathcal{S}}({\mathbf{R}}^{2+1})$.\\
(a): {\it{Product type estimate}}
\begin{equation}\nonumber\begin{split}
||P_{k}Q_{j}(P_{k_{1}}\phi_{1}P_{k_{2}}\phi_{2})||_{\dot{X}_{k}^{0,\frac{1}{2},\infty}}\lesssim
&2^{\min\{k_{1},k_{2}\}}2^{\min\{\frac{j-\min\{k,k_{1},k_{2}\}}{4+},0\}}\\&2^{\min\{\frac{\max\{k_{1},k_{2}\}-j}{2},0\}}||P_{k_{1}}\phi_{1}||_{S[k_{1}]}||P_{k_{2}}\phi_{2}||_{S[k_{2}]}\\
\end{split}\end{equation}
(b): {\it{Bilinear null-form estimate}}: let $R_{\nu}$ denote the
Riesz type operator
$R_{\nu}=\frac{\partial_{\nu}}{\sqrt{-\triangle}}$. For $0\leq
p<\frac{1}{4}$, we have\footnote{Moreover, applying an operator
$Q_{l}$ in front, where $l>>k$, we can include factors
$2^{\frac{\min\{l-k_{2},0\}}{2}}$ on the right hand side.}
\begin{equation}\label{eqn.Q(v,j)}\begin{split}
&||P_{k}[R_{1}P_{k_{1}}\psi_{1}R_{2}P_{k_{2}}\psi_{2}-R_{2}P_{k_{1}}\psi_{1}R_{2}P_{k_{2}}\psi_{2}]||_{\dot{X}_{0}^{p,-p,2}}
\\&\hspace{6cm}\leq
C_{p}2^{\frac{\min\{k_{1},k_{2},k\}}{2}}\prod_{i=1,2}||P_{k_{i}}\psi_{i}||_{S[k_{i}]}\\
\end{split}\end{equation}
\begin{equation}\nonumber\begin{split}
&||P_{k}[R_{1}P_{k_{1}}\psi_{1}R_{\nu}P_{k_{2}}\psi_{2}-R_{2}P_{k_{1}}\psi_{1}R_{\nu}P_{k_{2}}\psi_{2}]||_{\dot{X}_{0}^{p,-p,2}}
\\&\hspace{6cm}\leq
C_{p}|k-k_{1}|^{2}2^{\frac{\min\{k_{1},k_{2},k\}}{2}}\prod_{i=1,2}||P_{k_{i}}\psi_{i}||_{S[k_{i}]}\\
\end{split}\end{equation}
(c): {\it{Trilinear null-form estimates}}: these arise upon
formulating the derivative wave equations in the Coulomb Gauge and
applying Hodge type decompositions, as explained below: let
$I=\sum_{k\in{\mathbf{Z}}}P_{k}Q_{<k+100}$:
\begin{equation}\label{null-form1}\begin{split}
&||\partial^{\beta}P_{0}[R_{\alpha}P_{k_{1}}\psi_{1}\triangle^{-1}\sum_{j=1}^{2}\partial_{j}I[R_{\beta}P_{k_{2}}\psi_{2}R_{j}P_{k_{3}}\psi_{3}-R_{\beta}P_{k_{3}}\psi_{3}R_{j}P_{k_{2}}\psi_{2}]]
\\&+\partial_{\alpha}P_{0}[R_{\beta}P_{k_{1}}\psi_{1}\triangle^{-1}\sum_{j=1}^{2}\partial_{j}I[R^{\beta}P_{k_{2}}\psi_{2}R_{j}P_{k_{3}}\psi_{3}-R_{j}P_{k_{2}}\psi_{2}R^{\beta}P_{k_{3}}\psi_{3}]]||_{N[0]}\\
&\leq
C2^{\delta_{1}\min\{-\min\{k_{1},k_{2},k_{3}\},0\}}\prod_{i}2^{\delta_{2}\min\{\max_{j\neq i}\{k_{i},k_{i}-k_{j}\},0\}}\prod_{l}||P_{k_{l}}\psi_{l}||_{S[k_{l}]},\\
\end{split}\end{equation}
\begin{equation}\label{null-form2}\begin{split}
&||P_{0}\partial^{\beta}[R_{\beta}P_{k_{1}}\psi_{1}\triangle^{-1}\sum_{j}\partial_{j}I[R_{\alpha}P_{k_{2}}\psi_{2}R_{j}P_{k_{3}}\psi_{3}-R_{j}P_{k_{2}}\psi_{2}R_{\alpha}P_{k_{3}}\psi_{3}]]||_{N[0]}\\
&\leq
C2^{\delta_{1}\min\{-\min\{k_{1},k_{2},k_{3}\},0\}}\prod_{i}2^{\delta_{2}\min\{\max_{j\neq i}\{k_{i},k_{i}-k_{j}\},0\}}\prod_{l}||P_{k_{l}}\psi_{l}||_{S[k_{l}]}.\\
\end{split}\end{equation}
Of course one may rescale these, i. e. replace $P_{0}$ by $P_{k}$,
$k\in{\mathbf{Z}}$. Then one needs to replace
$\min\{-\min\{k_{1},k_{2},k_{3}\},0\}$ by
$\min\{-\min\{k_{1}-k,k_{2}-k,k_{3}-k\},0\}$ etc.\\
(d): {\it{'energy inequality'}} The following relates the spaces
$S[k]$ and $N[k]$:
\begin{equation}\label{energy}\begin{split}
&||P_{k}\psi||_{S[k]([-T,T]\times{\mathbf{R}}^{2})}\lesssim
\inf_{0<T_{0}\leq T}[\min\{2^{k}T_{0},1\}^{-\frac{1}{2}}||\Box
P_{k}\psi||_{N[k]([-T,T]\times{\mathbf{R}}^{2})}\\&\hspace{6cm}+\sup_{t_{0}\in[-T_{0},T_{0}]}||P_{k}\psi[t_{0}]||_{L^{2}\times \dot{H}^{-1}}].\\
\end{split}\end{equation}
In this inequality, one can leave out the factor
$\min\{2^{k}T_{0},1\}^{-\frac{1}{2}}$  and replace $S[k]$ with the
following stronger norm
\begin{equation}\nonumber\begin{split}
||\psi||_{S'[k]}:=&2^{-k}||\nabla_{x,t}\psi||_{L_{t}^{\infty}L_{x}^{2}}+||\nabla_{x,t}\psi||_{\dot{X}_{k}^{-1,\frac{1}{2},\infty}}\\
&+\sup_{\pm}\sup_{l<-10}\sup_{-10\geq\lambda\geq
 l}|\lambda|^{-1}(\sum_{\kappa\in
K_{l}}\sum_{R\in C_{k,\kappa,\lambda}}||\tilde{P}_{R}Q^{\pm}_{<k+2l}\psi||^{2}_{S[k,\pm\kappa]})^{\frac{1}{2}},\\
\end{split}\end{equation}
provided the norm $N[k]$ is
replaced by $L_{t}^{1}\dot{H}^{-1}$ or $\dot{X}_{k}^{-1,-\frac{1}{2},1}$.\\
(e): {\it{Relation to Strichartz type spaces:}} Let $p,\,q$
satisfy $\frac{1}{p}+\frac{1}{2q}<\frac{1}{4}$. Then we have
\begin{equation}\nonumber
||P_{0}\psi||_{L_{t}^{p}L_{x}^{q}}\leq C_{p,q}||P_{0}\psi||_{S[0]}
\end{equation}
(f): {\it{Relation to improved Strichartz type estimates:}} for
$l<-10$, let $C_{0,l}$ be a covering of the frequency region
$|\xi|\sim 1$ by uniformly finitely overlapping discs of radius
$\sim 2^{l}$. Let $P_{0,c}$ localize the Fourier support to the
disc $c$, such that $\sum_{c\in C_{0,l}}P_{c}=P_{0}$. Let $8\geq
p>4$. Then we have
\begin{equation}\nonumber
(\sum_{c\in
C_{0,l}}||P_{0,c}\psi||_{L_{t}^{p}L_{x}^{q}}^{2})^{\frac{1}{2}}\leq
C_{p}2^{(\frac{3}{4+}-\frac{2}{p})l}||\psi||_{S[0]}
\end{equation}
(g): {\it{Bilinear inequality relating the $S[k]$, $N[k]$:}} Let
$F,\,\psi\in {\mathcal{S}}({\mathbf{R}}^{2+1})$. Then for $j\leq
\min\{k_{1,2}\}+O(1)$ we have
\begin{equation}\nonumber\begin{split}
&||P_{k}[P_{k_{1}}\psi P_{k_{2}}Q_{j}F]||_{N[k]}\\&\lesssim
2^{\delta(j-\min\{k_{1},k_{2}\})}2^{\min\{k_{1}-k_{2},0\}}||P_{k_{1}}\psi_{1}||_{S[k_{1}]}||P_{k_{2}}F||_{\dot{X}_{k_{2}}^{0,-\frac{1}{2},\infty}}\\
\end{split}\end{equation}
If $k_{1}=k_{2}+O(1)$ $j\leq r+k$ for some
$r\leq 0$, we have
\begin{equation}\nonumber
||P_{k}Q_{<r+k}[P_{k_{1}}Q_{j}F
P_{k_{2}}Q_{<2k+r-k_{1}}\psi]||_{N[k]}\lesssim 2^{\delta
r}||F||_{\dot{X}_{k_{1}}^{0,-\frac{1}{2},\infty}}||P_{k_{2}}\psi_{2}||_{S[k_{2}]}.
\end{equation}

\subsection{A modification of the spaces; Moser type estimates}

In spite of the above properties, the spaces $S[k]$ don't appear
flexible enough to handle Moser type estimates of the kind we
shall need. More precisely, the property
\begin{equation}\nonumber
||\psi_{1} A(\nabla^{-1}\psi_{2})||_{S}\lesssim
C(||\psi_{1}||_{S},\,||\psi_{2}||_{S})
\end{equation}
where $\psi_{1,2}\in{\mathcal{S}}({\mathbf{R}}^{2})$, $A(.)\in
C^{\infty}({\mathbf{R}})$ with bounded derivatives, and
$\nabla^{-1}$ schematic notation for linear combinations of
operators of the form $\triangle^{-1}\partial_{j}$, appears
violated. This is a consequence of the fact that the product
estimate 3.4(a) does not allow us to recover enough exponential
gains in the difference $k-k_{1}$ if $k_{1}>>k$, the high-high
interaction case. \\
One way around this would be to re-engineer the way functions get
subdivided into 'free wave parts' and 'elliptic parts'. Indeed,
one has better product estimates than 3.4(a) for free waves, see
Klainerman-Foschi\footnote{It appears, however, that even for free
waves, a high-high-low interaction resulting in an elliptic
product (Fourier support very far from the light cone) does not
lead to the desired exponential gain in the frequencies.}
\cite{Kl-F}. Our way here around this shall be to 'enlarge the
space $S[k]$', shrinking the norm $||.||_{S[k]}$ suitably. More
precisely, we analyze the 'bad high-high' frequency interactions
and observe that by virtue of the spherical symmetry assumption,
these cases are actually favorable in some sense. Indeed, we shall
be able to exploit the well-known fact (e. g. \cite{So},
\cite{Ste}) that the range of admissible Strichartz estimates is
significantly improved in this situation:
\begin{theorem}\label{spherical Strichartz} Let $\psi\in
{\mathcal{S}}(\mathbf{R}^{2+1})$ be invariant under rotations.
Then we have the inequality
\begin{equation}\nonumber
||P_{0}\psi||_{L_{t}^{p}L_{x}^{q}}\lesssim
||\psi[0]||_{L_{t}^{\infty}L_{x}^{2}}+||\Box\psi||_{L_{t}^{1}L_{x}^{2}}
\end{equation}
provided the condition $\frac{1}{p}+\frac{1}{q}<\frac{1}{2}$
holds.
\end{theorem}
We note that theorem~\ref{spherical Strichartz} implies easily the
following corollary:
\begin{corollary} The derivative components $\psi_{\nu}$ satisfy
the estimates
\begin{equation}\nonumber
||P_{k}\psi_{\nu}||_{L_{t}^{p}L_{x}^{q}}\lesssim c_{k}
\end{equation}
for a system of numbers $\{c_{k}\}_{k\in{\mathbf{Z}}}$ ('frequency
envelope') which satisfies
$\sum_{k\in{\mathbf{Z}}}2^{\delta|k|}c_{k}<\infty$, provided the
condition $\frac{1}{p}+\frac{1}{q}<\frac{1}{2}$ holds and
$\delta>0$ is sufficiently small.
\end{corollary}
\begin{proof} This follows by applying a simple frequency
trichotomy to the frequency localized expression for $\psi_{\nu}$
in terms of $\frac{\partial_{\nu}{\bf{x}}}{{\bf{y}}}$,
$\frac{\partial_{\nu}{\bf{y}}}{{\bf{y}}}$. The latter are
controlled by application of theorem~\ref{spherical Strichartz} as
well as Corollary~\ref{crux}.
\end{proof}

\begin{definition}: We put
\begin{equation}\nonumber
||\psi||_{L}:=\sup_{(p,q)|\frac{1}{p}+\frac{1}{2q}\leq
\frac{1}{4}}2^{k(\frac{1}{p}+\frac{2}{q}-1)}\sup_{l\leq
0}2^{-\max\{\frac{1}{2+}-\frac{1}{q},0\}l}(\sum_{c\in
C_{k,l}}||P_{c}\psi||_{L_{t}^{p}L_{x}^{q}}^{2})^{\frac{1}{2}}
\end{equation}
Here $C_{k,l}$ is a finitely overlapping cover of the frequency
region $|\xi|\sim 2^{k}$ by discs of radius $2^{k+l}$, with
associated Fourier localizers $P_{c},\,c\in C_{k,l}$. Also, put
$2+=2+\frac{1}{1000}$, say, and let $\mu$ be a small positive
number, say $\frac{1}{1000}$. We let ${\mathcal{S}}[k]$ be the
atomic Banach space whose atoms satisfy one of the following:
\begin{enumerate}
\item{{\it{Type 1 atoms}}: Fix $\delta_{0}$ small, say
$\delta_{0}=\frac{1}{1000}$. These are functions
$\psi\in{\mathcal{S}}({\mathbf{R}}^{2+1})$ satisfying
\begin{equation}\nonumber\begin{split}
&||\psi||_{A[k]}=||\psi||_{L}+||\psi||_{\dot{X}_{k}^{0,\frac{1}{2},\infty}}+||\psi||_{\dot{X}_{k}^{-\frac{1}{2},1,2}}
+\sup_{\frac{1}{p}+\frac{1}{q}<\frac{1}{2}-\delta_{0}}2^{(\frac{1}{p}+\frac{2}{q}-1)k}||\psi||_{L_{t}^{p}L_{x}^{q}}\\&+\sup_{\pm}\sup_{l<-10}\sup_{-10\geq\lambda\geq
 l}|\lambda|^{-1}(\sum_{\kappa\in
K_{l}}\sum_{R\in C_{k,\kappa,\lambda}}||\tilde{P}_{R}Q^{\pm}_{<k+2l}\psi||^{2}_{S[k,\pm\kappa]})^{\frac{1}{2}}\leq 1\\
\end{split}\end{equation}}
\item{{\it{Atoms of the 2nd type}}: Let
$I=\sum_{k\in{\mathbf{Z}}}P_{k}Q_{<k+100}$. Then
$\psi\in{\mathcal{S}}({\mathbf{R}}^{2+1})$ is of the 2nd type
provided
\begin{equation}\nonumber\begin{split}
&||\psi||_{B[k]}:=\sup_{\frac{1}{p}+\frac{1}{q}<1-10\delta_{0}}2^{k[\frac{1}{p}+\frac{2}{q}-1]}||\psi||_{L_{t}^{p}L_{x}^{q}}
+||\psi||_{L}\\&\hspace{6cm}+||(1-I)\psi||_{\dot{X}_{k}^{-(\frac{1}{2}-\mu),1-\mu,1}}\leq
1,
\end{split}\end{equation}
The range of Lebesgue exponents $(p,q)$ includes the
pairs$(1+,\infty)$, $(\infty,1+)$.\\}
\end{enumerate}
\end{definition}
Note that any function $\psi\in {\mathcal{S}}({\mathbf{R}}^{2+1})$
may be decomposed into two pieces $\psi=\alpha+\beta$ satisfying
\begin{equation}\nonumber
||\alpha||_{A[k]}\lesssim
||\psi||_{{\mathcal{S}}[k]},\,||\beta||_{B[k]}\lesssim
||\psi||_{{\mathcal{S}}[k]}
\end{equation}
We call $\alpha$ 'of first type' and $\beta$ 'of 2nd type'. We let
\begin{equation}\nonumber
||\psi||_{{\mathcal{S}}}:=(\sum_{k\in{\mathbf{Z}}}||P_{k}\psi||_{{\mathcal{S}}[k]}^{2})
\end{equation}
Unfortunately, these norms do not quite suffice to close all the
estimates. The following theorem contains some bilinear estimates,
which we were unable to build into a linear framework. These have
to be proved independently. To state the theorem, we use the
concept of {\it{frequency envelope:}} following Tao \cite{Tao 1},
we call a sequence of positive numbers
$\{c_{k}\}_{k\in{\mathbf{Z}}}$ a frequency envelope provided
$C^{-1}c_{b}2^{-\sigma|a-b|}\leq c_{a}\leq C2^{\sigma|a-b|}c_{b}$
for some $\sigma>0$ and $C>>1$.

\begin{theorem}\label{Moser} Let
$\psi_{1,2}\in{\mathcal{S}}({\mathbf{R}}^{2+1})$. Assume that
$||P_{k}\psi_{1,2}||_{{\mathcal{S}}[k]}\lesssim c_{k}$ for a
sufficiently flat frequency envelope $\{c_{l}\}$. Let
\begin{equation}\nonumber
||\nabla^{-1}\psi_{2}||_{L_{t}^{\infty}L_{x}^{\infty}}\leq C_{0}
\end{equation}
More precisely, assume that for each $k\in{\mathbf{Z}}$ one may
split $P_{k}\psi_{2}=\alpha_{k}+\beta_{k}$ into functions of first
and 2nd type, respectively, such that
$||\nabla^{-1}\sum_{k}\alpha_{k}||_{L^{\infty}}\leq C_{0}$,
$||\nabla^{-1}\sum_{k}\beta_{k}||_{L^{\infty}}\leq C_{0}$. Also
assume that $\psi_{2}$ satisfies the bilinear estimates stated
further below\footnote{Substitute $\psi_{2}$ instead of
$\psi_{1}A(\nabla^{-1}\psi_{2})$.}. Then we can
conclude\footnote{The implied constants in the statements below
depend on $C_{0}$, the constant $C$ in the bilinear estimates
below as well as the decay of the frequency envelope and the
constants chosen in the definition of ${\mathcal{S}}[k]$.}
$\forall k\in {\mathbf{Z}}$
\begin{equation}\nonumber
||P_{k}[\psi_{1}\nabla^{-1}\psi_{2}]||_{{\mathcal{S}}[k]}\lesssim
c_{k}
\end{equation}
In particular, if $A(.):{\mathbf{R}}\rightarrow {\mathbf{C}}$ is
real analytic with bounded derivatives of arbitrary order, and
$\psi_{2}$ real valued, we can conclude that $\forall k\in
{\mathbf{Z}}$
\begin{equation}\nonumber
||P_{k}[\psi_{1}A(\nabla^{-1}\psi_{2})]||_{{\mathcal{S}}[k]}\lesssim
c_{k}
\end{equation}
A more precise version is as follows: for suitable $\delta>0$,
\begin{equation}\nonumber\begin{split}
&||P_{k}[P_{k_{1}}\psi_{1}\nabla^{-1}P_{k_{2}}\psi_{2}]||_{{\mathcal{S}}[k]}
\lesssim 2^{-\delta|k-k_{1}|}||P_{k_{1}}\psi_{1}||_{{\mathcal{S}}[k_{1}]}||P_{k_{2}}\psi_{2}||_{{\mathcal{S}}[k_{2}]}\\
\end{split}\end{equation}
Now assume that
$||\psi_{1,2}||_{S}+||\psi_{1,2}||_{\dot{B}_{2}^{0,1}}<C_{0}$.
Then we can conclude
\begin{equation}\nonumber\begin{split}
&||P_{k}[\psi_{1}A(\nabla^{-1}\psi_{2})]||_{\dot{X}_{k}^{0,\frac{1}{2},\infty}}\leq
C,\\
\end{split}\end{equation}
for a suitable $C=C(C_{0})$. More precisely, decomposing
\begin{equation}\nonumber
P_{k}[\psi_{1}A(\nabla^{-1}\psi_{2})]=\alpha+\beta
\end{equation}
into functions of first and 2nd type, respectively, we may assume
that
\begin{equation}\nonumber
||R_{0}(\alpha+\beta)||_{L_{t}^{\infty}L_{x}^{2}}\leq
C,\,||\beta||_{\dot{X}_{k}^{0,\frac{1}{2},1}}\leq C
\end{equation}
for
$C=C(C_{0},\sup_{k\in{\mathbf{Z}}}||R_{0}P_{k}\psi_{1}||_{L_{t}^{\infty}L_{x}^{2}},\sup_{k\in{\mathbf{Z}}}||R_{0}P_{k}\psi_{2}||_{L_{t}^{\infty}L_{x}^{2}})$.
Next, we have the bilinear estimates
\begin{equation}\nonumber\begin{split}
&\sup_{\phi|\,||\phi||_{S[k_{3}]}\leq
1}\sup_{k_{1,2,3}}\sup_{l<-10}2^{-\frac{\min\{k_{1},k_{3}\}}{2}}2^{-\delta_{1}l}
\\&\hspace{4cm}(\sum_{c\in
C_{k_{3},l}}||P_{k_{1}}[P_{k_{2}}R_{0}[\psi_{1}A(\nabla^{-1}\psi_{2})]P_{c}\phi]||_{L_{t}^{2}L_{x}^{2}}^{2})^{\frac{1}{2}}\lesssim
C\\
\end{split}\end{equation}
\begin{equation}\nonumber\begin{split}
&\sup_{\phi|\,||\phi||_{S[k_{3}]}\leq
1}\sup_{k_{1,2}\in{\mathbf{Z}}}2^{\frac{k_{2}(1-3\mu)}{2+2\mu}}||P_{k_{1}}[P_{k_{2}}R_{0}[\psi_{1}A(\nabla^{-1}\psi_{2})]P_{k_{2}+O(1)}\phi]||_{L_{t}^{2}L_{x}^{1+\mu}}
\leq C,\\
\end{split}\end{equation}
where
$C=C(||\psi_{1}||_{S},||\psi_{2}||_{S},||\psi_{2}||_{\dot{B}_{2}^{0,1}})$,
recall the definition in 3.4. The same estimates hold provided one
replaces $A(\nabla^{-1}\psi_{2})$ by
$A(\nabla^{-1}(\psi_{2}\nabla^{-1}\psi_{3}))$, where $\psi_{3}$
satisfies similar estimates as $\psi_{2}$.
\end{theorem}
The proof of this is a long calculation deferred to an appendix.

\subsection{Proof of the Theorem~\ref{Main}}

To show existence, it suffices to show that some subcritical norm
$\sum_{\nu=0}^{2}||\delta\psi_{\nu}(t)||_{H^{\delta}}$ is globally
bounded in $t$. Indeed, reasoning exactly as in \cite{Kr-4}, one
deduces that for every finite time interval $[-T,T]$, one has with
$\tilde{u}=(\tilde{{\bf{x}}},\,\tilde{{\bf{y}}})$ and
$\tilde{\delta}<\delta$:
\begin{equation}\nonumber
\sup_{t\in[-T,T]}\sum_{\nu=0}^{2}[||\tilde{{\bf{x}}}(t)||_{H^{1+\tilde{\delta}}}+||\tilde{{\bf{y}}}(t)||_{H^{1+\tilde{\delta}}}+||\partial_{t}\tilde{{\bf{x}}}(t)||_{H^{\tilde{\delta}}}+||\partial_{t}\tilde{{\bf{y}}}(t)||_{H^{\tilde{\delta}}}]<\infty
\end{equation}
Using the subcritical result of\\ Klainerman-Machedon \cite{Kl-M
1}, one deduces from here that there can't be breakdown of
smoothness
after finite time. \\
Global boundedness of a subcritical norm in turn shall follow from
the following Bootstrap Proposition: to formulate it, we shall
need time-localized versions of the spaces $S[k]$: for $\psi\in
C_{0}^{\infty}([-T,T]\times{\mathbf{R}}^{2})$ define
\begin{equation}\nonumber
||\psi||_{S[k]([-T,T]\times{\mathbf{R}}^{2})}:=\inf_{\tilde{\psi}\in{\mathcal{S}}({\mathbf{R}}^{2+1})|\tilde{\psi}|_{[-T,T]\times{\mathbf{R}}^{2}}=\psi}||\tilde{\psi}||_{S[k]}
\end{equation}
We use similar definitions for
$||.||_{S([-T,T]\times{\mathbf{R}}^{2})}$
$||\psi||_{S([-T,T]\times{\mathbf{R}}^{2+1})}$ etc. and also
different time intervals $[T_{1},T_{2}]$ etc.
\begin{proposition}\label{bootstrap}
In the situation of Theorem~\ref{Main}, let the smooth Wave Map
extending $\tilde{u}[0]$ exist on the time interval $[-T,T]$.
There exists $T_{1}>0$ such that for $T\geq T_{1}$ and every $K>0$
sufficiently large, there exists $\epsilon>0$ such that the
following conclusion applies: Introduce the frequency envelope
\begin{equation}\nonumber
\tilde{c}_{k}:=\sup_{t\in
[-T_{1},T_{1}]}\sum_{k_{1}\in{\mathbf{Z}}}2^{-\sigma|k-k_{1}|}||P_{k_{1}}\delta\psi_{\nu}(t)||_{L_{x}^{2}}
+\epsilon c_{k},
\end{equation}
where $\{c_{k}\}$ is as in the proof of lemma~\ref{auxiliary
envelope}, and assume \\ $ \sup_{t\in
[-T_{1},T_{1}]}\sum_{\nu=0,1,2}||\delta\psi_{\nu}(t)||_{L_{x}^{2}}\leq\epsilon$.
Then for any $T\geq T_{1}$ we have
\begin{equation}\nonumber\begin{split}
&\sup_{\nu}||P_{k}\delta\psi_{\nu}||_{S[k]([T_{1},T]\times{\mathbf{R}}^{2})}<K\tilde{c}_{k}
\Rightarrow ||P_{k}\delta\psi_{\nu}||_{S[k]([T_{1},T]\times{\mathbf{R}}^{2})}<\frac{K}{2}\tilde{c}_{k}.\\
\end{split}\end{equation}
A similar inequality holds by replacing $T,\,T_{1}$ by
$-T,\,-T_{1}$.

\end{proposition}
Assuming this for now, we continue with the proof of
Theorem~\ref{Main}. We claim that by local well-posedness of
\eqref{Euler1}, \eqref{Euler2}, there exists
$\epsilon=\epsilon(\mu)>0$ such that (using terminology of
Theorem~\ref{Main}) $||u[0]-\tilde{u}[0]||_{H^{1+\mu}\times
H^{\mu}}<\epsilon$ implies that $\tilde{u}$ extends smoothly to
$[-T_{1},T_{1}]$, where $T_{1}$ is as in the preceding
Proposition. To see that this is possible, we shall apply an
inequality of Klainerman-Selberg to the equation satisfied by the
differences ${\bf{\tilde{x}}}-{\bf{x}}$,
${\bf{\tilde{y}}}-{\bf{y}}$ of the coordinate representations of
the perturbed and the spherically symmetric Wave Map,
$\tilde{u}=({\bf{\tilde{x}}},{\bf{\tilde{y}}})$ and
$u=({\bf{x}},{\bf{y}})$. Subdivide the interval $[-T_{1},T_{1}]$
into small subintervals $I_{i}$, for which\footnote{The
$\theta>\frac{1}{2}$ is chosen in dependence of $\mu$, see
Klainerman-Selberg \cite{Kl-Se}.}
$||({\bf{x}},{\bf{y}})||_{{\mathcal{X}}^{1+\mu,\theta}(I_{i}\times{\mathbf{R}}^{2})}\lesssim
1$. This is possible by Corollary~\ref{crux} and local
well-posedness of \eqref{Euler1}, \eqref{Euler2} in $H^{1+\mu}$.
Note that
\begin{equation}\nonumber
\Box(\ln{\bf{\tilde{y}}}-\ln{\bf{y}})
=(\frac{\partial_{\nu}{\bf{\tilde{x}}}}{{\bf{\tilde{y}}}}-\frac{\partial_{\nu}{\bf{x}}}{{\bf{y}}})\frac{\partial^{\nu}{\bf{\tilde{x}}}}{\bf{\tilde{y}}}
+\frac{\partial^{\nu}{\bf{x}}}{{\bf{y}}}(\frac{\partial_{\nu}{\bf{\tilde{x}}}}{{\bf{\tilde{y}}}}-\frac{\partial_{\nu}{\bf{x}}}{{\bf{y}}})
\end{equation}
with a similar equation holding for
$\Box(\frac{{\bf{\tilde{x}}}}{{\bf{\tilde{y}}}}-\frac{{\bf{x}}}{{\bf{y}}})$.
We deduce that upon denoting $I_{i}=[a_{i},a_{i+1}]$, we have
\begin{equation}\nonumber\begin{split}
&\sum_{\nu}||\frac{\partial_{\nu}{\bf{\tilde{y}}}}{{\bf{\tilde{y}}}}-\frac{\partial_{\nu}{\bf{y}}}{{\bf{y}}}||_{{\mathcal{X}}^{\mu,\theta}(I_{i}\times{\mathbf{R}}^{2})}+||\partial_{\nu}(\frac{{\bf{\tilde{x}}}}{{\bf{\tilde{y}}}})-\partial_{\nu}(\frac{{\bf{\tilde{x}}}}{{\bf{\tilde{y}}}})||_{{\mathcal{X}}^{\mu,\theta}(I_{i}\times{\mathbf{R}}^{2})}
\\&\lesssim\sum_{\nu}||\frac{\partial_{\nu}{\bf{\tilde{y}}}}{{\bf{\tilde{y}}}}-\frac{\partial_{\nu}{\bf{y}}}{{\bf{y}}}||_{H^{\mu}(a_{i})}+||\partial_{\nu}(\frac{{\bf{\tilde{x}}}}{{\bf{\tilde{y}}}})-\partial_{\nu}(\frac{{\bf{\tilde{x}}}}{{\bf{\tilde{y}}}})||_{H^{\mu}(a_{i})}
\\&+|I_{i}|^{\epsilon(\mu)}[||\ln{\bf{\tilde{y}}}-\ln{\bf{y}}||_{{\mathcal{X}}^{1+\mu,\theta}(I_{i}\times{\mathbf{R}}^{2})}+||\frac{{\bf{\tilde{x}}}}{{\bf{\tilde{y}}}}-\frac{{\bf{x}}}{{\bf{y}}}||_{{\mathcal{X}}^{1+\mu,\theta}(I_{i}\times{\mathbf{R}}^{2})}]\\
\\&+[|I_{i}|^{\epsilon(\mu)}||\ln{\bf{\tilde{y}}}-\ln{\bf{y}}||_{{\mathcal{X}}^{1+\mu,\theta}(I_{i}\times{\mathbf{R}}^{2})}^{3}+||\frac{{\bf{\tilde{x}}}}{{\bf{\tilde{y}}}}-\frac{{\bf{x}}}{{\bf{y}}}||_{{\mathcal{X}}^{1+\mu,\theta}(I_{i}\times{\mathbf{R}}^{2})}^{3}]\\
\end{split}\end{equation}
We have used here the fact, due to Klainerman-Selberg, that
\begin{equation}\nonumber
||\phi||_{{\mathcal{X}}^{s,\theta}([-T,T]\times{\mathbf{R}}^{2})}\lesssim
||\phi[0]||_{H^{s}\times
H^{s-1}}+T^{\epsilon}||\Box\phi||_{X^{s-1,\theta-1}([-T,T]\times{\mathbf{R}}^{2})},\,\theta>\frac{1}{2},s>1
\end{equation}
as well as the following inequality of Klainerman-Machedon
\begin{equation}\nonumber
||\partial_{\nu}u_{1}\partial^{\nu}u_{2}||_{X^{s-1,\theta-1}}\lesssim
||u_{1}||_{{\mathcal{X}}^{s,\theta}}||u_{2}||_{{\mathcal{X}}^{s,\theta}}
\end{equation}
Refining the subdivision $[-T_{1},T_{1}]=\bigcup_{i=1}^{N} I_{i}$,
$N=N(u,\mu)$, if necessary, we see that
\begin{equation}\nonumber\begin{split}
&||\sum_{\nu}||\frac{\partial_{\nu}{\bf{\tilde{y}}}}{{\bf{\tilde{y}}}}-\frac{\partial_{\nu}{\bf{y}}}{{\bf{y}}}||_{H^{\mu}(a_{i+1})}
+||\partial_{\nu}(\frac{{\bf{\tilde{x}}}}{{\bf{\tilde{y}}}})-\partial_{\nu}(\frac{{\bf{\tilde{x}}}}{{\bf{\tilde{y}}}})||_{H^{\mu}(a_{i+1})}
\\&\lesssim
2[||\sum_{\nu}||\frac{\partial_{\nu}{\bf{\tilde{y}}}}{{\bf{\tilde{y}}}}-\frac{\partial_{\nu}{\bf{y}}}{{\bf{y}}}||_{H^{\mu}(a_{i})}
+||\partial_{\nu}(\frac{{\bf{\tilde{x}}}}{{\bf{\tilde{y}}}})-\partial_{\nu}(\frac{{\bf{\tilde{x}}}}{{\bf{\tilde{y}}}})||_{H^{\mu}(a_{i})}],\\
\end{split}\end{equation}
provided the quantity on the right hand side is less than some
constant $c$. Thus if we choose $\epsilon<\frac{c}{2^{N}}$, we see
that the Wave Map $\tilde{u}$ satisfying
$||\tilde{u}[0]-u[0]||_{H^{1+\mu}\times H^{\mu}}<<\epsilon$ will
exist and be smooth on the interval $[-T_{1},T_{1}]$. It follows
from the argument just given and a simple algebra type estimate
that by possibly shrinking the size of
$||u[0]-\tilde{u}[0]||_{H^{1+\mu}\times H^{\mu}}$ we can ensure
that
\begin{equation}\nonumber
||\delta\psi_{\nu}||_{L_{t}^{\infty}H^{\lambda}([-T_{1},T_{1}]\times{\mathbf{R}}^{2})}\leq
\epsilon,\,0<\lambda<\mu.
\end{equation}
Now assume that the perturbed Wave Map $\tilde{u}$ breaks down at
some time $T>T_{1}$. We claim that $\sup_{T_{1}\leq
t<T}\sup_{k\in{\mathbf{Z}}}\tilde{c}_{k}^{-1}||P_{k}\delta\psi_{\nu}||_{S[k]([T_{1},t]\times{\mathbf{R}^{2}})}<\infty$.
Indeed, in the opposite case, choosing $K$ large enough(and if
necessary shrinking $\epsilon$), and using the continuity of the
function $t\rightarrow
\sup_{k}\tilde{c}_{k}^{-1}||P_{k}\delta\psi_{\nu}||_{S[k]([T_{1},t]\times{\mathbf{R}^{2}})}$
for $t\in [T_{1},T)$, see e. g. \cite{Kr-4}, it follows that there
exists $T'$ satisfying the properties
\begin{equation}\nonumber
\sup_{k\in{\mathbf{Z}}}\tilde{c}_{k}^{-1}||P_{k}\delta\psi_{\nu}||_{S[k]([T_{1},T']\times{\mathbf{R}}^{2})}=K,\,T_{1}<T'<T
\end{equation}
This, however, contradicts Proposition~\ref{bootstrap}. This then
implies that\\
$\sup_{t<T}||P_{k}\delta\psi_{\nu}||_{S[k]([-t,t]\times{\mathbf{R}}^{2})}\lesssim\tilde{c}_{k}$.
But by definition of $\tilde{c}_{k}$, $||.||_{S[k]}$, this implies
that some subcritical norm
$||\delta\psi_{\nu}||_{L_{t}^{\infty}H^{\mu}([-T,T]\times{\mathbf{R}}^{2})}<\infty$,
$\mu>0$, which in turn implies that some
$||\tilde{u}||_{L_{t}^{\infty}H^{1+\mu'}([-T,T]\times{\mathbf{R}}^{2})}<\infty$.
This in turn contradicts breakdown by the result of
Klainerman-Machedon \cite{Kl-M 1}. Of course, the preceding
argument entails the bound
$||\delta\psi_{\nu}||_{L_{t}^{\infty}L_{x}^{2}}\lesssim \epsilon$.
Indeed, one obtains that some range of subcritical norms
$||.||_{H^{\delta}}$ satisfy that estimate.

\subsection{Proof of Proposition~\ref{bootstrap}}
We first recall the following theorem from \cite{Kr-4}: let
$N(\{\psi_{\nu}\})$ be the nonlinearity on the right hand side of
\eqref{Wave Equation}. Then we have
\begin{theorem}\label{small data}\cite{Kr-4} Let $\psi_{\nu}\in
C_{0}^{\infty}([-T,T]\times{\mathbf{R}}^{2})$ solve the combined
system \eqref{Div-Curl1}, \eqref{Div-Curl2}. Then provided
\begin{equation}\nonumber
\sup_{\nu}||\psi_{\nu}||_{S([-T,T]\times{\mathbf{R}}^{2})}<K\epsilon
\end{equation}
and $\epsilon$ is sufficiently small in relation to $K$, we have
\begin{equation}\nonumber
||N(\{\psi_{\nu}\})||_{N([-T,T]\times{\mathbf{R}}^{2})}\lesssim
K^{3}\epsilon^{3}
\end{equation}
More precisely, there exists a number $\delta_{1}>0$ such that we
have
\begin{equation}\nonumber
||P_{k}N(\{\psi_{\nu}\})||_{N[k]([-T,T]\times{\mathbf{R}}^{2})}\lesssim
(\sup_{\nu}\sup_{k_{1}\in{\mathbf{Z}}}2^{-\delta_{1}|k-k_{1}|}||P_{k_{1}}\psi_{\nu}||_{S[k_{1}]([-T,T]\times{\mathbf{R}}^{2})})K^{2}\epsilon^{2}
\end{equation}
\end{theorem}
The proof in \cite{Kr-4} of this relied on introduction of
null-form structure into the nonlinearities by means of Hodge type
decompositions, as briefly outlined in subsection 3.1. Thus
writing $\psi_{\nu}=R_{\nu}\psi+\chi_{\nu}$ and requiring
$\sum_{i=1,2}\partial_{i}\chi_{i}=0$ results in
\begin{equation}\label{elliptic}
\chi_{\nu} =
i\sum_{i,j=1}^{2}\partial_{i}\triangle^{-1}(\psi_{\nu}\triangle^{-1}\partial_{j}(\psi^{1}_{i}\psi^{2}_{j}-\psi^{1}_{j}\psi^{2}_{i})-\psi_{i}\triangle^{-1}\partial_{j}(\psi^{1}_{\nu}\psi^{2}_{j}-\psi^{1}_{j}\psi^{2}_{\nu})).
\end{equation}
\begin{equation}\label{hyperbolic}
\psi=-\sum_{i=1,2}R_{i}\psi_{i}
\end{equation}
One now writes the nonlinearities $N(\{\psi_{\nu}\})$ as sums of
various terms which are gotten by substituting either gradient
components $R_{\nu}\psi$ or elliptic components $\chi_{\nu}$ in
place of $\psi_{\nu}$, substituting\footnote{One reexpresses
$\psi$, $\chi_{\nu}$ in terms of $\psi_{\nu}$ via
\eqref{elliptic}, \eqref{hyperbolic}.} the Schwartz extensions
$\rho_{\nu}$ for $\psi_{\nu}$ which satisfy
\begin{equation}\nonumber
||\rho_{\nu}||_{S([-T,T]\times{\mathbf{R}}^{2})}\lesssim
K\epsilon,
\end{equation}
and further microlocalizing constituents of the expressions thus
obtained. One thereby obtains trilinear null-forms of the types
recorded in 3.4(c). Substituting elliptic components $\chi_{\nu}$
results in terms at least quintilinear in the variables
$\psi_{\nu}$, which are more elementary to estimate, but still
appear to require null-form structure, which is obtained upon
reiterating the Hodge type decomposition. One keeps going like
this until the error terms obtained can be estimated without using
null-structures, based only on Strichartz type estimates.
Summarizing, we have
\begin{theorem}\label{decomposition}\cite{Kr-4} Under the hypotheses of
Theorem~\ref{small data}, we can construct a function
\begin{equation}\nonumber
\tilde{N}(\{\rho_{\nu}\})\in{\mathcal{S}}({\mathbf{R}}^{2+1})
\end{equation}
which is expressible as a sum of terms trilinear, quadrilinear
etc. up to degree 11 in the $\rho_{\nu}$, and satisfies
\begin{equation}\nonumber
\tilde{N}(\{\rho_{\nu}\})|_{[-T,T]}=N(\psi_{\nu})
\end{equation}
\begin{equation}\nonumber
||P_{k}\tilde{N}(\{\rho_{\nu}\})||_{S[k]}\lesssim
(\sup_{\nu}\sup_{k_{1}\in{\mathbf{Z}}}2^{-\delta|k-k_{1}|}||P_{k_{1}}\rho_{\nu}||_{S[k_{1}]})\max\{||\rho_{\nu}||_{S}^{2},\,||\rho_{\nu}||_{S}^{10}\}
\end{equation}
\end{theorem}
We shall apply this theorem to our situation. The complication
that arises here has to do with the fact that the estimates for
$\{\psi_{\nu}\}$, the derivative components of the spherically
symmetric Wave Map, are not with respect to $||.||_{S}$, but
rather $||.||_{{\mathcal{S}}}$, in view of theorem~\ref{Moser}. We
state here
\begin{lemma}\label{auxiliary envelope} For any
$\sigma>0$ sufficiently small, there exists a frequency envelope
$\{c_{l}\}_{l\in{\mathbf{Z}}}$ with exponent $\sigma$ and
$\sum_{l\in{\mathbf{Z}}}c_{l}^{2}\lesssim 1$ such that $\forall
T>0$ we have
\begin{equation}\nonumber
||P_{k}\psi_{\nu}||_{{\mathcal{S}}[k]([-T,T]\times{\mathbf{R}}^{2})}\lesssim
c_{k}
\end{equation}
We can also assume $\sum_{l}2^{\mu |l|}c_{l}\lesssim 1$ for
$\mu>0$ sufficiently small. Moreover, choosing a Schwartz
extension $\widetilde{P_{k}\psi_{\nu}}$ of
$P_{k}\psi_{\nu}|_{[-T,T]}$ satisfying the above estimates, we may
decompose each $\widetilde{P_{k}\psi_{\nu}}$ into functions of
first and 2nd type,
$\widetilde{P_{k}\psi_{\nu}}=\alpha_{\nu}+\beta_{\nu}$, such that
the following properties hold:
\begin{equation}\nonumber
||\alpha_{\nu}||_{A[k]}\lesssim
c_{k},\,||\beta_{\nu}||_{B[k]}\lesssim
c_{k},\,||\beta_{\nu}||_{\dot{X}_{k}^{0,\frac{1}{2},1}}\leq
C,\,||R_{0}(\alpha_{\nu}+\beta_{\nu})||_{L_{t}^{\infty}L_{x}^{2}}\leq
C
\end{equation}
where $C$ depends on
$||\frac{\bf{x}}{{\bf{y}}}||_{\dot{B}_{2}^{1,1}}+||\ln{\bf{y}}||_{\dot{B}_{2}^{1,1}}$.
Moreover, the bilinear inequalities enunciated in
theorem~\ref{Moser} hold for $\psi_{\nu}$ in place of
$[\psi_{1}A(\nabla^{-1}\psi_{2})]$ there.
\end{lemma}
\begin{proof} We define
\begin{equation}\nonumber\begin{split}
&c_{k}:=\sum_{k_{1}\in{\mathbf{Z}}}2^{-\sigma|k-k_{1}|}[||P_{k_{1}}N(...)||_{L_{t}^{1}L_{x}^{2}({\mathbf{R}}^{1+2})}+||P_{k_{1}}N(...)||_{L_{t}^{2}\dot{H}^{-\frac{1}{2}}}]
\\&+\sum_{k_{1}\in{\mathbf{Z}}}2^{-\sigma|k-k_{1}|}||P_{k_{1}}\frac{\partial_{\nu}{\bf{y}}}{{\bf{y}}}(0)||_{L_{x}^{2}}
+\sup_{\nu=0,1,2}\sum_{k_{1}\in{\mathbf{Z}}}2^{-\sigma|k-k_{1}|}||P_{k_{1}}\partial_{\nu}(\frac{{\bf{x}}}{{\bf{y}}})(0)||_{L_{x}^{2}},\\
\end{split}\end{equation}
where
$N(...)=N(\nabla{\bf{x}},\,\nabla{\bf{y}},\,{\bf{x}},\,{\bf{y}})$
runs over the nonlinearities in \eqref{Euler1}, \eqref{Euler2}.
That this is indeed a frequency envelope with the desired
properties follows from Corollary~\ref{crux} as well as
lemma~\ref{technical25}. We need to exercise some care to get good
enough control over the elliptic portions of $\psi_{\nu}$. For
this, truncate
$N(\nabla{\bf{x}},\,\nabla{\bf{y}},\,{\bf{x}},\,{\bf{y}})$ past
some time $T_{0}>>\max\{2^{-k}, T\}$, and (committing abuse of
notation) decompose the nonlinearity
\begin{equation}\nonumber
P_{k}N(\nabla{\bf{x}},\,\nabla{\bf{y}},\,{\bf{x}},\,{\bf{y}})=P_{k}Q_{<k}N(\nabla{\bf{x}},\,\nabla{\bf{y}},\,{\bf{x}},\,{\bf{y}})
+P_{k}Q_{\geq
k}N(\nabla{\bf{x}},\,\nabla{\bf{y}},\,{\bf{x}},\,{\bf{y}}).
\end{equation}
Then consider
\begin{equation}\nonumber
\Box^{-1}P_{k}Q_{\geq
k}N(\nabla{\bf{x}},\,\nabla{\bf{y}},\,{\bf{x}},\,{\bf{y}}),
\end{equation}
where the operator $\Box^{-1}$ is division by the symbol
$(\tau^{2}-|\xi|^{2})$ on the space-time Fourier side. Clearly,
from definition we have
\begin{equation}\nonumber
||\Box^{-1}P_{k}Q_{\geq
k}\nabla_{x,t}N(\nabla{\bf{x}},\,\nabla{\bf{y}},\,{\bf{x}},\,{\bf{y}})||_{L_{t}^{1}L_{x}^{2}}\lesssim
2^{-k}c_{k}
\end{equation}
Thus there exists a time $t_{0}<T_{0}$ with the property
\begin{equation}\nonumber
||\Box^{-1}P_{k}Q_{\geq
k}\nabla_{x,t}N(\nabla{\bf{x}},\,\nabla{\bf{y}},\,{\bf{x}},\,{\bf{y}})(t_{0})||_{L_{x}^{2}}\lesssim
c_{k}
\end{equation}
We easily check that (for $C$ independent of $k$)
\begin{equation}\nonumber
||\Box^{-1}R_{0}P_{k}Q_{\geq
k}\nabla_{x,t}N(\nabla{\bf{x}},\,\nabla{\bf{y}},\,{\bf{x}},\,{\bf{y}})||_{L_{t}^{\infty}L_{x}^{2}}\leq
C,
\end{equation}
while also (using the wave equation)
\begin{equation}\nonumber
||R_{0}\nabla_{x,t}P_{k}(\frac{{\bf{x}}}{{\bf{y}}})||_{L_{t}^{\infty}L_{x}^{2}}+||R_{0}\nabla_{x,t}P_{k}\ln{\bf{y}}||_{L_{t}^{\infty}L_{x}^{2}}
\leq C
\end{equation}
Now construct a free wave $a$ with the properties
\begin{equation}\nonumber
a(t_{0})=P_{k}(\frac{{\bf{x}}}{{\bf{y}}})(t_{0})-\Box^{-1}P_{k}Q_{\geq
k}N(\nabla{\bf{x}},\,\nabla{\bf{y}},\,{\bf{x}},\,{\bf{y}})(t_{0})
\end{equation}
\begin{equation}\nonumber
\partial_{t}a(t_{0})=P_{k}\partial_{t}(\frac{{\bf{x}}}{{\bf{y}}})(t_{0})-\Box^{-1}\partial_{t}P_{k}Q_{\geq
k}N(\nabla{\bf{x}},\,\nabla{\bf{y}},\,{\bf{x}},\,{\bf{y}})(t_{0})
\end{equation}
and similarly for $\ln{\bf{y}}$. It follows that the quantity
$\frac{{\bf{x}}}{{\bf{y}}}-a-\Box^{-1}P_{k}Q_{\geq k}N(...)$
satisfies
\begin{equation}\nonumber\begin{split}
&\Box(\frac{{\bf{x}}}{{\bf{y}}}-a-\Box^{-1}P_{k}Q_{\geq
k}N(...))=P_{k}Q_{<k}N(\nabla{\bf{x}},\,\nabla{\bf{y}},\,{\bf{x}},\,{\bf{y}}),\\&||R_{0}\nabla_{x,t}P_{k}(\frac{{\bf{x}}}{{\bf{y}}}-a-\Box^{-1}P_{k}Q_{\geq
k}N(...))||_{L_{t}^{\infty}L_{x}^{2}}\leq C\\
\end{split}\end{equation}
as well as
$||\nabla_{x,t}(\frac{{\bf{x}}}{{\bf{y}}}-a-\Box^{-1}P_{k}Q_{\geq
k}N(...))(t_{0})||_{L_{x}^{2}}\lesssim c_{k}$. One also verifies
that
\begin{equation}\nonumber
||\Box^{-1}\nabla_{x,t}P_{k}Q_{\geq
k}N(\nabla{\bf{x}},\,\nabla{\bf{y}},\,{\bf{x}},\,{\bf{y}})||_{\dot{X}_{k}^{-\frac{1}{2},1,2}}\lesssim
c_{k},
\end{equation}
which by Sobolev's inequality also implies control over
$||.||_{L}$ as well as the Strichartz type norms
$||.||_{L_{t}^{p}L_{x}^{q}}$,
$\frac{1}{p}+\frac{1}{q}<\frac{1}{2}-\delta_{0}$, of this
expression. Now one solves the wave equation for
$\frac{{\bf{x}}}{{\bf{y}}}-a-\Box^{-1}P_{k}Q_{\geq k}N(...)$ with
initial data given at time $t_{0}$. Using 3.4(d) (which in turn
relies on a truncated Duhamel's formula, see \cite{Kr-4}), one
constructs Schwartz extensions
$\widetilde{P_{k}(\frac{\bf{x}}{{\bf{y}}})}$,
$\widetilde{P_{k}\ln{\bf{y}}}$ of
$P_{k}(\frac{{\bf{x}}}{{\bf{y}}})|_{[-T,T]}$,
$P_{k}\ln{\bf{y}}|_{[-T,T]}$, respectively, with the properties
\begin{equation}\nonumber
||\nabla_{x,t}\widetilde{P_{k}(\frac{\bf{x}}{{\bf{y}}})}||_{A[k]}\lesssim
c_{k},\,||\nabla_{x,t}\widetilde{P_{k}\ln{\bf{y}}}||_{A[k]}\lesssim
c_{k}
\end{equation}
as well as
\begin{equation}\nonumber
||R_{0}\nabla_{x,t}[\widetilde{P_{k}(\frac{\bf{x}}{{\bf{y}}})}]||_{L_{t}^{\infty}L_{x}^{2}}\leq
C,
||R_{0}\nabla_{x,t}\widetilde{P_{k}\ln{\bf{y}}}||_{L_{t}^{\infty}L_{x}^{2}}\leq
C
\end{equation}
where $C$ is independent of $k$. Using a partition of unity, one
glues these extensions together to get Schwartz extensions
$\widetilde{\frac{{\bf{x}}}{{\bf{y}}}}$, $\widetilde{\ln{\bf{y}}}$
of $\frac{{\bf{x}}}{{\bf{y}}}|_{[-T,T]}$, $\ln{\bf{y}}|_{[-T,T]}$
which satisfy
\begin{equation}\nonumber
||P_{k}(\widetilde{\frac{{\bf{x}}}{{\bf{y}}}})||_{A[k]}\lesssim
c_{k}
\end{equation}
etc. Now one recalls that
\begin{equation}\nonumber
\psi_{\nu}=(\frac{\partial_{\nu}{\bf{x}}}{{\bf{y}}}+i\frac{\partial_{\nu}{\bf{y}}}{{\bf{y}}})e^{i\sum_{j=1,2}\triangle^{-1}\partial_{j}
(\frac{\partial_{j}{\bf{x}}}{{\bf{y}}})},
\end{equation}
plugs in the Schwartz extensions of $\frac{\bf{x}}{{\bf{y}}}$ etc.
and uses theorem~\ref{Moser} to obtain the desired conclusion.
\end{proof}

Continuing with the proof of Proposition~\ref{bootstrap}, our
strategy now will be to analyze the wave equation satisfied by
$\delta\psi_{\nu}=\tilde{\psi}_{\nu}-\psi_{\nu}$. Using
\eqref{Wave Equation} for both $\psi_{\nu}$, $\tilde{\psi}_{\nu}$,
and subtracting, we obtain a first version. We eliminate
$\tilde{\psi}_{\nu}$ by substituting
$\delta\psi_{\nu}+\psi_{\nu}$. One thereby obtains a sum of
products of components $\delta\psi_{\nu}$, $\psi_{\nu}$ which are
at least linear in $\delta\psi_{\nu}$. Proceeding as in the
previous description, we decompose the $\delta\psi_{\nu}$,
$\psi_{\nu}$ into gradient and elliptic parts. For the
$\delta\psi_{\nu}$, this is obtained by applying the procedure to
$\tilde{\psi}_{\nu}$, $\psi_{\nu}$ and forming the difference,
resulting in
\begin{equation}\nonumber
\delta\psi_{\nu}=R_{\nu}(\delta\psi)+\delta\chi_{\nu},\,\delta\psi=-\sum_{i=1,2}R_{i}\psi_{i}
\end{equation}
\begin{equation}\nonumber\begin{split}
&\delta\chi_{\nu}=
i\sum_{i,j=1}^{2}\partial_{i}\triangle^{-1}(\tilde{\psi}_{\nu}\triangle^{-1}\partial_{j}(\tilde{\psi}^{1}_{i}\tilde{\psi}^{2}_{j}-\tilde{\psi}^{1}_{j}\tilde{\psi}^{2}_{i})-\tilde{\psi}_{i}\triangle^{-1}\partial_{j}(\tilde{\psi}^{1}_{\nu}\tilde{\psi}^{2}_{j}-\tilde{\psi}^{1}_{j}\tilde{\psi}^{2}_{\nu}))\\
&\hspace{0.5cm}-
i\sum_{i,j=1}^{2}\partial_{i}\triangle^{-1}(\psi_{\nu}\triangle^{-1}\partial_{j}(\psi^{1}_{i}\psi^{2}_{j}-\psi^{1}_{j}\psi^{2}_{i})-\psi_{i}\triangle^{-1}\partial_{j}(\psi^{1}_{\nu}\psi^{2}_{j}-\psi^{1}_{j}\psi^{2}_{\nu})).\\
\end{split}\end{equation}
Clearly one can reexpress the latter difference as a sum of terms
linear, quadratic and cubic in the $\delta\psi_{\nu}$, eliminating
the $\tilde{\psi}_{\nu}=\delta\psi_{\nu}+\psi_{\nu}$. In order to
demonstrate Proposition~\ref{bootstrap}, we shall rely on the
following refined
\begin{proposition}\label{refined bootstrap} Let
\begin{equation}\nonumber
\Box\delta\psi_{\alpha}=N_{\alpha}(\delta\psi_{\nu},\,\psi_{\nu})
\end{equation}
on $[-T,T]$. Proceeding as above, express the nonlinearity as a
sum of trilinear null-forms (substituting the gradient components
for $\delta\psi_{\nu}$, $\psi_{\nu}$), as well as error terms 'at
least quintilinear' in $\delta\psi_{\nu}$, $\psi_{\nu}$ (which
arise upon substituting $\delta\chi_{\nu}$, $\chi_{\nu}$). Denote
the sum of terms which are linear in $\delta\psi_{\nu}$ by
$N_{1\alpha}(\delta\psi_{\nu},\,\psi_{\nu})$. Then for any
$\zeta>0$ there exists $T_{0}>0$, such that for any fixed smooth
function $\chi(t)\in C^{\infty}({\mathbf{R}})$ with
$\text{supp}(\chi)\subset [-1,1]^{c}$, $\chi|_{[-2,2]^{c}}=1$, we
have
\begin{equation}\nonumber
||P_{k}(\delta\psi_{\nu})||_{S[k]([-T,T]\times{\mathbf{R}}^{2})}\leq
K\tilde{c}_{k},\,T>\tilde{T}_{0}\geq T_{0}
\end{equation}
\begin{equation}\nonumber
\Rightarrow||\chi(\frac{t}{\tilde{T}_{0}})P_{k}N_{1\alpha}(\delta\psi_{\nu},\,\psi_{\nu})||_{N[k]([-T,T]\times{\mathbf{R}}^{2})}\lesssim
\zeta K\tilde{c}_{k}
\end{equation}
Here $\{\tilde{c}_{k}\}$ is associated with $\tilde{T}_{0}$ as in
Proposition~\ref{bootstrap} (substitute $\tilde{T}_{0}$ for
$T_{1}$). Moreover, denoting the terms at least quadratic in
$\delta\psi_{\nu}$ by $N_{2}(\delta\psi_{\nu},\psi_{\nu})$, and
letting $\epsilon,\tilde{c}_{k}$ be as in the statement of
Proposition~\ref{bootstrap}, the following conclusion holds
provided $\epsilon$ is small enough and $\{\tilde{c}_{k}\}$ 'flat
enough':
\begin{equation}\nonumber
||P_{k}(\delta\psi_{\nu})||_{S[k]([-T,T]\times{\mathbf{R}}^{2})}\leq
K\tilde{c}_{k}\Rightarrow
||P_{k}N_{2\alpha}(\delta\psi_{\nu},\psi_{\nu})||_{N[k]([-T,T]\times{\mathbf{R}}^{2+1})}\lesssim
\epsilon K^{2}c_{k}.
\end{equation}
\end{proposition}

Deferring the proof of this for the moment, we continue with the
proof of Proposition~\ref{bootstrap}. Let $\zeta <1/C$ for some
$C>>1$, and construct $T_{0}$ as in Proposition~\ref{refined
bootstrap}; Define $T_{1}:=2T_{0}$. Now assume we have the
situation in the statement of Proposition~\ref{bootstrap}. We
intend to use the energy inequality 3.4(d). Fix
$k\in{\mathbf{Z}}$, and consider $P_{k}\delta\psi_{\nu}$. We
distinguish between the cases $T-T_{1}<\frac{2^{-k}}{C}$ and the
opposite. In the former case, the wave equation becomes useless,
and we use the divergence-curl system directly: observe that by
virtue of \eqref{Div-Curl1} we have for $i=1,2$, $T_{1}<t\leq T$
\begin{equation}\nonumber
P_{k}\delta\psi_{i}(t)-P_{k}\delta\psi_{i}(T_{1})=\int_{T_{1}}^{t}\partial_{i}P_{k}\delta\psi_{i}dt
+\int_{T_{1}}^{t}P_{k}[N(\psi,\delta\psi)]dt
\end{equation}
In this equation, by abuse of notation, $N(\psi,\delta\psi)$ is a
linear combination of terms of the schematic form
$\delta\psi\nabla^{-1}(\psi^{2})$,
$\delta\psi\nabla^{-1}(\psi\delta\psi)$ etc. Let's put
$N(\psi,\delta\psi)=\delta\psi\nabla^{-1}(\psi^{2})$, the other
terms being treated along the same lines (but also requiring
$\epsilon$ to be small enough). We note that
\begin{equation}\nonumber
||P_{k}[\delta\psi\nabla^{-1}(\psi^{2})]||_{L_{t}^{B}L_{x}^{2}}\lesssim
K 2^{(1-\frac{1}{B})k}\tilde{c}_{k},
\end{equation}
where $B$ is an arbitrarily large number (the implied constants
will depend on it). This follows from a simple frequency
trichotomy and the bootstrap assumption. Now using Holder's
inequality, we deduce
\begin{equation}\nonumber
||\int_{T_{1}}^{t}P_{k}[N(\psi,\delta\psi)]dt||_{L_{x}^{2}}\lesssim
|t-T_{1}|^{1-\frac{1}{B}}2^{(1-\frac{1}{B})k}K\tilde{c}_{k}\lesssim
\frac{1}{C^{1-\frac{1}{B}}}K\tilde{c}_{k}
\end{equation}
Clearly we also have
\begin{equation}\nonumber
||\int_{T_{1}}^{t}\partial_{i}P_{k}\delta\psi_{i}dt||_{L_{x}^{2}}\lesssim
|t-T_{1}|2^{k} K\tilde{c}_{k}\lesssim \frac{1}{C}K\tilde{c}_{k}.
\end{equation}
Therefore, we infer that
\begin{equation}\nonumber
||P_{k}\delta\psi_{i}||_{L_{t}^{\infty}L_{x}^{2}}\lesssim
(1+\frac{K}{C^{1-\frac{1}{B}}})\tilde{c}_{k}\lesssim
\frac{K}{100}\tilde{c}_{k}
\end{equation}
provided $K,\,C$ are chosen\footnote{Of course $C$ is chosen
independently of $K$.} large enough (in relation to $B$). Arguing
similarly, one deduces as well that
\begin{equation}\nonumber
2^{\frac{k}{2}}||P_{k}\delta\psi_{i}||_{L_{t}^{2}L_{x}^{2}}
+2^{-\frac{k}{2}}||P_{k}\partial_{t}\delta\psi_{i}||_{L_{t}^{2}L_{x}^{2}}\lesssim
\frac{K}{C^{\frac{1}{2}}}\tilde{c}_{k}.
\end{equation}
Using the fact that (see e. g. \cite{Kr-4})
\begin{equation}\nonumber\begin{split}
&||P_{k}\psi||_{S[k]([-T,T]\times{{\mathbf{R}}^{2}})}\lesssim
||P_{k}\psi||_{L_{t}^{\infty}L_{x}^{2}([-T,T]\times{{\mathbf{R}}^{2}})}+2^{\frac{k}{2}}||P_{k}\psi||_{L_{t}^{2}L_{x}^{2}([-T,T]\times{{\mathbf{R}}^{2}})}\\&\hspace{7cm}+2^{-\frac{k}{2}}||P_{k}\partial_{t}\psi||_{L_{t}^{2}L_{x}^{2}([-T,T]\times{{\mathbf{R}}^{2}})}
\\
\end{split}\end{equation}
and choosing $K,\,C$ large enough, one deduces from this that
\begin{equation}\nonumber
||P_{k}\delta\psi_{i}||_{S[k]([T_{1},T]\times{\mathbf{R}}^{2})}\lesssim
\frac{K}{2}\tilde{c}_{k},
\end{equation}
which is the desired conclusion for $P_{k}\delta\psi_{i}$. The
argument for $\delta\psi_{0}$ is similar using \eqref{Div-Curl2}.
Thus we see that we may assume $|T-T_{1}|\geq \frac{2^{-k}}{C}$.
Moreover, reiterating the preceding argument, and choosing $K$
large enough, we conclude that
$||P_{k}\delta\psi_{\nu}||_{L_{t}^{\infty}L_{x}^{2}([T_{1},T_{1}+\frac{2^{-k}}{C}]\times{\mathbf{R}}^{2})}\leq
\frac{K}{100C}\tilde{c}_{k}$. Now revert to the old notation
\begin{equation}\nonumber
\Box\delta\psi_{\nu}=N(\delta\psi_{\nu},\psi_{\nu})=N_{1}(\delta\psi_{\nu})+N_{2}(\delta\psi_{\nu},\psi_{\nu})
\end{equation}
as in Proposition~\ref{refined bootstrap}. Clearly, we have
\begin{equation}\nonumber
||P_{k}N_{1}(\delta\psi_{\nu},\psi_{\nu})||_{N[k]([T_{1},T]\times{\mathbf{R}}^{2})}\leq
||P_{k}\chi(\frac{t}{T_{0}})N_{1}(\delta\psi_{\nu},\psi_{\nu})||_{N[k]([T_{1},T]\times{\mathbf{R}}^{2})}
\end{equation}
\begin{equation}\nonumber
||P_{k}N_{2}(\delta\psi_{\nu},\psi_{\nu})||_{N[k]([T_{1},T]\times{\mathbf{R}}^{2})}\lesssim
||P_{k}N_{2}(\delta\psi_{\nu},\psi_{\nu})||_{N[k]([-T,T]\times{\mathbf{R}}^{2})}
\end{equation}
Using 3.4(d) as well as time translation invariance, we can now
infer that
\begin{equation}\nonumber
||P_{k}\delta\psi_{\nu}||_{S[k]([T_{1},T]\times{\mathbf{R}}^{2})}\lesssim
\frac{K}{100}\tilde{c}_{k}+\zeta
K\tilde{c}_{k}+\epsilon^{2}K^{2}\tilde{c}_{k}^{2}\lesssim
\frac{K}{2}\tilde{c}_{k}
\end{equation}
provided $\epsilon$, $\zeta$ are small enough. This yields the
desired conclusion.
\subsection{The proof of theorem~\ref{Main'}}
This is basically identical to the proof of theorem~\ref{Main}.
Control over some subcritical norm
$||u||_{L_{t}^{\infty}H^{1+\epsilon}}$ follows from standard Moser
estimates instead of Corollary~\ref{crux}.
\end{section}

\begin{section}{The proof of Proposition 3.17.}\markboth{Joachim
 Krieger}{The proof of Proposition 3.17}

We have thus reduced the proof of theorem~\ref{Main} to the
verification of Proposition~\ref{refined bootstrap} in addition to
the technical Moser type estimates allowing estimation of
$||P_{k}\psi||_{{\mathcal{S}}[k]}$. The proof of this Proposition
is divided into the part dealing with expressions linear in
$\delta\psi_{\nu}$, as well as those of higher degree of
linearity. We commence by spelling out in detail the decomposition
$N_{\alpha}(\delta\psi_{\nu},\psi_{\nu})=N_{1\alpha}(\delta\psi_{\nu},\psi_{\nu})+N_{2\alpha}(\delta\psi_{\nu},\psi_{\nu})$,
where
$\Box\delta\psi_{\alpha}=N_{\alpha}(\delta\psi_{\nu},\psi_{\nu})$.
As in \cite{Kr-4}, this decomposition requires extreme care in
order to avoid too many time derivatives. Recalling \eqref{Wave
Equation}, we define
\begin{equation}\nonumber\begin{split}
N_{1}(\delta\psi_{\nu},\psi_{\nu}):=\sum_{i=1}^{3}A_{i\alpha}(\delta\psi_{\nu},\psi_{\nu})+\sum_{i=1}^{5}B_{i\alpha}(\delta\psi_{\nu},\psi_{\nu})+\sum_{i=1}^{5}C_{i\alpha}\\
\end{split}\end{equation}
where
\begin{equation}\begin{split}
A_{1}(\delta\psi_{\nu},\psi_{\nu})=
&i\partial^{\beta}[\delta\psi_{\alpha}\triangle^{-1}\sum_{j=1}^{2}\partial_{j}[R_{\beta}\psi^{1}R_{j}\psi^{2}-R_{\beta}\psi^{2}R_{j}\psi^{1}]]\\
&-i\partial^{\beta}[\delta\psi_{\beta}\triangle^{-1}\sum_{j=1}^{2}\partial_{j}[R_{\alpha}\psi^{1}R_{j}\psi^{2}-R_{\alpha}\psi^{2}R_{j}\psi^{1}]]\\
&+i\partial_{\alpha}[\delta\psi_{\nu}\triangle^{-1}\sum_{j=1}^{2}\partial_{j}[R^{\nu}\psi^{1}R_{j}\psi^{2}-R^{\nu}\psi^{2}R_{j}\psi^{1}]]\\
\end{split}\end{equation}
\begin{equation}\nonumber\begin{split}
A_{1}(\delta\psi_{\nu},\psi_{\nu})=
&+i\partial^{\beta}[\psi_{\alpha}\triangle^{-1}\sum_{j=1}^{2}\partial_{j}[R_{\beta}\delta\psi^{1}R_{j}\psi^{2}-R_{\beta}\psi^{2}R_{j}\delta\psi^{1}]]\\
&-i\partial^{\beta}[\psi_{\beta}\triangle^{-1}\sum_{j=1}^{2}\partial_{j}[R_{\alpha}\delta\psi^{1}R_{j}\psi^{2}-R_{\alpha}\psi^{2}R_{j}\delta\psi^{1}]]\\
&+i\partial_{\alpha}[\psi_{\nu}\triangle^{-1}\sum_{j=1}^{2}\partial_{j}[R^{\nu}\delta\psi^{1}R_{j}\psi^{2}-R^{\nu}\psi^{2}R_{j}\delta\psi^{1}]]\\
\end{split}\end{equation}
\begin{equation}\nonumber\begin{split}
A_{3}(\delta\psi_{\nu},\psi_{\nu})=
&+i\partial^{\beta}[\psi_{\alpha}\triangle^{-1}\sum_{j=1}^{2}\partial_{j}[R_{\beta}\psi^{1}R_{j}\delta\psi^{2}-R_{\beta}\delta\psi^{2}R_{j}\psi^{1}]]\\
&-i\partial^{\beta}[\psi_{\beta}\triangle^{-1}\sum_{j=1}^{2}\partial_{j}[R_{\alpha}\psi^{1}R_{j}\delta\psi^{2}-R_{\alpha}\delta\psi^{2}R_{j}\psi^{1}]]\\
&+i\partial_{\alpha}[\psi_{\nu}\triangle^{-1}\sum_{j=1}^{2}\partial_{j}[R^{\nu}\psi^{1}R_{j}\delta\psi^{2}-R^{\nu}\delta\psi^{2}R_{j}\psi^{1}]].\\
\end{split}\end{equation}
\begin{equation}\nonumber
B_{1},C_{1}(\delta\psi_{\nu},\psi_{\nu})=\nabla_{x,t}[\delta\psi\nabla^{-1}[\psi\nabla^{-1}[\psi\nabla^{-1}(\psi^{2})]]]
\end{equation}
\begin{equation}\nonumber
B_{2},C_{2}(\delta\psi_{\nu},\psi_{\nu})=\nabla_{x,t}[\psi\nabla^{-1}[\delta\psi\nabla^{-1}[\psi\nabla^{-1}(\psi^{2})]]]
\end{equation}
\begin{equation}\nonumber
B_{3},C_{3}(\delta\psi_{\nu},\psi_{\nu})=\nabla_{x,t}[\psi\nabla^{-1}[\psi\nabla^{-1}[\delta\psi\nabla^{-1}(\psi^{2})]]]
\end{equation}
\begin{equation}\nonumber
B_{4,5},C_{4,5}(\delta\psi_{\nu},\psi_{\nu})=\nabla_{x,t}[\psi\nabla^{-1}[\psi\nabla^{-1}[\psi\nabla^{-1}(\delta\psi\psi)]]]
\end{equation}
Of course we have used schematic notation for the $B,\,C$'s, as
their fine structure won't matter. They are obtained by
substituting one $\chi_{\nu}$ instead of the corresponding entry
$\psi_{\nu}$ in the inner square bracket expressions on the right
hand side of \eqref{Wave Equation}, where $\chi_{\nu}$ is the
'elliptic component' of the spherically symmetric $\psi_{\nu}$ in
the decomposition $\psi_{\nu}=R_{\nu}\psi+\chi_{\nu}$. We recall
$\psi_{\nu}=\psi^{1}_{\nu}+i\psi^{2}_{\nu}$,
$\delta\psi_{\nu}=\delta\psi^{1}_{\nu}+i\delta\psi^{2}_{\nu}$,
$\psi=-\sum_{j=1,2}R_{j}\psi_{j}$,
$\delta\psi=-\sum_{j=1,2}R_{j}\delta\psi_{j}$. One can then define
$N_{2\alpha}(\delta\psi_{\nu},\psi_{\nu})=N_{\alpha}(\delta\psi_{\nu},\psi_{\nu})-N_{1\alpha}(\delta\psi_{\nu},\psi_{\nu})$.
The quintilinear terms above shall be relatively simple to
estimate on account of the strong Strichartz type estimates
satisfied by the $\psi_{\nu}$, see theorem~\ref{spherical
Strichartz} as well as the definition of
$||.||_{{\mathcal{S}}[k]}$. Unfortunately, the latter norm falls
short of controlling $||.||_{L_{t}^{2}L_{x}^{\infty}}$, which
appears necessary in order to grant an elementary estimation of
the trilinear terms $A_{i\alpha}$. We shall instead have to revert
to the inherent null-structure in these terms as was done already
in \cite{Kr-4}, in addition to the more complicated ingredients in
$||.||_{{\mathcal{S}}}$. The main new difficulty over the
estimates in \cite{Kr-4} has to do with the fact that we need to
gain explicitly in time in these estimates. This would be
relatively straightforward if we were working with Lebesgue type
spaces; however, we shall work with null-frame spaces of type
$L_{t_{\omega}}^{2}L_{x_{\omega}}^{\infty}$, which considerably
complicates obtaining gains in time. The main novelty
here(lemma~\ref{key}) shall be a special type of decomposition of
the spherical components $\psi_{\nu}$ into pieces which have
well-defined physical as well as frequency localization
properties. More precisely, we shall be able to physically
localize $\psi_{\nu}$ closely to the light cone. This part will
then be written as a sum of two components, the first of which can
be written as a sum of pieces which propagate in a direction
essentially opposite to their physical support. Thus the first
component is obtained by first localizing $\psi_{\nu}$ to an
angular sector in Fourier space, then multiplying with a physical
cutoff localizing to an {\it{opposite or identical\footnote{This
depends on whether the space-time Fourier support is contained in
the upper half-space $\tau>0$ or lower half-space $\tau<0$.}}}
angular sector, and finally summing over all sectors. The size of
the angular sectors shall essentially be dictated by the $\zeta$
in the statement of Proposition~\ref{refined bootstrap}. While the
first component is exactly the part which fails to decay in
$L_{t_{\omega}}^{2}L_{x_{\omega}}^{\infty}$ as
$t\rightarrow\infty$, it does lead to improved trilinear null-form
estimates due to the dual localization properties. The 2nd
component in turn will decay like a standard Lebesgue norm as
$t\rightarrow \infty$. The next subsection contains the core
estimates. As the estimates are rather technical, we briefly
explain the strategy of the proof, which is conceptually simple:
\\

{\bf{(1)}}: First, upon localizing the nonlinearity to a time
interval $t\sim 2^{i}$, one tries to reduce the frequencies of all
functions occuring inside the nonlinearity to absolute size
$<<2^{\delta i}$, for some small $\delta>0$. The idea here is that
far apart frequencies should interact little. But this in addition
to the refined control over the frequency modes of the spherically
symmetric components should suffice to get control over the cases
when extremely small or large frequencies are present. The tool to
achieve this are the refined trilinear estimates in 3.4(c).
Unfortunately, these estimates aren't quite good enough to control
certain high-high interactions, which accounts for a number of
extra cases that need to be considered.
\\

{\bf{(2)}} Having controlled the cases when the frequencies are
very small or large in relation to the time interval one works on,
one now tries to exploit the pointwise estimates provided by
Christodoulou-Tahvildar-Zadeh, since one has gained some room to
lose in the frequencies. The device here is the decomposition of
the spherically symmetric components referred to in the preceding
paragraph, which is a direct consequence of the pointwise decay
estimates. This allows one to decompose these components into
pieces that disperse quickly enough, as well as other pieces that
interact very weakly. Of course one exploits the trilinear
structure of the nonlinearity to make this work.

\subsection{Estimating the trilinear null-forms.}

We use the operator\\
$I=\sum_{k\in{\mathbf{Z}}}P_{k}Q_{<k+100}$ as before and employ
the schematic decomposition
\begin{equation}\nonumber
\nabla_{x,t}[\psi_{1}\nabla^{-1}[\psi_{2}\psi_{3}]]=\nabla_{x,t}[\psi_{1}I\nabla^{-1}[\psi_{2}\psi_{3}]]+\nabla_{x,t}[\psi_{1}(1-I)\nabla^{-1}[\psi_{2}\psi_{3}]]
\end{equation}
for each of the $A_{i\alpha}$'s. In order to make sense of this,
one needs to substitute Schwartz extensions for the inputs
$\delta\psi^{1,2}|_{[T_{1},T]},\,\psi^{1,2}|_{[T_{1},T]}$ of the
inner square brackets, in accordance with the bootstrap assumption
in Proposition~\ref{refined bootstrap}. In the following we shall
localize the frequency localized nonlinearities $P_{k}N(...)$ to a
dyadic time interval $t\sim 2^{i}$ and strive for an estimate of
the form
\begin{equation}\nonumber
||\chi_{i}(t)P_{k}N(...)||_{N[k]}\lesssim 2^{-\mu i}\tilde{c}_{k}
\end{equation}
One can then sum over $i$ large enough to obtain the estimate in
Proposition~\ref{bootstrap}.
\\

{\bf{(A): The large modulation case.}} Estimating the terms
\begin{equation}\nonumber
(I):\,\nabla_{x,t}[\delta\psi\triangle^{-1}\sum_{j=1,2}\partial_{j}(1-I)[R_{\beta}\psi_{2}R_{j}\psi_{3}-R_{j}\psi_{2}R_{\beta}\psi_{3}]]
\end{equation}
\begin{equation}\nonumber
(II):\,\nabla_{x,t}[\psi\triangle^{-1}\sum_{j=1,2}\partial_{j}(1-I)[R_{\beta}\delta\psi_{2}R_{j}\psi_{3}-R_{j}\delta\psi_{2}R_{\beta}\psi_{3}]]
\end{equation}
(I): {\it{The first term.}} We use the decomposition
\begin{equation}\nonumber\begin{split}
&\chi(\frac{t}{\tilde{T}_{0}})\nabla_{x,t}[\delta\psi\triangle^{-1}\sum_{j=1,2}\partial_{j}(1-I)[R_{\beta}\psi_{2}R_{j}\psi_{3}-R_{j}\psi_{2}R_{\beta}\psi_{3}]]\\
&=\sum_{i\geq\log_{2}\tilde{T}_{0}}\chi_{i}(t)\nabla_{x,t}[\delta\psi\triangle^{-1}\sum_{j=1,2}\partial_{j}(1-I)[R_{\beta}\psi_{2}R_{j}\psi_{3}-R_{j}\psi_{2}R_{\beta}\psi_{3}]],\\
\end{split}\end{equation}
where $\chi_{i}(t)$ smoothly localizes to the interval $t\sim
2^{i}$. Then we localize the frequencies and freeze
$i\in{\mathbf{Z}}$, arriving at an expression
\begin{equation}\nonumber
\chi_{i}(t)\nabla_{x,t}P_{k_{0}}[P_{k_{1}}\delta\psi\triangle^{-1}\sum_{j=1,2}\partial_{j}(1-I)P_{k}[R_{\beta}P_{k_{2}}\psi_{2}R_{j}P_{k_{3}}\psi_{3}-R_{j}P_{k_{2}}\psi_{2}R_{\beta}P_{k_{3}}\psi_{3}]]
\end{equation}
We distinguish between the following cases:\\
(I.a): {\it{One of the following options hold:}} $i\lesssim
|k_{2}|$, $i\lesssim |k_{3}|$, $i\lesssim |k_{0}-k_{1}|$,
$i\lesssim\min\{|k-k_{1}|,|k-k_{2}|\}$. This case is handled by
means of lemma~\ref{auxiliary envelope} as well as the following
lemma, provided $P_{k_{2,3}}\psi_{2,3}$ are of the first type:
\begin{lemma}\label{simple trilinear} \cite{Kr-4} Let $\psi_{1,2,3}\in{\mathbf{R}}^{2+1}$. Then, for integers
$k_{1,2,3}$ and suitable $\delta_{1,2}>0$, the following
inequality holds:
\\
\begin{equation}\nonumber\begin{split}
&||\nabla_{x,t}P_{0}[P_{k_{1}}\psi_{1}\nabla^{-1}(1-I)[R_{\nu}P_{k_{2}}\psi_{2}R_{j}P_{k_{3}}\psi_{3}-R_{j}P_{k_{2}}\psi_{2}R_{\nu}P_{k_{3}}\psi_{3}]]||_{N[0]}
\\&\lesssim
2^{\delta_{1}\min\{-\min\{k_{1},k_{2},k_{3}\},0\}}\prod_{i}2^{\delta_{2}\min\{\max_{j\neq
i}\{k_{i},k_{i}-k_{j}\},0\}}\prod_{l}||P_{k_{l}}\psi_{l}||_{S[k_{l}]}\\
\end{split}\end{equation}
\end{lemma}
Indeed, observe that if $i\lesssim \max\{|k_{2}|,|k_{3}|\}$, we
obtain from lemma~\ref{auxiliary envelope} that\\
$\min\{||P_{k_{2}}\psi_{2}||_{S[k_{2}]},||P_{k_{3}}\psi_{3}||_{S[k_{3}]}\}\lesssim
2^{-\mu i}$. Carrying out the summations over $k,k_{i}$ satisfying
these assumptions, we arrive at the upper bound $\lesssim 2^{-\mu
i}K\tilde{c}_{k_{0}}$. Summing over $i\geq \log_{2}\tilde{T}_{0}$
results in a small exponential gain in $T_{0}\leq \tilde{T}_{0}$.
If one of the other cases occurs, one gets an exponential gain
$2^{-\min\{\delta_{1},\delta_{2}\}i}$ from the above lemma. We are
fudging a bit since we have thrown the localizer $\chi_{i}(t)$ in
front, and this may affect the space-time Fourier support of the
expression, hence its norm $||.||_{N[k_{0}]}$. However, this is
detrimental only if the modulation (i. e. distance of the
space-time Fourier support to the light cone) is $\lesssim
2^{-i}$, and only affects those parts estimated with respect to
$||.||_{\dot{X}_{k_{0}}^{-1,-\frac{1}{2},1}}$, as null-frame
spaces aren't needed yet, see the proof in \cite{Kr-4}. Assuming
$Q_{<-i+O(1)}(\text{Output})$ to be a
$\dot{X}_{k_{0}}^{-1,-\frac{1}{2},1}$-atom, we estimate
\begin{equation}\nonumber\begin{split}
&||\chi_{i}(t)P_{k_{0}}Q_{<-i+O(1)}(\text{Output})||_{L_{t}^{1}\dot{H}^{-1}}\lesssim
||\chi_{i}(t)||_{L_{t}^{2}L_{x}^{\infty}}||P_{k_{0}}Q_{<-i+O(1)}(\text{Output})||_{L_{t}^{2}\dot{H}^{-1}}\\
&\lesssim \sum_{a<0}2^{\frac{i}{2}}2^{\frac{-i+a}{2}}\lesssim 1\\
\end{split}\end{equation}
Thus the cutoff is irrelevant.\\
Now assume at least one of $P_{k_{2,3}}$ is of the 2nd type. We
need the following lemma
\begin{lemma}\label{technical1} Let
$\psi_{2,3}\in{\mathcal{S}}({\mathbf{R}}^{2+1})$. Assume also that
$||P_{k_{2,3}}\psi_{2,3}||_{{\mathcal{S}}[k_{2,3}]}\lesssim
\frac{\tilde{c}_{k_{2,3}}}{\epsilon}$ with a frequency envelope
$\tilde{c}_{k}$ as in the preceding.  Assume that
$P_{k_{2}}\psi_{2}$ is of the 2nd type, and $P_{k_{3}}\psi_{3}$
admits a decomposition into functions of first and 2nd type as
enunciated in theorem~\ref{Moser}. Then we have for suitable
$\delta_{1,2}>0$
\begin{equation}\nonumber\begin{split}
&||P_{k}[R_{\beta}P_{k_{2}}\psi_{2}R_{j}P_{k_{3}}\psi_{3}-R_{j}P_{k_{2}}\psi_{2}R_{\beta}P_{k_{3}}\psi_{3}]||_{L_{t}^{2}L_{x}^{2}}\\
&\lesssim
2^{\frac{\min\{k_{3},k\}}{2}}2^{\delta_{1}\min\{k_{2}-k_{3},0\}}2^{\delta_{2}\min\{k-k_{2}\}}[\frac{\tilde{c}_{k_{2}}}{\epsilon}+\frac{\tilde{c}_{k_{3}}}{\epsilon}]\\
\end{split}\end{equation}
\end{lemma}
\begin{proof} First assume that $P_{k_{3}}\psi_{3}$ is of
the first type. Using the definition of
$||.||_{{\mathcal{S}}[k]}$, we infer the desired estimate for the
contributions of
\begin{equation}\nonumber
P_{k}[R_{\beta}P_{k_{2}}(1-I)\psi_{2}R_{j}P_{k_{3}}\psi_{3}-R_{j}P_{k_{2}}(1-I)\psi_{2}R_{\beta}P_{k_{3}}\psi_{3}]
\end{equation}
and similarly for
\begin{equation}\nonumber
P_{k}[R_{\beta}P_{k_{2}}\psi_{2}R_{j}P_{k_{3}}(1-I)\psi_{3}-R_{j}P_{k_{2}}\psi_{2}R_{\beta}(1-I)P_{k_{3}}\psi_{3}]
\end{equation}
Take the first expression: first consider the case
$k_{2}=k_{3}+O(1)$. We estimate, using theorem~\ref{Moser}
\begin{equation}\nonumber\begin{split}
&||P_{k}[R_{\beta}P_{k_{2}}(1-I)\psi_{2}R_{j}P_{k_{3}}\psi_{3}-R_{j}P_{k_{2}}(1-I)\psi_{2}R_{\beta}P_{k_{3}}\psi_{3}]||_{L_{t}^{2}L_{x}^{2}}\\
&\lesssim 2^{2k(\frac{1}{1+\mu}-\frac{1}{2})}||P_{k}[R_{\beta}P_{k_{2}}(1-I)\psi_{2}R_{j}P_{k_{3}}\psi_{3}-R_{j}P_{k_{2}}(1-I)\psi_{2}R_{\beta}P_{k_{3}}\psi_{3}]||_{L_{t}^{2}L_{x}^{1+\mu}}\\
&\lesssim 2^{\frac{k}{2}}2^{\frac{(k-k_{2})(1-3\mu)}{2+2\mu}}\frac{\tilde{c}_{k_{3}}}{\epsilon}\\
\end{split}\end{equation}
Next, in case $k_{2}<<k_{3}$, we estimate
\begin{equation}\nonumber\begin{split}
&||P_{k}[R_{\beta}P_{k_{2}}(1-I)\psi_{2}R_{j}P_{k_{3}}\psi_{3}]||_{L_{t}^{2}L_{x}^{2}}\\
&\lesssim (\sum_{c\in
C_{k_{3},k_{2}-k_{3}}}||P_{k}[R_{\beta}P_{k_{2}}(1-I)\psi_{2}R_{j}P_{c}\psi_{3}]||_{L_{t}^{2}L_{x}^{2}}^{2})^{\frac{1}{2}}\lesssim 2^{\frac{k_{3}}{2}}2^{\delta_{1}(k_{2}-k_{3})}\frac{\tilde{c}_{k_{3}}}{\epsilon}\\
\end{split}\end{equation}
The remaining term is estimated similarly, as is the case when
$k_{2}>>k_{3}$. Further, if for example $k_{2}=k_{3}+O(1)$, we can
estimate
\begin{equation}\nonumber\begin{split}
&||P_{k}[R_{\beta}P_{k_{2}}I\psi_{2}R_{j}P_{k_{3}}\psi_{3}-R_{j}P_{k_{2}}I\psi_{2}R_{\beta}P_{k_{3}}\psi_{3}]||_{L_{t}^{2}L_{x}^{2}}\\
&\lesssim
2^{(1-\epsilon)k}||P_{k_{2}}\psi_{2}||_{L_{t}^{2}L_{x}^{2+}}||P_{k_{3}}\psi_{3}||_{L_{t}^{\infty}L_{x}^{2}}
\lesssim
2^{(1-\epsilon)k-(\frac{1}{2}-\epsilon)k_{2}}\frac{\tilde{c}_{k_{2}}}{\epsilon}\\
\end{split}\end{equation}
The remaining cases $k=k_{2}+O(1)$, $k=k_{3}+O(1)$ are handled
similarly. Now assume that both $P_{k_{2}}\psi_{2}$ and
$P_{k_{3}}\psi_{3}$ are of 2nd type. In that case, if
$k_{1}=k_{2}+O(1)$, estimate
\begin{equation}\nonumber\begin{split}
&||P_{k}[R_{\beta}P_{k_{1}}\psi_{1}P_{k_{2}}R_{j}\psi_{2}]||_{L_{t}^{2}L_{x}^{2}}\lesssim
2^{(1-\epsilon)k}||R_{\beta}P_{k_{1}}\psi_{1}||_{L_{t}^{\infty}L_{x}^{2}}||P_{k_{2}}\psi_{2}||_{L_{t}^{2}L_{x}^{2+}}\\
&\lesssim
2^{(1-\epsilon)(k-k_{2})}2^{\frac{k}{2}}\frac{\tilde{c}_{k_{2}}}{\epsilon}
\end{split}\end{equation}
The remaining frequency interactions are treated similarly.
\end{proof}
 Returning to case (I.a) when at least one of $P_{k_{2,3}}\psi_{2,3}$ is of 2nd type, we claim that we
 have the estimate
 \begin{equation}\nonumber\begin{split}
 &||\nabla_{x,t}P_{k_{0}}[P_{k_{1}}\delta\psi_{1}\nabla^{-1}(1-I)P_{k}[R_{\nu}P_{k_{2}}\psi_{2}R_{j}P_{k_{3}}\psi_{3}-
 R_{j}P_{k_{2}}\psi_{2}R_{\nu}P_{k_{3}}\psi_{3}]]||_{N[k_{0}]}\\
&\lesssim
2^{\delta_{1}[\min_{i=2,3}\{k,k_{i}\}-\max_{i=2,3}\{k,k_{i}\}]}2^{-\delta_{2}(|k-k_{1}|)}[\frac{\tilde{c}_{k_{2}}}{\epsilon}+\frac{\tilde{c}_{k_{3}}}{\epsilon}]\tilde{c}_{k_{0}}\\
\end{split}\end{equation}
One could then sum over all frequency parameters (except $k_{0}$)
and obtain the required exponential gain in $i$ under the
hypotheses of case (I.a)\footnote{The cutoff $\chi_{i}(t)$ in
front is handled as before.}. To verify this estimate, we may
assume $k_{0}=0$. One needs to distinguish between $k_{1}\in
[-10,10]$, $k_{1}>10$, $k_{1}<-10$. These are similar, so we treat
the first case: we have
\begin{equation}\nonumber\begin{split}
&||\nabla_{x,t}P_{0}Q_{>10}[P_{k_{1}}\delta\psi_{1}\nabla^{-1}(1-I)P_{k}[R_{\nu}P_{k_{2}}\psi_{2}R_{j}P_{k_{3}}\psi_{3}
-R_{j}P_{k_{2}}\psi_{2}R_{\nu}P_{k_{3}}\psi_{3}]]||_{N[0]}\\
&\lesssim
||P_{0}Q_{>10}[P_{k_{1}}\delta\psi_{1}\nabla^{-1}(1-I)P_{k}[R_{\nu}P_{k_{2}}\psi_{2}R_{j}P_{k_{3}}\psi_{3}
-R_{j}P_{k_{2}}\psi_{2}R_{\nu}P_{k_{3}}\psi_{3}]]||_{L_{t}^{2}L_{x}^{2}}\\
&\lesssim ||P_{k_{1}}\delta\psi_{1}||_{L_{t}^{\infty}L_{x}^{2}}
2^{\frac{\min\{k,k_{3}\}}{2}}2^{\delta_{1}\min\{k_{2}-k_{3},0\}}2^{\delta_{2}\min\{k-k_{2},0\}}[\frac{\tilde{c}_{k_{2}}}{\epsilon}+\frac{\tilde{c}_{k_{3}}}{\epsilon}]\tilde{c}_{0}\\
\end{split}\end{equation}
One checks that this verifies the claim, with a lot to spare.
Next, we can estimate
\begin{equation}\nonumber\begin{split}
&||\nabla_{x,t}P_{0}Q_{<10}[P_{k_{1}}Q_{<k-100}\delta\psi_{1}\\
&\hspace{3cm}\nabla^{-1}(1-I)P_{k}[R_{\nu}P_{k_{2}}\psi_{2}R_{j}P_{k_{3}}\psi_{3}-R_{j}P_{k_{2}}\psi_{2}R_{\nu}P_{k_{3}}\psi_{3}]]||_{N[0]}\\
&\lesssim||\nabla_{x,t}P_{0}Q_{<10}[P_{k_{1}}Q_{<k-100}\delta\psi_{1}\\
&\hspace{3cm}\nabla^{-1}(1-I)P_{k}[R_{\nu}P_{k_{2}}\psi_{2}R_{j}P_{k_{3}}\psi_{3}-R_{j}P_{k_{2}}\psi_{2}R_{\nu}P_{k_{3}}\psi_{3}]]||_{\dot{X}_{0}^{-1,-\frac{1}{2},1}}\\
&\lesssim
2^{-\frac{k}{2}}2^{\frac{\min\{k,k_{3}\}}{2}}2^{\delta_{1}\min\{k_{2}-k_{3},0\}}2^{\delta_{2}\min\{k-k_{2},0\}}\tilde{c}_{0}[\frac{\tilde{c}_{k_{2}}}{\epsilon}
+\frac{\tilde{c}_{k_{3}}}{\epsilon}]\\
\end{split}\end{equation}
Again this verifies the claim. Finally, we have the estimate
\begin{equation}\nonumber\begin{split}
&||\nabla_{x,t}P_{0}Q_{<10}[P_{k_{1}}Q_{\geq k-100}\delta\psi_{1}\\
&\hspace{3cm}\nabla^{-1}(1-I)P_{k}[R_{\nu}P_{k_{2}}\psi_{2}R_{j}P_{k_{3}}\psi_{3}-R_{j}P_{k_{2}}\psi_{2}R_{\nu}P_{k_{3}}\psi_{3}]]||_{N[0]}\\
&\lesssim||\nabla_{x,t}P_{0}Q_{<10}[P_{k_{1}}Q_{\geq k-100}\delta\psi_{1}\\
&\hspace{3cm}\nabla^{-1}(1-I)P_{k}[R_{\nu}P_{k_{2}}\psi_{2}R_{j}P_{k_{3}}\psi_{3}-R_{j}P_{k_{2}}\psi_{2}R_{\nu}P_{k_{3}}\psi_{3}]]||_{L_{t}^{1}\dot{H}^{-1}}\\
&\lesssim ||P_{k_{1}}Q_{\geq
k-100}\delta\psi_{1}||_{L_{t}^{2}L_{x}^{2}}\\&\hspace{3cm}||\nabla^{-1}(1-I)P_{k}[R_{\nu}P_{k_{2}}\psi_{2}R_{j}P_{k_{3}}\psi_{3}-R_{j}P_{k_{2}}\psi_{2}R_{\nu}P_{k_{3}}\psi_{3}]||_{L_{t}^{2}L_{x}^{\infty}}\\
&\lesssim2^{-\frac{k}{2}}2^{\frac{\min\{k,k_{3}\}}{2}}2^{\delta_{1}\min\{k_{2}-k_{3},0\}}2^{\delta_{2}\min\{k-k_{2},0\}}\tilde{c}_{0}[\frac{\tilde{c}_{k_{2}}}{\epsilon}
+\frac{\tilde{c}_{k_{3}}}{\epsilon}]\\
\end{split}\end{equation}
as in the preceding estimate. This concludes case (I.a). \\
(I.b): {\it{$i\lesssim |k|$, $i\lesssim |k_{1}|$, and none of the
properties in (I.a) hold.}} Thus in this case, we have $k\lesssim
-i$, $k_{1}\lesssim -i$, $|k-k_{1}|<<i$; we may treat the last
difference as $O(1)$\footnote{We do this in order to avoid
carrying too many small constants around; this is legitimate since
the exponential gains obtained later are independent.}. In this
case we have to work harder to obtain the exponential gain in $i$,
since the previous trilinear estimates won't suffice. Observe that
we only need to worry about the case $\nu=0$, though, since
otherwise one can pull a derivative out of the inner square
bracket expression. Also, we may easily reduce the Fourier support
of $P_{k_{2,3}}\psi_{2,3}$ to the hyperbolic regime\footnote{We
shall not include the localizers $Q_{<k_{2,3}}$ everywhere in
order to streamline notation.} (distance to light cone at most
comparable to frequency). First, consider the case $|k|>(1+\mu)i$,
for some small $\mu>0$. In that case, we have
\begin{equation}\nonumber\begin{split}
 &||\nabla_{x,t}P_{k_{0}}\chi_{i}(t)[P_{k_{1}}\delta\psi_{1}\nabla^{-1}(1-I)P_{k}[R_{\nu}P_{k_{2}}\psi_{2}R_{j}P_{k_{3}}\psi_{3}-
 R_{j}P_{k_{2}}\psi_{2}R_{\nu}P_{k_{3}}\psi_{3}]]||_{N[k_{0}]}\\
 &\lesssim ||\nabla_{x,t}P_{k_{0}}\chi_{i}(t)[P_{k_{1}}Q_{\geq k-100}\delta\psi_{1}\\&\hspace{3.5cm}\nabla^{-1}(1-I)P_{k}[R_{\nu}P_{k_{2}}\psi_{2}R_{j}P_{k_{3}}\psi_{3}-
 R_{j}P_{k_{2}}\psi_{2}R_{\nu}P_{k_{3}}\psi_{3}]]||_{N[k_{0}]}\\
&+||\nabla_{x,t}P_{k_{0}}\chi_{i}(t)[P_{k_{1}}Q_{<k-100}\delta\psi_{1}\\&\hspace{3.5cm}\nabla^{-1}(1-I)P_{k}[R_{\nu}P_{k_{2}}\psi_{2}R_{j}P_{k_{3}}\psi_{3}-
 R_{j}P_{k_{2}}\psi_{2}R_{\nu}P_{k_{3}}\psi_{3}]]||_{N[k_{0}]}\\
\end{split}\end{equation}
The first summand is further decomposed as follows:
\begin{equation}\nonumber\begin{split}
 &||\nabla_{x,t}P_{k_{0}}\chi_{i}(t)[P_{k_{1}}Q_{\geq k-100}\delta\psi_{1}\\&\hspace{3.5cm}\nabla^{-1}(1-I)P_{k}[R_{\nu}P_{k_{2}}\psi_{2}R_{j}P_{k_{3}}\psi_{3}-
 R_{j}P_{k_{2}}\psi_{2}R_{\nu}P_{k_{3}}\psi_{3}]]||_{N[k_{0}]}\\
 &\lesssim ||\nabla_{x,t}P_{k_{0}}Q_{>k_{0}}\chi_{i}(t)[P_{k_{1}}Q_{\geq k-100}\delta\psi_{1}\\&\hspace{3cm}\nabla^{-1}(1-I)P_{k}[R_{\nu}P_{k_{2}}\psi_{2}R_{j}P_{k_{3}}\psi_{3}-
 R_{j}P_{k_{2}}\psi_{2}R_{\nu}P_{k_{3}}\psi_{3}]]||_{\dot{X}_{k_{0}}^{-\frac{1}{2},-1,2}}\\
&+||\nabla_{x,t}P_{k_{0}}Q_{<k_{0}}\chi_{i}(t)[P_{k_{1}}Q_{\geq
k-100}\delta\psi_{1}\\&\hspace{3cm}\nabla^{-1}(1-I)P_{k}[R_{\nu}P_{k_{2}}\psi_{2}R_{j}P_{k_{3}}\psi_{3}-
 R_{j}P_{k_{2}}\psi_{2}R_{\nu}P_{k_{3}}\psi_{3}]]||_{L_{t}^{1}\dot{H}^{-1}}\\
 \end{split}\end{equation}
We then estimate
\begin{equation}\nonumber\begin{split}
&||\nabla_{x,t}P_{k_{0}}Q_{>k_{0}}\chi_{i}(t)[P_{k_{1}}Q_{\geq
k-100}\delta\psi_{1}\\&\hspace{3cm}\nabla^{-1}(1-I)P_{k}[R_{\nu}P_{k_{2}}\psi_{2}R_{j}P_{k_{3}}\psi_{3}-
 R_{j}P_{k_{2}}\psi_{2}R_{\nu}P_{k_{3}}\psi_{3}]]||_{\dot{X}_{k_{0}}^{-\frac{1}{2},-1,2}}\\
&\lesssim
2^{\frac{k_{0}}{2}}||\chi_{i}(t)||_{L_{t}^{2}}||P_{k_{1}}Q_{\geq
k-100}\delta\psi_{1}||_{L_{t}^{\infty}L_{x}^{2}}\\&\hspace{3cm}||\nabla^{-1}(1-I)P_{k}[R_{\nu}P_{k_{2}}\psi_{2}R_{j}P_{k_{3}}\psi_{3}-
 R_{j}P_{k_{2}}\psi_{2}R_{\nu}P_{k_{3}}\psi_{3}]||_{L_{t}^{\infty}L_{x}^{2}}\\
 &\lesssim 2^{\frac{k_{0}}{2}+\frac{i}{2}}\tilde{c}_{k_{1}}\lesssim
 2^{-\frac{\mu}{2}
 i}\tilde{c}_{k_{1}}\\
 \end{split}\end{equation}
Similarly, we have
\begin{equation}\nonumber\begin{split}
&||\nabla_{x,t}P_{k_{0}}Q_{<k_{0}}\chi_{i}(t)[P_{k_{1}}Q_{\geq
k-100}\delta\psi_{1}\\&\hspace{3cm}\nabla^{-1}(1-I)P_{k}[R_{\nu}P_{k_{2}}\psi_{2}R_{j}P_{k_{3}}\psi_{3}-
 R_{j}P_{k_{2}}\psi_{2}R_{\nu}P_{k_{3}}\psi_{3}]]||_{L_{t}^{1}\dot{H}^{-1}}\\
 &\lesssim 2^{k_{0}}||\chi_{i}(t)||_{L_{t}^{2}}||P_{k_{1}}Q_{\geq
k-100}\delta\psi_{1}||_{L_{t}^{2}L_{x}^{2}}\\&\hspace{3cm}||\nabla^{-1}(1-I)P_{k}[R_{\nu}P_{k_{2}}\psi_{2}R_{j}P_{k_{3}}\psi_{3}-
 R_{j}P_{k_{2}}\psi_{2}R_{\nu}P_{k_{3}}\psi_{3}]||_{L_{t}^{\infty}L_{x}^{2}}, \\
\end{split}\end{equation}
and this is controlled by
$2^{k_{0}-\frac{k}{2}+\frac{i}{2}}\tilde{c}_{k_{1}}\lesssim
2^{-\frac{\mu }{2}i}\tilde{c}_{k_{1}}$ as desired. The remaining
terms are handled similarly. Thus we now assume that $i\lesssim
|k|\leq(1+\mu)|i|$. We then claim that we may replace the operator
$(1-I)$ by $Q_{>\frac{k}{2}}$. Indeed, we have
\begin{equation}\nonumber\begin{split}
&||\nabla_{x,t}P_{k_{0}}\chi_{i}(t)[P_{k_{1}}\delta\psi_{1}\\&\hspace{2cm}\nabla^{-1}(1-I)P_{k}Q_{<\frac{k}{2}}[R_{\nu}P_{k_{2}}\psi_{2}R_{j}P_{k_{3}}\psi_{3}-
 R_{j}P_{k_{2}}\psi_{2}R_{\nu}P_{k_{3}}\psi_{3}]]||_{N[k_{0}]}\\
&\lesssim
\sum_{k+100<a<\frac{k}{2}}||\nabla_{x,t}P_{k_{0}}\chi_{i}(t)[P_{k_{1}}Q_{<a-100}\delta\psi_{1}\\&\hspace{2cm}\nabla^{-1}(1-I)P_{k}Q_{a}[R_{\nu}P_{k_{2}}\psi_{2}R_{j}P_{k_{3}}\psi_{3}-
 R_{j}P_{k_{2}}\psi_{2}R_{\nu}P_{k_{3}}\psi_{3}]]||_{N[k_{0}]}\\
&+\sum_{k+100<a<\frac{k}{2}}||\nabla_{x,t}P_{k_{0}}\chi_{i}(t)[P_{k_{1}}Q_{\geq
a-100}\delta\psi_{1}\\&\hspace{2cm}\nabla^{-1}(1-I)P_{k}Q_{a}[R_{\nu}P_{k_{2}}\psi_{2}R_{j}P_{k_{3}}\psi_{3}-
 R_{j}P_{k_{2}}\psi_{2}R_{\nu}P_{k_{3}}\psi_{3}]]||_{N[k_{0}]}\\
 \end{split}\end{equation}
We treat the first summand, the 2nd being similar. We have
\begin{equation}\nonumber\begin{split}
&||\nabla_{x,t}P_{k_{0}}[\chi_{i}(t)[P_{k_{1}}Q_{<a-100}\delta\psi_{1}\\&\hspace{2cm}\nabla^{-1}(1-I)P_{k}Q_{a}[R_{\nu}P_{k_{2}}\psi_{2}R_{j}P_{k_{3}}\psi_{3}-
 R_{j}P_{k_{2}}\psi_{2}R_{\nu}P_{k_{3}}\psi_{3}]]]||_{N[k_{0}]}\\
 &\lesssim 2^{-\frac{k_{0}}{2}}||P_{k_{0}}Q_{a+O(1)}[\chi_{i}(t)[P_{k_{1}}Q_{<a-100}\delta\psi_{1}\\&\hspace{2cm}\nabla^{-1}(1-I)P_{k}Q_{a}[R_{\nu}P_{k_{2}}\psi_{2}R_{j}P_{k_{3}}\psi_{3}-
 R_{j}P_{k_{2}}\psi_{2}R_{\nu}P_{k_{3}}\psi_{3}]]]||_{L_{t}^{2}L_{x}^{2}}\\&
\lesssim 2^{\frac{a+k}{2}}\tilde{c}_{k_{1}}\\
 \end{split}\end{equation}
This is clearly acceptable. We now notice the identity
\begin{equation}\nonumber
-2(\partial_{t}\psi_{2}\partial_{r}\psi_{3}-\partial_{r}\psi_{2}\partial_{t}\psi_{3})=(\partial_{t}+\partial_{r})\psi_{2}
(\partial_{t}-\partial_{r})\psi_{3}-(\partial_{t}-\partial_{r})\psi_{2}(\partial_{t}+\partial_{r})\psi_{3}
\end{equation}
Applying this to our frequency localized situation, we have the
identity\footnote{Recall the suppressed localizations, see
previous footnote.}(recall that $\psi_{1,2}$ are radial)
\begin{equation}\nonumber\begin{split}
&R_{0}P_{k_{2}}Q_{<k_{2}}\psi_{2}R_{i}P_{k_{3}}Q_{<k_{3}}\psi_{3}-R_{i}P_{k_{2}}Q_{<k_{2}}\psi_{2}R_{0}P_{k_{3}}Q_{<k_{3}}\psi_{3}\\&=\frac{x_{i}}{r}[(\partial_{t}+\partial_{r})\nabla^{-1}P_{k_{2}}Q_{<k_{2}}\psi_{2}
(\partial_{t}-\partial_{r})\nabla^{-1}P_{k_{3}}Q_{<k_{3}}\psi_{3}\\&-(\partial_{t}-\partial_{r})\nabla^{-1}P_{k_{2}}Q_{<k_{2}}\psi_{2}(\partial_{t}+\partial_{r})\nabla^{-1}P_{k_{3}}Q_{<k_{3}}\psi_{3}]\\
\end{split}\end{equation}
Now let $\phi\in C_{0}^{\infty}({\mathbf{R}})$ be a smooth cutoff
and use the decomposition\\
$\psi_{2,3}=\phi_{2^{-\frac{i}{2+}}}(u)\psi_{2,3}+(1-\phi_{2^{-\frac{i}{2+}}}(u))\psi_{2,3}$,
where $u=t-r$ and $\phi_{\lambda}(u)=\phi(\frac{u}{\lambda})$. Now
Proposition~\ref{Chr1} implies that
\begin{equation}\nonumber
||\chi_{i}(t)\phi_{2^{-\frac{i}{2+}}}(u)\psi_{2,3}||_{L_{x}^{2}}\lesssim
\sqrt{t^{-1}\times t\times 2^{-\frac{i}{2+}}}\lesssim
2^{-\frac{i}{4+}}
\end{equation}
Now let $\frac{2}{4+}+\frac{1}{M}=\frac{1}{2}$ and return to the
full expression. We first estimate the large-modulation
contribution:
\begin{equation}\nonumber\begin{split}
&||P_{k_{0}}Q_{>k_{0}}\nabla_{x,t}\chi_{i}(t)[P_{k_{1}}\delta\psi_{1}
\\&\nabla^{-1}P_{k}Q_{>\frac{k}{2}}[R_{0}P_{k_{2}}Q_{<k_{2}}[\phi_{2^{-\frac{i}{2+}}}(u)\psi_{2}]R_{i}P_{k_{3}}Q_{<k_{3}}\psi_{3}\\
&\hspace{4.3cm}-R_{i}P_{k_{2}}Q_{<k_{2}}[\phi_{2^{-\frac{i}{2+}}}(u)\psi_{2}]R_{0}P_{k_{3}}Q_{<k_{3}}\psi_{3}||_{N[k_{0}]}\\
&\lesssim
||P_{k_{0}}Q_{>k_{0}}\nabla_{x,t}\chi_{i}(t)[P_{k_{1}}\delta\psi_{1}
\\&\hspace{2cm}\nabla^{-1}P_{k}Q_{\geq \frac{k}{2}}[R_{0}P_{k_{2}}Q_{<k_{2}}[\phi_{2^{-\frac{i}{2+}}}(u)\psi_{2}]R_{i}P_{k_{3}}Q_{<k_{3}}\psi_{3}\\
&\hspace{3.8cm}-R_{i}P_{k_{2}}Q_{<k_{2}}[\phi_{2^{-\frac{i}{2+}}}(u)\psi_{2}]R_{0}P_{k_{3}}Q_{<k_{3}}\psi_{3}||_{\dot{X}_{k_{0}}^{-\frac{1}{2},-1,2}}\\
&\lesssim
2^{-\frac{k_{0}}{2}}2^{(1-\epsilon)k_{0}-k_{1}}||\chi_{i}(t)||_{L_{t}^{4+}}||P_{k_{1}}\delta\psi_{1}||_{L_{t}^{4+}L_{x}^{\infty}}\\&
\hspace{4.1cm}||R_{0}P_{k_{2}}Q_{<k_{2}}[\phi_{2^{-\frac{i}{2+}}}(u)\psi_{2}]||_{L_{t}^{\infty}L_{x}^{2}}
||P_{k_{3}}\psi_{3}||_{L_{t}^{M}L_{x}^{2+}}\\
&\lesssim
2^{(\frac{1}{2}-\epsilon)k_{0}-\frac{k_{1}}{4+}}\tilde{c}_{k_{1}}\frac{\tilde{c}_{k_{3}}}{\epsilon}\\
\end{split}\end{equation}
Keeping our assumptions on the frequencies in mind, this is more
than what we need. Next, restricting the expression to modulation
$\leq 2^{k_{0}}$, we have
\begin{equation}\nonumber\begin{split}
&||P_{k_{0}}Q_{\leq
k_{0}}\nabla_{x,t}\chi_{i}(t)[P_{k_{1}}\delta\psi_{1}
\\&\nabla^{-1}P_{k}Q_{\geq \frac{k}{2}}[R_{0}P_{k_{2}}Q_{<k_{2}}[\phi_{2^{-\frac{i}{2+}}}(u)\psi_{2}]R_{i}P_{k_{3}}Q_{<k_{3}}\psi_{3}\\
&\hspace{4.3cm}-R_{i}P_{k_{2}}Q_{<k_{2}}[\phi_{2^{-\frac{i}{2+}}}(u)\psi_{2}]R_{0}P_{k_{3}}Q_{<k_{3}}\psi_{3}||_{N[k_{0}]}\\
&\lesssim \sum_{k_{2}+O(1)>a\geq\frac{k}{2}} ||P_{k_{0}}Q_{\leq
k_{0}}\nabla_{x,t}\chi_{i}(t)[P_{k_{1}}Q_{a+O(1)}\delta\psi_{1}
\\&\nabla^{-1}P_{k}Q_{a}[R_{0}P_{k_{2}}Q_{<k_{2}}[\phi_{2^{-\frac{i}{2+}}}(u)\psi_{2}]R_{i}P_{k_{3}}Q_{<k_{3}}\psi_{3}\\
&\hspace{4.3cm}-R_{i}P_{k_{2}}Q_{<k_{2}}[\phi_{2^{-\frac{i}{2+}}}(u)\psi_{2}]R_{0}P_{k_{3}}Q_{<k_{3}}\psi_{3}||_{L_{t}^{1}\dot{H}^{-1}}\\
&\lesssim\sum_{k_{2}+O(1)>a>\frac{k}{2}}
2^{a-k}2^{\frac{k}{2+}}||\chi_{i}(t)||_{L_{t}^{4+}}||P_{k_{1}}Q_{a+O(1)}\delta\psi_{1}||_{L_{t}^{2}L_{x}^{\infty}}\\&\hspace{5cm}||\phi_{2^{-\frac{i}{2+}}}(u)\psi_{2}||_{L_{t}^{\infty}L_{x}^{2}}
||P_{k_{3}}Q_{<k_{3}}\psi_{3}||_{L_{t}^{4-}L_{x}^{4+}}\\
\end{split}\end{equation}
We can bound this by
\begin{equation}\nonumber
2^{a-k}2^{\frac{i}{4+}}2^{\frac{k}{2+}}2^{k_{1}}2^{-\frac{i}{4+}}2^{-\frac{a}{2}}\tilde{c}_{k_{1}}
\end{equation}
Our assumptions ensure that we may sum over
$\frac{k}{2}<a<k_{2}+O(1)$, resulting in an exponential gain in
$i$. An identical argument may be used when
$P_{k_{3}}Q_{<k_{3}}\psi_{3}$ is replaced by
$P_{k_{3}}Q_{<k_{3}}[\phi_{2^{-\frac{i}{2+}}}(u)\psi_{3}]$, so we
may replace both $P_{k_{2,3}}Q_{<k_{2,3}}\psi_{2,3}$ by
$P_{k_{2,3}}Q_{<k_{2,3}}[(1-\phi_{2^{-\frac{i}{2+}}}(u))\psi_{3}]$.
In that case we utilize the null-form identity recorded earlier:
use
\begin{equation}\nonumber\begin{split}
&R_{0}P_{k_{2}}Q_{<k_{2}}[(1-\phi_{2^{-\frac{i}{2+}}})\psi_{2}]R_{i}P_{k_{3}}Q_{<k_{3}}[(1-\phi_{2^{-\frac{i}{2+}}})\psi_{3}]\\
&\hspace{3cm}-R_{i}P_{k_{2}}Q_{<k_{2}}[(1-\phi_{2^{-\frac{i}{2+}}})\psi_{2}]R_{0}P_{k_{3}}Q_{<k_{3}}[(1-\phi_{2^{-\frac{i}{2+}}})\psi_{3}]\\
&=\frac{x_{i}}{r}[(\partial_{t}+\partial_{r})\nabla^{-1}P_{k_{2}}Q_{<k_{2}}[(1-\phi_{2^{-\frac{i}{2+}}})\psi_{2}]
(\partial_{t}-\partial_{r})\nabla^{-1}P_{k_{3}}Q_{<k_{3}}[(1-\phi_{2^{-\frac{i}{2+}}})\psi_{3}]\\
&-\frac{x_{i}}{r}[(\partial_{t}-\partial_{r})\nabla^{-1}P_{k_{2}}Q_{<k_{2}}[(1-\phi_{2^{-\frac{i}{2+}}})\psi_{2}]
(\partial_{t}+\partial_{r})\nabla^{-1}P_{k_{3}}Q_{<k_{3}}[(1-\phi_{2^{-\frac{i}{2+}}})\psi_{3}]\\
\end{split}\end{equation}
We can now exploit the fact that
$|\frac{k}{2}|<\frac{(1+\mu)i}{2}<i(1-\delta)$ for $\mu$ small
enough, as well as $T_{0}$ large enough. Thus we may move the
multiplier $\chi_{i}(t)$ past the Fourier multiplier
$\nabla^{-1}P_{k}Q_{>\frac{k}{2}}$ while trading in errors
exponentially decreasing \footnote{More precisely, these errors
behave like $[1+2^{-(1-\delta)i}(2^{i-a}-t)]^{-N}$ for $t<2^{i-a}$
and $[1+2^{-(1-\delta)i}(2^{i+b}-t)]^{-N}$ for $t>2^{i+b}$,
provided $\chi_{i}(t)$ is supported in
$[2^{i-a},2^{i+b}]$.}outside of $t\sim 2^{i}$. In other words,
under the present assumptions on the frequencies, we may write
schematically
\begin{equation}\nonumber\begin{split}
&\chi_{i}(t)\nabla_{x,t}P_{k_{0}}[P_{k_{1}}\delta\psi_{1}\nabla^{-1}P_{k}Q_{>\frac{k}{2}}[,]]
\\&\hspace{3cm}=\chi_{i}(t)\nabla_{x,t}P_{k_{0}}[P_{k_{1}}\delta\psi_{1}\nabla^{-1}P_{k}Q_{>\frac{k}{2}}(\chi_{i1}(t)+\chi_{i2}(t))[,]],\\
\end{split}\end{equation}
where $\chi_{i1}(t)$ is supported on $t\sim 2^{i}$ while
$|\chi_{i2}(t)|\lesssim 2^{-Ni}$ for $t<2^{i+O(1)}$, as well as
$|\chi_{i2}(t)|<t^{-N}$ for $t>>2^{i}$. It is then easy to verify
that this leads to acceptable terms, so we may focus on the
contribution of $\chi_{1i}(t)$. We shall want to move the operator
$\partial_{t}+\partial_{r}$ past the Fourier localizer
$P_{k_{2}}Q_{<k_{2}}$. We write
\begin{equation}\nonumber\begin{split}
&P_{k_{2}}Q_{<k_{2}}[(1-\phi_{2^{-\frac{i}{2+}}}(u))\psi_{2}]\\&=\rho_{i}(r)P_{k_{2}}Q_{<k_{2}}[(1-\phi_{2^{-\frac{i}{2+}}}(u))\psi_{2}]
+(1-\rho_{i}(r))P_{k_{2}}Q_{<k_{2}}[(1-\phi_{2^{-\frac{i}{2+}}}(u))\psi_{2}],\\
\end{split}\end{equation}
where $\rho_{i}(r)$ localizes smoothly to a disc of radius $\sim
2^{i-10}$ around the origin. Now on account of the fact that
$\nabla^{-1}P_{k_{2}}$ is given by a convolution kernel which
decays rapidly outside of a disc of radius $<2^{\frac{i}{C}}$, we
see by means of Proposition~\ref{Chr1} that
\begin{equation}\nonumber
||\chi_{i1}(t)\rho_{i}(r)\nabla^{-1}P_{k_{2}}Q_{<k_{2}}[(1-\phi_{2^{-\frac{i}{2+}}}(u))\psi_{2}]||_{L_{x}^{\infty}}\lesssim
2^{-\frac{3i}{2+}}
\end{equation}
One  then easily concludes that the contribution of this term is
negligible: indeed, plugging it into the inner bracket instead of
$P_{k_{2}}\psi_{2}$ and using schematic notation, we can estimate
for example
\begin{equation}\nonumber\begin{split}
&||\chi_{i}(t)\nabla_{x,t}P_{k_{0}}[P_{k_{1}}Q_{<k_{1}-100}\delta\psi_{1}\nabla^{-1}P_{k}Q_{>\frac{k}{2}}(\chi_{i1}(t)[,])]||_{N[k_{0}]}\\
&\lesssim||\chi_{i}(t)\nabla_{x,t}P_{k_{0}}[P_{k_{1}}Q_{<k_{1}-100}\delta\psi_{1}\nabla^{-1}P_{k}Q_{>\frac{k}{2}}(\chi_{i1}(t)[,])]||_{\dot{X}_{k_{0}}^{-\frac{1}{2},-1,2}}\\
&\lesssim
2^{+\frac{k_{0}}{2}-k}||P_{k_{1}}Q_{<k_{1}-100}\delta\psi_{1}||_{L_{t}^{\infty}L_{x}^{2}}||P_{k}Q_{>\frac{k}{2}}(\chi_{i1}(t)[,])||_{L_{t}^{2}L_{x}^{2}},\\
\end{split}\end{equation}
where
\begin{equation}\nonumber\begin{split}
&||P_{k}Q_{>\frac{k}{2}}(\chi_{i1}(t)[,])||_{L_{t}^{2}L_{x}^{2}}\\&\lesssim
||\chi_{i1}(t)(\partial_{t}+\partial_{r})[\rho_{i}(r)\nabla^{-1}P_{k_{2}}Q_{<k_{2}}[(1-\phi_{2^{-\frac{i}{2+}}}(u))\psi_{2}]]||_{L_{t}^{2}L_{x}^{\infty}}
\\&\hspace{4cm}||(\partial_{t}-\partial_{r})\nabla^{-1}P_{k_{3}}Q_{<k_{3}}[(1-\phi_{2^{-\frac{i}{2+}}}(u))\psi_{2}]||_{L_{t}^{\infty}L_{x}^{2}}\\
&\lesssim 2^{-(1-\epsilon)i}\\
\end{split}\end{equation}
Our assumption $|k|<(1+\mu)i$ shows that putting these estimates
together gives an acceptable bound. The contribution of
$P_{k_{1}}Q_{\geq k_{1}-100}\delta\psi_{1}$ is handled similarly.
Now consider the contribution of the term with $\rho_{i}(r)$
replaced by $(1-\rho_{i}(r))$. We use the fact (see e. g.
\cite{Tao 1}) that
\begin{equation}\nonumber
[\frac{x_{i}}{r}(1-\rho_{i}(r)),P_{k_{2}}]"="2^{-k_{2}}\nabla(\frac{x_{i}}{r}(1-\rho_{i}(r)))
\end{equation}
where the latter expression stands for a weighted average of
translates of the derivatives of $\frac{x_{i}}{r}(1-\rho_{i}(r))$.
Notice that
\begin{equation}\nonumber
||\nabla_{x}(\frac{x_{i}}{r}(1-\rho_{i}(r)))||_{L_{x}^{\infty}}\lesssim
2^{-i},
\end{equation}
hence the contribution of the commutator is treated exactly as the
contribution of the term
$\rho_{i}(r)P_{k_{2}}Q_{<k_{2}}[(1-\phi_{2^{-\frac{i}{2+}}}(u))\psi_{2}]$.
This finally allows us to move the operator
$\partial_{t}+\partial_{r}$ past the operator
$\nabla^{-1}P_{k_{2}}Q_{<k_{2}}$. Arguing as before, one can also
move the operator $\chi_{i1}(t)$ past the localizer
$\nabla^{-1}P_{k_{2}}Q_{<k_{2}}$, generating acceptable error
terms. Using lemma~\ref{Chr2}, we can now estimate
\begin{equation}\nonumber\begin{split}
&||P_{k_{2}}Q_{<k_{2}}\nabla^{-1}\chi_{i1}(t)(1-\rho_{i}(r))(\partial_{t}+\partial_{r})[(1-\phi_{2^{-\frac{i}{2+}}}(u))\psi_{2}]
\\&\hspace{5.5cm}(\partial_{t}-\partial_{r})\nabla^{-1}P_{k_{3}}Q_{<k_{3}}[(1-\phi_{2^{-\frac{i}{2+}}}(u))\psi_{2}]||_{L_{t}^{2}L_{x}^{2}}\\
&\lesssim
||P_{k_{2}}Q_{<k_{2}}\nabla^{-1}\chi_{i1}(t)(1-\rho_{i}(r))(\partial_{t}+\partial_{r})[(1-\phi_{2^{-\frac{i}{2+}}}(u))\psi_{2}]||_{L_{t}^{2}L_{x}^{\infty}}
\\&\hspace{5cm}||(\partial_{t}-\partial_{r})\nabla^{-1}P_{k_{3}}Q_{<k_{3}}[(1-\phi_{2^{-\frac{i}{2+}}}(u))\psi_{2}]||_{L_{t}^{\infty}L_{x}^{2}}\\
&\lesssim 2^{-\frac{i}{2-}}\\
\end{split}\end{equation}
Proceeding as before, one deduces from this (and an analogous
estimate for the term with $(\partial_{t}\pm\partial_{r})$
interchanged) the following estimate for the full expression:
\begin{equation}\nonumber
||\chi_{i}(t)\nabla_{x,t}P_{k_{0}}[P_{k_{1}}\delta\psi_{1}\nabla^{-1}P_{k}Q_{>\frac{k}{2}}[,]]||_{N[k_{0}]}
\lesssim 2^{\frac{k_{0}}{2}-k}2^{-\frac{i}{2-}}\tilde{c}_{k_{1}}
\end{equation}
Since $|k|<(1+\mu)i$ and by assumption $k_{0}\leq k+O(1)$,
choosing $\mu$ small enough allows us to
get an exponential gain in $i$. This concludes case (I.b).\\
(I.c): {\it{None of (I.a), (I.b) hold, and $i\lesssim |k_{1}|$.}}
This then implies $|k_{2,3}|,|k|<<i$, and we shall treat these as
$O(1)$. Note that also necessarily $i\lesssim k_{1}=k_{0}+O(1)$.
One first reduces $P_{k_{2,3}}\psi_{2,3}$ to modulation
$<2^{\delta i}$, where $\delta$ is very small but such that
$2^{\delta i}>>\max\{|k_{2,3}|,|k|\}$. We shall treat the latter
quantities as $O(1)$. To achieve this, one estimates for example
\begin{equation}\nonumber\begin{split}
&||P_{k_{0}}Q_{<k_{0}}\nabla_{x,t}\chi_{i}(t)[P_{k_{1}}Q_{<\delta
i-10}\delta\psi_{1}\nabla^{-1}(1-I)[P_{k_{2}}Q_{>\delta
i}\psi_{2}, P_{k_{3}}Q_{<\delta
i-10}\psi_{3}]]||_{\dot{X}_{k_{0}}^{-1,-\frac{1}{2},1}}\\
&\lesssim
2^{-\frac{\delta}{2}i}||P_{k_{1}}\delta\psi_{1}||_{L_{t}^{\infty}L_{x}^{2}}||P_{k_{2}}Q_{>\delta
i}R_{\nu}\psi_{2}||_{L_{t}^{2}L_{x}^{2}}||P_{k_{3}}\psi_{3}||_{L_{t}^{\infty}L_{x}^{\infty}}\lesssim 2^{-\frac{\delta}{2+}i}\tilde{c}_{k_{1}}\\
\end{split}\end{equation}
Similarly, one has
\begin{equation}\nonumber\begin{split}
&||P_{k_{0}}Q_{<k_{0}}\nabla_{x,t}\chi_{i}(t)[P_{k_{1}}Q_{\geq
\delta i-10}\delta\psi_{1}\nabla^{-1}(1-I)[P_{k_{2}}Q_{>\delta
i}\psi_{2}, P_{k_{3}}Q_{<\delta
i-10}\psi_{3}]]||_{L_{t}^{1}\dot{H}^{-1}}\\
&\lesssim ||P_{k_{1}}Q_{\geq \delta
i-10}\delta\psi_{1}||_{L_{t}^{2}L_{x}^{2}}||R_{\nu}P_{k_{2}}Q_{>\delta
i}\psi_{2}||_{L_{t}^{2}L_{x}^{\infty}}||P_{k_{3}}\psi_{3}||_{L_{t}^{\infty}L_{x}^{\infty}}+\text{etc}\\
&\lesssim 2^{-\frac{\delta i}{2+}}\tilde{c}_{k_{1}}\\
\end{split}\end{equation}
The estimate when $P_{k_{3}}Q_{<\delta i-10}\psi_{3}$ gets
replaced by $P_{k_{3}}Q_{\geq \delta i-10}\psi_{3}$ is more of the
same. Moreover, we have by assumption
\begin{equation}\nonumber\begin{split}
&||P_{k_{0}}Q_{\geq
k_{0}}\nabla_{x,t}\chi_{i}(t)[P_{k_{1}}\delta\psi_{1}\nabla^{-1}(1-I)[P_{k_{2}}Q_{>\delta
i}\psi_{2}, P_{k_{3}}\psi_{3}]]||_{\dot{X}_{k_{0}}^{-\frac{1}{2},-1,1}}\\
&\lesssim
2^{-\frac{k_{0}}{2}}||P_{k_{1}}\delta\psi_{1}||_{L_{t}^{M}L_{x}^{2+}}||R_{\nu}P_{k_{2}}Q_{>\delta
i}\psi_{2}||_{L_{t}^{\infty}L_{x}^{M}}||P_{k_{3}}\psi_{3}||_{L_{t}^{2+}L_{x}^{\infty}},\\
\end{split}\end{equation}
which leads to an acceptable estimate. Now we estimate
\begin{equation}\nonumber\begin{split}
&||P_{k_{0}}Q_{<k_{0}}\nabla_{x,t}\chi_{i}(t)[P_{k_{1}}Q_{<k-10}\delta\psi_{1}\\&\hspace{3cm}\nabla^{-1}(1-I)P_{k}[P_{k_{2}}Q_{<\delta
i}\psi_{2}, P_{k_{3}}Q_{<\delta
i-10}\psi_{3}]]||_{\dot{X}_{k_{0}}^{-1,-\frac{1}{2},1}}\\
&\lesssim  2^{2\delta
i}||P_{k_{1}}\delta\psi_{1}||_{L_{t}^{\infty}L_{x}^{2}}||P_{k_{2}}Q_{<\delta
i}\chi_{i}(t)\psi_{2}||_{L_{t}^{M}L_{x}^{\infty}}||P_{k_{3}}Q_{<\delta
i-10}\psi_{3}||_{L_{t}^{2+}L_{x}^{\infty}}\\
&\lesssim 2^{(2\delta -\frac{1}{2+})i}\tilde{c}_{k_{1}}\\
\end{split}\end{equation}
(I.d): {\it{None of (I.a), (I.b), (I.c) hold.}} In this case, we
may treat all of $k,k_{i}$, $i=1,2,3$ as $O(1)$. The exponential
gain in $i$ is again obtained as in the preceding case. This
concludes the treatment of case (I).\\
(II): {\it{The 2nd term.}} This term is significantly simpler than
the preceding one: note that if at least one of $|k|$,$|k_{i}|$,
$i=0,\ldots 3$, is of size at least comparable to $\log_{2}T_{0}$,
one gets an exponential gain in $T_{0}$ from lemma~\ref{auxiliary
envelope} in conjunction with lemma~\ref{simple trilinear} and the
calculations in (I), provided $P_{k_{1}}\psi$ is of first type. If
it is of 2nd type, and $P_{k_{3}}\psi_{3}$ of 2nd type as well,
one also argues as in (I). If $P_{k_{1}}\psi_{1}$ is of 2nd type,
but $P_{k_{3}}\psi_{3}$ of first type, one uses the estimate
\begin{equation}\nonumber\begin{split}
&||P_{k}[R_{\beta}P_{k_{2}}Q_{<k_{2}}\delta\psi_{2}R_{j}P_{k_{3}}Q_{<k_{3}}\psi_{3}-R_{j}P_{k_{2}}Q_{<k_{2}}\delta\psi_{2}R_{\beta}P_{k_{3}}Q_{<k_{3}}\psi_{3}]||_{L_{t}^{\infty}L_{x}^{2}}
\lesssim 2^{k}\\
\end{split}\end{equation}
and (with a similar estimate when $Q_{<k_{3}}$ is replaced by
$Q_{\geq k_{3}}$)
\begin{equation}\label{technical6}\begin{split}
&||P_{k}[R_{\beta}P_{k_{2}}Q_{\geq
k_{2}}\delta\psi_{2}R_{j}P_{k_{3}}Q_{<k_{3}}\psi_{3}-R_{j}P_{k_{2}}Q_{\geq
k_{2}}\delta\psi_{2}R_{\beta}P_{k_{3}}Q_{<k_{3}}\psi_{3}]||_{L_{t}^{2}L_{x}^{2}}
\\&\lesssim
2^{\frac{\min\{k,k_{2},k_{3}\}}{2}}2^{-\delta|k_{2}-k_{3}|}\frac{\tilde{c}_{k_{2}}}{\epsilon}\tilde{c}_{k_{3}}\\
\end{split}\end{equation}
Using the first of these, one gets for example when
$k_{0}=k_{1}+O(1)$
\begin{equation}\nonumber\begin{split}
&||\nabla_{x,t}P_{k_{0}}Q_{>k_{0}}\chi(\frac{t}{T_{0}})[P_{k_{1}}\psi_{1}\triangle^{-1}\sum_{j=1,2}\partial_{j}(1-I)P_{k}[R_{\beta}P_{k_{2}}Q_{<k_{2}}\delta\psi_{2}R_{j}P_{k_{3}}Q_{<k_{3}}\psi_{3}\\&\hspace{5.5cm}-R_{j}P_{k_{2}}Q_{<k_{2}}\delta\psi_{2}R_{\beta}P_{k_{3}}Q_{<k_{3}}\psi_{3}]]||_{\dot{X}_{k_{0}}^{-\frac{1}{2},-1,2}}
\\&\lesssim
2^{-\frac{k_{0}}{2}}||P_{k_{1}}\psi_{1}||_{L_{t}^{2}L_{x}^{2+}}||\triangle^{-1}\sum_{j=1,2}\partial_{j}(1-I)P_{k}[R_{\beta}P_{k_{2}}Q_{<k_{2}}\delta\psi_{2}R_{j}P_{k_{3}}Q_{<k_{3}}\psi_{3}\\&\hspace{6.5cm}-R_{j}P_{k_{2}}Q_{<k_{2}}\delta\psi_{2}R_{\beta}P_{k_{3}}Q_{<k_{3}}\psi_{3}]||_{L_{t}^{\infty}L_{x}^{M}}\\
\end{split}\end{equation}
which is bounded by
\begin{equation}\nonumber
2^{(k-k_{0})(1-\epsilon)}2^{\delta[\min\{k_{2},k_{3}\}-\max\{k_{2},k_{3}\}]}\tilde{c}_{k_{1}}\frac{\tilde{c}_{k_{2}}}{\epsilon}\frac{\tilde{c}_{k_{3}}}{\epsilon}
\end{equation}
Assuming $\max\{|k|,|k_{i}\}$ to be at least comparable to $\delta
\log_{2}T_{0}$ and summing over $k,k_{2,3}$ results thus in the
estimate  $\tilde{c}_{k_{1}}T_{0}^{-\mu}$. Using the 2nd of the
above inequalities, \eqref{technical6}, and placing
$P_{k_{1}}\psi_{1}$ into $L_{t}^{\infty}L_{x}^{2}$ results in a
similar estimate provided $k_{0}=k_{1}+O(1)$. The remaining
frequency interactions $k_{0}<<k_{1}$ etc. are handled similarly,
as well as the contribution when $Q_{>k_{0}}$ is replaced by
$Q_{\leq k_{0}}$. Thus assume now that all the occuring
frequencies $|k|,\,|k_{i}|$, $i=0,\ldots,3$ are of size $<\delta
\log_{2}T_{0}$. Then we use Proposition~\ref{Chr1} directly in
conjunction with 3.4(b) resp. lemma~\ref{technical1}, to get
\begin{equation}\nonumber\begin{split}
&||\nabla_{x,t}P_{k_{0}}Q_{>k_{0}}\chi(\frac{t}{T_{0}})[P_{k_{1}}\psi_{1}\triangle^{-1}\sum_{j=1,2}\partial_{j}(1-I)P_{k}[R_{\beta}P_{k_{2}}\delta\psi_{2}R_{j}P_{k_{3}}\psi_{3}\\&\hspace{7.5cm}-R_{j}P_{k_{2}}\delta\psi_{2}R_{\beta}P_{k_{3}}\psi_{3}]]||_{\dot{X}_{k_{0}}^{-\frac{1}{2},-1,2}}\\&
\lesssim
||P_{k_{1}}\chi(\frac{t}{T_{0}})\psi_{1}||_{L_{t}^{\infty}L_{x}^{\infty}}||\triangle^{-1}\sum_{j=1,2}\partial_{j}(1-I)P_{k}[R_{\beta}P_{k_{2}}\delta\psi_{2}R_{j}P_{k_{3}}\psi_{3}\\&\hspace{7.5cm}-R_{j}P_{k_{2}}\delta\psi_{2}R_{\beta}P_{k_{3}}\psi_{3}]||_{L_{t}^{2}L_{x}^{2}}\\
&\lesssim
T_{0}^{-\frac{1}{2+}}[\tilde{c}_{k_{3}}+\tilde{c}_{k_{2}}]\\
\end{split}\end{equation}
But using the definition of frequency envelope we have
$\tilde{c}_{k_{2,3}}\lesssim
T_{0}^{\sigma\delta}\tilde{c}_{k_{1}}$, so we arrive at an
acceptable estimate upon summing over the admissible frequency
ranges. If we replace $Q_{>k_{0}}$ by $Q_{\leq k_{0}}$, we can for
example first reduce $P_{k_{1}}\psi_{1}$ to modulation
$<2\delta\log_{2}T_{0}$, then reduce both of
$P_{k_{2,3}}\psi_{2,3}$ to modulation $<3\delta\log_{2}T_{0}$, and
finally estimate
\begin{equation}\nonumber\begin{split}
&||\nabla_{x,t}P_{k_{0}}Q_{\leq
k_{0}}\chi(\frac{t}{T_{0}})[P_{k_{1}}Q_{k-10<.<2\delta\log_{2}T_{0}}\psi_{1}\\&\triangle^{-1}\sum_{j=1,2}\partial_{j}(1-I)P_{k}[R_{\beta}P_{k_{2}}Q_{<3\delta\log_{2}T_{0}}\delta\psi_{2}R_{j}P_{k_{3}}Q_{<3\delta\log_{2}T_{0}}\psi_{3}\\&\hspace{4.5cm}-R_{j}P_{k_{2}}Q_{<3\delta\log_{2}T_{0}}\delta\psi_{2}R_{\beta}P_{k_{3}}Q_{<3\delta\log_{2}T_{0}}\psi_{3}]]||_{L_{t}^{1}\dot{H}^{-1}}\\&
\lesssim
||P_{k_{1}}Q_{k-10<.<2\delta\log_{2}T_{0}}\psi_{1}||_{L_{t}^{2}L_{x}^{M}}||\nabla_{x,t}\nabla^{-1}P_{k_{2}}Q_{<3\delta\log_{2}T_{0}}\delta\psi_{2}||_{L_{t}^{M}L_{x}^{2+}}\\&\hspace{7cm}||\nabla_{x,t}\nabla^{-1}P_{k_{3}}Q_{<3\delta\log_{2}T_{0}}\psi_{3}||_{L_{t}^{2+}L_{x}^{\infty}}\\
\end{split}\end{equation}
Arranging that $\frac{1}{2}-\frac{1}{2+}>>\delta$ and using the
usual properties of the frequency envelope easily results in the
desired bound. If one replaces
$P_{k_{1}}Q_{k-10<.<2\delta\log_{2}T_{0}}\psi_{1}$ by
$P_{k_{1}}Q_{\leq k}\psi_{1}$, one can estimate the output with
respect to $||.||_{\dot{X}^{-1,-\frac{1}{2},1}}$ in the same
manner. This finishes case (II) and thereby the large modulation
case {\bf{(A)}}.
\\

{\bf{(B): The small modulation case.}} We now study the
expressions
\begin{equation}\nonumber
(I):\,
\partial^{\beta}[\delta\psi_{\alpha}\triangle^{-1}\sum_{j=1,2}\partial_{j}I[R_{\beta}\psi_{2}R_{j}\psi_{3}-R_{j}\psi_{2}R_{\beta}\psi_{3}]]
\end{equation}
\begin{equation}\nonumber
(II):\,\partial^{\beta}[\psi_{\alpha}\triangle^{-1}\sum_{j=1,2}\partial_{j}I[R_{\beta}\delta\psi_{2}R_{j}\psi_{3}-R_{j}\delta\psi_{2}R_{\beta}\psi_{3}]],
\end{equation}
as well as the analogous expressions
\begin{equation}\nonumber
(III):\,
\partial^{\beta}[\delta\psi_{\beta}\triangle^{-1}\sum_{j=1,2}\partial_{j}I[R_{\beta}\psi_{2}R_{j}\psi_{3}-R_{j}\psi_{2}R_{\beta}\psi_{3}]]
\end{equation}
\begin{equation}\nonumber
(IV):\,
\partial^{\beta}[\psi_{\beta}\triangle^{-1}\sum_{j=1,2}\partial_{j}I[R_{\beta}\delta\psi_{2}R_{j}\psi_{3}-R_{j}\delta\psi_{2}R_{\beta}\psi_{3}]]
\end{equation}
\begin{equation}\nonumber
(V):\,
\partial_{\alpha}[\delta\psi^{\nu}\triangle^{-1}\sum_{j=1,2}\partial_{j}I[R_{\nu}\psi_{2}R_{j}\psi_{3}-R_{j}\psi_{2}R_{\nu}\psi_{3}]]
\end{equation}
\begin{equation}\nonumber
(VI):\,
\partial_{\alpha}[\psi^{\nu}\triangle^{-1}\sum_{j=1,2}\partial_{j}I[R_{\nu}\delta\psi_{2}R_{j}\psi_{3}-R_{j}\delta\psi_{2}R_{\nu}\psi_{3}]]
\end{equation}
As is, these terms cannot yet be well estimated, and we need to
further decompose the first input into a gradient part and
elliptic error term: thus for example we write in term (I)
\begin{equation}\nonumber
\delta\psi_{\alpha}=R_{\alpha}\delta\psi+\delta\chi_{\alpha}
\end{equation}
Relegating the error terms involving $\chi_{\alpha}$ until later,
we substitute $R_{\alpha}\delta\psi$ for $\delta\psi_{\alpha}$,
and similarly for the other terms (II)-(VI).
We commence with the sum of first and fifth term in the list:\\
(I): {\bf{(I+V):}} As in the large modulation case, we shall have
to consider various types of frequency interactions. We also
reiterate the decomposition
\begin{equation}\nonumber\begin{split}
&\chi(\frac{t}{T_{0}})\partial^{\beta}[R_{\alpha}\delta\psi\triangle^{-1}\sum_{j=1,2}\partial_{j}I[R_{\beta}\psi_{2}R_{j}\psi_{3}-R_{j}\psi_{2}R_{\beta}\psi_{3}]]\\
&+\chi(\frac{t}{T_{0}})\partial_{\alpha}[R^{\nu}\delta\psi\triangle^{-1}\sum_{j=1,2}\partial_{j}I[R_{\nu}\psi_{2}R_{j}\psi_{3}-R_{j}\psi_{2}R_{\nu}\psi_{3}]]\\
&=\sum_{i>\log_{2}T_{0}}\chi_{i}(t)[\partial^{\beta}[R_{\alpha}\delta\psi\triangle^{-1}\sum_{j=1,2}\partial_{j}I[R_{\beta}\psi_{2}R_{j}\psi_{3}-R_{j}\psi_{2}R_{\beta}\psi_{3}]]\\
&+\partial_{\alpha}[R^{\nu}\delta\psi\triangle^{-1}\sum_{j=1,2}\partial_{j}I[R_{\nu}\psi_{2}R_{j}\psi_{3}-R_{j}\psi_{2}R_{\nu}\psi_{3}]]]\\
\end{split}\end{equation}
We frequency-localize this to obtain the following expression:
\begin{equation}\nonumber\begin{split}
&\chi_{i}(t)P_{k_{0}}[\partial^{\beta}[R_{\alpha}P_{k_{1}}\delta\psi_{1}P_{k}\triangle^{-1}\sum_{j=1,2}\partial_{j}I[R_{\beta}P_{k_{2}}\psi_{2}R_{j}P_{k_{3}}\psi_{3}-R_{j}P_{k_{2}}\psi_{2}R_{\beta}P_{k_{3}}\psi_{3}]]\\
&+\partial_{\alpha}[R^{\nu}P_{k_{1}}\delta\psi_{1}P_{k}\triangle^{-1}\sum_{j=1,2}\partial_{j}I[R_{\nu}P_{k_{2}}\psi_{2}R_{j}P_{k_{3}}\psi_{3}-R_{j}P_{k_{2}}\psi_{2}R_{\nu}P_{k_{3}}\psi_{3}]]]\\
\end{split}\end{equation}
Now we subdivide into the following possibilities:\\
(I.a): {\it{One of the following options hold:}} $i\lesssim
|k_{2}|$, $i\lesssim |k_{3}|$, $i\lesssim |k_{0}-k_{1}|$,
$i\lesssim\min\{|k-k_{1}|,|k-k_{2}|\}$. In this case, we obtain
the desired estimate involving an exponential gain in $i$ from
3.4(c) as well as lemma~\ref{auxiliary envelope} if both
$P_{k_{2,3}}\psi_{2,3}$ are of the first type. If at least one of
them is of the 2nd type, this is again straightforward due to the
strong estimates satisfied by these: then we have
\begin{equation}\nonumber\begin{split}
&||P_{k}\triangle^{-1}\sum_{j=1,2}\partial_{j}I[R_{\beta}P_{k_{2}}\psi_{2}R_{j}P_{k_{3}}\psi_{3}-R_{j}P_{k_{2}}\psi_{2}R_{\beta}P_{k_{3}}\psi_{3}]||_{L_{t}^{1}L_{x}^{\infty}}\\&\hspace{6cm}\lesssim
2^{\delta[\min\{k,k_{2},k_{3}\}-\max\{k,k_{2},k_{3}\}]}\frac{\tilde{c}_{k_{2}}}{\epsilon}\frac{\tilde{c}_{k_{3}}}{\epsilon}\\
\end{split}\end{equation}
\begin{equation}\nonumber\begin{split}
&||P_{k}\triangle^{-1}\sum_{j=1,2}\partial_{j}I[R_{\beta}P_{k_{2}}\psi_{2}R_{j}P_{k_{3}}\psi_{3}-R_{j}P_{k_{2}}\psi_{2}R_{\beta}P_{k_{3}}\psi_{3}]||_{L_{t}^{1}L_{x}^{2+}}\\&\hspace{5cm}\lesssim
2^{-k(1-\epsilon)}2^{\delta[\min\{k,k_{2},k_{3}\}-\max\{k,k_{2},k_{3}\}]}\frac{\tilde{c}_{k_{2}}}{\epsilon}\frac{\tilde{c}_{k_{3}}}{\epsilon}\\
\end{split}\end{equation}
From this we get
\begin{equation}\nonumber\begin{split}
&||P_{k_{0}}Q_{<k_{0}}[\partial^{\beta}[R_{\alpha}P_{k_{1}}\delta\psi_{1}\\&\hspace{2cm}P_{k}\triangle^{-1}\sum_{j=1,2}\partial_{j}I[R_{\beta}P_{k_{2}}\psi_{2}R_{j}P_{k_{3}}\psi_{3}-R_{j}P_{k_{2}}\psi_{2}R_{\beta}P_{k_{3}}\psi_{3}]]||_{L_{t}^{1}\dot{H}^{-1}}\\
&\hspace{5cm}\lesssim
2^{\delta[\min\{k,k_{0},\ldots,k_{3}\}-\max\{k,k_{0},\ldots,k_{3}\}]}\frac{1}{\epsilon^{2}}\prod_{i=1,2,3}\tilde{c}_{k_{i}}\\
\end{split}\end{equation}
The contribution when one has $P_{k_{0}}Q_{\geq k_{0}}$ in front is even simpler and left out.
Term (V) is treated by exact analogy.\\
(I.b): {\it{$i\lesssim |k|$, $i\lesssim |k_{1}|$, and none of the
properties in (I.a) hold.}} This implies $k,k_{1}\lesssim -i$. We
may and shall assume $k-k_{1}=O(1)$, $k_{0}-k_{1}=O(1)$. Also, we
may and shall assume that $|k_{2,3}|=O(1)$. We shall again treat
term (I), term (V) being treated analogously. We start out by
observing that we may assume $k>-i(1+\mu)$ for any $\mu>0$.
Indeed, assume the opposite. Again considering term (I), we have
\begin{equation}\nonumber\begin{split}
&||P_{k_{0}}Q_{<k_{0}}\chi_{i}(t)[\partial^{\beta}Q_{>k_{0}+100}[R_{\alpha}P_{k_{1}}\delta\psi_{1}P_{k}\\&\hspace{3cm}\triangle^{-1}\sum_{j=1,2}\partial_{j}I[R_{\beta}P_{k_{2}}\psi_{2}R_{j}P_{k_{3}}\psi_{3}-R_{j}P_{k_{2}}\psi_{2}R_{\beta}P_{k_{3}}\psi_{3}]]||_{L_{t}^{1}\dot{H}^{-1}}\\
&\lesssim \sum_{k_{0}+100\leq j\leq -i+O(1)}||P_{k_{0}}Q_{<k_{0}}\chi_{i}(t)[\partial^{\beta}Q_{j}[R_{\alpha}P_{k_{1}}Q_{j+O(1)}\delta\psi_{1}P_{k}\\&\hspace{3cm}\triangle^{-1}\sum_{j=1,2}\partial_{j}I[R_{\beta}P_{k_{2}}\psi_{2}R_{j}P_{k_{3}}\psi_{3}-R_{j}P_{k_{2}}\psi_{2}R_{\beta}P_{k_{3}}\psi_{3}]]||_{L_{t}^{1}\dot{H}^{-1}}\\
&\lesssim
\sum_{k_{0}+100\leq j\leq -i+O(1)}2^{j}||\chi_{i}(t)||_{L_{t}^{2}}||R_{\alpha}P_{k_{1}}Q_{j+O(1)}\delta\psi_{1}||_{L_{t}^{2}L_{x}^{2}}\\&\hspace{7cm}||P_{k}[\nabla^{-1}P_{k_{2}}\psi_{2}P_{k_{3}}\psi_{3}]||_{L_{t}^{\infty}L_{x}^{\infty}}\\
\end{split}\end{equation}
One checks that this is estimated by $\lesssim
\tilde{c}_{k_{1}}2^{k_{0}}$, which is acceptable. Similarly, we
estimate
\begin{equation}\nonumber\begin{split}
&||P_{k_{0}}Q_{<k_{0}}\chi_{i}(t)[\partial^{\beta}Q_{\leq k_{0}+100}[R_{\alpha}P_{k_{1}}\delta\psi_{1}\\&\hspace{2cm}P_{k}\triangle^{-1}\sum_{j=1,2}\partial_{j}I[R_{\beta}P_{k_{2}}\psi_{2}R_{j}P_{k_{3}}\psi_{3}-R_{j}P_{k_{2}}\psi_{2}R_{\beta}P_{k_{3}}\psi_{3}]]||_{L_{t}^{1}\dot{H}^{-1}}\\
&\lesssim||\chi_{i}(t)||_{L_{t}^{2}}||P_{k_{1}}\delta\psi_{1}||_{L_{t}^{\infty}L_{x}^{\infty}}\\&\hspace{2cm}||P_{k}\triangle^{-1}\sum_{j=1,2}\partial_{j}I[R_{\beta}P_{k_{2}}\psi_{2}R_{j}P_{k_{3}}\psi_{3}-R_{j}P_{k_{2}}\psi_{2}R_{\beta}P_{k_{3}}\psi_{3}]||_{L_{t}^{2}L_{x}^{2}},\\
\end{split}\end{equation}
which, upon using 3.4(b) as well as lemma~\ref{technical1}, can be
estimated by $\lesssim \tilde{c}_{k_{0}}2^{\frac{k_{0}+i}{2}}$,
which is acceptable.\\
The estimates when $Q_{<k_{0}}$ is replaced by $Q_{\geq k_{0}}$
are similar. Thus we now assume that $|k|<i(1+\mu)$. Arguing as in
case {\bf{(A)}}(I.b), we may replace $P_{k}I$ be
$P_{k}Q_{<\frac{k}{2}}$ while only generating acceptable error
terms. We may then move the multiplier $\chi_{i}(t)$ past the
operator $P_{k}\nabla^{-1}Q_{<\frac{k}{2}}$, and transform the
latter back into $\nabla^{-1}P_{k}I$ innocuously. Proceeding as in
\cite{Kr-4}, we intend to exploit the null-structure of the
expression. Before being able to do so, we need to effect a few
more reductions: we need to reduce the modulation of the output
and first input to size $<2^{k_{0,1}}$, respectively. For this,
note that
\begin{equation}\nonumber\begin{split}
&||P_{k_{0}}Q_{\geq k_{0}+100}[\partial^{\beta}[R_{\alpha}P_{k_{1}}\delta\psi_{1}P_{k}\\&\hspace{1.5cm}\triangle^{-1}\sum_{j=1,2}\partial_{j}I[R_{\beta}P_{k_{2}}\chi_{i}(t)\psi_{2}R_{j}P_{k_{3}}\psi_{3}-R_{j}P_{k_{2}}\chi_{i}(t)\psi_{2}R_{\beta}P_{k_{3}}\psi_{3}]]||_{\dot{X}_{k_{0}}^{-\frac{1}{2},-1,2}}\\
&\lesssim ||P_{k_{0}}Q_{\geq k_{0}+100}[\partial^{\beta}[R_{\alpha}P_{k_{1}}Q_{\geq k_{1}}\delta\psi_{1}P_{k}\\&\hspace{1.5cm}\triangle^{-1}\sum_{j=1,2}\partial_{j}I[R_{\beta}P_{k_{2}}\chi_{i}(t)\psi_{2}R_{j}P_{k_{3}}\psi_{3}-R_{j}P_{k_{2}}\chi_{i}(t)\psi_{2}R_{\beta}P_{k_{3}}\psi_{3}]]||_{\dot{X}_{k_{0}}^{-\frac{1}{2},-1,2}}\\
&\lesssim 2^{-\frac{k_{0}}{2}}||R_{\alpha}P_{k_{1}}Q_{\geq
k_{1}}\delta\psi_{1}||_{L_{t}^{2}L_{x}^{\infty}}||P_{k}[\nabla^{-1}P_{k_{2}}\psi_{2}P_{k_{3}}R_{\beta}\psi_{3}]||_{L_{t}^{\infty}L_{x}^{2}}\\
&\lesssim
2^{k}\tilde{c}_{k_{1}}[\frac{\tilde{c}_{k_{2}}}{\epsilon}+\frac{\tilde{c}_{k_{3}}}{\epsilon}].\\
\end{split}\end{equation}
We can also reduce $P_{k_{2,3}}\psi_{2,3}$ to modulation
$<2^{\max\{k_{2},k_{3}\}+O(1)}$. For this, note that
\begin{equation}\nonumber\begin{split}
&||P_{k}I[\nabla^{-1}P_{k_{2}}\psi_{2}P_{k_{3}}Q_{>\max\{k_{2},k_{3}\}+100}\psi_{3}]||_{L_{t}^{1}L_{x}^{2}}\\
&\lesssim
\sum_{j>\max\{k_{2},k_{3}\}+100}2^{k}||\nabla^{-1}P_{k_{2}}Q_{j+O(1)}\psi_{2}||_{L_{t}^{2}L_{x}^{2}}||P_{k_{3}}Q_{j}\psi_{3}||_{L_{t}^{2}L_{x}^{2}}\\
&\lesssim\sum_{j>\max\{k_{2},k_{3}\}+100}2^{k}2^{-(2-2\mu)j}2^{(1-2\mu)\max\{k_{2},k_{3}\}}\\&\hspace{4cm}||\nabla^{-1}P_{k_{2}}\psi_{2}||_{\dot{X}_{k_{2}}^{-(\frac{1}{2}-\mu),1-\mu,1}}||P_{k_{3}}\psi_{3}||_{\dot{X}_{k_{3}}^{-(\frac{1}{2}-\mu),1-\mu,1}}\\
&\lesssim
2^{-k_{2}}\frac{\tilde{c}_{k_{2}}}{\epsilon}\frac{\tilde{c}_{k_{3}}}{\epsilon}\\
\end{split}\end{equation}
From this one deduces that for $a=\max\{k_{2},k_{3}\}+100$
\begin{equation}\nonumber\begin{split}
&||P_{k_{0}}Q_{<k_{0}+O(1)}[\partial^{\beta}[R_{\alpha}P_{k_{1}}\delta\psi_{1}\\&P_{k}\triangle^{-1}\sum_{j=1,2}\partial_{j}I[R_{\beta}P_{k_{2}}\chi_{i}(t)\psi_{2}R_{j}P_{k_{3}}Q_{>a}\psi_{3}-R_{j}P_{k_{2}}\chi_{i}(t)\psi_{2}R_{\beta}P_{k_{3}}Q_{>a}\psi_{3}]]||_{L_{t}^{1}\dot{H}^{-1}}\\
&\lesssim
||R_{\alpha}P_{k_{1}}Q_{<k_{1}+O(1)}\delta\psi_{1}||_{L_{t}^{\infty}L_{x}^{\infty}}||P_{k}I[\nabla^{-1}P_{k_{2}}\psi_{2}P_{k_{3}}Q_{>\max\{k_{2},k_{3}\}+100}\psi_{3}]||_{L_{t}^{1}L_{x}^{2}}\\
&\lesssim 2^{k-k_{2}}\tilde{c}_{k_{1}}\frac{\tilde{c}_{k_{2}}}{\epsilon}\frac{\tilde{c}_{k_{3}}}{\epsilon}\\
\end{split}\end{equation}
One can sum over the appropriate range of $k_{2,3}$, deducing the
desired estimate. We shall always assume these reductions of
modulation, but sometimes omit them to simplify notation. Now we
expand the null-structure as in \cite{Kr-3}, \cite{Kr-4}:
schematically we have
\begin{equation}\label{null-form1}\begin{split}
&2\sum_{j=1}^{2}\triangle^{-1}\partial_{j}[R_{\nu}fR_{j}g-R_{j}f
R_{\nu}g]\partial^{\nu}h \\
&=\sum_{j=1}^{2}\Box[\triangle^{-1}\partial_{j}[\nabla^{-1}fR_{j}g]h]-\sum_{j=1}^{2}\Box\triangle^{-1}\partial_{j}[\nabla^{-1}fR_{j}g]h\\
&-\sum_{j=1}^{2}\triangle^{-1}\partial_{j}[\nabla^{-1}fR_{j}g]\Box
h-\nabla^{-1}f\Box((\nabla^{-1}g)
h)\\&+\nabla^{-1}f\Box(\nabla^{-1}g)
h+\nabla^{-1}f(\nabla^{-1}g)\Box h,\\
\end{split}\end{equation}
\begin{equation}\label{null-form2}\begin{split}
&\sum_{j=1}^{2}\triangle^{-1}\partial_{j}\partial^{\nu}[R_{\nu}f
R_{j}g-R_{j}f R_{\nu}g]h\\
&=\sum_{j=1}^{2}[\triangle^{-1}\partial_{j}\Box[\nabla^{-1}f
g]h-\frac{1}{2}\Box[\nabla^{-1}f\nabla^{-1}g]h\\&+\frac{1}{2}\Box\nabla^{-1}f
\nabla^{-1}g h-\frac{1}{2}\nabla^{-1}f\Box g h]\\
\end{split}\end{equation}
The first of these identities is useful when the outer derivative
$\partial^{\beta}$ falls on the first input
$R_{\alpha}\delta\psi_{1}$. The 2nd is useful provided the outer
derivative lands on the inner square bracket. We shall treat each
of these terms. Clearly, the terms in the 2nd expansion are almost identical to the ones in the first.
We treat the first in detail, the 2nd being treated similarly.\\
(I.b.1) {\it{The first term in the expansion.}} This is the
expression
\begin{equation}\nonumber\begin{split}
&P_{k_{0}}Q_{<k_{0}}\Box
[R_{\alpha}P_{k_{1}}Q_{<k_{0}}\delta\psi_{1}P_{k}\triangle^{-1}\sum_{j=1,2}\partial_{j}I[\nabla^{-1}P_{k_{2}}\chi_{i}(t)\psi_{2}R_{j}P_{k_{3}}\psi_{3}]]\\
\end{split}\end{equation}
This is straightforward to estimate: we have
\begin{equation}\nonumber\begin{split}
&||P_{k_{0}}Q_{<k_{0}}\Box
[R_{\alpha}P_{k_{1}}Q_{<k_{0}}\delta\psi_{1}P_{k}\triangle^{-1}\sum_{j=1,2}\partial_{j}I[\nabla^{-1}P_{k_{2}}\chi_{i}(t)\psi_{2}R_{j}P_{k_{3}}\psi_{3}]]||_{\dot{X}_{k_{0}}^{-1-\frac{1}{2},1}}\\
&\lesssim
2^{-\frac{k_{0}}{2}}||R_{\alpha}P_{k_{1}}Q_{<k_{0}}\delta\psi_{1}||_{L_{t}^{\infty}L_{x}^{M}}||\nabla^{-1}P_{k_{2}}\chi_{i}(t)\psi_{2}||_{L_{t}^{2+}L_{x}^{\infty}}||P_{k_{3}}\psi_{3}||_{L_{t}^{M}L_{x}^{2+}}
\lesssim 2^{\frac{k_{1}}{2+}}\tilde{c}_{k_{1}}\\
\end{split}\end{equation}
(I.b.2) {\it{The 2nd term in the expansion.}} This is the
expression
\begin{equation}\nonumber\begin{split}
P_{k_{0}}Q_{<k_{0}}
[R_{\alpha}P_{k_{1}}Q_{<k_{0}}\delta\psi_{1}P_{k}\triangle^{-1}\sum_{j=1,2}\partial_{j}\Box I[\nabla^{-1}P_{k_{2}}\chi_{i}(t)\psi_{2}R_{j}P_{k_{3}}\psi_{3}]]\\
\end{split}\end{equation}
This turns out to be significantly more complicated. The reason
for this is that we need to exploit the bilinear inequality
3.4(g); using Strichartz type norms here appears to result in a
loss in the low frequencies, or in $i$. The only way we can
possibly squeeze out a small gain in $i$ is to exploit the
temporal cutoff $\chi_{i}(t)$ applied to $P_{k_{2}}\psi_{2}$. This
is a non-trivial task on account of the fact that the only way to
place the inner bracket $[,]$ into $L_{t}^{2}L_{x}^{2}$ appears to
involve null-frame spaces. Our main tool for this is the following
lemma:
\begin{lemma}\label{key} Let $C$ be a sufficiently large number.
The following limits hold: for any $k$ with
$|k|<\frac{\epsilon}{1000 C}i$, and arbitrary $\epsilon>0$,
\begin{equation}\nonumber
\lim_{i\rightarrow\infty}||P_{k}Q_{[k-(1-\epsilon)i-C,k]}[\chi_{i}(t)\psi_{\nu}]||_{\dot{X}_{k}^{0,\frac{1}{2},\infty}}=0
\end{equation}
More precisely, for appropriate $\mu(\epsilon)>0$, we have
\begin{equation}\nonumber
||P_{k}Q_{[k-(1-\epsilon)i-C,k]}[\chi_{i}(t)\psi_{\nu}]||_{\dot{X}_{k}^{0,\frac{1}{2},\infty}}\lesssim
2^{-\mu(\epsilon)i}
\end{equation}
Next, denote by $\chi_{i,\kappa}^{c}$ a smooth bump function which
localizes to the complement in the $2^{\epsilon i}/C$-neighborhood
of the (physical) light cone of a slab of length $2^{i-1}$
centered at time $t=2^{i}$, angular opening $2\kappa$ with
$|\kappa|\sim 2^{i(\frac{\epsilon-1}{2})}$ and distance $\lesssim
2^{\epsilon i}/C$ from the light cone. Then with the same
assumption on $k$,
\begin{equation}\nonumber\begin{split}
&\lim_{i\rightarrow\infty}\sum_{\pm}\sup_{l\in
[\frac{\epsilon-1}{2}i,-10]}2^{-\frac{l}{2}}(\sum_{\tilde{\kappa}\in
K_{l}}||P_{k,\tilde{\kappa}}\sum_{\kappa\in
K_{\frac{\epsilon-1}{2}i}}\chi_{i,\mp\kappa}^{c}(t,x)\\&\hspace{4cm}P_{k,\kappa}Q^{\pm}_{<k+i(\epsilon-1)-C}[\chi_{i}(t)\psi_{\nu}]||_{PW[\pm\tilde{\kappa}]}^{2})^{\frac{1}{2}}=0
\end{split}\end{equation}
More precisely, this quantity decays like $2^{-\mu i}$ for
suitable $\mu>0$.
\end{lemma}
\begin{proof}: We first estimate
\begin{equation}\nonumber
||\nabla_{x,t}P_{k}Q_{>k-(1-\epsilon)i}[\chi_{i}(t)\frac{{\bf{x}}}{{\bf{y}}}]||_{A[k]},\,
||\nabla_{x,t}P_{k}Q_{>k-(1-\epsilon)i}[\chi_{i}(t)\ln{\bf{y}}]||_{A[k]}
\end{equation}
for arbitrary $k\in{\mathbf{Z}}$ with $|k|<\frac{\epsilon}{C}i$.
Observe that
\begin{equation}\nonumber
\Box[\chi_{i}(t)\frac{{\bf{x}}}{{\bf{y}}}]=\chi_{i}''(t)\frac{{\bf{x}}}{{\bf{y}}}
+\chi_{i}'(t)\partial_{t}(\frac{{\bf{x}}}{{\bf{y}}})+\chi_{i}(t)\Box[\frac{{\bf{x}}}{{\bf{y}}}]
\end{equation}
Therefore, we obtain for $j\in [k-(1-\epsilon)i,k]$:
\begin{equation}\nonumber\begin{split}
&||P_{k}Q_{j}\nabla_{x}\chi_{i}(t)\frac{{\bf{x}}}{{\bf{y}}}||_{A[k]}\lesssim
2^{-\frac{j}{2}}[||P_{k}(\chi_{i}''(t)\frac{{\bf{x}}}{{\bf{y}}})||_{L_{t}^{2}L_{x}^{2}}+||P_{k}(\chi_{i}'(t)\partial_{t}(\frac{{\bf{x}}}{{\bf{y}}}))||_{L_{t}^{2}L_{x}^{2}}
\\&\hspace{8cm}+||P_{k}Q_{j}(\chi_{i}(t)\Box[\frac{{\bf{x}}}{{\bf{y}}}])||_{L_{t}^{2}L_{x}^{2}}]\\
\end{split}\end{equation}
The first two terms on the right hand side are elementary to
estimate, using Holder's inequality, finite propagation speed and
the energy inequality:
\begin{equation}\nonumber\begin{split}
&||P_{k}(\chi_{i}''(t)\frac{{\bf{x}}}{{\bf{y}}})||_{L_{t}^{2}L_{x}^{2}}\lesssim
||\chi_{i}''(t)||_{L_{t}^{2}}||P_{k}(\frac{{\bf{x}}}{{\bf{y}}})||_{L_{t}^{\infty}L_{x}^{2}}\lesssim
2^{\frac{i}{2}}2^{-2i}2^{i}\lesssim 2^{-\frac{i}{2}}\\
\end{split}\end{equation}
\begin{equation}\nonumber
||P_{k}(\chi_{i}'(t)\partial_{t}(\frac{{\bf{x}}}{{\bf{y}}}))||_{L_{t}^{2}L_{x}^{2}}\lesssim
||\chi_{i}'(t)||_{L_{t}^{2}}||\partial_{t}(\frac{{\bf{x}}}{{\bf{y}}})||_{L_{t}^{\infty}L_{x}^{2}}.
\lesssim 2^{-\frac{i}{2}}
\end{equation}
We conclude that the contribution from these terms is at most
\begin{equation}\nonumber
\lesssim
\sum_{j>k-i(1-\epsilon)}2^{-\frac{j}{2}}2^{-\frac{i}{2}}\lesssim
2^{-\frac{\epsilon}{2+}i}.
\end{equation}
Proceeding to the last term above, we use the null-structure in
\eqref{Euler2} as well as Proposition~\ref{Chr1} to get
\begin{equation}\nonumber\begin{split}
&||P_{k}(\chi_{i}(t)\Box(\frac{{\bf{x}}}{{\bf{y}}})||_{L_{t}^{2}L_{x}^{2}}\\&\lesssim
||\chi_{i}(t)\frac{(\partial_{t}+\partial_{r}){\bf{x}}}{{\bf{y}}}||_{L_{t}^{2}L_{x}^{\infty}}
||\chi_{i}(t)\frac{(\partial_{t}-\partial_{r}){\bf{x}}}{{\bf{y}}}||_{L_{t}^{\infty}L_{x}^{2}}+\text{etc}\lesssim
2^{-i}\\
\end{split}\end{equation}
This establishes the claim for $\frac{{\bf{x}}}{{\bf{y}}}$ because
of 3.4(d). As far as
$\partial_{t}P_{k}Q_{>k+(\epsilon-1)i}(\frac{{\bf{x}}}{{\bf{y}}})$
is concerned this only differs from the preceding as far as the
estimate for \\ $P_{k}Q_{\geq
k}R_{0}\partial_{t}(\frac{{\bf{x}}}{{\bf{y}}})$. This is treated
by using the equation for $\Box(\frac{{\bf{x}}}{{\bf{y}}})$ and
arguing as before. The estimate for $\ln{\bf{y}}$ is similar. We
now establish the 2nd inequality stated in the lemma provided
$\psi_{\nu}$ is replaced by
$\nabla_{x,t}\frac{{\bf{x}}}{{\bf{y}}}$, $\nabla_{x,t}\ln
{\bf{y}}$. First, let $\phi\in C^{\infty}({\mathbf{R}}^{2+1})$ be
a rotationally symmetric free wave with $\phi[0,x]=(0,g(x))$. We
use the representation formula
\begin{equation}\nonumber\begin{split}
&P_{k,\kappa}\phi(t,x)=c\int_{S^{1}}{\mathcal{F}}(\hat{g}(|\xi|)m_{k}(|\xi|)|\xi|^{2})a_{\kappa}(\omega)(-t+x\cdot\omega)d\omega
\\&\hspace{5cm}-c\int_{S^{1}}{\mathcal{F}}(\hat{g}(|\xi|)m_{k}(|\xi|)|\xi|^{2})a_{\kappa}(\omega)(t+x\cdot\omega)d\omega\\
\end{split}\end{equation}
where(committing abuse of notation) we wrote
$\hat{g}(\xi)=\hat{g}(|\xi|)$. Now assume that
$||(1+\nabla_{|\xi|}^{\alpha})\hat{g}(|\xi|)||_{L_{|\xi|}^{2}}<C$.
This implies that
\begin{equation}\nonumber
||\chi_{i}(.){\mathcal{F}}(\hat{g}(|\xi|))(.)||_{L^{2}}\lesssim
2^{-\alpha i}
\end{equation}
Now observe that on the support of $\chi^{c}_{i,\mp\kappa}$, we
have
\begin{equation}\nonumber
|t\pm x\cdot\omega|\geq -|t-|x||+||x|\pm x\cdot\omega|\geq
||t-|x||-|x|\sin<x,\mp\omega>^{2}|\geq c2^{i\epsilon}
\end{equation}
Next, observe that the multiplier $\chi^{c}_{i,\mp\kappa}(t,x)$
smears out the Fourier support in the angular direction by an
amount $\sim 2^{-\frac{(1+\epsilon)i}{2}}$, while $|\kappa|\sim
2^{\frac{i(\epsilon-1)}{2}}$ and $|k|<\frac{\epsilon}{C}i$. This
entails that one can include an operator $P'_{k,\kappa}$ with the
same properties as $P_{k,\kappa}$ in the expression and reason as
follows\footnote{Let $Q^{\pm}$ microlocalize to the upper or lower
half-space $\tau><0$.}:
\begin{equation}\nonumber\begin{split}
&2^{-\frac{l}{2}}||P_{k,\tilde{\kappa}}\sum_{\kappa\in
K_{\frac{(\epsilon-1)}{2}i}}\chi^{c}_{i,\mp\kappa}(t,x)P_{k,\kappa}Q^{\pm}\phi||_{PW[\pm\tilde{\kappa}]}
\\&=2^{-\frac{l}{2}}||P_{k,\tilde{\kappa}}\sum_{\kappa\in
K_{\frac{(\epsilon-1)}{2}i}}P'_{k,\kappa}\chi^{c}_{i,\mp\kappa}(t,x)P_{k,\kappa}Q^{\pm}\phi||_{PW[\pm\tilde{\kappa}]}\\
&\lesssim 2^{i\frac{1-\epsilon}{4}}(\sum_{\kappa\in
K_{\frac{(\epsilon-1)i}{2}},\kappa\cap 2\tilde{\kappa}\neq
\emptyset}||\chi^{c}_{i,\mp\kappa}(t,x)P_{k,\kappa}Q^{\pm}\phi||_{PW[\pm\kappa]}^{2})^{\frac{1}{2}}\\
\end{split}\end{equation}
Using Plancherel's theorem, we can further estimate this as
\begin{equation}\nonumber\begin{split}
&\lesssim 2^{i\frac{1-\epsilon}{4}}(\sum_{\kappa\in
K_{\frac{(\epsilon-1)i}{2}},\kappa\cap 2\tilde{\kappa}\neq
\emptyset}[\int_{\omega\in\kappa}||\chi_{t\pm
x\cdot\omega>2^{\epsilon
i}}{\mathcal{F}}(\hat{g}(|\xi|))(t\pm x\cdot\omega)||_{L^{2}}a_{\kappa}(\omega)d\omega]^{2})^{\frac{1}{2}}\\
&\lesssim 2^{-\epsilon\alpha i}(\sum_{\kappa\in
K_{\frac{(\epsilon-1)i}{2}},\kappa\cap 2\tilde{\kappa}\neq
\emptyset}[\int_{\omega\in\kappa}||(1+\nabla_{|\xi|}^{\alpha})\hat{g}(|\xi|)||_{L_{|\xi|}^{2}}a_{\kappa}(\omega)d\omega]^{2})^{\frac{1}{2}}\\
\end{split}\end{equation}
Using Cauchy Schwarz' inequality, this leads to the estimate
\begin{equation}\nonumber\begin{split}
&2^{-\frac{l}{2}}\sum_{\pm}(\sum_{\tilde{\kappa}\in
K_{l}}||P_{k,\tilde{\kappa}}\sum_{\kappa\in
K_{\frac{(\epsilon-1)}{2}i}}\chi^{c}_{i,\mp\kappa}(t,x)P_{k,\kappa}Q^{\pm}\phi||_{PW[\pm\tilde{\kappa}]}^{2})^{\frac{1}{2}}
\\&\hspace{7cm}\lesssim 2^{-\alpha\epsilon
i}||(1+\nabla_{|\xi|}^{\alpha})\hat{g}(|\xi|)||_{L_{|\xi|}^{2}}\\
\end{split}\end{equation}
Now we proceed to the inhomogeneous situation at
hand: let $S(t)$ denote the free wave propagator, i. e.
$\Box[S(t)(f,g)]=0$, $S(0)(f,g)=f$, $\partial_{t}[S(0)(f,g)]=g$.
Also, let $U(t)g=S(t)(0,f)$. Then we can write
\begin{equation}\nonumber
\frac{{\bf{x}}}{{\bf{y}}}(t,.)=
S(t)(\frac{{\bf{x}}}{{\bf{y}}}(0,.),\partial_{t}[\frac{{\bf{x}}}{{\bf{y}}}](0,.))+\int_{0}^{t}U(t-s)\Box[\frac{{\bf{x}}}{{\bf{y}}}](s,.)ds
\end{equation}
Reasoning as above, we immediately get the desired estimate for
the linear part. As concerns the inhomogeneity, we decompose the
integral as
\begin{equation}\nonumber
\int_{0}^{t}U(t-s)\Box[\frac{{\bf{x}}}{{\bf{y}}}](s,.)ds=\int_{0}^{2^{\frac{\epsilon}{C}i}}U(t-s)\Box[\frac{{\bf{x}}}{{\bf{y}}}](s,.)ds
+\int_{2^{\frac{\epsilon}{C}i}}^{t}U(t-s)\Box[\frac{{\bf{x}}}{{\bf{y}}}](s,.)ds
\end{equation}
Then we observe that from the argument given above we have for
$l\in [\frac{(\epsilon-1)i}{2},-10]$
\begin{equation}\nonumber\begin{split}
&2^{-\frac{l}{2}}(\sum_{\tilde{\kappa}\in
K_{l}}||P_{k,\tilde{\kappa}}\sum_{\kappa\in
K_{\frac{\epsilon-1}{2}i}}\chi^{c}_{i,\mp\kappa}P_{k,\kappa}Q^{\pm}[\int_{0}^{2^{\frac{\epsilon}{C}i}}U(t-s)\Box[\frac{{\bf{x}}}{{\bf{y}}}](s,.)ds]||_{PW[\pm\tilde{\kappa}]}^{2})^{\frac{1}{2}}\\
&\lesssim 2^{-\epsilon\alpha
i}2^{-k}||x^{\alpha}\Box[\frac{{\bf{x}}}{{\bf{y}}}]||_{L_{t}^{1}L_{x}^{2}}\lesssim 2^{-\frac{\epsilon}{2}\alpha i},\\
\end{split}\end{equation}
provided $\frac{1}{C}<<\alpha$ and we choose $\alpha>0$ small
enough, by the proof of Corollary~\ref{crux} and finite
propagation speed. Next, we can estimate by using 3.4(d):
\begin{equation}\nonumber\begin{split}
&2^{\frac{1-\epsilon}{4}i}(\sum_{\kappa\in
K_{\frac{\epsilon-1}{2}i}}||P_{k,\pm\kappa}Q^{\pm}_{<k+(\epsilon-1)i}\chi_{i}(t)[\int_{2^{\frac{\epsilon}{C}i}}^{t}U(t-s)\Box[\frac{{\bf{x}}}{{\bf{y}}}](s,.)ds]||_{PW[\tilde{\kappa}]}^{2})^{\frac{1}{2}}\\
&\lesssim
||\chi_{t>2^{\frac{\epsilon}{C}i}}\Box[\frac{{\bf{x}}}{{\bf{y}}}]||_{L_{t}^{1}L_{x}^{2}}\lesssim
2^{-\frac{\epsilon}{C}\alpha i}\\
\end{split}\end{equation}
Call the integral in the preceding $\psi$. From this we get
control over the 2nd more complicated norm in the lemma: using
Cauchy-Schwartz and the fact that $||.||_{PW[\kappa]}$ is
essentially unaffected by multiplication with bounded functions,
we get
\begin{equation}\label{calculation}\begin{split}
&2^{-\frac{l}{2}}(\sum_{\tilde{\kappa}\in
K_{l}}||P_{k,\tilde{\kappa}}\sum_{\kappa\in
K_{\frac{\epsilon-1}{2}i}}\chi^{c}_{i,\mp\kappa}P_{k,\kappa}Q^{\pm}_{<k+(\epsilon-1)i}\psi||_{PW[\pm\tilde{\kappa}]}^{2})^{\frac{1}{2}}\\
&\lesssim 2^{\frac{1-\epsilon}{4}i}(\sum_{\tilde{\kappa}\in
K_{l}}\sum_{\kappa\in K_{\frac{\epsilon-1}{2}i,\kappa\subset
2\tilde{\kappa}}}||P_{k,\kappa}Q^{\pm}_{<k+(\epsilon-1)i}\psi||_{PW[\pm\kappa]}^{2})^{\frac{1}{2}},\\
\end{split}\end{equation}
and we just bounded this expression. The estimate for
$\ln{\bf{y}}$ is of course analogous. Let
$\mu_{1}(\epsilon)=\frac{\epsilon}{C}\alpha $. Now we need to
transfer these statements to $\psi_{\nu}$. We recall the identity
\begin{equation}\nonumber
\psi_{\nu}=(\frac{\partial_{\nu}{\bf{x}}}{{\bf{y}}}+i\frac{\partial_{\nu}{\bf{y}}}{{\bf{y}}})e^{i\sum_{j=1,2}\triangle^{-1}\partial_{j}
(\frac{\partial_{j}{\bf{x}}}{{\bf{y}}})}
\end{equation}
We observe as usual that
$\frac{\partial_{\nu}{\bf{x}}}{{\bf{y}}}=\partial_{\nu}(\frac{{\bf{x}}}{{\bf{y}}})-\frac{{\bf{x}}}{{\bf{y}}}\frac{\partial_{\nu}{\bf{y}}}{{\bf{y}}}$.
For example consider the term
$\frac{{\bf{x}}}{{\bf{y}}}\frac{\partial_{\nu}{\bf{y}}}{{\bf{y}}}e^{i\sum_{j=1,2}\triangle^{-1}\partial_{j}
(\frac{\partial_{j}{\bf{x}}}{{\bf{y}}})}$, the other terms being
treated similarly. One expands the exponential in a Taylor series,
which results in schematic terms of the form
\begin{equation}\nonumber
a_{k}\psi\nabla^{-1}\psi_{1}\nabla^{-1}\psi_{2}\ldots\nabla^{-1}\psi_{k},
\end{equation}
where $\nabla^{-1}\psi$ stands for expressions like
$\frac{{\bf{x}}}{{\bf{y}}}$,
$\nabla^{-1}(\frac{{\bf{x}}}{{\bf{y}}}\frac{\partial_{i}{\bf{y}}}{{\bf{y}}})$,
and $\psi$ stands for either
$\frac{\partial_{\nu}{\bf{y}}}{{\bf{y}}}$ or
$\frac{\partial_{\nu}{\bf{x}}}{{\bf{y}}}$. Also note that the
coefficients of these expressions decay faster than exponentially.
Now apply a localizer $P_{k_{0}}$, $|k_{0}|<\frac{\epsilon
i}{1000C}$ in front. We redefine $C$ so large that
$\frac{\epsilon}{C}<<\mu_{1}(\epsilon)$, the latter coming from
the preceding computation. We claim that if one of the input
frequencies has absolute value greater than $i\delta$ for suitable
$\delta>\frac{\epsilon}{500 C}$, one obtains an exponential gain
in $i$ for the norms in the statement of the lemma. This is done
inductively: write the expression under consideration as
\begin{equation}\nonumber
P_{k}[P_{k_{1}}(\psi_{1}\nabla^{-1}\psi_{2}\ldots\nabla^{-1}\psi_{k-1})\nabla^{-1}P_{k_{2}}\psi_{k}]
\end{equation}
If $k_{1}>\frac{\epsilon i}{1000C}+5$, so is $k_{2}$. If
$P_{k_{1}}(...)$, $P_{k_{2}}\psi_{2}$ are both of first type, one
estimates using 3.4(a)
\begin{equation}\nonumber\begin{split}
&||P_{k}[P_{k_{1}}(\psi_{1}\nabla^{-1}\psi_{2}\ldots\nabla^{-1}\psi_{k-1})\nabla^{-1}P_{k_{2}}\psi_{k}]||_{\dot{X}_{k}^{0,\frac{1}{2},1}}\\
&\lesssim
||P_{k_{1}}(\psi_{1}\nabla^{-1}\psi_{2}\ldots\nabla^{-1}\psi_{k-1})||_{A[k_{1}]}||P_{k_{2}}\psi_{2}||_{A[k_{2}]}\lesssim
\frac{\tilde{c}_{k_{1}}}{\epsilon}\frac{\tilde{c}_{k_{2}}}{\epsilon}\\
\end{split}\end{equation}
If one of them is of 2nd type, one places this into
$L_{t}^{2}L_{x}^{2+}$ to control the portion of the output at
modulation $<2^{k+100}$ (the other portion being controlled by
theorem~\ref{Moser}.). One obtains the same bound, and our decay
assumptions on the frequency envelope yield the claim, provided
one shows that the 2nd more complicated norm in the statement of
the lemma is controlled by
$||.||_{\dot{X}_{k}^{0,\frac{1}{2},1}}$. This follows from
\eqref{tame} as well as the preceding computation
\eqref{calculation}. In case $k_{2}<-\frac{\epsilon i}{1000 C}-5$,
first assume $P_{k_{2}}\psi_{2}$ to be of 2nd type. If $k_{2}\geq
j-10$, we estimate
\begin{equation}\nonumber
||P_{k}Q_{j}[P_{k_{1}}\psi_{1}\nabla^{-1}P_{k_{2}}\psi_{2}]||_{\dot{X}_{k}^{0,\frac{1}{2},\infty}}
\lesssim
2^{\frac{j}{2}}||P_{k_{1}}\psi_{1}||_{L_{t}^{\infty}L_{x}^{2}}||\nabla^{-1}P_{k_{2}}\psi_{2}||_{L_{t}^{2}L_{x}^{\infty}}
\lesssim 2^{\frac{j-k_{2}}{2}}\tilde{c}_{k_{2}}
\end{equation}
If $k_{2}<j-10$, we estimate
\begin{equation}\nonumber\begin{split}
&||P_{k}Q_{j}[P_{k_{1}}\psi_{1}\nabla^{-1}P_{k_{2}}\psi_{2}]||_{\dot{X}_{k}^{0,\frac{1}{2},\infty}}
\lesssim
||P_{k}Q_{j}[P_{k_{1}}\psi_{1}\nabla^{-1}P_{k_{2}}Q_{\geq j-10}\psi_{2}]||_{\dot{X}_{k}^{0,\frac{1}{2},\infty}}\\
&\hspace{4cm}+||P_{k}Q_{j}[P_{k_{1}}Q_{\geq
j-10}\psi_{1}\nabla^{-1}P_{k_{2}}Q_{<j-10}\psi_{2}]||_{\dot{X}_{k}^{0,\frac{1}{2},\infty}}\\
\end{split}\end{equation}
Then we have
\begin{equation}\nonumber\begin{split}
&||P_{k}Q_{j}[P_{k_{1}}\psi_{1}\nabla^{-1}P_{k_{2}}Q_{\geq
j-10}\psi_{2}]||_{\dot{X}_{k}^{0,\frac{1}{2},\infty}}\lesssim
||P_{k_{1}}\psi_{1}||_{L_{t}^{\infty}L_{x}^{2}}||\nabla^{-1}P_{k_{2}}Q_{\geq
j-10}\psi_{2}||_{L_{t}^{2}L_{x}^{\infty}}\\
&\lesssim 2^{\frac{1}{2+}(k_{2}-j)}\tilde{c}_{k_{2}}\\
\end{split}\end{equation}
Also, we have
\begin{equation}\nonumber\begin{split}
&||P_{k}Q_{j}[P_{k_{1}}Q_{\geq
j-10}\psi_{1}\nabla^{-1}P_{k_{2}}Q_{<j-10}\psi_{2}]||_{\dot{X}_{k}^{0,\frac{1}{2},\infty}}\\&\lesssim
2^{\frac{j}{2}}||P_{k_{1}}Q_{\geq
j-10}\psi_{1}||_{L_{t}^{2}L_{x}^{2}}||\nabla^{-1}P_{k_{2}}Q_{<j-10}\psi_{2}||_{L_{t}^{\infty}L_{x}^{\infty}}
\lesssim \tilde{c}_{k_{2}}\\
\end{split}\end{equation}
Control over the more complicated norm in the lemma is a bit more
difficult, and follows from a computation performed in the
appendix, in addition to the calculation \eqref{calculation}
performed above: one again obtains the bound $\lesssim
\tilde{c}_{k_{2}}$, and the claim follows by summing over
$k_{2}<-\frac{\epsilon i}{1000 C}-5$. If, on the other hand,
$P_{k_{2}}\psi_{2}$ is of first type, then one estimates
\begin{equation}\nonumber\begin{split}
&||P_{k}[P_{k_{1}}\psi_{1}\nabla^{-1}P_{k_{2}}\psi_{2}]||_{\dot{X}_{k}^{0,\frac{1}{2},\infty}}\lesssim
||P_{k_{1}}\psi_{1}||_{A[k_{1}]+\dot{X}_{k_{1}}^{0,\frac{1}{2},1}}||P_{k_{2}}\psi_{2}||_{A[k_{2}]}
\lesssim \tilde{c}_{k_{2}}\\
\end{split}\end{equation}
with a similar estimate (using \eqref{calculation}) controlling
the 2nd more complicated norm in the lemma. Summing over
$k_{2}<-\frac{\epsilon i}{1000 C}-5$ results in the desired
exponential gain in $i$. Now assume $k_{1}<-\frac{\epsilon i}{1000
C}-5$. In that case, the claim follows from the estimate
\begin{equation}\nonumber
||P_{k}[P_{k_{1}}\psi_{1}\nabla^{-1}P_{k_{2}}\psi_{2}]||_{\dot{X}_{k}^{0,\frac{1}{2},1}}\lesssim
2^{\delta_{1}(k_{1}-k)}|k_{1}-k_{2}|||P_{k_{1}}\psi_{1}||_{{\mathcal{S}}[k_{1}]}||P_{k_{2}}\psi_{2}||_{A[k_{2}]+\dot{X}_{k_{2}}^{0,\frac{1}{2},1}},
\end{equation}
by a calculation similar to the immediately preceding. Summing
over $k_{1}<-\frac{\epsilon i}{1000 C}-5$ results in the desired
exponential gain. Having thus shown that we may assume
$|k_{1}|<\frac{\epsilon i}{1000 C}+5 $, we apply the same
procedure to the expression\\
$P_{k_{1}}[\psi_{1}\nabla^{-1}\psi_{2}\ldots\nabla^{-1}\psi_{k-1}]$.
Of course the frequency bounds will grow the further one continues
this process, but this is counteracted by the rapidly decreasing
coefficients coming from the Taylor expansion\footnote{More
precisely, one needs to go to expressions of length up to
$\frac{\epsilon i}{10000 C}$.}. Arguing inductively, we see that
it suffices to show that provided two functions
$\psi_{1}\in{\mathcal{S}}({\mathbf{R}}^{2+1})$,
$\psi_{2}\in{\mathcal{S}}({\mathbf{R}}^{2+1})$ satisfy the
assertions in the lemma, then we get the same assertions for the
expression
\begin{equation}\nonumber
P_{k}\chi_{i}(t)[P_{k_{1}}\psi_{1}\nabla^{-1}P_{k_{2}}\psi_{2}]
\end{equation}
where $|k|,|k_{1}|,|k_{2}|<\delta i$, $\delta<<\epsilon$. First,
we estimate
\begin{equation}\nonumber\begin{split}
&||P_{k}[P_{k_{1}}Q_{>k_{1}-(1-\epsilon)i}[\chi_{i}(t)\psi_{1}]\nabla^{-1}P_{k_{2}}\psi_{2}]||_{\dot{X}_{k}^{0,\frac{1}{2},1}}
\\&\lesssim
(|k|+|k_{1}|+|k_{2}|)||P_{k_{1}}Q_{>k_{1}-(1-\epsilon)i}[\chi_{i}(t)\psi_{1}]||_{\dot{X}_{k_{1}}^{0,\frac{1}{2},1}}
||\nabla^{-1}P_{k_{2}}\psi_{2}||_{\dot{X}_{k_{2}}^{0,\frac{1}{2},1}+A[k_{2}]}\\
&\lesssim |i|2^{-\mu(\epsilon)i}.\\
\end{split}\end{equation}
We have used the bound on $P_{k_{2}}\psi_{2}$ from
theorem~\ref{Moser}, as well as the easily verified fact that
\begin{equation}\nonumber
||P_{k_{1}}Q_{>k_{1}-(1-\epsilon)i}[\chi_{i}(t)\psi_{1}]||_{\dot{X}_{k_{1}}^{0,\frac{1}{2},1}}\lesssim
2^{-\mu(\epsilon)i}
\end{equation}
One proceeds similarly when $P_{k_{2}}\psi_{2}$ is replaced by
$P_{k_{2}}Q_{>k_{2}-(1-\epsilon)i}\psi_{2}$. Thus it suffices to
consider the expression
\begin{equation}\nonumber
P_{k}\chi_{i}(t)[P_{k_{1}}Q_{<k_{1}-(1-\epsilon)i}\phi\nabla^{-1}P_{k_{2}}Q_{<k_{2}-(1-\epsilon)i}\psi]
\end{equation}
We may always reduce the modulation of the output to size
$>2^{k-\frac{5\epsilon i}{C}}$, as follows easily from 3.4(a). Now
we carry out the following decomposition
\begin{equation}\label{decompose21}\begin{split}
&P_{k}\chi_{i}(t)[P_{k_{1}}Q^{\pm}_{<k_{1}-(1-\epsilon)i}\phi\nabla^{-1}P_{k_{2}}Q_{<k_{2}-(1-\epsilon)i}\psi]
\\&=P_{k}\chi_{i}(t)[\sum_{\kappa\in
K_{\frac{\epsilon-1}{2}i}}\chi^{c}_{i,\mp\kappa}P_{k_{1},\kappa}Q^{\pm}_{<k_{1}-(1-\epsilon)i}\phi\nabla^{-1}P_{k_{2}}Q_{<k_{2}-(1-\epsilon)i}\psi]
\\&+P_{k}\chi_{i}(t)[\sum_{\kappa\in
K_{\frac{\epsilon-1}{2}i}}\chi_{i,\mp\kappa}P_{k_{1},\kappa}Q^{\pm}_{<k_{1}-(1-\epsilon)i}\phi\nabla^{-1}P_{k_{2}}Q_{<k_{2}-(1-\epsilon)i}\psi]
\\&+P_{k}\chi_{i}(t)[\phi_{i}(t,x)P_{k_{1}}Q^{\pm}_{<k_{1}-(1-\epsilon)i}\phi\nabla^{-1}P_{k_{2}}Q_{<k_{2}-(1-\epsilon)i}\psi]\\
\end{split}\end{equation}
In the immediately preceding we let $\chi_{i,\pm\kappa}(t,x)$
localize to two slabs aligned or anti-aligned with $\kappa$ of
length $\sim 2^{i}$ and thickness $\sim 2^{\epsilon i}$ (the
complement of $\chi^{c}_{i,\pm\kappa}$ within the
$\frac{2^{\epsilon i}}{C}$-neighborhood of the physical light
cone), and we let $\phi_{i}(t,x)$ smoothly localize to the
intersection of the complement of the $\frac{2^{\epsilon i}}{C}$
neighborhood of the light cone and a box of dimensions $\sim
2^{i}\times2^{i}\times 2^{i}$(finite propagation speed and
properties of the Fourier multipliers). We estimate
\begin{equation}\nonumber\begin{split}
&||P_{k}\chi_{i}(t)[\phi_{i}(t,x)P_{k_{1}}Q_{<k_{1}-(1-\epsilon)i}\phi\nabla^{-1}P_{k_{2}}Q_{<k_{2}-(1-\epsilon)i}\psi]||_{\dot{X}_{k}^{0,\frac{1}{2},1}}
\\&\lesssim
2^{\frac{\max\{k,k_{1},k_{2}\}}{2}}||\phi_{i}(t,x)\chi_{i}(t)\nabla^{-1}P_{k_{2}}Q_{<k_{2}}\psi||_{L_{t}^{2}L_{x}^{\infty}}||P_{k_{1}}Q_{<k_{1}-(1-\epsilon)i}\phi||_{L_{t}^{\infty}L_{x}^{2}}
\\&\hspace{1cm}+||P_{k}\chi_{i}(t)[\phi_{i}(t,x)P_{k_{1}}Q_{<k_{1}-(1-\epsilon)i}\phi\nabla^{-1}P_{k_{2}}Q_{[k_{2}+(\epsilon-1)i,k_{2}]}\psi]||_{\dot{X}_{k}^{0,\frac{1}{2},1}}\\
\end{split}\end{equation}
Now the first summand on the right hand side is immediately
controlled\footnote{One may have to shrink the size of $|k|$,
$|k_{1,2}|$ if necessary.} from Proposition~\ref{Chr1}, while we
estimate the 2nd as follows(modulo error terms of order of
magnitude $2^{-Ni}$):
\begin{equation}\nonumber\begin{split}
&||P_{k}\chi_{i}(t)[\phi_{i}(t,x)P_{k_{1}}Q_{<k_{1}-(1-\epsilon)i}\phi\nabla^{-1}P_{k_{2}}Q_{[
k_{2}+(\epsilon-1)i,k_{2}]}\psi]||_{\dot{X}_{k}^{0,\frac{1}{2},1}}\\&\lesssim
2^{\frac{\max\{k,k_{1},k_{2}\}}{2}}||P_{k}\chi_{i}(t)[\phi_{i}(t,x)P_{k_{1}}Q_{<k_{1}-(1-\epsilon)i}\phi\nabla^{-1}P_{k_{2}}Q_{[
k_{2}+(\epsilon-1)i,k_{2}]}\psi]||_{L_{t}^{2}L_{x}^{2}}\\&\lesssim
2^{\frac{\max\{k,k_{1},k_{2}\}}{2}}||\phi_{i}(t,x)||_{L_{t}^{\infty}L_{x}^{M}}||P_{k_{1}}Q_{<k_{1}-(1-\epsilon)i}\phi||_{L_{t}^{M}L_{x}^{2+}}
\\&\hspace{7cm}||\nabla^{-1}P_{k_{2}}Q_{[
k_{2}+(\epsilon-1)i,k_{2}]}\chi_{i}(t)\psi||_{L_{t}^{2+}L_{x}^{\infty}}\\&\lesssim 2^{\frac{\max\{k,k_{1},k_{2}\}}{2}}C(M)2^{(\frac{1}{M}-\mu)i}\\
\end{split}\end{equation}
This is acceptable if we choose $M$ large enough. Now we proceed
to the other two terms in \eqref{decompose21}. Consider the first.
We have for $\min\{k_{1},k_{2},k\}+O(1)>j>k-\frac{5\epsilon i}{C}$
\begin{equation}\nonumber\begin{split}
&||P_{k}Q_{j}\chi_{i}(t)[\sum_{\kappa\in
K_{\frac{\epsilon-1}{2}i}}\chi^{c}_{i,\mp\kappa}P_{k_{1},\kappa}Q^{\pm}_{<k_{1}-(1-\epsilon)i}\phi\nabla^{-1}P_{k_{2}}Q^{\pm}_{<k_{2}-(1-\epsilon)i}\psi]||_{\dot{X}_{k}^{0,\frac{1}{2},1}}\\
&\lesssim
2^{\frac{\min\{k,k_{1},k_{2}\}-j}{4}}(\sum_{\tilde{\kappa}\in
K_{\frac{j-\min\{k,k_{1},k_{2}\}}{2}}}||P_{k_{1},\tilde{\kappa}}\sum_{\kappa\in
K_{\frac{\epsilon-1}{2}i}}\chi^{c}_{i,\mp\kappa}P_{k_{1},\kappa}Q^{\pm}_{<k_{1}-(1-\epsilon)i}\phi||_{PW[\pm\tilde{\kappa}]}^{2})^{\frac{1}{2}}
\\&\hspace{4cm}(\sum_{\tilde{\kappa}\in
K_{\frac{j-\min\{k,k_{1},k_{2}\}}{2}}}||\nabla^{-1}P_{k_{2},\tilde{\kappa}}Q^{\pm}_{<k_{2}-(1-\epsilon)i}\psi||_{NFA^{*}[\pm\tilde{\kappa}]}^{2})^{\frac{1}{2}},\\
\end{split}\end{equation}
which is bounded by $2^{-\mu i}$. One can easily sum over $j\in
[k-ai, k+bi]$, obtaining the desired estimate. If
$j>\min\{k_{1},k_{2},k\}+100$, say, then necessarily
$k_{1}=k_{2}+O(1)=j+O(1)$, since the inputs
$\chi^{c}_{i,\mp\kappa}P_{k_{1},\kappa}Q^{\pm}_{<k_{1}-(1-\epsilon)i}\psi_{2}$
etc have modulation $<C2^{-\epsilon
i}<<\max\{2^{k},2^{k_{1,2}}\}$. Then one can concurrently
microlocalize the inputs to caps $\kappa_{1,2}\in K_{-100}$ of
distance $\sim 1$, and argue just as before. The 2nd term in
\eqref{decompose21} has to be decomposed further as follows:
\begin{equation}\nonumber\begin{split}
&P_{k}\chi_{i}(t)[\sum_{\kappa\in
K_{\frac{\epsilon-1}{2}i}}\chi_{i,\mp\kappa}P_{k_{1},\kappa}Q^{\pm}_{<k_{1}-(1-\epsilon)i}\phi\nabla^{-1}P_{k_{2}}Q_{<k_{2}-(1-\epsilon)i}\psi]\\
&=P_{k}\chi_{i}(t)[\sum_{\kappa\in
K_{\frac{\epsilon-1}{2}i}}\chi_{i,\mp\kappa}P_{k_{1},\kappa}Q^{\pm}_{<k_{1}-(1-\epsilon)i}\phi\nabla^{-1}\sum_{\kappa\in K_{\frac{\epsilon-1}{2}i}}\chi^{c}_{i,\mp\kappa}P_{k_{2},\kappa}Q^{\pm}_{<k_{2}-(1-\epsilon)i}\psi]]\\
&+P_{k}\chi_{i}(t)[\sum_{\kappa\in
K_{\frac{\epsilon-1}{2}i}}\chi_{i,\mp\kappa}P_{k_{1},\kappa}Q^{\pm}_{<k_{1}-(1-\epsilon)i}\phi\nabla^{-1}\sum_{\kappa\in K_{\frac{\epsilon-1}{2}i}}\chi_{i,\mp\kappa}P_{k_{2},\kappa}Q^{\pm}_{<k_{2}-(1-\epsilon)i}\psi]]\\
\end{split}\end{equation}
For the first of the immediately preceding terms, we can
estimate\footnote{Exploit the fact that $\chi_{i,\kappa}$ only
smears out the Fourier support by about $2^{-\epsilon i}$, up to
negligible error terms.} provided $j<\min\{k,k_{1},k_{2}\}+O(1)$
\begin{equation}\nonumber\begin{split}
&2^{\frac{j}{2}}||P_{k}Q_{j}\chi_{i}(t)[\sum_{\kappa\in
K_{\frac{\epsilon-1}{2}i}}\chi_{i,\mp\kappa}P_{k_{1},\kappa}Q^{\pm}_{<k_{1}-(1-\epsilon)i}\phi\\&\hspace{6cm}\nabla^{-1}\sum_{\kappa\in
K_{\frac{\epsilon-1}{2}i}}\chi^{c}_{i,\mp\kappa}P_{k_{2},\kappa}Q^{\pm}_{<k_{2}-(1-\epsilon)i}\psi]]||_{L_{t}^{2}L_{x}^{2}}\\
&\lesssim (\sum_{\kappa\in
K_{\frac{\epsilon-1}{2}i}}||P_{k}Q_{j}[P_{k_{1},\kappa}Q^{\pm}_{<k_{1}-(1-\epsilon)i}\phi\nabla^{-1}\sum_{\kappa\in
K_{\frac{\epsilon-1}{2}i}}\chi^{c}_{i,\mp\kappa}P_{k_{2},\kappa}Q^{\pm}_{<k_{2}-(1-\epsilon)i}\psi]||_{L_{t}^{2}L_{x}^{2}}^{2})^{\frac{1}{2}}\\
\end{split}\end{equation}
Using elementary geometry and the definition of $PW[\kappa]$,
$NFA^{*}[\kappa]$ etc., this in turn is estimated by
\begin{equation}\nonumber\begin{split}
&\lesssim 2^{-\frac{j-\min\{k,k_{1},k_{2}\}}{4}}(\sum_{\kappa\in
K_{\frac{\epsilon-1}{2}i}}||P_{k_{1},\kappa}Q^{\pm}_{<k_{1}-(1-\epsilon)i}\phi||_{NFA^{*}[\pm\kappa]}^{2})^{\frac{1}{2}}
\\&\hspace{1cm}(\sum_{\tilde{\kappa}\in
K_{\frac{j-\min\{k,k_{1},k_{2}\}}{2}}}||P_{k,\tilde{\kappa}}\sum_{\kappa\in
K_{\frac{\epsilon-1}{2}i}}\chi^{c}_{i,\mp\kappa}P_{k_{2},\kappa}Q^{\pm}_{<k_{2}-(1-\epsilon)i}\psi||_{PW[\pm\tilde{\kappa}]}^{2})^{\frac{1}{2}}\\
\end{split}\end{equation}
This can again be estimated by $\lesssim 2^{-\mu i}$, as desired.
The case $j>\min\{k,k_{1},k_{2}\}+O(1)$ is treated as before.
Finally, we have
\begin{equation}\nonumber\begin{split}
&||P_{k}Q_{j}\chi_{i}(t)[\sum_{\kappa\in
K_{\frac{\epsilon-1}{2}i}}\chi_{i,\mp\kappa}P_{k_{1},\kappa}Q^{\pm}_{<k_{1}-(1-\epsilon)i}\phi\\&\hspace{5cm}\nabla^{-1}\sum_{\kappa\in K_{\frac{\epsilon-1}{2}i}}\chi_{i,\mp\kappa}P_{k_{2},\kappa}Q^{\pm}_{<k_{2}-(1-\epsilon)i}\psi]]||_{L_{t}^{2}L_{x}^{2}}\\
&\lesssim 2^{-Ni}\\
\end{split}\end{equation}
on account of the support properties of the (Fourier)multipliers.
We deduce from this that the expression
$P_{k}\chi_{i}(t)[P_{k_{1}}\phi\nabla^{-1}P_{k_{2}}\psi]$ also
satisfies the 2nd property of the lemma: simply apply \eqref{tame}
in addition to \eqref{calculation}.
\end{proof}
We can now conclude case (I.b.2). The preceding proof in addition
to 3.4(g) imply that
\begin{equation}\nonumber\begin{split}
&||P_{k_{0}}Q_{<k_{0}}
[R_{\alpha}P_{k_{1}}Q_{<k_{0}}\delta\psi_{1}P_{k}\triangle^{-1}\sum_{j=1,2}\partial_{j}\Box I[\nabla^{-1}P_{k_{2}}\chi_{i}(t)\psi_{2}R_{j}P_{k_{3}}\psi_{3}]]||_{N[k_{0}]}\\
&\lesssim
||P_{k_{1}}\delta\psi_{1}||_{S[k_{1}]}||P_{k}I[\nabla^{-1}P_{k_{2}}\chi_{i}(t)\psi_{2}R_{j}P_{k_{3}}\psi_{3}]||_{\dot{X}_{k}^{0,\frac{1}{2},1}}\lesssim
2^{-\mu i}\tilde{c}_{k_{1}}\\
\end{split}\end{equation}
One can sum over the admissible frequency ranges here (picking up
factors $O(i)$), which yields the desired result.\\
(I.b.3): {\it{The third term in the expansion:}}
 \begin{equation}\nonumber\begin{split}
&P_{k_{0}}Q_{<k_{0}}
[R_{\alpha}P_{k_{1}}Q_{<k_{0}}\Box\delta\psi_{1}P_{k}\triangle^{-1}\sum_{j=1,2}\partial_{j}I[\nabla^{-1}P_{k_{2}}\chi_{i}(t)\psi_{2}R_{j}P_{k_{3}}\psi_{3}]]\\
\end{split}\end{equation}
This is much simpler to estimate, on account of the strong
Strichartz type estimates available for $\psi_{\nu}$: we have
\begin{equation}\nonumber\begin{split}
&||P_{k_{0}}Q_{<k_{0}}
[R_{\alpha}P_{k_{1}}Q_{<k_{0}}\Box\delta\psi_{1}P_{k}\triangle^{-1}\sum_{j=1,2}\partial_{j}I[\nabla^{-1}P_{k_{2}}\chi_{i}(t)\psi_{2}R_{j}P_{k_{3}}\psi_{3}]]||_{L_{t}^{1}\dot{H}^{-1}}\\
&\lesssim||P_{k_{1}}Q_{<k_{0}}\Box\delta\psi_{1}||_{L_{t}^{2}L_{x}^{M}}||\nabla^{-1}P_{k_{2}}\chi_{i}(t)\psi_{2}||_{L_{t}^{4}L_{x}^{4+}}||P_{k_{3}}\psi_{3}||_{L_{t}^{4}L_{x}^{4+}}
\lesssim 2^{-\mu i}\tilde{c}_{k_{1}}\\
\end{split}\end{equation}
(I.b.4) {\it{The fourth term of the expansion:}}
 \begin{equation}\nonumber\begin{split}
&P_{k_{0}}Q_{<k_{0}}
[P_{k_{3}}Q_{<k_{3}}\nabla^{-1}\psi_{3}IP_{k_{2}+O(1)}\Box[\nabla^{-1}P_{k_{2}}Q_{<k_{2}}\chi_{i}(t)\psi_{2}R_{\alpha}P_{k_{1}}Q_{<k_{0}}\delta\psi_{1}]]\\
\end{split}\end{equation}
This is straightforward by means of 3.4(a), 3.4(g): one
estimates\footnote{Recall the definition of ${\mathcal{S}}[k]$ via
$||.||_{A[k]}$, $||.||_{B[k]}$. Also recall that we may replace
$P_{k_{2}}\psi_{2}$ by $P_{k_{2}}Q_{<k_{2}}\psi_{2}$ etc.}
\begin{equation}\nonumber\begin{split}
&||P_{k_{0}}Q_{<k_{0}}
[P_{k_{3}}Q_{<k_{3}}\nabla^{-1}\psi_{3}P_{k_{2}+O(1)}I\Box[\nabla^{-1}P_{k_{2}}Q_{<k_{2}}\chi_{i}(t)\psi_{2}R_{\alpha}P_{k_{1}}Q_{<k_{0}}\delta\psi_{1}]]||_{N[k_{0}]}\\
&\lesssim
|k_{0}-k_{3}| ||P_{k_{3}}\nabla^{-1}\psi_{3}||_{A[k_{3}]+\dot{X}_{k_{3}}^{0,\frac{1}{2},1}}\\&\hspace{5cm}||[\nabla^{-1}P_{k_{2}}Q_{<k_{2}}\chi_{i}(t)\psi_{2}R_{\alpha}P_{k_{1}}Q_{<k_{0}}\delta\psi_{1}]||_{\dot{X}_{k_{2}}^{1,\frac{1}{2},1}}\\
&\lesssim |k_{0}-k_{3}|
||P_{k_{3}}\nabla^{-1}\psi_{3}||_{\dot{X}_{k_{3}}^{0,\frac{1}{2},1}+A[k_{3}]}
\\&\hspace{4cm}2^{k_{1}}||P_{k_{1}}\delta\psi_{1}||_{S[k_{1}]}||\nabla^{-1}P_{k_{2}}Q_{<k_{2}}\chi_{i}(t)\psi_{2}||_{\dot{X}_{k_{2}}^{0,\frac{1}{2},1}+A[k_{2}]}\\&
\lesssim 2^{-\mu i}\tilde{c}_{k_{1}}\\
\end{split}\end{equation}
(I.b.5) {\it{The fifth term of the expansion.}}
\begin{equation}\nonumber\begin{split}
&P_{k_{0}}Q_{<k_{0}}
[P_{k_{3}}Q_{<k_{3}}\nabla^{-1}\psi_{3}P_{k_{3}}I[\nabla^{-1}\Box P_{k_{2}}Q_{<k_{2}+O(1)}\chi_{i}(t)\psi_{2}R_{\alpha}P_{k_{1}}Q_{<k_{0}}\delta\psi_{1}]]\\
\end{split}\end{equation}
This can be estimated by means of lemma~\ref{key}: observe that we
may throw an operator $Q_{[k_{2}+(\epsilon-1)i,k_{2}]}$ in front
of $P_{k_{2}}Q_{<k_{2}+O(1)}\chi_{i}(t)\psi_{2}$, since in the
opposite case we have
\begin{equation}\nonumber\begin{split}
&||P_{k_{0}}Q_{<k_{0}}
[P_{k_{3}}Q_{<k_{3}}\nabla^{-1}\psi_{3}P_{k_{3}}I[\nabla^{-1}\Box P_{k_{2}}Q_{<k_{2}+(\epsilon-1)i}\chi_{i}(t)\psi_{2}R_{\alpha}P_{k_{1}}Q_{<k_{0}}\delta\psi_{1}]]||_{L_{t}^{1}\dot{H}^{-1}}\\
&\lesssim
2^{-k_{0}}||P_{k_{3}}Q_{<k_{3}}\nabla^{-1}\psi_{3}||_{L_{t}^{2+}L_{x}^{\infty}}||\nabla^{-1}\Box
P_{k_{2}}Q_{<k_{2}+(\epsilon-1)i}\chi_{i}(t)\psi_{2}||_{L_{t}^{2}L_{x}^{2}}||P_{k_{1}}\delta\psi_{1}||_{L_{t}^{M}L_{x}^{\infty}}\\
&\lesssim
2^{-k_{0}}2^{\frac{\epsilon-1}{2}i}2^{(1-\frac{1}{M})k_{1}}\tilde{c}_{k_{1}}\lesssim
2^{-\mu i}\tilde{c}_{k_{1}}\\
\end{split}\end{equation}
for very large $M$, on account of the assumptions on $k_{1},
k_{0}$. Using lemma~\ref{key} we have
\begin{equation}\nonumber\begin{split}
&||P_{k_{0}}Q_{<k_{0}}
[R_{\alpha}P_{k_{3}}Q_{<k_{3}}\nabla^{-1}\psi_{3}P_{k_{3}}I[\nabla^{-1}\Box P_{k_{2}}Q_{[k_{2}+(\epsilon-1)i,k_{2}]}\chi_{i}(t)\psi_{2}R_{j}P_{k_{1}}\delta\psi_{1}]]||_{L_{t}^{1}\dot{H}^{-1}}\\
&\lesssim
2^{-k_{0}}||P_{k_{3}}\nabla^{-1}\psi_{3}||_{L_{t}^{2+}L_{x}^{\infty}}||P_{k_{2}}Q_{[k_{2}+(\epsilon-1)i,k_{2}]}\chi_{i}(t)\psi_{2}||_{\dot{X}_{k_{2}}^{0,\frac{1}{2},1}}||P_{k_{1}}\delta\psi_{1}||_{L_{t}^{M}L_{x}^{\infty}}\\
&\lesssim 2^{-\mu i}2^{-\frac{k_{1}}{M}}\tilde{c}_{k_{1}}\\
\end{split}\end{equation}
Choosing $M$ large enough results in the desired gain in $i$.\\
(I.b.6): {\it{The sixth term of the expansion.}} This is similar
to the third and hence omitted. The terms in the expansion
\eqref{null-form2} are simple variations of the same kind of
reasoning and hence omitted. This concludes case (I.b).\\
(I.c): {\it{None of (I.a), (I.b) hold but $i\lesssim |k_{1}|$.}}
This implies $|k_{2,3}|, |k|<<i$, and we may treat these as
$O(1)$. We also have $i\lesssim k_{1}$, and $k_{1}=k_{0}+O(1)$.
Now we proceed in close analogy to the immediately preceding case.
We may pull the multiplier $\chi_{i}(t)$ past the operator
$P_{k}QI\nabla^{-1}$ right in front of $P_{k_{2}}\psi_{2}$. Also,
we may reduce the output and first input $P_{k_{1}}\delta\psi_{1}$
to modulation $<O(1)$: for example, one estimates
\begin{equation}\nonumber\begin{split}
&||P_{k_{0}}Q_{\geq
0}\nabla_{x,t}[P_{k_{1}}Q_{<k_{1}}R_{\alpha}\delta\psi_{1}\nabla^{-1}P_{k}I[P_{k_{2}}\chi_{i}(t)\psi_{2}P_{k_{3}}\psi_{3}]]||_{\dot{X}_{k_{0}}^{-1,-\frac{1}{2},1}}\\
&\lesssim
||P_{k_{1}}\delta\psi_{1}||_{L_{t}^{\infty}L_{x}^{2}}||\nabla^{-1}P_{k}I[P_{k_{2}}\chi_{i}(t)\psi_{2}P_{k_{3}}\psi_{3}]||_{L_{t}^{2}L_{x}^{\infty}}\\
&\lesssim 2^{-\mu i}\tilde{c}_{k_{1}}\\
\end{split}\end{equation}
on account of the proof of lemma~\ref{key}, and the fact that by
theorem~\ref{spherical Strichartz} we have
\begin{equation}\nonumber\begin{split}
&||\nabla^{-1}P_{k}I[P_{k_{2}}\chi_{i}(t)\psi_{2}P_{k_{3}}\psi_{3}]||_{L_{t}^{2}L_{x}^{\infty}}\lesssim
||\nabla^{-1}P_{k}I[P_{k_{2}}\chi_{i}(t)\psi_{2}P_{k_{3}}\psi_{3}]||_{\dot{X}_{k}^{\epsilon,\frac{1}{2+},1}}\\
\end{split}\end{equation}
The estimate for when $P_{k_{1}}R_{\alpha}\delta\psi_{1}$ is
replaced by $P_{k_{1}}Q_{\geq k_{1}}R_{\alpha}\psi_{1}$ is of
course similar, since one may place this into $L_{t}^{2}L_{x}^{2}$
and gains $2^{-\frac{k_{1}}{2}}$. Then we resort to the null-form
identities \eqref{null-form1}, \eqref{null-form2}. We wind up
having to estimate the same types of expressions as in (I.b.1-6),
with our changed assumptions on
the frequencies. We shall do this in a brisk manner here:\\
(I.c.1): we need to exert some care not to lose in $k_{0}$: use
the proof of lemma~\ref{key} as above to obtain
\begin{equation}\nonumber\begin{split}
&||P_{k_{0}}Q_{<O(1)}\Box
[R_{\alpha}P_{k_{1}}Q_{<O(1)}\delta\psi_{1}P_{k}\triangle^{-1}\sum_{j=1,2}\partial_{j}I[\nabla^{-1}P_{k_{2}}\chi_{i}(t)\psi_{2}R_{j}P_{k_{3}}\psi_{3}]]||_{\dot{X}_{k_{0}}^{-1,-\frac{1}{2},1}}\\
&\lesssim
||P_{k_{1}}Q_{<O(1)}\delta\psi_{1}||_{L_{t}^{\infty}L_{x}^{2}}||P_{k}\triangle^{-1}\sum_{j=1,2}\partial_{j}I[\nabla^{-1}P_{k_{2}}\chi_{i}(t)\psi_{2}R_{j}P_{k_{3}}\psi_{3}]||_{L_{t}^{2}L_{x}^{\infty}}\\
&\lesssim ||P_{k_{1}}Q_{<O(1)}\delta\psi_{1}||_{L_{t}^{\infty}L_{x}^{2}}||P_{k}\triangle^{-1}\sum_{j=1,2}\partial_{j}I[\nabla^{-1}P_{k_{2}}\chi_{i}(t)\psi_{2}R_{j}P_{k_{3}}\psi_{3}]||_{\dot{X}_{k}^{\epsilon,\frac{1}{2+},1}}\\
&\lesssim 2^{-\mu i}\tilde{c}_{k_{0}}\\
\end{split}\end{equation}
(I.c.2): This is estimated exactly like (I.b.2).\\
(I.c.3): The argument is here is slightly more complicated than in
(I.b.3); we argue as in (I.c.2):
 \begin{equation}\nonumber\begin{split}
&||P_{k_{0}}Q_{<O(1)}
[R_{\alpha}P_{k_{1}}Q_{<O(1)}\Box\delta\psi_{1}P_{k}\triangle^{-1}\sum_{j=1,2}\partial_{j}I[\nabla^{-1}P_{k_{2}}\chi_{i}(t)\psi_{2}R_{j}P_{k_{3}}\psi_{3}]]||_{L_{t}^{1}\dot{H}^{-1}}\\
&\lesssim
2^{-k_{0}}||P_{k_{1}}Q_{<O(1)}\Box\delta\psi_{1}||_{L_{t}^{2}L_{x}^{2}}||P_{k}\triangle^{-1}\sum_{j=1,2}\partial_{j}I[\nabla^{-1}P_{k_{2}}\chi_{i}(t)\psi_{2}R_{j}P_{k_{3}}\psi_{3}]||_{L_{t}^{2}L_{x}^{\infty}}\\
&\lesssim
\tilde{c}_{k_{0}}||P_{k}\triangle^{-1}\sum_{j=1,2}\partial_{j}I[\nabla^{-1}P_{k_{2}}\chi_{i}(t)\psi_{2}R_{j}P_{k_{3}}\psi_{3}]||_{\dot{X}_{k}^{\epsilon,\frac{1}{2+},1}}
\lesssim 2^{-\mu i}\tilde{c}_{k_{0}}\\
\end{split}\end{equation}
(I.c.4): The argument here is more complicated than in (I.b.4) on
account of the fact that we don't get the desired exponential gain
from simple application of 3.4(a). Instead, we shall have to
resort to lemma~\ref{key}. We split the expression into several
manageable pieces: first, using 3.4(g), we can
estimate\footnote{It is easy to see that we may restrict
$P_{k_{1}}R_{\alpha}\delta\psi_{1}$ to modulation $<O(1)$, so that
we don't lose from $R_{\alpha}$. We shall omit this to streamline
the formulae.}
\begin{equation}\nonumber\begin{split}
&||P_{k_{0}}Q_{<O(1)}
[P_{k_{3}}Q_{<k_{3}}\nabla^{-1}\psi_{3}\\&\hspace{2cm}I\Box[\nabla^{-1}P_{k_{2}}Q_{[k_{2}+(\epsilon-1)i,k_{2}+O(1)]}(\chi_{i}(t)\psi_{2})R_{j}P_{k_{1}}R_{\alpha}\delta\psi_{1}]]||_{N[k_{0}]}\\
&\lesssim
||P_{k_{3}}\psi_{3}\nabla^{-1}\psi_{3}||_{\dot{X}_{k_{3}}^{0,\frac{1}{2},1}+A[k_{3}]}\\&\hspace{1.5cm}||\nabla^{-1}P_{k_{2}}Q_{[k_{2}+(\epsilon-1)i,k_{2}+O(1)]}(\chi_{i}(t)\psi_{2})||_{\dot{X}_{k_{2}}^{0,\frac{1}{2},1}}
||P_{k_{1}}R_{\alpha}\delta\psi_{1}||_{S[k_{3}]}\lesssim 2^{-\mu i}\tilde{c}_{k_{0}}\\
\end{split}\end{equation}
Next we reduce the large-frequency input $P_{k_{1}}\delta\psi_{1}$
to modulation $<2^{-i\delta}$. For this, throw an operator
$Q_{\geq-i\delta}$ in front of $P_{k_{1}}\delta\psi_{1}$, let the
operator $\Box$ fall inside the inner square bracket and use the
proof of lemma~\ref{key}, placing the products
$P_{k_{3}}Q_{<k_{3}}\nabla^{-1}\psi_{3}\nabla^{-1}P_{k_{2}}Q_{<k_{2}+(\epsilon-1)i}(\chi_{i}(t)\psi_{2})$
etc. into $L_{t}^{2}L_{x}^{\infty}$ and $P_{k_{1}}Q_{\geq
-i\delta}\delta\psi_{1}$ into $L_{t}^{2}L_{x}^{2}$ and the output
into $L_{t}^{1}\dot{H}^{-1}$. Next, we claim that we may apply an
operator $Q_{\geq -i\delta+C}$ in front of the inner square
bracket. Indeed, if we apply an operator $Q_{<-i\delta+C}$, we
estimate using 3.4(a) that
\begin{equation}\nonumber\begin{split}
&||P_{k_{0}}Q_{<O(1)}
[P_{k_{3}}Q_{<k_{3}}\nabla^{-1}\psi_{3}\\&\hspace{2cm}I\Box Q_{<-i\delta+C}[\nabla^{-1}P_{k_{2}}Q_{<k_{2}+(\epsilon-1)i}(\chi_{i}(t)\psi_{2})R_{\alpha}P_{k_{1}}Q_{<-i\delta}\delta\psi_{1}]]||_{N[k_{0}]}\\
&\lesssim
||P_{k_{3}}\psi_{3}||_{\dot{X}_{k_{3}}^{0,\frac{1}{2},1}+A[k]}\\&\hspace{2cm}||\Box
Q_{<-i\delta+C}[\nabla^{-1}P_{k_{2}}Q_{<k_{2}+(\epsilon-1)i}(\chi_{i}(t)\psi_{2})R_{j}P_{k_{1}}Q_{<-i\delta}\delta\psi_{1}]||_{\dot{X}_{k_{1}}^{0,-\frac{1}{2},1}}\\
&\lesssim 2^{\frac{-i\delta-k_{2}}{4+}}\tilde{c}_{k_{1}},\\
\end{split}\end{equation}
which is acceptable provided $|k_{2}|$ is small enough in relation
to $i\delta$, which we may always arrange. Thus we now reduce to
estimating
\begin{equation}\nonumber\begin{split}
&P_{k_{0}}Q_{<O(1)}
[P_{k_{3}}Q_{<k_{3}}\nabla^{-1}\psi_{3}\\&\hspace{2cm}I\Box
Q_{\geq-i\delta+C}[\nabla^{-1}P_{k_{2}}Q_{<k_{2}+(\epsilon-1)i}(\chi_{i}(t)\psi_{2})R_{\alpha}P_{k_{1}}Q_{<-i\delta}\delta\psi_{1}]]\\
\end{split}\end{equation}
We note that on account of lemma~\ref{key}, we may
estimate\footnote{We need to assume that $\delta<<\epsilon$, which
we may always arrange.}
\begin{equation}\nonumber\begin{split}
&||P_{k_{0}}Q_{<O(1)}
[P_{k_{3}}Q_{<k_{3}}\nabla^{-1}\psi_{3}\\&I\Box
Q_{\geq-i\delta+C}[\nabla^{-1}\sum_{\kappa\in K_{\frac{\epsilon-1}{2}i}}\chi^{c}_{i,\mp\kappa}P_{k_{2},\kappa}Q^{\pm}_{<k_{2}+(\epsilon-1)i}(\chi_{i}(t)\psi_{2})R_{\alpha}P_{k_{1}}Q_{<-i\delta}\delta\psi_{1}]]||_{N[k_{0}]}\\
&\lesssim ||P_{k_{3}}\psi_{3}||_{S[k_{3}]}||I\Box
Q_{\geq-i\delta+C}[\nabla^{-1}\sum_{\kappa\in
K_{\frac{\epsilon-1}{2}i}}\chi^{c}_{i,\mp\kappa}P_{k_{2},\kappa}Q^{\pm}_{<k_{2}+(\epsilon-1)i}(\chi_{i}(t)\psi_{2})\\&\hspace{8.5cm}R_{\alpha}P_{k_{1}}Q_{<-i\delta}\delta\psi_{1}]]||_{\dot{X}_{k_{1}}^{0,-\frac{1}{2},1}},\\
\end{split}\end{equation}
and one has (see the arguments in \cite{Kr-4})
\begin{equation}\nonumber\begin{split}
&||P_{k_{1}+O(1)}Q_{\geq-i\delta+C}[\nabla^{-1}\sum_{\kappa\in
K_{\frac{\epsilon-1}{2}i}}\chi^{c}_{i,\mp\kappa}P_{k_{2},\kappa}Q^{\pm}_{<k_{2}+(\epsilon-1)i}(\chi_{i}(t)\psi_{2})\\&\hspace{8.5cm}R_{j}P_{k_{1}}Q_{<-i\delta}\delta\psi_{1}]||_{\dot{X}_{k_{1}}^{0,\frac{1}{2},1}},\\
&\lesssim(\sum_{\tilde{\kappa}\in
K_{\frac{-i\delta+C-k_{2}}{2}}}||\nabla^{-1}P_{k_{2},\tilde{\kappa}}\sum_{\kappa\in
K_{\frac{\epsilon-1}{2}i}}\chi^{c}_{i,\mp\kappa}P_{k_{2},\kappa}Q^{\pm}_{<k_{2}+(\epsilon-1)i}(\chi_{i}(t)\psi_{2})||_{PW[\pm\tilde{\kappa}]}^{2})^{\frac{1}{2}}
\\&\hspace{4.5cm} 2^{\frac{-\delta i-k_{2}}{4}}(\sum_{\tilde{\kappa}\in K_{\frac{-i\delta+C-k_{2}}{2}}}
||P_{k_{1},\tilde{\kappa}}Q^{\pm}_{<-i\delta}\delta\psi_{1}||_{NFA^{*}[\pm\tilde{\kappa}]}^{2})^{\frac{1}{2}}\\
&\lesssim 2^{-\mu i}\tilde{c}_{k_{1}}\\
\end{split}\end{equation}
Also, arguing as in the proof of lemma~\ref{key} and using the
same notation, we can estimate
\begin{equation}\nonumber\begin{split}
&||P_{k_{0}}Q_{<O(1)}
[P_{k_{3}}Q_{<k_{3}}\nabla^{-1}\psi_{3}\\&\hspace{1.5cm}I\Box
Q_{\geq-i\delta+C}[\phi_{i}(t,x)\nabla^{-1}P_{k_{2}}Q_{<k_{2}+(\epsilon-1)i}(\chi_{i}(t)\psi_{2})R_{\alpha}P_{k_{1}}Q_{<-i\delta}\delta\psi_{1}]]||_{L_{t}^{1}\dot{H}^{-1}}\\
&\lesssim
||\phi_{i}(t,x)P_{k_{3}}Q_{<k_{3}}\nabla^{-1}\psi_{3}||_{L_{t}^{2}L_{x}^{\infty}}||\phi_{i}(t,x)\nabla^{-1}P_{k_{2}}Q_{<k_{2}}\chi_{i}(t)\psi_{2}||_{L_{t}^{2}L_{x}^{\infty}}||P_{k_{1}}\delta\psi_{1}||_{L_{t}^{\infty}L_{x}^{2}}\\
&+||\phi_{i}(t,x)R_{\alpha}P_{k_{3}}Q_{<k_{3}}\nabla^{-1}\psi_{3}||_{L_{t}^{2}L_{x}^{\infty}}||\phi_{i}(t,x)||_{L_{t}^{M}L_{x}^{\infty}}\\&\hspace{4.5cm}||\nabla^{-1}P_{k_{2}}Q_{[k_{2}+(\epsilon-1)i,k_{2}]}(\chi_{i}(t)\psi_{2})||_{L_{t}^{2+}L_{x}^{\infty}}||P_{k_{1}}\delta\psi_{1}||_{L_{t}^{\infty}L_{x}^{2}}\\
&\lesssim 2^{-\mu i}\tilde{c}_{k_{1}}\\
\end{split}\end{equation}
This means we only need to estimate the following expression:
\begin{equation}\nonumber\begin{split}
&P_{k_{0}}Q_{<O(1)}
[P_{k_{3}}Q_{<k_{3}}\nabla^{-1}\psi_{3}\\&\hspace{1cm}I\Box
Q_{\geq-i\delta+C}[\nabla^{-1}\sum_{\kappa\in
K_{\frac{\epsilon-1}{2}i}}\chi_{i,\mp\kappa}P_{k_{2},\kappa}Q^{\pm}_{<k_{2}+(\epsilon-1)i}(\chi_{i}(t)\psi_{2})R_{\alpha}P_{k_{1}}Q_{<-i\delta}\delta\psi_{1}]]\\
\end{split}\end{equation}
Using 3.4(g), it suffices to estimate
\begin{equation}\nonumber
||I\Box Q_{\geq-i\delta+C}[\nabla^{-1}\sum_{\kappa\in
K_{\frac{\epsilon-1}{2}i}}\chi_{i,\mp\kappa}P_{k_{2},\kappa}Q^{\pm}_{<k_{2}+(\epsilon-1)i}(\chi_{i}(t)\psi_{2})R_{\alpha}P_{k_{1}}Q_{<-i\delta}\delta\psi_{1}]||_{\dot{X}_{k_{1}}^{0,-\frac{1}{2},1}}
\end{equation}
We have the identity\footnote{Again this requires that
$\delta<<\epsilon$.}
\begin{equation}\nonumber\begin{split}
&Q_{\geq-i\delta+C}[\nabla^{-1}\sum_{\kappa\in
K_{\frac{\epsilon-1}{2}i}}\chi_{i,\mp\kappa}\nabla^{-1}P_{k_{2},\kappa}Q^{\pm}_{<k_{2}+(\epsilon-1)i}(\psi_{2})R_{\alpha}P_{k_{1}}Q^{\pm}_{<-i\delta}\delta\psi_{1}]\\
&=\sum_{\tilde{\kappa}_{1,2}\in K_{\frac{-\delta
i-k_{2}}{2}},\text{dist}(\pm\kappa_{1},\pm\kappa_{2})\sim
2^{\frac{-i\delta-k_{2}}{2}}}\\&Q_{\geq-i\delta+C}[P_{k_{2},\tilde{\kappa_{1}}}\nabla^{-1}\sum_{\kappa\in
K_{\frac{\epsilon-1}{2}i}}\chi_{i,\mp\kappa}\nabla^{-1}P_{k_{2},\kappa}Q^{\pm}_{<k_{2}+(\epsilon-1)i}(\psi_{2})R_{\alpha}P_{k_{1},\tilde{\kappa_{2}}}Q^{\pm}_{<-i\delta}\delta\psi_{1}]
\end{split}\end{equation}
One can move the multiplier $\chi_{i,\kappa}$ to the outside of
this expression while generating error terms rapidly decaying in
$i$ and hence negligible for the argument. Hence exploiting
orthogonality we may now estimate
\begin{equation}\nonumber\begin{split}
&||Q_{\geq-i\delta+C}[\sum_{\kappa\in
K_{\frac{\epsilon-1}{2}i}}\chi_{i,\pm\kappa}\nabla^{-1}P_{k_{2},\kappa}Q^{\pm}_{<k_{2}+(\epsilon-1)i}(\psi_{2})R_{\alpha}P_{k_{1}}Q^{\pm}_{<-i\delta}\delta\psi_{1}]||_{L_{t}^{2}L_{x}^{2}}\\
&\lesssim (\sum_{\kappa\in
K_{\frac{\epsilon-1}{2}i}}[\sum_{\tilde{\kappa}_{1,2}\in
K_{\frac{-\delta
i-k_{2}}{2}},\text{dist}(\pm\tilde{\kappa}_{1},\pm\tilde{\kappa}_{2})\sim
2^{\frac{-i\delta-k_{2}}{2}},\kappa\cap\tilde{\kappa}_{1}\neq\emptyset}\\&||Q_{\geq-i\delta+C}[P_{k_{2},\pm\tilde{\kappa_{1}}}
\nabla^{-1}P_{k_{2},\kappa}Q^{\pm}_{<k_{2}+(\epsilon-1)i}\psi_{2}R_{\alpha}P_{k_{1},\tilde{\kappa_{2}}}Q^{\pm}_{<-i\delta}\delta\psi_{1}]||_{L_{t}^{2}L_{x}^{2}}]^{2})^{\frac{1}{2}}\\
\end{split}\end{equation}
On the other hand, using \eqref{bilinear1}, we can easily estimate
that
\begin{equation}\nonumber\begin{split}
&||Q_{\geq-i\delta+C}[P_{k_{2},\tilde{\kappa_{1}}}
\nabla^{-1}P_{k_{2},\kappa}Q^{\pm}_{<k_{2}+(\epsilon-1)i}\psi_{2}R_{\alpha}P_{k_{1},\tilde{\kappa_{2}}}Q^{\pm}_{<-i\delta}\delta\psi_{1}]||_{L_{t}^{2}L_{x}^{2}}\\
&\lesssim
2^{-\frac{-i\delta-k_{2}}{4}}||\nabla^{-1}P_{k_{2},\kappa}Q^{\pm}_{<k_{2}+(\epsilon-1)i}\psi_{2}||_{PW[\pm\kappa]}||P_{k_{1},\tilde{\kappa_{2}}}Q^{\pm}_{<-i\delta}\delta\psi_{1}||_{NFA[\pm\tilde{\kappa}_{2}]}\\
\end{split}\end{equation}
Plugging this back into the preceding inequality and using the
definition of $S[k,\kappa]$, we arrive at the following:
\begin{equation}\nonumber\begin{split}
&||Q_{\geq-i\delta+C}[\sum_{\kappa\in
K_{\frac{\epsilon-1}{2}i}}\chi_{i,\mp\kappa}\nabla^{-1}P_{k_{2},\kappa}Q^{\pm}_{<k_{2}+(\epsilon-1)i}(\psi_{2})R_{\alpha}P_{k_{1}}Q^{\pm}_{<-i\delta}\delta\psi_{1}]||_{L_{t}^{2}L_{x}^{2}}\\
&\lesssim
2^{\frac{\epsilon-1}{4}i}2^{-\frac{-i\delta-k_{2}}{4}}(\sum_{\kappa\in
K_{\frac{\epsilon-1}{2}i}}||\nabla^{-1}P_{k_{2},\kappa}Q^{\pm}_{<k_{2}+(\epsilon-1)i}\psi_{2}||_{S[k_{2},\pm\kappa_{2}]}^{2})^{\frac{1}{2}}
\\&\hspace{6cm}(\sum_{\tilde{\kappa}\in K_{\frac{-\delta i-k_{2}}{2}}}||P_{k_{1},\tilde{\kappa}}Q^{\pm}_{<-i\delta}\delta\psi_{1}||_{S[k_{1},\pm\tilde{\kappa}]}^{2})^{\frac{1}{2}}\\
\end{split}\end{equation}
It follows easily that this expression is $\lesssim 2^{-\mu
i}\tilde{c}_{k_{1}}$, which is as desired. \\
(I.c.5): This is much simpler since a derivative has been moved
from a high-frequency term to a low-frequency term. We use finite
propagation speed and Hoelder's inequality to conclude
\begin{equation}\nonumber\begin{split}
&||P_{k_{0}}Q_{<O(1)}[P_{k_{3}}Q_{<k_{3}}\nabla^{-1}\psi_{3}\nabla^{-1}\Box[P_{k_{2}}Q_{<k_{2}}\psi_{2}]R_{\alpha}P_{k_{1}}\delta\psi_{1}]||_{L_{t}^{1}\dot{H}^{-1}}\\
&\lesssim 2^{-k_{0}}2^{\delta
i}||P_{k_{3}}Q_{<k_{3}}\nabla^{-1}\psi_{3}||_{L_{t}^{2+}L_{x}^{\infty}}||\nabla^{-1}\Box[P_{k_{2}}Q_{<k_{2}}\psi_{2}]||_{L_{t}^{2}L_{x}^{\infty}}
||R_{\alpha}P_{k_{1}}\delta\psi_{1}||_{L_{t}^{M}L_{x}^{2+}}\\
&\lesssim 2^{\delta i-k_{0}}\tilde{c}_{k_{1}}\\
\end{split}\end{equation}
One can choose $\delta=\delta(M)$ arbitrarily small so one gets
the desired exponential gain in $i$.\\
(I.c.6): This is the expression
\begin{equation}\nonumber
\nabla^{-1}P_{k_{2}}Q_{<k_{2}+O(1)}\psi_{2}\nabla^{-1}P_{k_{3}}\psi_{3}\Box
P_{k_{1}}Q_{<O(1)}\delta\psi_{1}
\end{equation}
One estimates this using 3.4(g) as well as lemma~\ref{key}:
\begin{equation}\nonumber\begin{split}
&||P_{k_{0}}Q_{<O(1)}[\nabla^{-1}P_{k_{2}}Q_{<k_{2}+O(1)}\psi_{2}\nabla^{-1}P_{k_{3}}\psi_{3}\Box
P_{k_{1}}Q_{<O(1)}\delta\psi_{1}]||_{N[k_{0}]}\\
&\lesssim
2^{k-k_{1}}\sum_{k<\max\{k_{2},k_{3}\}+O(1)}||P_{k}(\nabla^{-1}P_{k_{2}}Q_{<k_{2}+O(1)}\psi_{2}\nabla^{-1}P_{k_{3}}\psi_{3})||_{\dot{X}_{k}^{0,\frac{1}{2},1}}\\&
\hspace{8cm}||\Box
P_{k_{1}}Q_{<O(1)}\delta\psi_{1}||_{\dot{X}_{k_{1}}^{0,-\frac{1}{2},1}}\\
&\lesssim 2^{-\mu i}\tilde{c}_{k_{1}}\\
\end{split}\end{equation}
The terms on the right-hand side of \eqref{null-form2} are treated
analogously. This concludes case (I.c).\\
(I.d): {\it{None of (I.a), (I.b), (I.c) hold.}} In that case we
may assume $|k_{1,2,3}|=O(1)$, $|k|=O(1)$. One proceeds exactly as
in the immediately preceding case (I.c), but has to argue
differently for case (I.d.5): In that case, we split the
expression into two manageable pieces:
\begin{equation}\nonumber\begin{split}
&||P_{k_{0}}
[P_{k_{3}}Q_{<k_{3}}\nabla^{-1}\psi_{3}\nabla^{-1}\Box[P_{k_{2}}Q_{[k_{2}+(\epsilon-1)i,k_{2}]}\chi_{i}(t)\psi_{2}]R_{\alpha}P_{k_{1}}Q_{<O(1)}\delta\psi_{1}]||_{N[k_{0}]}\\
&\lesssim
\sum_{k<\max\{k_{3},k_{1}\}+O(1)}2^{k-k_{2}}||R_{\alpha}P_{k_{3}}Q_{<k_{3}}\nabla^{-1}\psi_{3}R_{\alpha}P_{k_{1}}Q_{<k_{1}}\delta\psi_{1}||_{\dot{X}_{k}^{0,\frac{1}{2},1}}\\&\hspace{6.5cm}||\Box[P_{k_{2}}Q_{[k_{2}+(\epsilon-1)i,k_{2}]}\chi_{i}(t)\psi_{2}]||_{\dot{X}_{k_{2}}^{0,-\frac{1}{2},1}}\\
&\lesssim 2^{-\mu i}\tilde{c}_{k_{1}}\\
\end{split}\end{equation}
We have used here 3.4(g), 3.4(a), as well as lemma~\ref{key}.
Next, once $P_{k_{2}}\psi_{2}$ is reduced to very low modulation,
we can estimate
\begin{equation}\nonumber\begin{split}
&||P_{k_{0}}
[P_{k_{3}}Q_{<k_{3}}\nabla^{-1}\psi_{3}\nabla^{-1}\Box[P_{k_{2}}Q_{<k_{2}+(\epsilon-1)i}\chi_{i}(t)\psi_{2}]R_{\alpha}P_{k_{1}}Q_{<O(1)}\delta\psi_{1}]||_{L_{t}^{1}\dot{H}^{-1}}\\
&\lesssim
||P_{k_{3}}\nabla^{-1}\psi_{3}||_{L_{t}^{2+}L_{x}^{\infty}}||\nabla^{-1}\Box[P_{k_{2}}Q_{<k_{2}+(\epsilon-1)i}\chi_{i}(t)\psi_{2}]||_{L_{t}^{2}L_{x}^{2}}||P_{k_{1}}\delta\psi_{1}||_{L_{t}^{M}L_{x}^{2+}}\\
&\lesssim 2^{-\mu i}\tilde{c}_{k_{1}}\\
\end{split}\end{equation}
This concludes the treatment of case (I). \\
(II): {\it{The term {\bf{(III)}}.}} This is treated analogously to
case (I). One invokes the identity
\begin{equation}\nonumber\begin{split}
&\partial^{\nu}[R_{\nu}h\sum_{j=1,2}\triangle^{-1}\partial_{j}[R_{i}f
R_{j}g-R_{j}f R_{i}g]]\\
&=\frac{1}{2}\Box\nabla^{-1}h\sum_{j=1,2}\triangle^{-1}\partial_{j}[R_{i}f
R_{j}g-R_{j}f
R_{i}g]\\&+\frac{1}{2}\Box[\nabla^{-1}h\sum_{j=1,2}\triangle^{-1}\partial_{j}[R_{i}f
R_{j}g-R_{j}f
R_{i}g]]\\&-\frac{1}{2}[\nabla^{-1}h\sum_{j=1,2}\triangle^{-1}\partial_{j}\Box[R_{i}f
R_{j}g-R_{j}f R_{i}g]]\\
\end{split}\end{equation}
(III): {\it{The sum of terms {\bf{(II+VI)}}.}} We frequency
localize the expression and its inputs as in the preceding. If
both $P_{k_{1}}\psi_{1}$, $P_{k_{3}}\psi_{3}$ are of the first
type, we may assume that $|k_{j}|<<i$, $|k|<<i$, due to
lemma~\ref{auxiliary envelope} as well as the trilinear estimates
3.4(c). Now assume at least one of them is of the 2nd kind. Then
the estimate is straightforward on account of the strong Lebesgue
type estimates available: for example, assume that
$P_{k_{1}}\psi_{1}$ is of the 2nd type. We need to estimate the
schematically written expression\footnote{It is to be kept in mind
that the inner square bracket stands for a null-form of type
$Q_{\nu j}$.}
\begin{equation}\nonumber
\nabla_{x,t}P_{k_{0}}[P_{k_{1}}R_{\nu}\psi_{1}\nabla^{-1}P_{k}I[P_{k_{2}}\delta\psi_{2}P_{k_{3}}\psi_{3}]]
\end{equation}
In this case one has the estimate
\begin{equation}\nonumber\begin{split}
&||\nabla_{x,t}P_{k_{0}}[P_{k_{1}}R_{\nu}\psi_{1}\nabla^{-1}P_{k}I[P_{k_{2}}\delta\psi_{2}P_{k_{3}}\psi_{3}]]||_{N[k_{0}]}
\\&\lesssim
2^{-\delta_{1}|k_{1}-k_{0}|}2^{\delta_{2}[\min\{k,k_{2,3}\}-\max\{k,k_{2,3}\}]}\tilde{c}_{k_{0}}
\end{split}\end{equation}
This is a simple exercise with the exception of
the case when $\nu=0$ and  $P_{k_{1}}R_{0}\psi_{1}$ has large
modulation, i. e. we replace this expression by
$P_{k_{1}}Q_{>k_{1}+100}R_{0}\psi_{1}$. Assume that $k_{1}>>k$,
the other situations being similar or simpler. Then we have
\begin{equation}\nonumber\begin{split}
&||\nabla_{x,t}P_{k_{0}}[P_{k_{1}}Q_{>k_{1}+100}R_{0}\psi_{1}\nabla^{-1}P_{k}I[P_{k_{2}}\delta\psi_{2}P_{k_{3}}\psi_{3}]]||_{N[k_{0}]}
\\&\lesssim
2^{-\frac{k_{1}}{2}}||P_{k_{0}}[P_{k_{1}}Q_{>k_{1}+100}R_{0}\psi_{1}\nabla^{-1}P_{k}I[P_{k_{2}}\delta\psi_{2}P_{k_{3}}\psi_{3}]]||_{L_{t}^{2}L_{x}^{2}}\\
&\lesssim
2^{-\frac{k_{1}}{2}}||P_{k_{0}}[P_{k_{1}}Q_{>k_{1}+100}R_{0}\psi_{1}||_{L_{t}^{\infty}L_{x}^{2}}||\nabla^{-1}P_{k}I[P_{k_{2}}\delta\psi_{2}P_{k_{3}}\psi_{3}]||_{L_{t}^{2}L_{x}^{\infty}}
\\&\lesssim
2^{\frac{k-k_{1}}{2}}\epsilon 2^{\frac{\min\{k,k_{2},k_{3}\}-\max\{k,k_{2},k_{3}\}}{2}}\frac{\tilde{c}_{k_{2}}}{\epsilon}\frac{\tilde{c}_{k_{3}}}{\epsilon}.\\
\end{split}\end{equation}
Thus in that case, too, we may assume that
$\max_{j=0,\ldots,3}\{|k|,|k_{j}|\}<<i$. Next, we claim that we
may also reduce the modulations of all inputs to size $<1$.
Indeed, for example assume that the first input
$P_{k_{1}}\psi_{1}$ has modulation $>1$. Then we may estimate
(using schematic notation)
\begin{equation}\nonumber\begin{split}
&||P_{k_{0}}Q_{<k_{0}}\nabla_{x,t}\chi_{i}(t)[P_{k_{1}}Q_{>0}R_{\nu}\psi_{1}\nabla^{-1}P_{k}I[P_{k_{2}}\delta\psi_{2},P_{k_{3}}\psi_{3}]]||_{L_{t}^{1}\dot{H}^{-1}}\\
&\lesssim
||P_{k_{1}}Q_{[0,\max\{k_{2},k_{3}\}]}\psi_{1}||_{\dot{X}_{k_{1}}^{-\frac{1}{2},1,2}}||P_{k_{2}}\delta\psi_{2}||_{L_{t}^{M}L_{x}^{2+}}||P_{k_{3}}\chi_{i}(t)\psi_{3}||_{L_{t}^{2+}L_{x}^{\infty}}
\lesssim 2^{-\mu i}\tilde{c}_{k_{0}}\\
\end{split}\end{equation}
If, on the other hand, we replace $P_{k_{0}}Q_{<k_{0}}$ by
$P_{k_{0}}Q_{\geq k_{0}}$, we have to be careful when
$P_{k_{1}}\psi_{1}$ is of the 2nd kind. Then we first reduce both
$P_{k_{2}}\delta\psi_{2}$, $P_{k_{3}}\psi_{3}$ to modulation
$<2^{\delta i+O(1)}$ (which is straightforward), and estimate
\begin{equation}\nonumber\begin{split}
&||P_{k_{0}}Q_{\geq k_{0}}\nabla_{x,t}\chi_{i}(t)[P_{k_{1}}Q_{>0}R_{\nu}\psi_{1}\nabla^{-1}P_{k}I[P_{k_{2}}Q_{<\delta i}\delta\psi_{2},P_{k_{3}}Q_{<\delta i}\psi_{3}]]||_{\dot{X}_{k_{0}}^{-\frac{1}{2},-1,2}}\\
&\lesssim
2^{-k}||P_{k_{1}}Q_{>0}R_{\nu}\psi_{1}||_{L_{t}^{\infty}L_{x}^{2}}||\nabla_{x,t}\nabla^{-1}P_{k_{2}}Q_{<\delta
i}\delta\psi_{2}||_{L_{t}^{M}L_{x}^{2+}}\\&\hspace{7cm}||\chi_{i}(t)P_{k_{3}}Q_{<\delta
i}\nabla_{x,t}\nabla^{-1}\psi_{3}||_{L_{t}^{2+}L_{x}^{\infty}}\\
\end{split}\end{equation}
Assuming (as we may) that all the absolute frequencies
$|k|,\,|k_{i}|$ as well as $\delta i$ are much less than
$(\frac{1}{2}-\frac{1}{2+})i$ and summing over these frequency
ranges, one obtains from this and the definition of frequency
envelope the upper bound $\lesssim \tilde{c}_{k_{0}}2^{-\mu i}$,
which is as desired. The remaining cases are handled similarly.
Now one reverts to the expansions \eqref{null-form1},
\eqref{null-form2}. The only term requiring different treatment
than in case (I) is the 2nd term of the expansion. We record this
with the appropriate microlocalizations as follows:
\begin{equation}\nonumber
P_{k_{0}}\chi_{i}(t)[P_{k_{1}}Q_{<0}\psi_{1}\sum_{j=1,2}\triangle^{-1}\partial_{j}\Box
I P_{k}[\nabla^{-1}P_{k_{2}}\delta\psi_{2}P_{k_{3}}\psi_{3}]]
\end{equation}
As before, we may innocuously move the localizer $\chi_{i}(t)$
right in front of the inner square bracket $[,]$. Then we
decompose this term as a sum of manageable expressions, keeping in
mind lemma~\ref{key}, as well as the estimates in 3.4(a)-3.4(g):
first, we have
\begin{equation}\nonumber\begin{split}
&||P_{k_{0}}\chi_{i}(t)[P_{k_{1}}Q_{[k_{1}+(\epsilon-1)i,k_{1}]}\psi_{1}\sum_{j=1,2}\triangle^{-1}\partial_{j}\Box
I P_{k}[\nabla^{-1}P_{k_{2}}\delta\psi_{2}P_{k_{3}}\psi_{3}]]||_{N[k_{0}]}\\
&\lesssim
||P_{k_{1}}Q_{[k_{1}+(\epsilon-1)i,k_{1}]}\chi_{i}(t)\psi_{1}||_{\dot{X}_{k_{1}}^{0,\frac{1}{2},1}}
||I[\nabla^{-1}P_{k_{2}}\delta\psi_{2}P_{k_{3}}\psi_{3}]||_{\dot{X}_{k}^{0,\frac{1}{2},1}}
\lesssim 2^{-\mu i}\tilde{c}_{k_{0}}\\
\end{split}\end{equation}
Similarly, one estimates
\begin{equation}\nonumber\begin{split}
&||P_{k_{0}}\chi_{i}(t)[P_{k_{1}}Q_{<(\epsilon-1)i+k_{1}}\psi_{1}\\&\hspace{3cm}\sum_{j=1,2}\triangle^{-1}\partial_{j}\Box
I P_{k}[\nabla^{-1}P_{k_{2}}\delta\psi_{2}P_{k_{3}}Q_{[k_{3}+(\epsilon-1)i,k_{3}]}\psi_{3}]]||_{N[k_{0}]}\\
&\lesssim 2^{-\mu i}\tilde{c}_{k_{0}}\\
\end{split}\end{equation}
It is further easy to see that one may reduce the output to
modulation $<2^{-i\delta}$ for some very small\footnote{As usual,
we choose $\delta<<\epsilon$.} $\delta>0$ (in other words, throw
an operator $Q_{<-i\delta}$ in front of it), and also reduce the
inner square bracket $\Box I P_{k}[...]$ to modulation
$>2^{-i\delta +C}$. Further, one easily reduces the input
$P_{k_{2}}\delta\psi_{2}$ to modulation $<2^{-i\delta}$. One can
then rewrite the expression as follows:
\begin{equation}\label{technical100}\begin{split}
&P_{k_{0}}Q_{<-i\delta}[\chi_{i}(t)[P_{k_{1}}Q_{<(\epsilon-1)i+k_{1}}\psi_{1}\\&\hspace{2cm}\sum_{j=1,2}\triangle^{-1}\partial_{j}\Box
I
P_{k}Q_{>-i\delta+C}[\nabla^{-1}P_{k_{2}}Q_{<-i\delta}\delta\psi_{2}P_{k_{3}}Q_{<k_{3}+(\epsilon-1)i}\psi_{3}]]]\\
&=\sum_{\pm,\pm}\sum_{\kappa_{1,2}\in
K_{\frac{-i\delta-\min\{k_{1},k\}}{2}},\text{dist}(\pm\kappa_{1},\pm\kappa_{2})\sim
2^{\frac{-i\delta-\min\{k_{1},k\}}{2}}}\\&
P_{k_{0},\kappa_{1}}Q^{\pm}_{<-i\delta}\chi_{i}(t)[P_{k_{1},\kappa_{2}}Q^{\pm}_{<(\epsilon-1)i+k_{1}}\psi_{1}\\&\hspace{2cm}\sum_{j=1,2}\triangle^{-1}\partial_{j}\Box
I
P_{k}Q_{>-i\delta+C}[\nabla^{-1}P_{k_{2}}Q_{<-i\delta}\delta\psi_{2}P_{k_{3}}Q_{<k_{3}+(\epsilon-1)i}\psi_{3}]]\\
\end{split}\end{equation}
Now we have achieved the kind of situation in which the 2nd part
of lemma~\ref{key} becomes useful: indeed, we can estimate
\begin{equation}\nonumber\begin{split}
&||P_{k_{0}}Q_{<-i\delta}\chi_{i}(t)[\sum_{\kappa\in
K_{\frac{\epsilon-1}{2}i}}\chi^{c}_{i,\mp\kappa}P_{k_{1},\kappa}Q^{\pm}_{<(\epsilon-1)i+k_{1}}\psi_{1}\\&\hspace{1.6cm}\sum_{j=1,2}\triangle^{-1}\partial_{j}\Box
I
P_{k}Q_{>-i\delta+C}[\nabla^{-1}P_{k_{2}}Q_{<-\delta i}\delta\psi_{2}P_{k_{3}}Q_{<k_{3}+(\epsilon-1)i}\psi_{3}]]||_{N[k_{0}]}\\
&\lesssim (\sum_{\kappa_{1}\in
K_{\frac{-i\delta-\min\{k_{1},k\}}{2}}}||P_{k_{1},\kappa_{1}}\sum_{\kappa\in
K_{\frac{\epsilon-1}{2}i}}\chi^{c}_{i,\mp\kappa}P_{k_{1},\kappa}Q^{\pm}_{<(\epsilon-1)i+k_{1}}\psi_{1}||_{PW[\kappa_{1}]}^{2})^{\frac{1}{2}}
\\&\hspace{6cm}||P_{k_{2}}Q_{<-\delta i}\delta\psi_{2}||_{S[k_{2}]}||P_{k_{3}}\psi_{3}||_{A[k_{3}]+\dot{X}_{k_{3}}^{0,\frac{1}{2},1}}\\
\end{split}\end{equation}
This furnishes the desired estimate $\lesssim 2^{-\mu
i}\tilde{c}_{k_{0}}$. One similarly treats the contribution when
$P_{k_{3}}\psi_{3}$ is replaced by $\sum_{\pm}\sum_{\kappa\in
K_{\frac{\epsilon-1}{2}i}}\chi^{c}_{i,\mp\kappa}P_{k_{3},\kappa}Q^{\pm}_{<k_{3}+(\epsilon-1)i}\psi_{3}$.
Thus we now replace both $P_{k_{1,3}}\psi_{1,3}$ by
\begin{equation}\nonumber
\sum_{\pm}\sum_{\kappa\in
K_{\frac{\epsilon-1}{2}i}}\chi_{i,\mp\kappa}P_{k_{1,3}}Q^{\pm}_{<k_{1,3}+(\epsilon-1)i}\psi_{1,3},
\end{equation}
respectively. Observe that on account of the rapid decay
properties of the kernels of the multipliers $Q_{<>-i\delta}$
etc., we can rewrite our term (up to errors of order of magnitude
$2^{-Ni}$) as
\begin{equation}\nonumber\begin{split}
&\sum_{\pm,\pm}\sum_{\kappa_{1,2}\in
K_{\frac{\epsilon-1}{2}i},\text{dist}(\pm\kappa_{1},\pm\kappa_{2})\lesssim
2^{\frac{\epsilon-1}{2}i}}P_{k_{0},\kappa_{2}}Q^{\pm}_{<-i\delta}\chi_{i}(t)[\chi_{i,\mp\kappa_{1}}P_{k_{1},\kappa_{1}}Q^{\pm}_{<(\epsilon-1)i+k_{1}}\psi_{1}\\&\hspace{1.5cm}\sum_{j=1,2}\triangle^{-1}\partial_{j}\Box
I
P_{k}Q_{>-i\delta+C}[\nabla^{-1}P_{k_{2}}Q_{<-\delta i}\delta\psi_{2}\chi_{i,\mp\kappa_{2}}P_{k_{3},\kappa_{2}}Q^{\pm}_{<k_{3}+(\epsilon-1)i}\psi_{3}]]\\
\end{split}\end{equation}
We can now invoke \eqref{bilinear2} as well as
\eqref{technical100} in order to conclude that the preceding
expression is bounded with respect to $||.||_{N[k_{0}]}$ by
\begin{equation}\nonumber\begin{split}
&\lesssim 2^{\delta
i}2^{\frac{\min\{k_{i},k\}}{2}}(\sum_{\pm}\sum_{\kappa_{1}\in
K_{\frac{\epsilon-1}{2}i}}||\chi_{i,\mp\kappa_{1}}P_{k_{1},\kappa_{1}}Q^{\pm}_{<(\epsilon-1)i+k_{1}}\psi_{1}||_{PW[\pm\kappa_{1}]}^{2})^{\frac{1}{2}}
\\&\hspace{2cm}(\sum_{\pm}\sum_{\kappa_{3}\in
K_{\frac{\epsilon-1}{2}i}}||\chi_{i,\mp\kappa_{3}}P_{k_{1},\kappa_{3}}Q^{\pm}_{<(\epsilon-1)i+k_{1}}\psi_{1}||_{PW[\pm\kappa_{3}]}^{2})^{\frac{1}{2}}
||P_{k_{2}}\delta\psi_{2}||_{S[k_{2}]}\\
\end{split}\end{equation}
On the other hand, it is easily seen that
\begin{equation}\nonumber\begin{split}
&(\sum_{\kappa_{1}\in
K_{\frac{\epsilon-1}{2}i}}||\chi_{i,\mp\kappa_{1}}P_{k_{1},\kappa_{1}}Q^{\pm}_{<(\epsilon-1)i+k_{1}}\psi_{1}||_{PW[\pm\kappa_{1}]}^{2})^{\frac{1}{2}}
\lesssim 2^{\frac{\epsilon-1}{2}i}||P_{k_{1}}\psi_{1}||_{A[k_{1}]+\dot{X}_{k_{1}}^{0,\frac{1}{2},1}}.\\
\end{split}\end{equation}
This implies that the preceding expression may be bounded by
$\lesssim 2^{-\mu i}\tilde{c}_{k_{0}}$, as desired. This concludes
case (III).\\
(IV): {\it{The term {\bf{(IV)}}.}} This is similar to (II) and the
preceding case (III). Details are omitted. We are now done with
establishing Proposition~\ref{refined bootstrap} for the
trilinear null-forms linear in the perturbation $\delta\psi$.\\
\\

{\bf{(C): Estimating the trilinear null-forms at least quadratic
in $\delta\psi$.}} The claim here follows directly from 3.4(c) in
conjunction with lemma~\ref{auxiliary envelope}, provided any
function $P_{k}\psi_{\nu}$ present is of first type. The other
case is treated just like for the expressions linear in
$\delta\psi$. This completes the trilinear estimates.

\subsection{The quintilinear and higher order terms linear in the perturbation.}
 These turn out to be fairly simple to estimate on account of the
 favorable Strichartz type estimates available for the radial
 components $\psi_{\nu}$. We recall that these terms have the
 following schematic structure:
\begin{equation}\nonumber
{\bf{A}}(\delta\psi_{\nu},\psi_{\nu})=\nabla_{x,t}[\delta\psi\nabla^{-1}[R_{\nu}\psi\nabla^{-1}[\psi\nabla^{-1}(\psi^{2})]]]
\end{equation}
\begin{equation}\nonumber
{\bf{B}}(\delta\psi_{\nu},\psi_{\nu})=\nabla_{x,t}[\psi\nabla^{-1}[R_{\nu}\delta\psi\nabla^{-1}[\psi\nabla^{-1}(\psi^{2})]]]
\end{equation}
\begin{equation}\nonumber
{\bf{C}}(\delta\psi_{\nu},\psi_{\nu})=\nabla_{x,t}[\psi\nabla^{-1}[R_{\nu}\psi\nabla^{-1}[\delta\psi\nabla^{-1}(\psi^{2})]]]
\end{equation}
\begin{equation}\nonumber
{\bf{D}}(\delta\psi_{\nu},\psi_{\nu})=\nabla_{x,t}[\psi\nabla^{-1}[R_{\nu}\psi\nabla^{-1}[\psi\nabla^{-1}(\delta\psi\psi)]]]
\end{equation}
\begin{equation}\nonumber
{\bf{E}}(\delta\psi_{\nu},\psi_{\nu})=\nabla_{x,t}[\nabla^{-1}[\psi\nabla^{-1}(\psi^{2})]\nabla^{-1}IQ_{\nu
j}(\delta\psi,\psi)]\,\text{etc}
\end{equation}
\begin{equation}\nonumber
{\bf{F}}(\delta\psi_{\nu},\psi_{\nu})=\nabla_{x,t}[\delta\psi\nabla^{-1}[\nabla^{-1}(\psi\nabla^{-1}(\psi^{2}))\nabla^{-1}[\psi\nabla^{-1}(\psi^{2})]]]\,\text{etc}
\end{equation}
We have left out the other terms like ${\bf{E}}$, ${\bf{F}}$, in
which $\delta\psi$ gets shifted to different positions. We treat
here the first quintilinear term in the list, the other terms
being tedious reiterations of similar computations.
\\
{\bf{(A):}} We microlocalize this term as follows:
\begin{equation}\nonumber
\nabla_{x,t}P_{k_{0}}[P_{k_{1}}\delta\psi\nabla^{-1}P_{k_{2}}[P_{k_{3}}R_{\nu}\psi\nabla^{-1}P_{k_{4}}[P_{k_{5}}\psi\nabla^{-1}P_{k_{6}}(P_{k_{7}}\psi
P_{k_{8}}\psi)]]]
\end{equation}
We first dispose of the case when $\nu=0$ $R_{\nu}\psi$ has
elliptic microsupport\footnote{When $\nu\neq 0$, this distinction
is irrelevant.}, i. e. is replaced by $R_{0}P_{k_{3}}Q_{\geq
k_{3}+100}\psi$. Recall that the expression
\begin{equation}\nonumber
\nabla^{-1}P_{k_{4}}[P_{k_{5}}\psi\nabla^{-1}P_{k_{6}}(P_{k_{7}}\psi
P_{k_{8}}\psi)]
\end{equation}
arises upon substituting the elliptic part for an input of the
form $P_{k_{4}}\psi_{j}$, where $j\neq 0$. We then rewrite this
expression as $P_{k_{4}}\psi_{j}-R_{j}\sum_{i=1,2}R_{i}\psi_{i}$,
and are led to a term of the schematic form
\begin{equation}\nonumber
\nabla_{x,t}P_{k_{0}}[P_{k_{1}}\delta\psi\nabla^{-1}P_{k_{2}}[P_{k_{3}}Q_{\geq
k_{3}+100}R_{0}\psi P_{k_{4}}\psi]]
\end{equation}
We claim the estimate (under our bootstrap assumption)
\begin{equation}\nonumber\begin{split}
&\nabla_{x,t}P_{k_{0}}[P_{k_{1}}\delta\psi\nabla^{-1}P_{k_{2}}[P_{k_{3}}Q_{\geq
k_{3}+100}R_{0}\psi P_{k_{4}}\psi]]\\
&\lesssim
2^{-\delta_{1}|k_{1}-k_{0}|}2^{\delta_{2}[\min_{i=2,3,4}\{k_{i}\}-\max_{i=2,3,4}\{k_{i}\}]}\frac{\tilde{c}_{k_{4}}}{\epsilon}\tilde{c}_{k_{0}}
\end{split}\end{equation}
As usual, we may assume $k_{0}=0$. We first treat the case when
$k_{2}<-10$. Then consider the cases\\
{\bf{(a)}}: $k_{2}\in [k_{3}-10,k_{3}+10]$ whence $k_{4}\leq
k_{3}+O(1)$. Either $P_{k_{4}}\psi$ has modulation at least
comparable to that of $P_{k_{3}}Q_{\geq k_{3}}\psi$, or else
either $P_{k_{1}}\delta\psi$ or the output has modulation at least
comparable. These are all similar, so we treat the first
possibility. Write this contribution as
\begin{equation}\nonumber\begin{split}
&\sum_{a\geq
k_{3}+100}\nabla_{x,t}P_{k_{0}}[P_{k_{1}}\delta\psi\nabla^{-1}P_{k_{2}}[P_{k_{3}}Q_{a}R_{0}\psi
P_{k_{4}}Q_{\geq a+O(1)}\psi]]\\
\end{split}\end{equation}
Then we can estimate
\begin{equation}\nonumber\begin{split}
&||\sum_{a\geq
k_{3}+100}\nabla_{x,t}P_{k_{0}}Q_{<k_{0}}[P_{k_{1}}\delta\psi\nabla^{-1}P_{k_{2}}[P_{k_{3}}Q_{a}R_{0}\psi
P_{k_{4}}Q_{\geq a+O(1)}\psi]]||_{L_{t}^{1}\dot{H}^{-1}}\\
&\lesssim \sum_{a\geq
k_{3}+100}||P_{k_{1}}\delta\psi||_{L_{t}^{\infty}L_{x}^{2}}||P_{k_{3}}Q_{a}R_{0}\psi||_{L_{t}^{2}L_{x}^{2}}||P_{k_{4}}Q_{\geq
a+O(1)}\psi||_{L_{t}^{2}L_{x}^{\infty}}\\
&\lesssim \sum_{a\geq
k_{3}+100}2^{-\frac{k_{3}}{2}}2^{(1-\epsilon)(k_{4}-a)}2^{\frac{k_{4}}{2}}\tilde{c}_{k_{1}}\frac{\tilde{c}_{k_{3}}}{\epsilon}\frac{\tilde{c}_{k_{4}}}{\epsilon},\\
\end{split}\end{equation}
irrespective of whether $P_{k_{3,4}}\psi$ is of the first or 2nd
type. This is an acceptable bound. Similarly, we have
\begin{equation}\nonumber\begin{split}
&||\sum_{a\geq k_{3}+100}\nabla_{x,t}P_{k_{0}}Q_{\geq
k_{0}}[P_{k_{1}}\delta\psi\nabla^{-1}P_{k_{2}}[P_{k_{3}}Q_{a}R_{0}\psi
P_{k_{4}}Q_{\geq a+O(1)}\psi]]||_{\dot{X}_{k_{0}}^{-\frac{1}{2},-1,2}}\\
&\lesssim\sum_{a\geq k_{3}+100}
2^{-\frac{k_{0}}{2}}||P_{k_{1}}\delta\psi||_{L_{t}^{\infty}L_{x}^{2}}||P_{k_{3}}Q_{a}R_{0}\psi||_{L_{t}^{\infty}L_{x}^{2}}
||P_{k_{4}}Q_{\geq a+O(1)}\psi||_{L_{t}^{2}L_{x}^{\infty}}\\
&\lesssim \sum_{a\geq
k_{3}+100}2^{\frac{k_{4}-k_{0}}{2}}2^{(1-\epsilon)(k_{4}-a)}\frac{\tilde{c}_{k_{4}}}{\epsilon}\tilde{c}_{k_{0}},\\
\end{split}\end{equation}
again an acceptable bound.\\
{\bf{(b)}}: $k_{2}<k_{3}-10$, whence $k_{4}=k_{3}+O(1)$.  One
proceeds as in the preceding case. Again treating the case when
$P_{k_{4}}\psi_{4}$ is at modulation at least comparable to that
of $P_{k_{3}}\psi_{3}$, we can estimate
\begin{equation}\nonumber\begin{split}
&||\sum_{a\geq
k_{3}+100}\nabla_{x,t}P_{k_{0}}Q_{<k_{0}}[P_{k_{1}}\delta\psi\nabla^{-1}P_{k_{2}}[P_{k_{3}}Q_{a}R_{0}\psi
P_{k_{4}}Q_{\geq a+O(1)}\psi]]||_{L_{t}^{1}\dot{H}^{-1}}\\
&2^{k_{2}}||P_{k_{1}}\delta\psi||_{L_{t}^{\infty}L_{x}^{2}}||P_{k_{3}}Q_{a}R_{0}\psi||_{L_{t}^{2}L_{x}^{2}}||P_{k_{4}}Q_{\geq
a+O(1)}\psi||_{L_{t}^{2}L_{x}^{2}}\\
&\lesssim\sum_{a\geq k_{3}+100}
2^{k_{2}-k_{3}}2^{(1-\epsilon)(k_{3}-a)}\frac{\tilde{c}_{k_{3}}}{\epsilon}\frac{\tilde{c}_{k_{4}}}{\epsilon}\tilde{c}_{k_{0}},\\
\end{split}\end{equation}
which leads to an acceptable estimate. The contribution of $P_{k_{0}}Q_{\geq k_{0}}$ is treated similarly.\\
{\bf{(c)}}: The case $k_{3}\leq k_{2}-10$. Rewrite the term as
\begin{equation}\nonumber\begin{split}
&\nabla_{x,t}P_{k_{0}}[P_{k_{1}}\delta\psi_{1}\nabla^{-1}P_{k_{4}}\psi
P_{k_{3}}R_{0}(1-I)\psi]=\\&\hspace{5cm}\sum_{a>k_{3}+100}\nabla_{x,t}P_{k_{0}}[P_{k_{1}}\delta\psi_{1}\nabla^{-1}P_{k_{4}}\psi
P_{k_{3}}R_{0}Q_{a}\psi]\\
\end{split}\end{equation}
Then one estimates
\begin{equation}\nonumber\begin{split}
&||\nabla_{x,t}P_{k_{0}}Q_{<k_{0}}[Q_{\geq
a-10}[P_{k_{1}}\delta\psi_{1}\nabla^{-1}P_{k_{4}}\psi]
P_{k_{3}}R_{0}Q_{a}\psi]||_{L_{t}^{1}\dot{H}^{-1}}\\
&\lesssim
||Q_{\geq a-10}[P_{k_{1}}\delta\psi_{1}\nabla^{-1}P_{k_{4}}\psi]||_{L_{t}^{2}L_{x}^{2}}||P_{k_{3}}R_{0}Q_{a}\psi||_{L_{t}^{2}L_{x}^{\infty}}\\
&\lesssim
2^{-\frac{a}{2}}2^{\frac{\min\{a-k_{4},0\}}{4+}}||P_{k_{1}}\delta\psi_{1}||_{S[k_{1}]}||P_{k_{4}}\psi_{4}||_{\dot{X}_{k_{4}}^{0,\frac{1}{2},1}+A[k_{4}]}
\\&\hspace{6cm}2^{\frac{k_{3}}{2}}2^{\mu(a-k_{3})}||P_{k_{3}}R_{0}Q_{a}\psi||_{\dot{X}_{k_{4}}^{-(\frac{1}{2}-\mu),1-\mu,1}}\\
&\lesssim
2^{\min\{\frac{a-k_{4}}{4+},0\}}2^{\frac{k_{3}-a}{2+}}\tilde{c}_{k_{1}}\frac{\tilde{c}_{k_{4}}}{\epsilon}\\
\end{split}\end{equation}
One can sum over $a>k_{3}+100$ to obtain the desired estimate.
Similarly, we have
\begin{equation}\nonumber\begin{split}
&||\nabla_{x,t}P_{k_{0}}Q_{<k_{0}}[Q_{<a-10}[P_{k_{1}}\delta\psi_{1}\nabla^{-1}P_{k_{4}}\psi]
P_{k_{3}}R_{0}Q_{a}\psi]||_{\dot{X}_{k_{0}}^{-1,-\frac{1}{2},1}}\\
&\lesssim
2^{-\frac{a}{2}}||Q_{<a-10}[P_{k_{1}}\delta\psi_{1}\nabla^{-1}P_{k_{4}}\psi]||_{L_{t}^{\infty}L_{x}^{2}}
||P_{k_{3}}R_{0}Q_{a}\psi||_{L_{t}^{2}L_{x}^{\infty}}\\
&\lesssim
2^{\min\{\frac{a-k_{4}}{4+},0\}}2^{\frac{k_{3}-a}{2+}}\frac{\tilde{c}_{k_{3}}}{\epsilon}\tilde{c}_{k_{1}}\\
\end{split}\end{equation}
Next, as to the large modulation contribution of the output, we
estimate using theorem~\ref{Moser}
\begin{equation}\nonumber\begin{split}
&||\nabla_{x,t}P_{k_{0}}Q_{\geq
k_{0}}[P_{k_{1}}\delta\psi_{1}\nabla^{-1}P_{k_{4}}\psi
P_{k_{3}}R_{0}Q_{a}\psi]||_{\dot{X}_{k_{0}}^{-\frac{1}{2},-1,2}}\\
&\lesssim
2^{-\frac{k_{0}}{2}}||P_{k_{1}}\delta\psi_{1}||_{L_{t}^{\infty}L_{x}^{2}}||\nabla^{-1}P_{k_{4}}\psi
P_{k_{3}}R_{0}Q_{a}\psi||_{L_{t}^{2}L_{x}^{\infty}}\\
&\lesssim
2^{\delta(k_{3}-k_{4})}2^{\frac{k_{4}-k_{0}}{2}}\tilde{c}_{k_{1}}\frac{\tilde{c}_{k_{4}}}{\epsilon}\\
\end{split}\end{equation}
Now consider the case $k_{2}>10$, whence $k_{1}=k_{2}+O(1)$. Again
this only requires trivial modifications: for example, one
estimates in case $k_{3}=k_{4}+O(1)$
\begin{equation}\nonumber\begin{split}
&||\sum_{a\geq
k_{3}+100}\nabla_{x,t}P_{k_{0}}Q_{<k_{0}}[P_{k_{1}}\delta\psi\nabla^{-1}P_{k_{2}}[P_{k_{3}}Q_{a}R_{0}\psi
P_{k_{4}}Q_{\geq a+O(1)}\psi]]||_{L_{t}^{1}\dot{H}^{-1}}\\
&\lesssim\sum_{a\geq k_{3}+100}
2^{-k_{2}}||P_{k_{1}}\delta\psi||_{L_{t}^{\infty}L_{x}^{2}}||P_{k_{3}}Q_{a}R_{0}\psi||_{L_{t}^{2}L_{x}^{2}}||P_{k_{4}}Q_{\geq
a+O(1)}\psi||_{L_{t}^{2}L_{x}^{2}}\\
&\lesssim\sum_{a\geq
k_{3}+100}2^{-k_{3}}2^{(1-\epsilon)(k_{3}-a)}\tilde{c}_{k_{1}}\frac{\tilde{c}_{k_{4}}}{\epsilon}\lesssim 2^{-(1-\epsilon)k_{1}}2^{(1-\epsilon)(k_{1}-k_{3})}\tilde{c}_{k_{0}}\frac{\tilde{c}_{k_{4}}}{\epsilon},\\
\end{split}\end{equation}
and this is again an acceptable estimate. The remaining cases are
more of the same. We conclude from this\footnote{Applying the
multiplier $\chi_{i}(t)$ doesn't affect the estimate for reasons
explained earlier.} that the estimate
\begin{equation}\nonumber
||\chi_{i}(t)P_{k_{0}}[P_{k_{1}}\delta\psi
\nabla^{-1}P_{k_{2}}[P_{k_{3}}(1-I)R_{0}\psi
P_{k_{4}}\psi]]||_{N[k_{0}]}\lesssim 2^{-\mu i}\tilde{c}_{k_{0}}
\end{equation}
holds, provided we have $i\lesssim \max_{i=2,3,4}\{|k_{i}|\}$. Now
assume that $\max\{|k_{2,3,4}|\}<<i$. We shall treat this quantity
as $O(1)$. Also, assume $i\lesssim k_{1}$. We also omit the
operator $\nabla^{-1}P_{k_{2}}$ for simplicity's sake, and
estimate
\begin{equation}\nonumber\begin{split}
&||\nabla_{x,t}P_{k_{0}}Q_{>k_{0}}[P_{k_{1}}\delta\psi
P_{k_{4}}\psi R_{0}(1-I)\psi||_{N[k_{0}]}\\&\lesssim
2^{-\frac{k_{0}}{2}}||P_{k_{1}}\delta\psi||_{L_{t}^{\infty}L_{x}^{2}}||P_{k_{4}}\psi
R_{0}(1-I)\psi||_{L_{t}^{2}L_{x}^{2}}\lesssim
2^{-\frac{k_{0}}{2}}2^{\frac{k_{4}}{2}}\tilde{c}_{k_{1}}\frac{\tilde{c}_{k_{4}}}{\epsilon}\\
\end{split}\end{equation}
This yields the desired exponential gain in $i$. The contribution
of the hyperbolic part is unfortunately a bit more complicated.
First, observe that if $a>\delta i$, one estimates
\begin{equation}\nonumber\begin{split}
&||\nabla_{x,t}P_{k_{0}}Q_{\leq
k_{0}}[Q_{<a-10}[P_{k_{1}}\delta\psi P_{k_{4}}\psi]
R_{0}Q_{a}\psi||_{N[k_{0}]}\\&\lesssim
2^{-\frac{a}{2}}||P_{k_{1}}\delta\psi||_{L_{t}^{\infty}L_{x}^{2}}||P_{k_{4}}\psi||_{L_{t}^{\infty}L_{x}^{\infty}}||R_{0}Q_{a}P_{k_{3}}\psi||_{L_{t}^{2}L_{x}^{\infty}}
\lesssim
2^{-\frac{a}{2+}}\tilde{c}_{k_{1}}\frac{\tilde{c}_{k_{4}}}{\epsilon},
\end{split}\end{equation}
which upon summing over $a>\delta i$ results in the desired
exponential gain in $i$. The case when $Q_{<a-10}$ is replaced by
$Q_{\geq a-10}$ is similar (place the output into
$L_{t}^{1}\dot{H}^{-1}$). Now consider
\begin{equation}\nonumber
\nabla_{x,t}P_{k_{0}}Q_{\leq
k_{0}}\chi_{i}(t)[[P_{k_{1}}\delta\psi P_{k_{4}}\psi]
R_{0}(1-I)P_{k_{3}}Q_{<\delta i}\psi]
\end{equation}
We have
\begin{equation}\nonumber\begin{split}
&||\nabla_{x,t}P_{k_{0}}Q_{\leq
k_{0}}\chi_{i}(t)[Q_{<a-10}[P_{k_{1}}\delta\psi P_{k_{4}}\psi]
R_{0}(1-I)P_{k_{3}}Q_{a}\psi]||_{N[k_{0}]}\\
&\lesssim
2^{-\frac{a}{2}}||P_{k_{1}}\delta\psi||_{L_{t}^{\infty}L_{x}^{2}}||\chi_{i}(t)P_{k_{4}}\psi_{4}||_{L_{t}^{\infty}L_{x}^{\infty}}||R_{0}(1-I)P_{k_{3}}Q_{a}\psi||_{L_{t}^{2}L_{x}^{2}}\\
&\lesssim
2^{-\frac{i}{2+}}\tilde{c}_{k_{1}}\frac{\tilde{c}_{k_{3}}}{\epsilon}\\
\end{split}\end{equation}
Summing over $a<\delta i$ yields an acceptable estimate. Now
consider for $k_{3}+100<a<\delta i$
\begin{equation}\nonumber
\nabla_{x,t}P_{k_{0}}Q_{\leq k_{0}}\chi_{i}(t)[Q_{\geq
a-10}[P_{k_{1}}\delta\psi P_{k_{4}}\psi]
R_{0}(1-I)P_{k_{3}}Q_{a}\psi]
\end{equation}
We first reduce $P_{k_{1}}\delta\psi_{1}$ to modulation
$<2^{a-100}$, which is straightforward. We estimate
\begin{equation}\nonumber\begin{split}
&||\nabla_{x,t}P_{k_{0}}Q_{\leq k_{0}}[Q_{\geq
a-10}[P_{k_{1}}Q_{<a-100}\delta\psi
\\&\hspace{4cm}P_{k_{4}}Q_{>k+(\epsilon-1)i}(\chi_{i}(t)\psi)]
R_{0}(1-I)P_{k_{3}}Q_{a}\psi]||_{L_{t}^{1}\dot{H}^{-1}}\\
&\lesssim 2^{-\frac{a}{2}}||P_{k_{1}}Q_{<a-100}\delta\psi
P_{k_{4}}Q_{>k+(\epsilon-1)i}(\chi_{i}(t)\psi)]||_{\dot{X}_{k_{1}}^{0,\frac{1}{2},\infty}}||R_{0}(1-I)P_{k_{3}}Q_{a}\psi||_{L_{t}^{2}L_{x}^{2}}\\
&\lesssim
2^{-\frac{a}{2}}||P_{k_{1}}Q_{<a-100}\delta\psi||_{S[k_{1}]}||P_{k_{4}}Q_{>k+(\epsilon-1)i}(\chi_{i}(t)\psi)||_{\dot{X}_{k_{4}}^{0,\frac{1}{2},1}}
||R_{0}(1-I)P_{k_{3}}Q_{a}\psi||_{L_{t}^{2}L_{x}^{2}}\\&\lesssim 2^{-\frac{a}{2+}}2^{-\mu i}\tilde{c}_{k_{1}},\\
\end{split}\end{equation}
using lemma~\ref{key}. Next, borrowing notation form the proof of
lemma~\ref{key}, estimate
\begin{equation}\nonumber\begin{split}
&||\nabla_{x,t}P_{k_{0}}Q_{\leq k_{0}}[Q_{\geq
a-10}[P_{k_{1}}Q_{<a-100}\delta\psi
\phi_{i}(t,x)\\&\hspace{4cm}P_{k_{4}}Q_{<k+(\epsilon-1)i}(\chi_{i}(t)\psi)]
R_{0}(1-I)P_{k_{3}}Q_{a}\psi]||_{L_{t}^{1}\dot{H}^{-1}}\\
&\lesssim
||P_{k_{1}}Q_{<a-100}\delta\psi||_{L_{t}^{\infty}L_{x}^{2}}||\phi_{i}(t,x)P_{k_{4}}\psi||_{L_{t}^{2}L_{x}^{\infty}}||R_{0}(1-I)P_{k_{3}}Q_{a}\psi||_{L_{t}^{2}L_{x}^{2}}\\
&+||Q_{\geq a-10}[P_{k_{1}}Q_{<a-100}\delta\psi \phi_{i}(t,x)
Q_{>k+(\epsilon-1)i}\psi]||_{L_{t}^{2}L_{x}^{2}}||R_{0}(1-I)P_{k_{3}}Q_{a}\psi||_{L_{t}^{2}L_{x}^{2}}\\
\end{split}\end{equation}
The 2nd summand may be estimated as in the immediately preceding
since we may assume that $|a|, |k|$ etc are $<<\epsilon i$. As to
the first summand, we estimate it using
\begin{equation}\nonumber
||\phi_{i}(t)P_{k_{4}}\psi||_{L_{t}^{2}L_{x}^{\infty}}\lesssim
2^{-\mu i},
\end{equation}
from which the desired estimate follows. Now one decomposes
\begin{equation}\label{decomposition40}\begin{split}
&(1-\phi_{i}(t,x))P_{k_{4}}Q_{<k+(\epsilon-1)i}(\chi_{i}(t)\psi)
=\sum_{\pm}\sum_{\kappa\in
K_{\frac{\epsilon-1}{2}i}}\chi^{c}_{i,\mp\kappa}P_{k_{4},\kappa}Q^{\pm}_{<k+(\epsilon-1)i}(\chi_{i}(t)\psi)\\
&\hspace{4cm}+\sum_{\pm}\sum_{\kappa\in
K_{\frac{\epsilon-1}{2}i}}\chi_{i,\mp\kappa}P_{k_{4},\kappa}Q^{\pm}_{<k+(\epsilon-1)i}(\chi_{i}(t)\psi),\\
\end{split}\end{equation}
and proceeds as in the trilinear estimates: plugging in the first
summand on the right, we have
\begin{equation}\nonumber\begin{split}
&||Q_{\geq a-10}[P_{k_{1}}Q_{<a-100}\delta\psi
\sum_{\pm}\sum_{\kappa\in
K_{\frac{\epsilon-1}{2}i}}\chi^{c}_{i,\mp\kappa}P_{k_{4},\kappa}Q^{\pm}_{<k+(\epsilon-1)i}(\chi_{i}(t)\psi)]||_{\dot{X}_{k_{1}}^{0,\frac{1}{2},\infty}}\\
&\lesssim (\sum_{\tilde{\kappa_{1}}\in
K_{\frac{a-k_{4}}{2}}}||P_{k_{1},
\tilde{\kappa}_{1}}Q^{\pm}_{<a-100}\delta\psi||_{NFA^{*}[\pm\kappa_{1}]}^{2})^{\frac{1}{2}}
\\&\hspace{1cm}(\sum_{\tilde{\kappa}_{2}\in
K_{\frac{a-k_{4}}{2}}}||P_{k_{4},\tilde{\kappa_{4}}}\sum_{\kappa\in
K_{\frac{\epsilon-1}{2}i}}\chi^{c}_{i,\mp\kappa}P_{k_{4},\kappa}Q^{\pm}_{<k+(\epsilon-1)i}(\chi_{i}(t)\psi)||_{PW[\pm\tilde{\kappa}_{4}]}^{2})^{\frac{1}{2}}\\
\end{split}\end{equation}
This can be estimated by $\lesssim 2^{-\mu i}\tilde{c}_{k_{1}}$,
as desired. Plugging in the 2nd part of \eqref{decomposition40} is
handled as in the trilinear estimates, exploiting the
orthogonality of these pieces. This gives the desired estimate.
The case when $k_{1}<<i$ is more elementary. One may always assume
that $R_{0}\psi_{3}$ lives at modulation $<\delta i$ (argue as
before), whence one may estimate
\begin{equation}\nonumber\begin{split}
&||\nabla_{x,t}P_{k_{0}}Q_{\leq
k_{0}}\chi_{i}(t)[[P_{k_{1}}\delta\psi P_{k_{4}}\psi]
R_{0}(1-I)P_{k_{3}}Q_{<\delta i}\psi]||_{L_{t}^{1}\dot{H}^{-1}}\\
&\lesssim
||P_{k_{1}}\delta\psi||_{L_{t}^{M}L_{x}^{2+}}||P_{k_{4}}\chi_{i}(t)\psi||_{L_{t}^{2+}L_{x}^{\infty}}||R_{0}(1-I)P_{k_{3}}Q_{<\delta
i}\psi||_{L_{t}^{2}L_{x}^{M}}\\
\end{split}\end{equation}
Choosing $\frac{1}{2}-\frac{1}{2+}>>\delta$ results in the desired
estimate $\lesssim 2^{-\mu i}\tilde{c}_{k_{1}}$. Also, we have
\begin{equation}\nonumber\begin{split}
&||\nabla_{x,t}P_{k_{0}}Q_{>
k_{0}}\chi_{i}(t)[[P_{k_{1}}\delta\psi P_{k_{4}}\psi]
R_{0}(1-I)P_{k_{3}}Q_{<\delta i}\psi]||_{\dot{X}_{k_{0}}^{-\frac{1}{2},-1,2}}\\
&\lesssim
||P_{k_{1}}\delta\psi||_{L_{t}^{\infty}L_{x}^{2}}||\chi_{i}(t)P_{k_{4}}\psi||_{L_{t}^{\infty}L_{x}^{\infty}}
||R_{0}(1-I)P_{k_{3}}Q_{<\delta
i}\psi||_{L_{t}^{2}L_{x}^{2}}\lesssim
2^{-\frac{i}{2+}}\tilde{c}_{k_{1}}\\
\end{split}\end{equation}
which is as desired. This finally concludes estimating the
contribution when $P_{k_{3}}R_{\nu}\psi$ is replaced by
$P_{k_{3}}(1-I)R_{\nu}\psi$. Thus we may replace
$P_{k_{3}}R_{\nu}\psi$ by $P_{k_{3}}IR_{\nu}\psi$, whence we may
discard $IR_{\nu}$ for all intents and purposes. We revert to the
original formulation of this term given by
${\bf{A}}(\delta\psi,\psi)$, and decompose the innermost bracket
$(\psi^{2})$ into a $Q_{\nu j}$-type null-form as well as error
terms at least quadrilinear by means of our standard Hodge-type
decompositions, i. e. we write (using schematic notation)
\begin{equation}\nonumber\begin{split}
&(\psi^{2})=
R_{\nu}\psi^{1}R_{j}\psi^{2}-R_{j}\psi^{1}R_{\nu}\psi^{2}
+\nabla^{-1}(\psi\nabla^{-1}(\psi^{2}))R_{\nu}\psi\\&\hspace{6cm}+\nabla^{-1}(\psi\nabla^{-1}(\psi^{2}))\nabla^{-1}(\psi\nabla^{-1}(\psi^{2})).\\
\end{split}\end{equation}
Substituting the $Q_{\nu j}$-form for now, we claim that we have
the estimate
\begin{equation}\nonumber\begin{split}
&||\nabla_{x,t}P_{k_{0}}[P_{k_{1}}\delta\psi\nabla^{-1}P_{k_{2}}[P_{k_{3}}\psi\nabla^{-1}P_{k_{4}}[P_{k_{5}}\psi\nabla^{-1}P_{k_{6}}Q_{\nu
j}(P_{k_{7}}\psi P_{k_{8}}\psi)]]]||_{N[k_{0}]}\\&\lesssim
\tilde{c}_{k_{1}}\tilde{c}_{k_{3}}\tilde{c}_{k_{5}}(\tilde{c}_{k_{7}}+\tilde{c}_{k_{8}})\\&\hspace{3cm} 2^{\delta_{1}[\min\{k_{7},k_{8}\}-\max\{k_{7},k_{8}\}]}2^{\delta_{2}[\min_{i=2,\ldots,6}\{k_{i}\}-\max_{i=2,\ldots,6}]\{k_{i}\}}\\
\end{split}\end{equation}
To prove this, one needs to analyze the possible frequency
interactions. By scaling invariance, we may assume $k_{0}=0$. The
first step consists in reducing \\ $P_{k_{6}}Q_{\nu
j}(P_{k_{7}}\psi P_{k_{8}}\psi)$ to hyperbolic microsupport,
dealing with the contribution of $(1-I)P_{k_{6}}Q_{\nu
j}(P_{k_{7}}\psi P_{k_{8}}\psi)$. The argument for this is given
in the appendix\footnote{If one of the inputs of $Q_{\nu j}(.)$ is
of 2nd type, one use lemma~\ref{technical1} instead of 3.4(b).} of
\cite{Kr-4}. Then
we deal with the following cases in schematic fashion:\\
({\bf{A}}.a): $k_{2}<-10$, $k_{5}=k_{6}+O(1)$, $k_{3}=k_{4}+O(1)$.
We combine 3.4(b) as well as lemma~\ref{technical1} to conclude
that
\begin{equation}\nonumber\begin{split}
&||P_{k_{6}}Q_{\nu j}I(P_{k_{7}}\psi
P_{k_{8}}\psi)||_{L_{t}^{2}L_{x}^{2}}\lesssim
2^{\frac{\min\{k_{6,7,8}\}}{2}}(\tilde{c}_{k_{7}}+\tilde{c}_{k_{8}}).\\
\end{split}\end{equation}
Next, we have
\begin{equation}\nonumber\begin{split}
&||P_{k_{4}}[P_{k_{5}}\psi\nabla^{-1}P_{k_{6}}A]||_{L_{t}^{\frac{4}{3}}L_{x}^{2}}\lesssim
2^{-k_{6}}(\sum_{c\in
C_{k_{5},k_{4}-k_{5}}}||P_{c}\psi||_{L_{t}^{4}L_{x}^{\infty}}^{2})^{\frac{1}{2}}||P_{k_{6}}A||_{L_{t}^{2}L_{x}^{2}}\\
&\lesssim
2^{\frac{3k_{5}}{4}}2^{\frac{k_{4}-k_{5}}{2}}2^{-k_{6}}||P_{k_{5}}\psi||_{{\mathcal{S}}[k_{5}]}||P_{k_{6}}A||_{L_{t}^{2}L_{x}^{2}}\\
\end{split}\end{equation}
Finally, we have
\begin{equation}\nonumber\begin{split}
&||P_{k_{2}}[P_{k_{3}}\psi\nabla^{-1}P_{k_{4}}A]||_{L_{t}^{1}L_{x}^{2}}\lesssim
2^{-k_{4}}(\sum_{c\in
C_{k_{3},k_{2}-k_{3}}}||P_{c}\psi||_{L_{t}^{4}L_{x}^{\infty}}^{2})^{\frac{1}{2}}||P_{k_{4}}A||_{L_{t}^{\frac{4}{3}}L_{x}^{2}}\\
&\lesssim
2^{\frac{3k_{3}}{4}}2^{\frac{k_{2}-k_{3}}{2}}2^{-k_{4}}||P_{k_{3}}\psi||_{{\mathcal{S}}[k_{3}]}||P_{k_{4}}A||_{L_{t}^{\frac{4}{3}}L_{x}^{2}}\\
\end{split}\end{equation}
We have used here that $(\sum_{c\in
C_{k_{4},k_{2}-k_{4}}}||P_{c}\psi_{4}||_{L_{t}^{\frac{4}{3}}L_{x}^{2}}^{2})^{\frac{1}{2}}\lesssim
||P_{k_{4}}\psi_{4}||_{L_{t}^{\frac{4}{3}}L_{x}^{2}}$. Combining
these inequalities, we easily see that
\begin{equation}\nonumber\begin{split}
&||\nabla_{x,t}P_{k_{0}}Q_{<k_{0}}[P_{k_{1}}\delta\psi\nabla^{-1}P_{k_{2}}[P_{k_{3}}\psi\nabla^{-1}P_{k_{4}}[P_{k_{5}}\psi\nabla^{-1}P_{k_{6}}IQ_{\nu
j}(P_{k_{7}}\psi
P_{k_{8}}\psi)]]]||_{L_{t}^{1}\dot{H}^{-1}}\\&\lesssim
||P_{k_{1}}\delta\psi||_{L_{t}^{\infty}L_{x}^{2}}||P_{k_{2}}[P_{k_{3}}\psi\nabla^{-1}P_{k_{4}}[P_{k_{5}}\psi\nabla^{-1}P_{k_{6}}Q_{\nu
j}(P_{k_{7}}\psi P_{k_{8}}\psi)]]||_{L_{t}^{1}L_{x}^{2}}\\
&\lesssim
2^{\frac{3k_{3}}{4}+\frac{k_{2}-k_{3}}{2}-k_{4}+\frac{3k_{5}}{4}+\frac{k_{4}-k_{5}}{2}+\frac{\min\{k_{6,7,8}\}}{2}-k_{6}}\\&\hspace{5cm}||P_{k_{1}}\delta\psi||_{L_{t}^{\infty}L_{x}^{2}}||P_{k_{3}}\psi||_{{\mathcal{S}}[k_{3}]}||P_{k_{5}}\psi||_{{\mathcal{S}}[k_{5}]}
(\tilde{c}_{k_{7}}+\tilde{c}_{k_{8}}).\\
\end{split}\end{equation}
Moreover, we can estimate
\begin{equation}\nonumber\begin{split}
&||\nabla_{x,t}P_{k_{0}}Q_{\geq
k_{0}}[P_{k_{1}}\delta\psi\nabla^{-1}P_{k_{2}}[P_{k_{3}}\psi\nabla^{-1}P_{k_{4}}[P_{k_{5}}\psi\nabla^{-1}P_{k_{6}}IQ_{\nu
j}(P_{k_{7}}\psi
P_{k_{8}}\psi)]]]||_{\dot{X}_{k_{0}}^{-\frac{1}{2},-1,2}}\\
&\lesssim
2^{-\frac{k_{0}}{2}}2^{k_{2}}||P_{k_{1}}\delta\psi||_{L_{t}^{\infty}L_{x}^{2}}||P_{k_{2}}[P_{k_{3}}\psi\nabla^{-1}P_{k_{4}}[P_{k_{5}}\psi\nabla^{-1}P_{k_{6}}IQ_{\nu
j}(P_{k_{7}}\psi P_{k_{8}}\psi)]]||_{L_{t}^{2}L_{x}^{1}}\\
&\lesssim
2^{k_{2}-\frac{k_{0}}{2}}2^{\frac{\min\{k_{6,7,8}\}}{2}-k_{6}}||P_{k_{1}}\delta\psi||_{L_{t}^{\infty}L_{x}^{2}}||P_{k_{3}}\psi_{3}||_{L_{t}^{\infty}L_{x}^{2}}||P_{k_{5}}\psi_{5}||_{L_{t}^{\infty}L_{x}^{2}}
(\tilde{c}_{k_{7}}+\tilde{c}_{k_{8}}),\\
\end{split}\end{equation}
One easily verifies that these estimates verify the claim in the
case under consideration.\\
({\bf{A}}.b): $k_{2}+O(1)\leq k_{6}<<k_{5}$. In this case, we
rewrite the term under consideration as
\begin{equation}\nonumber
\nabla_{x,t}P_{k_{0}}[P_{k_{1}}\delta\psi\nabla^{-1}P_{k_{2}}[(P_{k_{3}}\psi\nabla^{-1}P_{k_{5}}\psi)\nabla^{-1}P_{k_{6}}IQ_{\nu
j}(P_{k_{7}}\psi P_{k_{8}}\psi)]]].
\end{equation}
Note that necessarily $k_{3}=k_{5}+O(1)$. We have
\begin{equation}\nonumber\begin{split}
&||\nabla_{x,t}P_{k_{0}}Q_{<k_{0}}[P_{k_{1}}\delta\psi\nabla^{-1}P_{k_{2}}[(P_{k_{3}}\psi\nabla^{-1}P_{k_{5}}\psi)\nabla^{-1}P_{k_{6}}IQ_{\nu
j}(P_{k_{7}}\psi P_{k_{8}}\psi)]]]||_{L_{t}^{1}\dot{H}^{-1}}\\
&\lesssim 2^{\frac{4 k_{2}}{M}-k_{6}}
||P_{k_{1}}\delta\psi||_{L_{t}^{\infty}L_{x}^{2}}||P_{\leq
k_{6}+O(1)}(P_{k_{3}}\psi\nabla^{-1}P_{k_{5}}\psi)||_{L_{t}^{2}L_{x}^{\frac{M}{2}}}
\\&\hspace{8cm}||P_{k_{6}}IQ_{\nu j}(P_{k_{7}}\psi
P_{k_{8}}\psi)||_{L_{t}^{2}L_{x}^{2}}\\
\end{split}\end{equation}
Then we estimate
\begin{equation}\nonumber\begin{split}
&||P_{\leq
k_{6}+O(1)}(P_{k_{3}}\psi\nabla^{-1}P_{k_{5}}\psi)||_{L_{t}^{2}L_{x}^{\frac{M}{2}}}
\lesssim
\sum_{a<k_{6}+O(1)}||P_{a}(P_{k_{3}}\psi\nabla^{-1}P_{k_{5}}\psi)||_{L_{t}^{2}L_{x}^{\frac{M}{2}}}\\
&\lesssim\sum_{a<k_{6}+O(1)}(\sum_{c\in
C_{k_{3},a-k_{3}}}||P_{c}\psi||_{L_{t}^{4}L_{x}^{M}}^{2})^{\frac{1}{2}}(\sum_{c\in
C_{k_{5},a-k_{5}}}||P_{c}\psi||_{L_{t}^{4}L_{x}^{M}}^{2})^{\frac{1}{2}}\\
&\lesssim
\sum_{a<k_{6}+O(1)}2^{\frac{a-k_{3}}{2+}}2^{\frac{a-k_{5}}{2+}}2^{(\frac{3}{4}-\frac{2}{M})(k_{3}+k_{5})}2^{-k_{3}}||P_{k_{3}}\psi||_{{\mathcal{S}}[k_{3}]}
||P_{k_{5}}\psi||_{{\mathcal{S}}[k_{5}]},\\
\end{split}\end{equation}
where $2+=2+(M)$. Putting this together with the above estimate
for \\ $IQ_{\nu j}(P_{k_{7}}\psi P_{k_{8}}\psi)$ easily results in
the claim for this case as well. Replacing the operator
$Q_{<k_{0}}$ by $Q_{\geq k_{0}}$, we have
\begin{equation}\nonumber\begin{split}
&||\nabla_{x,t}P_{k_{0}}Q_{\geq
k_{0}}[P_{k_{1}}\delta\psi\nabla^{-1}P_{k_{2}}[(P_{k_{3}}\psi\nabla^{-1}P_{k_{5}}\psi)\nabla^{-1}P_{k_{6}}IQ_{\nu
j}(P_{k_{7}}\psi P_{k_{8}}\psi)]]||_{\dot{X}_{k_{0}}^{-\frac{1}{2},-1,2}}\\
&\lesssim
2^{k_{2}}||P_{k_{1}}\delta\psi||_{L_{t}^{\infty}L_{x}^{2}}||P_{k_{2}}[(P_{k_{3}}\psi\nabla^{-1}P_{k_{5}}\psi)\nabla^{-1}P_{k_{6}}IQ_{\nu
j}(P_{k_{7}}\psi P_{k_{8}}\psi)]||_{L_{t}^{2}L_{x}^{1}}\\
&\lesssim 2^{\frac{4
(k_{2}-k_{5})}{M}}2^{\frac{1}{2+}(k_{6}-k_{5})}2^{\frac{\min\{k_{6,7,8}\}-k_{6}}{2}}\\&\hspace{4cm}||P_{k_{1}}\delta\psi||_{L_{t}^{\infty}L_{x}^{2}}||P_{k_{3}}\psi_{3}||_{L_{t}^{\infty}L_{x}^{2}}
||P_{k_{5}}\psi_{5}||_{L_{t}^{\infty}L_{x}^{2}}(\tilde{c}_{k_{7}}+\tilde{c}_{k_{8}})\\
\end{split}\end{equation}
This again verifies the claim.\\
({\bf{A}}.c): $k_{6}\leq k_{2}+O(1)<<k_{5}$. Consider the case
when the output is reduced to modulation $<2^{k_{0}}$, the
opposite case being treated similarly. We rewrite the term as
above and estimate
\begin{equation}\nonumber\begin{split}
&||\nabla_{x,t}P_{k_{0}}Q_{<k_{0}}[P_{k_{1}}\delta\psi\nabla^{-1}P_{k_{2}}[(P_{k_{3}}\psi\nabla^{-1}P_{k_{5}}\psi)\nabla^{-1}P_{k_{6}}IQ_{\nu
j}(P_{k_{7}}\psi P_{k_{8}}\psi)]]||_{L_{t}^{1}\dot{H}^{-1}}\\
&\lesssim
2^{(\frac{2}{M}-1)k_{2}}||P_{k_{1}}\delta\psi||_{L_{t}^{\infty}L_{x}^{2}}||(P_{k_{3}}\psi\nabla^{-1}P_{k_{5}}\psi)||_{L_{t}^{2}L_{x}^{\frac{M}{2}}}
||\nabla^{-1}P_{k_{6}}IQ_{\nu j}(P_{k_{7}}\psi
P_{k_{8}}\psi)||_{L_{t}^{2}L_{x}^{\infty}}\\
&\lesssim
2^{\frac{\min\{k_{7,8}\}-\max\{k_{7},k_{8}\}}{2}}2^{\frac{k_{6}-k_{2}}{2}}2^{\frac{k_{2}-k_{3}}{2+}}
\tilde{c}_{k_{1}}\tilde{c}_{k_{3}}\tilde{c}_{k_{5}}(\tilde{c}_{k_{7}}+\tilde{c}_{k_{8}})\\
\end{split}\end{equation}
This again yields the desired claim. \\
({\bf{A}}.d): $k_{5}\leq k_{6}+O(1)$. This case is treated
similarly and left out. The remaining frequency interactions are
also simple variations of this kind of reasoning and left out.
This establishes the claim for term ${\bf{(A)}}(\delta\psi,\psi)$.
We immediately deduce that if we reduce the expression to dyadic
time $\sim 2^{i}$, i. e. we apply a multiplier $\chi_{i}(t)$ in
front, we may reduce almost all (logarithmic) frequencies to norm
$<<i$, i. e. we may assume $|k_{2}|,\ldots, |k_{8}|<<i$. We shall
treat these frequencies as $O(1)$. In that case, though, we can
argue rather simply: observe that
\begin{equation}\nonumber\begin{split}
&||\nabla_{x,t}P_{k_{0}}Q_{<k_{0}}\chi_{i}(t)[P_{k_{1}}\delta\psi\nabla^{-1}P_{k_{2}}[(P_{k_{3}}\psi\nabla^{-1}P_{k_{5}}\psi)\nabla^{-1}P_{k_{6}}IQ_{\nu
j}(P_{k_{7}}\psi P_{k_{8}}\psi)]]]||_{L_{t}^{1}\dot{H}^{-1}}\\
&\lesssim
||P_{k_{1}}\delta\psi||_{L_{t}^{\infty}L_{x}^{2}}||P_{k_{3}}\psi||_{L_{t}^{2+}L_{x}^{\infty}}||P_{k_{5}}\psi||_{L_{t}^{M}L_{x}^{2+}}
||\nabla^{-1}P_{k_{6}}IQ_{\nu j}(P_{k_{7}}\psi P_{k_{8}}\psi)||_{L_{t}^{2}L_{x}^{2}}\\
&\lesssim 2^{-\mu i}\tilde{c}_{k_{0}}.\\
\end{split}\end{equation}
We have used 3.4(c) as well as lemma~\ref{technical1}. Similarly,
we can estimate
\begin{equation}\nonumber\begin{split}
&||\nabla_{x,t}P_{k_{0}}Q_{\geq
k_{0}}\chi_{i}(t)[P_{k_{1}}\delta\psi\nabla^{-1}P_{k_{2}}[(P_{k_{3}}\psi\nabla^{-1}P_{k_{5}}\psi)\\&\hspace{7cm}\nabla^{-1}P_{k_{6}}IQ_{\nu
j}(P_{k_{7}}\psi P_{k_{8}}\psi)]]]||_{\dot{X}_{k_{0}}^{-\frac{1}{2},-1,2}}\\
&\lesssim
2^{-\frac{k_{0}}{2}}||P_{k_{1}}\delta\psi||_{L_{t}^{\infty}L_{x}^{2}}||\chi_{i}(t)\psi_{3}||_{L_{t}^{\infty}L_{x}^{\infty}}
||P_{k_{5}}\psi||_{L_{t}^{\infty}L_{x}^{\infty}}||P_{k_{6}}IQ_{\nu
j}(P_{k_{7}}\psi P_{k_{8}}\psi)||_{L_{t}^{2}L_{x}^{2}}\\
&\lesssim 2^{-\mu i}\tilde{c}_{k_{0}}\\
\end{split}\end{equation}
We still need to consider the case when $(\psi^{2})$ gets replaced
by the error terms
$R_{\nu}\psi\nabla^{-1}(\psi\nabla^{-1}(\psi^{2}))$ etc. But this
is straightforward: first, if $\nu=0$, one reduces $R_{0}\psi$ to
$IR_{0}\psi$ arguing as before; the latter is morally equivalent
to $\psi$. One then winds up with the following schematic
expression:
\begin{equation}\nonumber
\nabla_{x,t}[\delta\psi\nabla^{-1}[\psi\nabla^{-1}[\psi\nabla^{-1}[\psi\nabla^{-1}[\psi\nabla^{-1}[\psi^{2}]]]]]]
\end{equation}
One can place two of the inputs $\psi$ into
$L_{t}^{2+}L_{x}^{\infty}$ and two others into
$L_{t}^{\infty}L_{x}^{2}$, $L_{t}^{M}L_{x}^{2+}$,
respectively\footnote{One avoids losses for high-high interactions
by interpolating with the improved
$L_{t}^{4}L_{x}^{\infty}$-norms.}. Details are tedious
reiterations of previously given arguments. This concludes the
estimates for term ${\bf{A}}(\delta\psi,\psi)$. The estimates for
the remaining terms ${\bf{B}}(\delta\psi,\psi)$ etc. are more of
the same. One invokes the weaker form of the improved Strichartz
norms available for the spaces $S[k]$ available by 3.4(f), which
is good enough when combined with the stronger version available
for the $\psi_{\nu}$.

\subsection{The higher order error terms at least quadratic in $\delta\psi$. Completing the proof of
Proposition~\ref{refined bootstrap}.} We need to estimate
expressions like ${\bf{A}}(\delta\psi,\psi)$ as in the preceding
section in which at least one additional $\psi$ has been replaced
by $\delta\psi$. Of course we no longer need to gain in time, but
we need to obtain an estimate like in theorem~\ref{decomposition}
upon frequency localizing the expression and its inputs. This will
eventually allow us to sum over all frequency interactions. The
degree of difficulty of these terms varies depending on how many
inputs $\delta\psi$ are present, due to the weaker nature of the
estimates satisfied by these. However, we shall try to treat all
of these terms in a uniform manner, as these estimates are not
interesting in and of themselves. Proceeding as in \cite{Kr-4},
one introduces a null-structure into these terms by splitting all
inputs into gradient and elliptic parts, just as before.
Substituting the elliptic parts results in error terms of higher
order, while the gradient parts contribute to the null-structure.
Of course, as before one may apply this splitting to both the
inputs $\delta\psi$ as well as $\psi$. Carrying out this process a
finite number of times results in the decomposition described in
theorem~\ref{decomposition}. This theorem then ensures that all
the resulting terms upon frequency localizing will satisfy the
desired estimates, provided all inputs $P_{k_{i}}\psi$ are of the
first type. Our concern here is whether the same kind of estimates
hold if they are of the 2nd type. Thus consider for example a
typical quintilinear term of the form
\begin{equation}\nonumber
\nabla_{x,t}[P_{k_{1}}\delta\psi_{1}P_{k_{2}}\nabla^{-1}[P_{k_{3}}R_{\nu}\delta\psi\nabla^{-1}[P_{k_{4}}\psi
P_{k_{5}}\psi]]]
\end{equation}
The first step here consists in reducing  $R_{\nu}\delta\psi$ to
the 'hyperbolic version ' $R_{\nu}I\delta\psi$, arguing as in the
preceding subsection. Next, one applies the usual Hodge-type
decomposition to the innermost square bracket $[\psi^{2}]$,
replacing this by the sum of a $Q_{\nu j}$-type null-form as well
as error terms at least quadrilinear. One treats the resulting
quintilinear null-form just as in \cite{Kr-4}, resulting in the
desired estimate, provided both inputs $P_{k_{4,5}}\psi_{4,5}$ are
of the first type. Now assume at least one is of the 2nd type. We
first consider the expression
\begin{equation}\nonumber
\nabla_{x,t}[P_{k_{1}}\delta\psi_{1}P_{k_{2}}\nabla^{-1}[P_{k_{3}}R_{\nu}I\delta\psi\nabla^{-1}P_{k_{4}}[P_{k_{5}}\psi\nabla^{-1}(1-I)Q_{\nu
j}[P_{k_{6}}\psi P_{k_{7}}\psi]]]
\end{equation}
If both $P_{k_{6,7}}\psi$ are of the first type, one argues here
as in the appendix of \cite{Kr-4},  resulting in the desired
estimate. If at least one of these inputs is of the 2nd type, one
substitutes lemma~\ref{technical1} instead of 3.4(b) in that same
argument. Thus we can now replace $Q_{\nu j}(P_{k_{6}}\psi,
P_{k_{7}}\psi)$ by $IQ_{\nu j}(P_{k_{6}}\psi, P_{k_{7}}\psi)$. If
both $P_{k_{6,7}}\psi$ are of first type, one again argues as in
\cite{Kr-4} (these quintilinear null-forms are part of the
expansion in theorem~\ref{decomposition}). If one of $P_{k_{6,7}}$
is of 2nd type, the estimate becomes quite simple due to the
strong estimates available: observe that then
\begin{equation}\nonumber
||P_{k_{5}}\nabla ^{-1}Q_{\nu j}I[P_{k_{6}}\psi,
P_{k_{7}}\psi]||_{L_{t}^{1}L_{x}^{2+}}\lesssim
2^{(1-\delta_{1})\max\{k_{i}\}}2^{\delta_{2}\min\{k_{i}\}}\tilde{c}_{k_{6}}\tilde{c}_{k_{7}}
\end{equation}
Repeating the usual frequency trichotomies, one gets from here
that
\begin{equation}\nonumber\begin{split}
&||\nabla_{x,t}P_{k_{0}}Q_{<k_{0}}[P_{k_{1}}\delta\psi_{1}P_{k_{2}}\nabla^{-1}[P_{k_{3}}R_{\nu}I\delta\psi\\&\hspace{4.5cm}\nabla^{-1}P_{k_{4}}[P_{k_{5}}\psi\nabla^{-1}(1-I)Q_{\nu
j}[P_{k_{6}}\psi P_{k_{7}}\psi]]]||_{L_{t}^{1}\dot{H}^{-1}}\\
&\lesssim
2^{-\delta_{1}|k_{1}-k_{0}|}2^{\delta_{2}[\min_{i=2,\ldots,
7}\{k_{i}\}-\max_{i=2,\ldots,7}\{k_{i}\}]}\tilde{c}_{k_{0}}\\
\end{split}\end{equation}
\begin{equation}\nonumber\begin{split}
&||\nabla_{x,t}P_{k_{0}}Q_{\geq
k_{0}}[P_{k_{1}}\delta\psi_{1}P_{k_{2}}\nabla^{-1}[P_{k_{3}}R_{\nu}I\delta\psi\\&\hspace{4cm}\nabla^{-1}P_{k_{4}}[P_{k_{5}}\psi\nabla^{-1}(1-I)Q_{\nu
j}[P_{k_{6}}\psi P_{k_{7}}\psi]]]||_{\dot{X}_{k_{0}}^{-\frac{1}{2},-1,2}}\\
&\lesssim
2^{-\delta_{1}|k_{1}-k_{0}|}2^{\delta_{2}[\min_{i=2,\ldots,
7}\{k_{i}\}-\max_{i=2,\ldots,7}\{k_{i}\}]}\tilde{c}_{k_{0}}\\
\end{split}\end{equation}
The remaining error terms are treated analogously.
\end{section}
\begin{section}{Appendix: Proof of theorem~\ref{Moser}.}

We first check the algebra type estimate. Thus let $\psi_{1,2}\in
{\mathcal{S}}({\mathbf{R}}^{2+1})$; we need to estimate
$||P_{k}[P_{k_{1}}\psi_{1}\nabla^{-1}P_{k_{2}}\psi_{2}]||_{{\mathcal{S}}[k]}$.
Of course we may assume that $k=0$. We decompose
\begin{equation}\nonumber
P_{0}[P_{k_{1}}\psi_{1}\nabla^{-1}P_{k_{2}}\psi_{2}]=P_{0}Q_{<100}[P_{k_{1}}\psi_{1}\nabla^{-1}P_{k_{2}}\psi_{2}]
+P_{0}Q_{\geq 100}[P_{k_{1}}\psi_{1}\nabla^{-1}P_{k_{2}}\psi_{2}]
\end{equation}
We first consider the large modulation case, i. e. the 2nd summand
on the right. Commence with the case $k_{1}>10$. We freeze the
modulation of the output to dyadic size $\sim 2^{l}$, and further
decompose into the following cases:
\begin{equation}\nonumber\begin{split}
&P_{0}Q_{l}[P_{k_{1}}\psi_{1}\nabla^{-1}P_{k_{2}}\psi_{2}]=P_{0}Q_{l}[P_{k_{1}}Q_{\geq l-10}\psi_{1}\nabla^{-1}P_{k_{2}}\psi_{2}]\\
&+P_{0}Q_{l}[P_{k_{1}}Q_{<l-10}\psi_{1}\nabla^{-1}P_{k_{2}}Q_{\geq l-10}\psi_{2}]+P_{0}Q_{l}[P_{k_{1}}Q_{<l-10}\psi_{1}\nabla^{-1}P_{k_{2}}Q_{<l-10}\psi_{2}]\\
\end{split}\end{equation}
We treat each of the summands on the right: for the first, observe
that irrespective of whether $P_{k_{1}}\psi_{1}$ is of first or
2nd type,
\begin{equation}\nonumber\begin{split}
&||P_{0}Q_{l}[P_{k_{1}}Q_{\geq
l-10}\psi_{1}\nabla^{-1}P_{k_{2}}\psi_{2}]||_{\dot{X}_{0}^{-(\frac{1}{2}-\mu),1-\mu,1}}
\\&\hspace{3cm}\lesssim 2^{(1-\mu)l}||P_{k_{1}}Q_{\geq
l-10}\psi_{1}||_{L_{t}^{2}L_{x}^{2}}||\nabla^{-1}P_{k_{2}}\psi_{2}||_{L_{t}^{\infty}L_{x}^{2}}\\
\end{split}\end{equation}
When $l\geq k_{1}$, one estimates this by
\begin{equation}\nonumber
\lesssim \sum_{a\geq
l-10}2^{(1-\mu)(l-a)}||P_{k_{1}}Q_{a}\psi_{1}||_{\dot{X}_{k_{1}}^{-(\frac{1}{2}-\mu),(1-\mu),1}}||\nabla^{-1}P_{k_{2}}\psi_{2}||_{L_{t}^{\infty}L_{x}^{2}}
\end{equation}
If, on the other hand, $l<k_{1}$, one estimates this by
\begin{equation}\nonumber
\lesssim
2^{(\frac{1}{2}-\mu)l}||P_{k_{1}}\psi_{1}||_{\dot{X}_{k_{1}}^{0,\frac{1}{2},\infty}}||\nabla^{-1}P_{k_{2}}\psi_{2}||_{L_{t}^{\infty}L_{x}^{2}}
\end{equation}
if $P_{k_{1}}\psi_{1}$ is of first type, and by
\begin{equation}\nonumber
\lesssim
2^{(1-\mu)l}||P_{k_{1}}\psi_{1}||_{L_{t}^{2}L_{x}^{2+}}||P_{k_{2}}\nabla^{-1}\psi_{2}||_{L_{t}^{\infty}L_{x}^{2}}
\end{equation}
if it is of 2nd type. Moreover, one estimates
\begin{equation}\nonumber
||P_{0}[P_{k_{1}}\psi_{1}\nabla^{-1}\psi_{2}]||_{L_{t}^{\infty}L_{x}^{1+}}\lesssim
2^{-k_{1}}||P_{k_{1}}\psi_{1}||_{L_{t}^{\infty}L_{x}^{2}}||P_{k_{2}}\psi_{2}||_{L_{t}^{\infty}L_{x}^{2}}
\end{equation}
Next, one estimates the output with respect to
$||.||_{L_{t}^{1+}L_{x}^{\infty}}$ by interpolating between a
crude estimate for $L_{t}^{1+}L_{x}^{\infty}$ gotten by placing
the inputs into $L_{t}^{2+}L_{x}^{\infty}$ and a refined estimate
for $||.||_{L_{t}^{2}L_{x}^{\infty}}$ by using improved
$L_{t}^{4}L_{x}^{\infty}$-Strichartz  norms for the inputs,
resulting again in a small exponential gain in $k_{1}$. This
yields the desired bound upon summing over
$k_{1}=k_{2}+O(1)>O(1)$, showing that this contribution to the
output is of 2nd type. Now consider the third summand in the above
trichotomy: we have $k_{1}=k_{2}+O(1)=l+O(1)$ in this case. First
assume both $P_{k_{1,2}}\psi_{1,2}$ are of first type. We can
decompose
\begin{equation}\nonumber\begin{split}
&P_{0}Q_{l}[P_{k_{1}}Q_{<l-10}\psi_{1}\nabla^{-1}P_{k_{2}}Q_{<l-10}\psi_{2}]
=\sum_{\pm}P_{0}Q_{l}[P_{k_{1}}Q^{\pm}_{<l-10}\psi_{1}\nabla^{-1}P_{k_{2}}Q^{\pm}_{<l-10}\psi_{2}]\\
&=\sum_{\pm}\sum_{\kappa_{1,2}\in K_{-k_{1}},\text{dist}(\kappa_{1},\kappa_{2})\sim 1}P_{0}Q_{l}[P_{k_{1},\kappa_{1}}Q^{\pm}_{<l-10}\psi_{1}\nabla^{-1}P_{k_{2},\kappa_{2}}Q^{\pm}_{<l-10}\psi_{2}]\\
\end{split}\end{equation}
Now one estimates, using the definition of $S[k,\kappa]$:
\begin{equation}\nonumber\begin{split}
&||\sum_{\pm}\sum_{\kappa_{1,2}\in K_{-k_{1}},\text{dist}(\kappa_{1},\kappa_{2})\sim 1}P_{0}Q_{l}[P_{k_{1},\kappa_{1}}Q^{\pm}_{<l-10}\psi_{1}\nabla^{-1}P_{k_{2},\kappa_{2}}Q^{\pm}_{<l-10}\psi_{2}]||_{L_{t}^{2}L_{x}^{2}}\\
&\lesssim\sum_{\pm}(\sum_{\kappa_{1}\in
K_{-k_{1}}}||P_{k_{1},\kappa_{1}}Q^{\pm}_{<l-10}\psi_{1}||_{PW[\pm\kappa_{1}]}^{2})^{\frac{1}{2}}(\sum_{\kappa_{2}\in
K_{-k_{1}}}||P_{k_{2},\kappa_{2}}Q^{\pm}_{<l-10}\nabla^{-1}\psi_{2}||_{NFA[\pm\kappa_{1}]}^{2})^{\frac{1}{2}}\\
\end{split}\end{equation}
This in turn can be bounded by
\begin{equation}\nonumber
\lesssim
l^{2}2^{-l}||P_{k_{1}}\psi_{1}||_{A[k_{1}]}||P_{k_{2}}\psi_{2}||_{A[k_{2}]}
\end{equation}
Multiplying by $2^{(1-\mu)l}$ results in an acceptable estimate.
Now assume $P_{k_{1}}\psi_{1}$ is of 2nd type, say. In that case,
we estimate
\begin{equation}\nonumber\begin{split}
&||P_{0}Q_{l}[P_{k_{1}}Q_{<l-10}\psi_{1}\nabla^{-1}P_{k_{2}}Q_{<l-10}\psi_{2}]||_{L_{t}^{2}L_{x}^{2}}
\\&\hspace{4cm}\lesssim
||P_{k_{1}}Q_{<l-10}\psi_{1}||_{L_{t}^{2}L_{x}^{2+}}||P_{k_{2}}\nabla^{-1}\psi_{2}||_{L_{t}^{\infty}L_{x}^{2}}\\
\end{split}\end{equation}
Multiplication with $2^{(1-\mu)l}$ again yields an acceptable
bound. One estimates the output with respect to
$||.||_{L_{t}^{\infty}L_{x}^{1}}$ as well as
$||.||_{L_{t}^{1+}L_{x}^{\infty}}$ as before, showing that this
contribution is of 2nd type as well. The 2nd term of the
trichotomy is treated like the first. Now consider the case
$k_{1}\in [-10,10]$. We decompose
\begin{equation}\nonumber\begin{split}
&P_{0}Q_{>100}\partial_{t}[P_{k_{1}}\psi_{1}\nabla^{-1}P_{k_{2}}\psi_{2}]
\\&=P_{0}Q_{>100}[P_{k_{1}}\partial_{t}\psi_{1}\nabla^{-1}P_{k_{2}}\psi_{2}]+P_{0}Q_{>100}[P_{k_{1}}\psi_{1}R_{0}P_{k_{2}}\psi_{2}]\\
\end{split}\end{equation}
We claim that if $P_{k_{1}}\psi_{1}$ is of first type, so is the
output. This is immediate when $P_{k_{2}}\psi_{2}$ is of 1st type.
Now assume that $P_{k_{2}}\psi_{2}$ is of 2nd type. If
$P_{k_{1}}\psi_{1}$ has modulation $>2^{10}$, this is again
immediate. In the opposite case, one calculates
\begin{equation}\nonumber\begin{split}
&||P_{0}Q_{>100}[P_{k_{1}}Q_{<10}\partial_{t}\psi_{1}\nabla^{-1}P_{k_{2}}\psi_{2}]||_{L_{t}^{2}L_{x}^{2}}\\&\lesssim
||P_{k_{1}}Q_{<10}\partial_{t}\psi_{1}||_{L_{t}^{\infty}L_{x}^{2}}||\nabla^{-1}P_{k_{2}}Q_{>10}\psi_{2}||_{L_{t}^{2}L_{x}^{\infty}}\\
\end{split}\end{equation}
Using the definition of $B[k_{2}]$, one checks that the summation
over $k_{2}$ can be carried out. Next, we use the first bilinear
property in theorem~\ref{Moser}(which will be proved later
independently of this) to calculate
\begin{equation}\nonumber\begin{split}
&||P_{0}Q_{>100}[P_{k_{1}}\psi_{1}R_{0}P_{k_{2}}\psi_{2}]||_{L_{t}^{2}L_{x}^{2}}
\\&\lesssim (\sum_{c\in
C_{k_{1},k_{2}-k_{1}}}||P_{0}Q_{>100}[P_{c}Q_{<10}\psi_{1}R_{0}Q_{>90}\psi_{2}]||_{L_{t}^{2}L_{x}^{2}}^{2})^{\frac{1}{2}}\\
&\lesssim
2^{\delta(k_{2}-k_{1})}||P_{k_{1}}\psi_{1}||_{A[k_{1}]}\\
\end{split}\end{equation}
One can sum over $k_{2}<15$, getting the desired bound. Control
over $||.||_{L}$ again follows via the Sobolev embedding. If
$P_{k_{1}}\psi_{1}$ is of 2nd type, so is the output. This is a
simple repetition of arguments before. The case $k_{1}<-10$ is
more of the same, which finishes the large modulation case. Now
consider
$P_{0}Q_{<100}[P_{k_{1}}\psi_{1}\nabla^{-1}P_{k_{2}}\psi_{2}]$.
First assume $k_{1}>10$, and both $P_{k_{1,2}}\psi_{1,2}$ of first
type. In that case, we claim that the output will be of 2nd type.
We need to check that it is controlled with respect to both
$||.||_{L_{t}^{\infty}L_{x}^{1+}}$ as well as
$||.||_{L_{t}^{1+}L_{x}^{\infty}}$. For the first norm, this is
immediate from Bernstein's inequality. For the 2nd norm, one
interpolates between a crude estimate for
$L_{t}^{1+}L_{x}^{\infty}$ gotten by placing both
$P_{k_{1,2}}\psi_{1,2}$ into $L_{t}^{2+}L_{x}^{\infty}$ and an
estimate for $||.||_{L_{t}^{2}L_{x}^{\infty}}$ gotten by using
improved Strichartz type norms for the inputs. The remaining
improved Strichartz type norms constituting $||.||_{L}$ are
controlled from Bernstein's inequality. Now assume
$P_{k_{1}}\psi_{1}$ is of 2nd type. Then the output will again be
of 2nd type, as is easily verified by placing $P_{k_{2}}\psi_{2}$
into $L_{t}^{\infty}L_{x}^{2}$. This concludes the case
$k_{1}>10$. Now assume $k_{1}\in [-10,10]$, in which case
$k_{2}<15$. First assume $P_{k_{1}}\psi_{1}$ is of first type. We
claim that then the output will be of first type, irrespective of
the type of $P_{k_{2}}\psi_{2}$. Commence with the case when
$P_{k_{2}}\psi_{2}$ is of first type. We need to check the various
parts constituting $||.||_{A[0]}$. First, we consider
$||.||_{\dot{X}_{0}^{0,\frac{1}{2},\infty}}$. Freeze the
modulation of the output and decompose
\begin{equation}\nonumber\begin{split}
&P_{0}Q_{j}[P_{k_{1}}\psi_{1}\nabla^{-1}P_{k_{2}}\psi_{2}]=P_{0}Q_{j}[P_{k_{1}}Q_{\geq j-10}\psi_{1}\nabla^{-1}P_{k_{2}}\psi_{2}]\\
&+P_{0}Q_{j}[P_{k_{1}}Q_{<j-10}\psi_{1}\nabla^{-1}P_{k_{2}}Q_{\geq j-10}\psi_{2}]+P_{0}Q_{j}[P_{k_{1}}Q_{<j-10}\psi_{1}\nabla^{-1}P_{k_{2}}Q_{<j-10}\psi_{2}]\\
\end{split}\end{equation}
We start with the first term on the righthand side: estimate
\begin{equation}\nonumber\begin{split}
&||P_{0}Q_{j}[P_{k_{1}}Q_{\geq
j-10}\psi_{1}\nabla^{-1}P_{k_{2}}\psi_{2}]||_{\dot{X}_{0}^{0,\frac{1}{2},\infty}}
\\&\hspace{5cm}\lesssim 2^{\frac{j}{2}}||P_{k_{1}}Q_{\geq
j-10}\psi_{1}||_{L_{t}^{2}L_{x}^{2}}||\nabla^{-1}P_{k_{2}}\psi_{2}||_{L_{t}^{\infty}L_{x}^{\infty}},\\
\end{split}\end{equation}
and one easily bounds this by
\begin{equation}\nonumber
||P_{k_{1}}\psi_{1}||_{A[k_{1}]}||P_{k_{2}}\psi_{2}||_{A[k_{2}]}
\end{equation}
Now consider the 2nd term in the above trichotomy. We can estimate
this as follows:
\begin{equation}\nonumber\begin{split}
&||P_{0}Q_{j}[P_{k_{1}}Q_{<j-10}\psi_{1}\nabla^{-1}P_{k_{2}}Q_{\geq
j-10}\psi_{2}]||_{\dot{X}_{0}^{0,\frac{1}{2},\infty}}\\&\lesssim
||P_{k_{1}}Q_{<j-10}\psi_{1}||_{L_{t}^{\infty}L_{x}^{2}}||\nabla^{-1}P_{k_{2}}Q_{\geq
j-10}\psi_{2}||_{L_{t}^{2}L_{x}^{\infty}}\\
&\lesssim
2^{\min\{(\frac{1}{2}-\mu)(k_{2}-j),0\}}2^{\min\{\frac{j-k_{2}}{4},0\}}||P_{k_{1}}\psi_{1}||_{A[k_{1}]}||P_{k_{2}}\psi_{2}||_{A[k_{2}]}\\
\end{split}\end{equation}
We have invoked the improved Bernstein's inequality \cite{Tao 2},
\cite{Kr-4}. Note that one may sum here over $k_{2}<15$ for fixed
$j<100$. Now assume that $P_{k_{2}}\psi_{2}$ is of 2nd type. Then
we don't use the preceding trichotomy, but divide into the cases
$k_{2}\geq j-100$ and $k_{2}<j-100$. In the first case, we
estimate
\begin{equation}\nonumber\begin{split}
&||P_{0}Q_{j}[P_{k_{1}}\psi_{1}\nabla^{-1}P_{k_{2}}\psi_{2}]||_{\dot{X}_{0}^{0,\frac{1}{2},\infty}}\lesssim
||P_{k_{1}}\psi_{1}||_{L_{t}^{\infty}L_{x}^{2}}||\nabla^{-1}P_{k_{2}}\psi_{2}||_{L_{t}^{2}L_{x}^{\infty}}\\&
\lesssim
2^{\frac{j-k_{2}}{2}}||P_{k_{1}}\psi_{1}||_{A[k_{1}]}||P_{k_{2}}\psi_{2}||_{B[k_{2}]}.\\
\end{split}\end{equation}
In the case $k_{2}<j-100$, we split
\begin{equation}\nonumber\begin{split}
&P_{0}Q_{j}[P_{k_{1}}\psi_{1}\nabla^{-1}P_{k_{2}}\psi_{2}]=P_{0}Q_{j}[P_{k_{1}}Q_{\geq
j-10}\psi_{1}\nabla^{-1}P_{k_{2}}\psi_{2}]\\&\hspace{6cm}+P_{0}Q_{j}[P_{k_{1}}Q_{<j-10}\psi_{1}\nabla^{-1}P_{k_{2}}Q_{\geq
j-10}\psi_{2}]\\
\end{split}\end{equation}
The first summand on the right is estimated just as before:
\begin{equation}\nonumber\begin{split}
&||P_{0}Q_{j}[P_{k_{1}}Q_{\geq
j-10}\psi_{1}\nabla^{-1}P_{k_{2}}\psi_{2}]||_{\dot{X}_{0}^{0,\frac{1}{2},\infty}}\\&\hspace{5cm}\lesssim
2^{\frac{j}{2}}||P_{k_{1}}Q_{\geq
j-10}\psi_{1}||_{L_{t}^{2}L_{x}^{2}}||\nabla^{-1}P_{k_{2}}\psi_{2}||_{L_{t}^{\infty}L_{x}^{\infty}}\\
\end{split}\end{equation}
as well as
\begin{equation}\nonumber\begin{split}
&||\sum_{k_{2}<j-100}P_{0}Q_{j}[P_{k_{1}}Q_{\geq
j-10}\psi_{1}\nabla^{-1}P_{k_{2}}\psi_{2}]||_{\dot{X}_{0}^{0,\frac{1}{2},\infty}}\\&\hspace{6cm}\lesssim
2^{\frac{j}{2}}||P_{k_{1}}Q_{\geq
j-10}\psi_{1}||_{L_{t}^{2}L_{x}^{2}}||\nabla^{-1}\psi_{2}||_{L_{t}^{\infty}L_{x}^{\infty}}\\
\end{split}\end{equation}
As for the 2nd summand, we have
\begin{equation}\nonumber\begin{split}
&||P_{0}Q_{j}[P_{k_{1}}Q_{<j-10}\psi_{1}\nabla^{-1}P_{k_{2}}Q_{\geq
j-10}\psi_{2}]||_{\dot{X}_{0}^{0,\frac{1}{2},\infty}}\\&\lesssim
2^{\frac{j}{2}}||P_{k_{1}}Q_{<j-10}\psi_{1}||_{L_{t}^{\infty}L_{x}^{2}}2^{(\frac{1}{2}-\mu)k_{2}}2^{-j(1-\mu)}||\nabla^{-1}P_{k_{2}}Q_{\geq
j-10}\psi_{2}||_{\dot{X}_{k_{2}}^{-(\frac{1}{2}-\mu),1-\mu,1}}\\
\end{split}\end{equation}
Thus we can estimate this contribution by
\begin{equation}\nonumber
\lesssim
2^{(\frac{1}{2}-\mu)(k_{2}-j)}||P_{k_{1}}\psi_{1}||_{A[k_{1}]}||P_{k_{2}}\psi_{2}||_{B[k_{2}]}
\end{equation}
Observe that we have obtained the gain
$2^{\frac{\min\{k_{2}-j,0\}}{2}}2^{(\frac{1}{2}-\mu)(j-k_{2})}$,
which allows us to sum over $k_{2}<15$. This concludes the
estimate for $||.||_{\dot{X}_{0}^{0,\frac{1}{2},\infty}}$. We
still need to control the complicated null-frame part, as well as
$||.||_{L}$. Thus fix $l<-10$, and $-10\geq \lambda\geq l$, and
consider an expression
\begin{equation}\nonumber
|\lambda|^{-1}(\sum_{\kappa\in K_{l}}\sum_{R\in
C_{0,\kappa,\lambda}}||\tilde{P}_{R}Q^{\pm}_{<2l}[P_{k_{1}}\psi_{1}\nabla^{-1}P_{k_{2}}\psi_{2}]||_{S[0,\pm\kappa]}^{2})^{\frac{1}{2}}
\end{equation}
The following estimates are irrespective of the type of
$P_{k_{2}}\psi_{2}$. Then we split this into a bunch of
contributions: Observe the identity
\begin{equation}\nonumber
\tilde{P}_{R}Q^{\pm}_{<2l}[P_{k_{1}}\psi_{1}\nabla^{-1}P_{<2l}Q_{<2l}\psi_{2}]
=\tilde{P}_{R}Q^{\pm}_{<2l}[\tilde{P}_{R_{1}}Q^{\pm}_{<2l+O(1)}\psi_{1}\nabla^{-1}P_{<2l}Q_{<2l}\psi_{2}],
\end{equation}
where $R_{1}=(1+\frac{1}{1000})R$. Next, note that from the
definition of $S[k,\kappa]$, we have
\begin{equation}\nonumber\begin{split}
&||\tilde{P}_{R_{1}}Q^{\pm}_{<2l+O(1)}\psi_{1}||_{S[k,\pm\kappa]}\lesssim
\sum_{\kappa_{1}\in K_{l-100},\kappa_{1}\subset
(1+\frac{1}{100})\kappa}||\tilde{P}_{R_{1}}P_{k_{1},\kappa_{1}}Q^{\pm}_{<2l+O(1)}\psi_{1}||_{S[k,\pm\kappa]}\\
&\lesssim\sum_{\kappa_{1}\in K_{l-100},\kappa_{1}\subset
(1+\frac{1}{100})\kappa}||\tilde{P}_{R_{1}}P_{k_{1},\kappa_{1}}Q^{\pm}_{<2l+O(1)}\psi_{1}||_{S[k,\pm\kappa_{1}]},\\
\end{split}\end{equation}
Note that
$\frac{11}{10}\kappa_{1}\subset\frac{11}{10}\kappa=\tilde{\kappa}$,
the latter as in the definition of $PW[\kappa]$. Hence the above
inequality. Putting these observations together, we get
\begin{equation}\nonumber\begin{split}
&|\lambda|^{-1}(\sum_{\kappa\in K_{l}}\sum_{R\in
C_{0,\kappa,\lambda}}||\tilde{P}_{R}Q^{\pm}_{<2l}[P_{k_{1}}\psi_{1}\nabla^{-1}P_{<2l}Q_{<2l}\psi_{2}]||_{S[0,\pm\kappa]}^{2})^{\frac{1}{2}}\\
&\lesssim |\lambda|^{-1}(\sum_{\kappa\in K_{l+O(1)}}\sum_{R\in
C_{0,\kappa,\lambda+O(1)}}||\tilde{P}_{R}Q^{\pm}_{<2l}\psi_{1}||_{S[k,\pm\kappa]}^{2})^{\frac{1}{2}}||\nabla^{-1}\psi_{2}||_{L_{t}^{\infty}L_{x}^{\infty}},\\
\end{split}\end{equation}
which leads to the desired estimate. Further, using
\begin{equation}\nonumber\begin{split}
&||P_{0}Q_{<2l}[P_{k_{1}}\psi_{1}\nabla^{-1}P_{<2l}Q_{\geq
2l}\psi_{2}]||_{\dot{X}_{0}^{0,\frac{1}{2},1}}\\&\hspace{3cm}+||P_{0}Q_{<2l}[P_{k_{1}}\psi_{1}\nabla^{-1}P_{\geq
2l}\psi_{2}]||_{\dot{X}_{0}^{0,\frac{1}{2},1}}\lesssim
||P_{k_{1}}\psi_{1}||_{A[k_{1}]}\\
\end{split}\end{equation}
We can estimate
\begin{equation}\nonumber\begin{split}
&|\lambda|^{-1}(\sum_{\kappa\in K_{l}}\sum_{R\in
C_{0,\kappa,\lambda}}||\tilde{P}_{R}Q^{\pm}_{<2l}[P_{k_{1}}\psi_{1}\nabla^{-1}P_{\geq
2l}\psi_{2}]||_{S[0,\pm\kappa]}^{2})^{\frac{1}{2}}
\\&\hspace{3cm}\lesssim
||P_{0}Q^{\pm}_{<2l}[P_{k_{1}}\psi_{1}\nabla^{-1}P_{\geq 2l}\psi_{2}]||_{\dot{X}_{0}^{0,\frac{1}{2},1}},\\
\end{split}\end{equation}
\begin{equation}\nonumber\begin{split}
&|\lambda|^{-1}(\sum_{\kappa\in K_{l}}\sum_{R\in
C_{0,\kappa,\lambda}}||\tilde{P}_{R}Q^{\pm}_{<2l}[P_{k_{1}}\psi_{1}\nabla^{-1}P_{<
2l}Q_{\geq 2l}\psi_{2}]||_{S[0,\pm\kappa]}^{2})^{\frac{1}{2}}
\\&\hspace{3cm}\lesssim
||P_{0}Q^{\pm}_{<2l}[P_{k_{1}}\psi_{1}\nabla^{-1}P_{<2l}Q_{\geq 2l}\psi_{2}]||_{\dot{X}_{0}^{0,\frac{1}{2},1}}\\
\end{split}\end{equation}
and so the desired estimate follows easily for the contributions
of $P_{k_{1}}\psi_{1}\nabla^{-1}P_{\geq 2l}\psi_{2}$,
$P_{k_{1}}\psi_{1}\nabla^{-1}P_{<2l}Q_{\geq 2l}\psi_{2}$. The
estimate for $||.||_{L}$ is quite similar. This concludes the
estimates when $P_{k_{1}}\psi_{1}$ is of first type. If it is of
2nd type, so will be the output. This is straightforward to check.
The case $k_{1}<-10$ is a tedious reiteration of similar estimates
and hence omitted. This concludes the proof for the assertions of
theorem~\ref{Moser} as far as the estimates concerning
$||.||_{{\mathcal{S}}}$ are concerned. We now proceed to the
assertions concerning the bilinear estimates, as well as the
estimates for
\begin{equation}\nonumber
||R_{0}[\psi_{1}A(\nabla^{-1}\psi_{2})]||_{L_{t}^{\infty}L_{x}^{2}},\,||P_{k}Q_{<k+O(1)}[\psi_{1}A(\nabla^{-1}\psi_{2})]||_{\dot{X}_{k}^{0,\frac{1}{2},\infty}}
\end{equation}
We start with the latter, which follows from the refined assertion
for the function $\beta$ in the decomposition
\begin{equation}\nonumber
P_{k}Q_{<k+O(1)}[\psi_{1}A(\nabla^{-1}\psi_{2})]=\alpha+\beta
\end{equation}
To understand this decomposition, one expands
$A(\nabla^{-1}\psi_{2})$ into a Taylor series (using real
analyticity). One winds up with schematic expressions of the form
$\psi\nabla^{-1}\psi$, $\psi\nabla^{-1}\psi\nabla^{-1}\psi$ etc.
For the first type of expression, $P_{k}[\psi\nabla^{-1}\psi]$,
the previous proof revealed that the only contribution of 2nd type
arises from high-high interactions. But for these it is
straightforward to verify that they are controlled with respect to
$||.||_{\dot{X}_{k}^{0,\frac{1}{2},1}}$, see 3.4(a). Now one
proceeds inductively, assuming the assertion to be true for both
$P_{k_{1,2}}\psi_{1,2}$ $\forall k_{1,2}\in{\mathbf{Z}}$, and
considering
$P_{k}[P_{k_{1}}\psi_{1}\nabla^{-1}P_{k_{2}}\psi_{2}]$. For
example, considering high-high interactions, when
$P_{k_{1}}\psi_{1}$ is of 2nd type, one estimates
\begin{equation}\nonumber\begin{split}
&||P_{k}Q_{<k+O(1)}[P_{k_{1}}\psi_{1}\nabla^{-1}P_{k_{2}}\psi_{2}]||_{\dot{X}_{k}^{0,\frac{1}{2},1}}
\lesssim
2^{(\frac{3}{2}-\epsilon)k}||P_{k_{1}}\psi_{1}||_{L_{t}^{2}L_{x}^{2+}}||\nabla^{-1}P_{k_{2}}\psi_{2}||_{L_{t}^{\infty}L_{x}^{2}}\\
&\lesssim
2^{(\frac{3}{2}-\epsilon)(k-k_{2})}||P_{k_{1}}\psi_{1}||_{A[k_{1}]}||P_{k_{2}}\psi_{2}||_{{\mathcal{S}}[k_{2}]}\\
\end{split}\end{equation}
One can sum here over $k_{1}=k_{2}+O(1)\geq k+O(1)$, obtaining the
desired estimate. In case of high-low interactions, one reasons as
follows: assume $P_{k_{1}}\psi_{1}$ is of 2nd type. Then for
$j<k+O(1)$, $k=k_{1}+O(1)$, we have
\begin{equation}\nonumber\begin{split}
&P_{k}Q_{j}[P_{k_{1}}\psi_{1}\nabla^{-1}\psi_{2}]=P_{k}Q_{j}[P_{k_{1}}Q_{\geq j-10}\psi_{1}\nabla^{-1}P_{<j-10}Q_{<j-10}\psi_{2}]\\
&+P_{k}Q_{j}[P_{k_{1}}\psi_{1}\nabla^{-1}P_{<j-10}Q_{\geq
j-10}\psi_{2}]+P_{k}Q_{j}[P_{k_{1}}\psi_{1}\nabla^{-1}P_{\geq
j-10}\psi_{2}]\\
\end{split}\end{equation}
We estimate each of the terms on the right: for the first, we have
\begin{equation}\nonumber\begin{split}
&||\sum_{j<k+O(1)}P_{k}Q_{j}[P_{k_{1}}Q_{\geq
j-10}\psi_{1}\nabla^{-1}P_{<j-10}Q_{<j-10}\psi_{2}]||_{\dot{X}_{k}^{0,\frac{1}{2},1}}\\&\lesssim
\sum_{j<k+O(1)}2^{\frac{j}{2}}\sum_{a\geq
j-10}||P_{k_{1}}Q_{a}\psi_{1}||_{L_{t}^{2}L_{x}^{2}}||\nabla^{-1}\psi_{2}||_{L_{t}^{\infty}L_{x}^{\infty}}\\
&\lesssim
||P_{k_{1}}\psi_{1}||_{\dot{X}_{k_{1}}^{0,\frac{1}{2},1}}||\nabla^{-1}\psi_{2}||_{L_{t}^{\infty}L_{x}^{\infty}}\\
\end{split}\end{equation}
Now consider the 2nd term on the right. We exploit the improved
range of Strichartz type norms available for $P_{k_{1}}\psi_{1}$:
\begin{equation}\nonumber\begin{split}
&||P_{k}Q_{j}[P_{k_{1}}\psi_{1}\nabla^{-1}P_{<j-10}Q_{\geq
j-10}\psi_{2}]||_{\dot{X}_{k}^{0,\frac{1}{2},\infty}}\\&\lesssim
2^{\frac{j}{2}}||P_{k_{1}}\psi_{1}||_{L_{t}^{2}L_{x}^{2+}}||\nabla^{-1}P_{<j-10}Q_{\geq
j-10}\psi_{2}||_{L_{t}^{\infty}L_{x}^{M}}\\
&\lesssim 2^{\frac{j-k_{1}}{2+}}||P_{k_{1}}\psi_{1}||_{A[k_{1}]}\\
\end{split}\end{equation}
One can sum over $j<k+O(1)$, resulting in the desired bound.
Finally, the last term in the above trichotomy is handled
similarly. One handles low-high interactions analogously. This
shows that the desired property for $\beta$ is inherited from one
stage of the expansion to the next.\\
Next, consider $P_{k}R_{0}[\psi_{1}A(\nabla^{-1}\psi_{2})]$. One
expands $A(\nabla^{-1}\psi_{2})$ into a Taylor series, and
proceeds inductively. Assume one has
$\sup_{k\in{\mathbf{Z}}}||P_{k}R_{0}\psi_{1,2}||_{L_{t}^{\infty}L_{x}^{2}}\leq
C$, and consider $P_{k}R_{0}[\psi_{1}\nabla^{-1}\psi_{2}]$. One
rewrites this as
\begin{equation}\nonumber\begin{split}
&P_{k}R_{0}[\psi_{1}\nabla^{-1}\psi_{2}]=\sum_{k_{1}>k+10,k_{1}=k_{2}+O(1)}P_{k}R_{0}[P_{k_{1}}\psi_{1}\nabla^{-1}P_{k_{2}}\psi_{2}]\\
&+P_{k}R_{0}[P_{[k-10,k+10]}\psi_{1}\nabla^{-1}\psi_{2}]+P_{k}R_{0}[P_{<k-10}\psi_{1}\nabla^{-1}P_{[k-10,k+10]}\psi_{2}]\\
\end{split}\end{equation}
One treats each of these terms separately: for the first, let the
derivative $\partial_{t}$ fall inside, replacing this by
\begin{equation}\nonumber\begin{split}
&||\sum_{k_{1}>k+10,k_{1}=k_{2}+O(1)}P_{k}\nabla^{-1}[P_{k_{1}}\partial_{t}\psi_{1}\nabla^{-1}P_{k_{2}}\psi_{2}+P_{k_{1}}\psi_{1}R_{0}\psi_{2}]||_{L_{t}^{\infty}L_{x}^{2}}\\
&\lesssim
\sum_{k_{1}=k_{2}+O(1)}[||P_{k_{1}}R_{0}\psi_{1}||_{L_{t}^{\infty}L_{x}^{2}}||P_{k_{2}}\psi_{2}||_{L_{t}^{\infty}L_{x}^{2}}+
||P_{k_{1}}\psi_{1}||_{L_{t}^{\infty}L_{x}^{2}}||R_{0}P_{k_{2}}\psi_{2}||_{L_{t}^{\infty}L_{x}^{2}}]\leq C\\
\end{split}\end{equation}
Next one estimates
\begin{equation}\nonumber\begin{split}
&||P_{k}R_{0}[P_{[k-10,k+10]}\psi_{1}\nabla^{-1}\psi_{2}]||_{L_{t}^{\infty}L_{x}^{2}}\\&\leq
||P_{k}\nabla^{-1}[P_{[k-10,k+10]}\partial_{t}\psi_{1}\nabla^{-1}\psi_{2}+P_{[k-10,k+10]}\psi_{1}R_{0}\psi_{2}]||_{L_{t}^{\infty}L_{x}^{2}}\\
&\lesssim
||P_{[k-10,k+10]}R_{0}\psi_{1}||_{L_{t}^{\infty}L_{x}^{2}}||\nabla^{-1}\psi_{2}||_{L_{t}^{\infty}L_{x}^{\infty}}\\&\hspace{4.5cm}+
||P_{[k-10,k+10]}\nabla^{-1}\psi_{1}||_{L_{t}^{\infty}L_{x}^{2}}||R_{0}P_{<k+15}\psi_{2}||_{L_{t}^{\infty}L_{x}^{\infty}}\\
\end{split}\end{equation}
Again one checks that this can be bounded by $C$. The third term
above is more of the same. The assertion follows from this. We
proceed to the bilinear estimates. For
this we expand the expression\\
$P_{k}[\psi_{1}A(\nabla^{-1}\psi_{2})]$ as a Taylor series,
obtaining terms of the schematic form $\psi$,
$\psi\nabla^{-1}\psi$, $\psi\nabla^{-1}\psi\nabla^{-1}\psi$ etc.
Consider a typical such term of the form
\begin{equation}\nonumber
P_{0}[\psi_{1}\nabla^{-1}\psi_{2}\nabla^{-1}\psi_{3}\ldots\nabla^{-1}\psi_{a}]
\end{equation}
This term will have a coefficient decaying like $(a!)^{-1}$.
Freeze the modulation of the output to dyadic size $\sim 2^{l}$.
If $\nabla^{-1}\psi_{a}$ has modulation $\geq 2^{l+10}$, we can
estimate
\begin{equation}\nonumber\begin{split}
&||P_{0}Q_{l}R_{0}[\psi_{1}\nabla^{-1}\psi_{2}\nabla^{-1}\psi_{3}\ldots\nabla^{-1}Q_{>l+10}\psi_{a}]||_{L_{t}^{2}L_{x}^{2}}\\
&\lesssim
2^{l}||Q_{>l+O(1)}[\psi_{1}\nabla^{-1}\psi_{2}\nabla^{-1}\psi_{3}\ldots\nabla^{-1}\psi_{a-1}]||_{L_{t}^{2}L_{x}^{2}}||\nabla^{-1}Q_{>l+10}\psi_{a}||_{L_{t}^{2}L_{x}^{2}+L_{t}^{2}L_{x}^{\infty}}\\
&\lesssim 2^{-\frac{l}{2+}}
\end{split}\end{equation}
Summing over $l>O(1)$ results in the desired upper bound.
Similarly, we have
\begin{equation}\nonumber\begin{split}
&||\sum_{l>O(1)}P_{0}Q_{l}R_{0}[\psi_{1}\nabla^{-1}\psi_{2}\nabla^{-1}\psi_{3}\ldots\nabla^{-1}Q_{[l-10,l+10]}\psi_{a}]||_{L_{t}^{2}L_{x}^{2}}\\
&\lesssim
(\sum_{l>O(1)}||Q_{[l-10,l+10]}R_{0}\psi_{a}||_{L_{t}^{2}L_{x}^{2}+L_{t}^{2}L_{x}^{\infty}}^{2})^{\frac{1}{2}}\lesssim C\\
\end{split}\end{equation}
which is again acceptable. Thus assume now that
$\nabla^{-1}\psi_{a}$ has modulation $<2^{l-10}$. If
$[\psi_{1}\ldots\nabla^{-1}\psi_{a-1}]$ has modulation $\geq
2^{l-10}$, repeat the same process with this expression instead of
the longer one. Continuing in this fashion, one either
eventually\footnote{We may assume $a<<l$, since otherwise one gets
a large gain in $l$ just from the Taylor coefficient.} forces
$\psi_{1}$ to be at modulation $>2^{l-10a}$, or else one arrives
at a situation of the following sort:
\begin{equation}\nonumber
P_{0}Q_{l}[Q_{>l-10(a-m-1)}[Q_{<l-10(a-m)}[\psi_{1}\ldots\nabla^{-1}\psi_{m}]\nabla^{-1}\psi_{m+1}]\ldots\nabla^{-1}\psi_{a}]
\end{equation}
If $\nabla^{-1}\psi_{m+1}$ has modulation $>2^{l-10(a-m)}$, one
argues just as before. One loses exponentially in $m$, which is
counteracted by the small coefficient eventually applied to the
expression from the Taylor expansion. Thus we may also apply an
operator $Q_{<l-10(a-m)}$ in front of $\nabla^{-1}\psi_{m+1}$. In
this case, however, we can write
\begin{equation}\nonumber\begin{split}
&P_{0}Q_{l}[Q_{>l-10(a-m-1)}[Q_{<l-10(a-m)}[\psi_{1}\ldots\nabla^{-1}\psi_{m}]\nabla^{-1}Q_{<l-10(a-m)}\psi_{m+1}]\ldots\nabla^{-1}\psi_{a}]\\
&=P_{0}Q_{l}[Q_{>l-10(a-m-1)}[P_{l-10(a-m)+O(1)}Q_{<l-10(a-m)}[\psi_{1}\ldots\nabla^{-1}\psi_{m}]\\&\hspace{4.5cm}\nabla^{-1}P_{l-10(a-m)+O(1)}Q_{<l-10(a-m)}\psi_{m+1}]\ldots\nabla^{-1}\psi_{a}]\\
\end{split}\end{equation}
As to the first bilinear inequality in Theorem~\ref{Moser}, we
decompose
\begin{equation}\nonumber
P_{k}\phi=P_{k}Q_{\geq k+100}\phi+\sum_{\pm}\sum_{\kappa\in
K_{-100}}P_{k,\kappa}Q^{\pm}_{<k+100}\phi
\end{equation}
Then we have
\begin{equation}\nonumber\begin{split}
&(\sum_{c\in C_{k,r}}||P_{0}Q_{l}R_{0}[Q_{>l-10(a-m-1)}[P_{l-10(a-m)+O(1)}Q_{<l-10(a-m)}[\psi_{1}\ldots\nabla^{-1}\psi_{m}]\\&\hspace{1cm}\nabla^{-1}P_{l-10(a-m)+O(1)}Q_{<l-10(a-m)}\psi_{m+1}]\ldots\nabla^{-1}\psi_{k}]P_{c}Q_{\geq k+100}\phi||_{L_{t}^{2}L_{x}^{2}}^{2})^{\frac{1}{2}}\\
&\lesssim
2^{\min\{k+r,0\}}||P_{0}Q_{l}R_{0}[Q_{>l-10(a-m-1)}[P_{l-10(a-m)+O(1)}Q_{<l-10(a-m)}[\psi_{1}\ldots\nabla^{-1}\psi_{m}]\\&\hspace{1cm}\nabla^{-1}P_{l-10(a-m)+O(1)}Q_{<l-10(a-m)}\psi_{m+1}]\ldots\nabla^{-1}\psi_{k}]||_{L_{t}^{\infty}L_{x}^{2}}
||P_{k}Q_{\geq k+100}\phi||_{L_{t}^{2}L_{x}^{2}}\\
&\lesssim 2^{\min\{\frac{k}{2},0\}}2^{r}2^{10m}\\
\end{split}\end{equation}
The loss in $m$ will be counteracted by the small Taylor
coefficients. Now consider the contribution of
\begin{equation}\nonumber
\sum_{\pm}\sum_{\kappa\in
K_{-100}}P_{k,\kappa}Q^{\pm}_{<k+100}\phi
\end{equation}
We have the identity
\begin{equation}\nonumber\begin{split}
&P_{0}Q_{l}[Q_{>l-10(a-m-1)}[P_{l-10(a-m)+O(1)}Q_{<l-10(a-m)}[\psi_{1}\ldots\nabla^{-1}\psi_{m}]\\&\hspace{4.5cm}\nabla^{-1}P_{l-10(a-m)+O(1)}Q_{<l-10(a-m)}\psi_{m+1}]\ldots\nabla^{-1}\psi_{k}]\\
&=\sum_{\pm}\sum_{\kappa_{1,2}\in
K_{-100},\,\text{dist}(\kappa_{1},\kappa_{2})\sim 1}
\\&P_{0}Q_{l}[Q_{>l-10(a-m-1)}[P_{l-10(a-m)+O(1),\kappa_{1}}Q^{\pm}_{<l-10(a-m)}[\psi_{1}\ldots\nabla^{-1}\psi_{m}]\\&\hspace{4cm}\nabla^{-1}P_{l-10(a-m)+O(1),\kappa_{2}}Q^{\pm}_{<l-10(a-m)}\psi_{m+1}]\ldots\nabla^{-1}\psi_{k}]\\
\end{split}\end{equation}
We may thus assume that either $\pm\kappa_{1}$ or $\pm\kappa_{2}$
has angular separation $\sim 1$ from $\pm\kappa$. Assume w. l. o.
g. that $\text{dist}(\pm\kappa_{1},\pm\kappa)\sim 1$. Then
estimate
\begin{equation}\nonumber\begin{split}
&(\sum_{c\in C_{k,r}}||P_{0}Q_{l}R_{0}[Q_{>l-10(a-m-1)}[P_{l-10(a-m)+O(1)}Q_{<l-10(a-m)}[\psi_{1}\ldots\nabla^{-1}\psi_{m}]\\&\hspace{1cm}\nabla^{-1}P_{l-10(a-m)+O(1)}Q_{<l-10(a-m)}\psi_{m+1}]\ldots\nabla^{-1}\psi_{k}]P_{c}P_{k,\kappa}Q^{\pm}_{<k+100}\phi||_{L_{t}^{2}L_{x}^{2}}^{2})^{\frac{1}{2}}\\
&\lesssim
2^{l}||[P_{l-10(a-m)+O(1),\kappa_{1}}Q^{\pm}_{<l-10(a-m)}[\psi_{1}\ldots\nabla^{-1}\psi_{m}]||_{NFA^{*}[\pm\kappa]}\\&||\nabla^{-1}P_{l-10(a-m)+O(1)}Q_{<l-10(a-m)}\psi_{m+1}||_{L_{t}^{\infty}L_{x}^{2}}
(\sum_{c\in
C_{k,r}}||P_{c}P_{k,\kappa}Q^{\pm}_{<k+100}\phi||_{PW[\pm\kappa]}^{2})^{\frac{1}{2}}\\
&\lesssim 2^{\frac{k}{2}}2^{\frac{r}{2+}}2^{l}||\nabla^{-1}P_{l-10(a-m)+O(1)}Q_{<l-10(a-m)}\psi_{m+1}||_{L_{t}^{\infty}L_{x}^{2}}\\
\end{split}\end{equation}
We have used that
\begin{equation}\nonumber\begin{split}
&||[P_{l-10(a-m)+O(1),\kappa_{1}}Q^{\pm}_{<l-10(a-m)}[\psi_{1}\ldots\nabla^{-1}\psi_{m}]||_{NFA^{*}[\pm\kappa]}\\
&\lesssim||[P_{l-10(a-m)+O(1),\kappa_{1}}Q^{\pm}_{<l-10(a-m)}[\psi_{1}\ldots\nabla^{-1}\psi_{m}]||_{\dot{X}_{l-10(a-m)}^{\frac{1}{2},1}+A[l-10(a-m)]}\leq
C\\
\end{split}\end{equation}
One can sum now over $l$, obtaining the desired estimate with a
loss $2^{10m}$, which is made up for by the small Taylor
coefficient in front. The 2nd bilinear inequality is proved
similarly, as is the version when $A(\nabla^{-1}\psi_{2})$ is
replaced by $A(\nabla^{-1}(\psi_{2}\nabla^{-1}\psi_{3}))$. This
completes the proof of theorem~\ref{Moser}.

\noindent

\medskip\noindent
\textsc{Harvard University, Dept. of Mathematics, Science Center, 1 Oxford Street, Cambridge, MA 02138, USA}\\
{\em email: }\textsf{jkrieger at math.harvard.edu}\\

\end{section}
\end{document}